\documentclass[onecolumn,10pt]{IEEEtran}
\usepackage[normalem]{ulem}
\RequirePackage{amsthm,amsmath,amsfonts,amssymb}
\RequirePackage[numbers,sort&compress]{natbib}

\theoremstyle{plain}
\newtheorem{thm}{Theorem}[section]
\newtheorem{lem}[thm]{Lemma}
\newtheorem{prop}[thm]{Proposition}

\theoremstyle{remark}
\newtheorem{rem}[thm]{Remark}

\newtheorem{exa}[thm]{Example}
\newtheorem{definition}[thm]{Definition}

\usepackage{mathrsfs} 

\usepackage{color}
\usepackage{enumitem}

\usepackage{mathtools}  
\usepackage{verbatim} 
\usepackage{bm} 
\usepackage{graphicx}
\usepackage{subcaption}
\graphicspath{{graphs/}}
\usepackage{pinlabel}

\usepackage[utf8]{inputenc}
\usepackage[T1]{fontenc}

\usepackage{pifont}

\usepackage{hyperref} 
\hypersetup{
	colorlinks=true,
	linkcolor=blue,
	filecolor=magenta,      
	urlcolor=blue,
	allcolors=blue,
}
\usepackage{xcolor} 
\usepackage{imakeidx} 
\makeindex

\def \P{\mathbb{P}}
\def \E{\mathbb{E}}
\def \R{\mathbb{R}}
\def \Sbb{\mathbb{S}}

\def \Ec{\mathcal{E}}
\def \Dc{\mathcal{D}}
\def \Fc{\mathcal{F}}
\def \Gc{\mathcal{G}}
\def \Hc{\mathcal{H}}

\def \Kc{\mathcal{K}}

\def \Nc{\mathcal{N}}

\def \Vc{\mathcal{V}}

\def \Er{\mathrm{E}}
\def \Vr{\mathrm{V}}

\def \Sum{\displaystyle\sum}

\def \lambdamin{\operatorname{\lambda_{min}}}
\def \lambdamax{\operatorname{\lambda_{max}}}
\def \smin{\text{S}_{\text{min}}}
\def \smax{\text{S}_{\text{max}}}

\def \qve{\bm{q}}
\def \uve{\bm{u}}
\def \vve{\bm{v}}
\def \wve{\bm{w}}
\def \wwve{\bm{w}}
\def \xve{\bm{x}}
\def \yve{\bm{y}}
\def \zve{\bm{z}}

\def \muve{\bm{\mu}}

\def \Ama{\bm{A}}
\def \Bma{\bm{B}}
\def \Dma{\bm{D}}
\def \Ema{\bm{E}}

\def \Hma{\bm{H}}
\def \Ima{\bm{I}}
\def \Pma{\bm{P}}
\def \Qma{\bm{Q}}
\def \Rma{\bm{R}}
\def \Sma{\bm{S}}
\def \Tma{\bm{T}}
\def \Uma{\bm{U}}

\def \Xma{\bm{X}}
\def \Yma{\bm{Y}}
\def \Zma{\bm{Z}}
\def \Phima{\bm{\Phi}}
\def \Psima{\bm{\Psi}}
\def \Sigmama{\bm{\Sigma}}

\def \nn{n}

\def \qq{\tau}
\def \varrho{\gamma}

\newif\iftag
\tagtrue
\newcommand\numberthis{\addtocounter{equation}{1}\tag{\theequation}}

\def\VV{V}
\def\EE{E}
\def\VVV{\breve{V}}

\def\Nedged{N_{\EE_\delta}}

\def\mylll{\delta}
\def\kkkk{\kappa}

\newcommand{\YW}[1]{{\color{blue} {#1}}}
\newcommand{\myeoe}{\hfill $\Diamond$}
\newcommand{\ONE}[1]{{ {#1}}}
\newcommand{\ONENEW}[1]{{\color{brown} {#1}}}
\newcommand{\aoh}[1]{{\color{magenta}[AOH]: {#1}}}

\def\nmd{\operatorname{{\bf NUV}}}
\def\pnmd{\operatorname{ {\bf PNUV} }}
\def \diag{\operatorname{{\bf diag } }}

\def \unif{\operatorname{{\bf unif}}}
\def \1{{1}}
\def \CP{\operatorname{CP}}
\def \SC{\operatorname{SC}}
\def \deg{\operatorname{deg}}
\def \TV{\operatorname{TV}}
\def \Pois{\operatorname{Pois}}
\def \Dirac{\operatorname{Dirac}}

\def\mgeo{{\delta+1}}
\def\mgeominus{\delta}
\def\mgeominusminus{\delta-1}
\def\Ngeo{n-2}
\def\Ngeoplus{n-1}
\def\Ngeoplusplus{n}

\allowdisplaybreaks

\begin{document}
\title{A unified framework for correlation mining 
\\ in ultra-high dimension}
\author{Yun Wei, Bala Rajaratnam and Alfred O. Hero 
\thanks{This work was partially supported by the US Army Research Office under grants W911NF-15-1-0479 and W911NF-19-1-0269,  and by the US National Science Foundation under grants NSF DMS 1976787, NSF DMS 1916787 and NSF CCF 1934568.}
\thanks{Yun Wei is with Department of Statistical Science, Duke University, Durham, NC, USA 27707 (Email: yun.wei@duke.edu) and SAMSI. He was with the Department of Mathematics, University of Michigan, Ann Arbor, MI, USA, 48109.}
\thanks{Bala Rajaratnam is with the Department of Statistics, University of California, Davis, CA, USA, 95616 (Email: brajaratnam@ucdavis.edu)}
\thanks{Alfred Hero is with the Department of Electrical Engineering and Computer Science, University of Michigan, Ann Arbor, MI, USA, 48109 (Email: hero@eecs.umich.edu).}
}

\maketitle

\begin{abstract}
Many applications benefit from theory relevant to the identification of variables having large correlations or partial correlations in high dimension. 
Recently there has been progress in the ultra-high dimensional setting when the sample size $n$ is fixed and the dimension $p$ tends to infinity. Despite these advances, the correlation screening framework suffers from 
practical, methodological and theoretical deficiencies. For instance, previous correlation screening theory requires that the population covariance matrix be sparse and block diagonal. This block sparsity assumption is however 
restrictive in 
practical applications. As a second example,  correlation and partial correlation screening  requires the estimation of dependence measures, which can be 
computationally prohibitive. In this paper, we propose a unifying  approach to correlation and partial correlation mining that is not restricted to  block diagonal correlation structure, thus yielding a methodology that is suitable for modern applications. By making connections to random geometric graphs, the number of highly correlated or partial correlated variables are shown to have 
compound Poisson finite-sample characterizations, which hold for both the finite $p$ case and when  $p$ tends to infinity. The unifying  framework also demonstrates 
a duality between correlation and partial correlation screening with 
theoretical and practical consequences.
\end{abstract}


\tableofcontents

\section{Introduction}
\label{sec:intro}

This paper considers the problem of identifying high correlations and partial correlations in modern ultra high dimensional setting. In particular we study the problem of screening $\nn$ identically distributed $p$-variate samples for variables that have high correlation or high partial correlation with at least one other variable when the sample size $\nn \leq C_0\ln p$ for some constant $C_0$. 
In the screening framework one applies a threshold to the sample correlation matrix or the sample partial correlation matrix to detect variables with at least one significant correlation. The threshold serves to separate signal from noise. Correlation and partial correlation screening in ultra-high dimensions with few samples arises frequently in applications where the per-sample cost of collecting high dimensional data is much more costly than the per-variable cost. 
For example, in genomic correlation screening the cost of high throughput RNAseq assays is decreasing faster than the cost of biological samples \cite{hero2016Large}. In such situations $p$ is much larger than $\nn$.

The ultra-high dimensional regime when $\nn\leq C_0\ln p$ is challenging since the number of samples is insufficient to guarantee reliability of commonly applied statistical methods. For example, one way to undertake partial correlation screening is to first estimate the population covariance matrix, then obtain the inverse, from which a partial correlation matrix can be estimated. However, to get a reliable estimate of a general covariance matrix, the number of samples $\nn$ must be at least 
$\Omega(p)$ as shown in \cite[Section 5.4.3]{vershynin2012introduction}. Even if the covariance matrix has a special structure like sparsity, covariance estimation requires the number of samples be of order $\Omega(\ln p)$ 
The reader is referred to  
\cite{dawid1993hyper, letac2007wishart, anandkumar2009detection, castro2014detection, tan2014learning, khare2015convex, firouzi2016two, dalal2017sparse, heydari2017quickest, tarzanagh2018estimation,   engle2019large, cao2019posterior, li2020adaptable, ledoit2020analytical} and the references therein for related work in modern high dimensional covariance selection and estimation. 

While estimating the correlation matrix or partial correlation matrix is challenging in ultra-high dimensions, recent work \cite{hero2011large,hero2012hub} has shown that it is possible to accurately test the number of highly (partial) correlated variables under a false positive probability constraint; in particular the probability that a variable is spuriously (partially) correlated with at least one other variable. While correlation screening finds variables that have a high marginal correlation with at least another variable, partial correlation screening identifies variables that have high conditional correlations with one other variable conditioned on the rest. 
In \cite{hero2011large}, the ultra-high dimensional correlation screening problem is studied under a row-sparsity assumption on the population covariance matrix. A phase transition in the number of false positive correlations was characterized as a function of the correlation threshold and the true covariance. In the case of block sparse covariance, the critical phase transition threshold becomes independent of the true covariance. In \cite{hero2012hub} the partial correlation screening problem was studied, and similar phase transition results as in correlation screening \cite{hero2011large} were obtained under the block-sparse assumption on the population covariance matrix. The survey \cite{hero2015foundational} reviews the correlation and partial correlation screening problems. 
A follow on work \cite{zhang2017spherical} of \cite{hero2011large} applies a similar framework to the spurious correlations problem and the low-rank detection problem. \ONE{The reader is referred to \cite{anandkumar2009detection,castro2014detection, tan2014learning, heydari2017quickest, tarzanagh2018estimation, arias2012detection, fan2018discoveries} and the references therein \ONE{for additional work related to high dimensional sample correlation matrices.} 
} 


Despite these advances in correlation and partial correlation screening, the screening framework proposed in \cite{hero2011large, hero2012hub} has some serious methodological, theoretical and practical shortcomings. For instance, results for partial correlation screening impose a 
restrictive block sparsity assumption on the true underlying correlation matrix. The block sparsity condition in \cite{hero2012hub} requires that only a small group of the variables have correlation within the blocks and have no correlations with variables outside the block. This assumption is severely restrictive for cases where variables have correlations within a group and also correlations with variables outside their respective groups. Furthermore, expressions for false discovery probabilities in \cite{hero2011large, hero2012hub} require estimating dependence functionals. Estimating such functionals lead to computationally prohibitive non-parametric estimation, rendering the screening methodology disconnected from the very setting it was designed for. 

In this paper we propose a  unifying  framework for correlation and partial correlation screening that delivers a practical and scalable variable selection 
in the ultra-high dimensional regime. By making novel connections to random 
geometric graphs \cite{penrose2003random}, we demonstrate that the distribution of the number of discoveries beyond a certain threshold is approximated by a compound Poisson distribution, with different parameters in the regimes when $p$ is finite and when $p$ approaches $\infty$. To the best of our knowledge, such characterization has not previously appeared in the literature. 
Furthermore, our results are proved in greater generality by relaxing the block-sparse assumption to a new sparsity condition, defined as $(\tau,\kappa)$ sparsity in Section \ref{sec:taukappaspa}, on the population covariance matrix. 
The block-sparse assumption is a special case of the $(\tau,\kappa)$ sparsity assumption. The characterizations established in this paper depend on the covariance matrix only through the $(\tau,\kappa)$ sparsity condition. 
Moreover, the assumptions of $(\tau,\kappa)$ sparsity allows us to formulate unified theorems, which covers both the cases of correlation and of partial correlation screening. 

The theory in this paper is directly relevant to hypothesis testing concerning the empirical degree distribution of a correlation graph. This topic arises in a wide spectrum of areas including graph mining, network science, social science, and natural sciences \cite{chakrabarti2006graph,kolaczyk2014statistical,kolaczykstatistical}. 
Variables having strong sample correlations will appear in the correlation graph as vertices having positive vertex degree. As one sweeps over vertex degree values, the histogram of vertex degrees specifies the empirical degree distribution of the graph. 
From this perspective, this paper provides compound Poisson characterizations of the empirical degree distribution for large correlation graphs under more realistic sparsity conditions on the population covariance. The expressions that are derived from our theorems also provide approximations to family-wise error rates associated with false discoveries of vertices of degree exceeding a specified fixed degree.

Finite sample  results for controlling the probability of discovering a
(false) partial correlation in high-dimensional thresholded covariance 
settings have been
elusive for the better part of the last decade and a half. Indeed,
evaluating expressions for such probabilities in the fixed $n$ setting
are known to be a notoriously difficult problem. This difficulty is in
part attributed to the dependence of such false 
positive probabilities on the
unknown covariance parameter. Previous work has instead provided
expressions for the probability that two distinct connectivity
components of the partial correlation graph are falsely joined (see
\cite{banerjee2008model} and the references therein for more detail).
Controlling such probabilities implicitly assumes that the covariance
parameter is block diagonal. Such an assumption is tantamount to
requiring that the  true partial correlation graph is not fully connected, a
restrictive assumption in many application areas. In contrast, the $(\tau,\kappa)$ sparsity condition introduced in this paper allows the underlying graph to be fully connected. 
Moreover, we provide finite sample results for
controlling the probability that a (partial) correlation is falsely
discovered when the population correlation matrix is $(\tau, \kappa)$ sparse.

\ONE{
\subsection{Contribution}
\label{sec:contribution}
\ONE{We summarize the principal contributions of the paper. As above $p$ denotes the number of variables and $n$ denotes the number of samples.  

\begin{enumerate}
    \item 
The paper presents a unified and complete asymptotic analysis of the star subgraph counts, the counts of vertices of a given degree and the counts of vertices above a given degree in the random graphs obtained by thresholding the sample correlation and sample partial correlation. This unification of different types of random counts represents an 
improvement over previous work \cite{hero2011large,hero2012hub} where only 
the counts of vertices above a given degree are studied.

\item 
We approximate the full distributions of the random counts for finite $p$ and as $p\to\infty$.  In addition, we characterize the first and second moments of these random counts. This is a 
generalization of previous results \cite{hero2011large, hero2012hub} that only established approximations for the mean number of random counts and for the probabilities that these counts were positive.   
    \item 
    We obtain a compound Poisson characterization of the distributions of the random counts.  
    The compound Poisson limit and approximation are well approximated by the standard Poisson limit when $n$ is moderately large (Section \ref{sec:limitingcompoundPoisson}). This result corrects and refines the claim in \cite{hero2012hub} that erroneously asserted a Poisson limit.

    \item 
    The theory in this paper is developed under a novel sparsity condition on the population dispersion matrix. This sparsity condition, called   $(\tau, \kappa)$ sparsity in Sec. \ref{sec:taukappaspa}, is 
    weaker than previously assumed conditions, which makes our theory more broadly applicable. Specifically, while the block sparsity condition in previous work  \cite{hero2012hub,fan2008sure,firouzi2016two} 
    imposes that correlation can only occur locally in small blocks of variables,  the $(\tau,\kappa)$ sparsity condition relaxes this condition to more general global correlation patterning. 

    \end{enumerate}
    
 
Some of the broader implications of the technical contributions of this paper are described below.

    Previous conditions on population partial correlation networks assume they are of lower dimension.
    In particular, conditions such as block sparsity do not allow for completely connected partial correlation graphs which involve all $p$ variables. Such restrictive assumptions are difficult to validate and rule out many realistic population (partial) correlation structures. Overcoming this hurdle has been an open problem for several years. The newly introduced  $(\tau,\kappa)$ sparsity condition on the population covariance matrix settles this longstanding problem by successfully allowing for completely connected partial correlation graphs over the entire set of features. 

Historically, the literature on correlation estimation and graphical models has separated the treatments of covariance graph models and undirected graphical models (or inverse covariance graph models) \cite{cox2014multivariate}. Unifying the two classes of statistical models has been an open problem for the better part of almost 3 decades. While this separate treatment may be appropriate in low dimensional settings when there are few variables, it is not immediately obvious which of the two frameworks is appropriate for a given data set in modern ultra-high dimensional regimes. To our knowledge, the framework in this paper is the first in the graphical model or correlation graph estimation literature to propose methodology which brings both approaches under one umbrella.  

The results in the paper also have 
relevance to applications. Recall that our Poisson and compound Poisson expressions effectively describe the number of false discoveries and hence allow us to obtain results for the familywise error rate (FWER) or k-FWER, that is the probability of obtaining k or more false discoveries. Note that we can also obtain the marginal distributions of correlation 
estimates, which in turn allow us to obtain expressions for p-values for testing correlation estimates. These marginal p-values allow us to establish FDR control too using either the Benjamini-Hochberg \cite{benjamini1995controlling} or Benjamini-Yakutelli procedures \cite{benjamini2001control}. In summary, the correlations screening framework is sufficiently rich that it allows us to undertake statistical error control in terms of FWER, k-FWER and FDR. This is one of the main strengths of our results: a rigorous inferential framework in the ultra-high dimensional setting.
}
}

The remainder of this paper is organized as follows. Section \ref{sec:mainresult} outlines the framework and presents our main theorem which characterizes the compound Poisson approximations when $p$ approaches $\infty$. In Section \ref{mainres} an approximating theorem when $p$ is finite is presented, based on which the main theorem follows. Section \ref{sec:convergenceofmoments} covers convergence of moments. Section \ref{sec:limitingcompoundPoisson} provides explicit expressions for the parameters of the compound Poisson characterizations. 
Notation and symbols used in this paper are collected in the Section \ref{sec:symbols} of the Appendix. 
Most of the technical proofs and auxiliary results are given in the Appendix. 

\section{Main results}
\label{sec:mainresult}

\subsection{Framework}
\label{Hubfra}
\noindent
Available is a data matrix consisting of multivariate samples

\begin{equation}
\Xma= [\xve^{(1)},\xve^{(2)},\cdots, \xve^{(\nn)}]^\top= [\xve_1,\xve_2,\cdots, \xve_p] 
\in \R^{\nn\times p}, \label{Xtildedef}
\end{equation}
where $\{\xve^{(i)}\}_{i=1}^{\nn} \subset \R^p$ are samples from a $p$-dimensional distribution. We assume that the $\nn\times p$ data matrix $\Xma$ follows a vector elliptically contoured distribution \cite{dawid1977spherical, anderson1990theory,anderson1992nonnormal}. 
A random matrix $\Xma\in \R^{\nn\times p}$ is vector-elliptical 
with positive definite covariance or dispersion parameter $\Sigmama\in \R^{p\times p}$ and location parameter $\muve$ if its density satisfies
\begin{equation}
f_{\Xma}(\Xma) = \text{det}(\Sigmama)^{-\nn/2}g(\text{tr}((\Xma-\bm{1}\muve^\top)\Sigmama^{-1}(\Xma^\top-\bm{1}\muve^\top))),  \label{eqn:ellden}
\end{equation}
for a shape function $g: \R\to [0,\infty)$ such that $\int f_{\Xma}(\Xma) =1$. In \eqref{eqn:ellden}, $\bm{1}$ is a column vector with all elements equal to $1$. 
We use the shorthand $\Xma\sim \mathcal{VE}(\bm{\mu},\Sigmama,g)$ to denote that $\Xma$ follows a vector elliptically contoured distribution with density \eqref{eqn:ellden}. Note that the rows $\{\xve^{(i)}\}_{i=1}^n $ of $\Xma$  are uncorrelated but not necessarily independent \cite{anderson1992nonnormal}. An example of a vector-elliptical distributed is the matrix normal distribution, for which the rows $\{\xve^{(i)}\}_{i=1}^{\nn}\subset \R^p$ are i.i.d. samples from $\Nc(\muve,\Sigmama)$. Specifically, the matrix normal density is obtained when, in \eqref{eqn:ellden}, $g(w)=g_0(w)=(2\pi)^{-\frac{\nn p}{2}}\exp(-\frac{1}{2}w)$ and in this case $\Xma\sim \mathcal{VE}(\muve,\Sigmama,g_0)$. 

Given a data matrix $\Xma\sim \mathcal{VE}(\muve,\Sigmama,g)$\footnote{\ONE{In previous work \cite{hero2011large,hero2012hub} it was assumed that the samples $\xve^{(i)}$ are i.i.d. elliptical contoured distributed. This condition is in fact insufficient and the stronger vector elliptical contoured distribution condition (\ref{eqn:ellden}) 
is required.} 
}, the sample mean $\bar{\xve}$ is given as a row vector
$$
\bar{\xve}=\frac{1}{\nn}\Sum_{i=1}^{\nn}\xve^{(i)} = \frac{1}{\nn}\Xma^\top\bm{1}.
$$
The sample covariance matrix $\Sma$ is
 \begin{equation}
 \Sma=\frac{1}{\nn-1}\sum_{i=1}^{\nn} (\xve^{(i)}-\bar{\xve})(\xve^{(i)}-\bar{\xve})^\top
 =\frac{1}{\nn-1}\Xma^\top\left(\Ima_n -\frac{1}{n}\bm{1}\bm{1}^\top\right) \Xma.\label{def:samcov}
 \end{equation}
The sample correlation matrix $\Rma$ is defined as:
\begin{equation}
\Rma=\diag(\Sma)^{-\frac{1}{2}}\Sma \diag(\Sma)^{-\frac{1}{2}}, \label{defR}
\end{equation}
where $\diag(\Ama)$ for a matrix $\Ama\in \R^{n\times n}$ is the diagonal part of $\Ama$ and $\Bma^{-1/2}$ for a diagonal matrix $\Bma$ is a diagonal matrix formed by raising every diagonal element of $\Bma$ to the power $-1/2$. Since $\Rma$ is not invertible, we define the sample partial correlation matrix  $\Pma$ by
\begin{equation}
\label{Pdef}
\Pma = \diag(\Rma^\dagger)^{-\frac{1}{2}}\Rma^\dagger \diag(\Rma^\dagger)^{-\frac{1}{2}}.
\end{equation}
where $\Rma^\dagger$ is the Moore-Penrose pseudo-inverse of $\Rma$.

  Let $\Psima=(\Psi_{ij})_{i,j\in [p]}$ be generic notation for a correlation-type matrix like $\Rma$ or $\Pma$. Given a threshold $\rho\in [0, 1)$   define the undirected graph induced by thresholding $\Psima$, denoted by $\Gc_\rho(\Psima)$,  as follows. The vertex set of graph $\Gc_\rho(\Psima)$ is $\Vr^{(\Psima)} = [p]:=\{1,2,\cdots,p\}$ and the edge set is $\Er^{(\Psima)} \subset \Vr^{(\Psima)}\times\Vr^{(\Psima)}$, with $(i,j) \in \Er^{(\Psima)}$ if $|\Psi_{ij} | \geq \rho$, where $(i,j)$ denotes an edge between $i$ and $j$ $(i\not =j)$.  We call $\Gc_\rho(\Psima)$ the empirical correlation graph and the empirical partial correlation graph, respectively, when $\Psima = \Rma$ and $\Psima = \Pma$. Let $\Phima^{(\Psima)}(\rho)$ be the adjacency matrices associated with the graph $\Gc_\rho(\Psima)$, with elements 
 $\Phi_{ij}^{(\Psima)}(\rho):=\1(|\Psi_{ij} | \geq \rho)$ for $i\not = j$, where $\1(\cdot)$ is the indicator function.
  The dependence of $\Phima^{(\Psima)}(\rho)$ and $\Phi_{ij}^{(\Psima)}(\rho)$ on $\rho$ will be suppressed when it is clear from the context. 
  
  The focus of this paper is correlation and partial screening, which counts the number of 
  vertices of prescribed degree, the number of star subgraphs, or the number of edges in  $\Gc_\rho(\Psima)$. 
  The objective is to characterize the  
  distributions of these counting statistics. 
  More specifically, for the graph $\Gc_\rho(\Psima)$ with $\Psima =\Rma$ or $\Psima =\Pma$, the degree of vertex $i$ is defined as $\sum_{j=1,j\not =i}^p \Phi_{ij}^{(\Psima)}(\rho)$.
  For $1\leq \delta \leq p-1$, the total number of vertices with degree exactly  $\delta$ (at least $\delta$), denoted by  $N_{\VVV_\delta}^{(\Psima)}$ ($N_{\VV_\delta}^{(\Psima)}$), is of particular interest. Note that if a vertex has degree exactly $\delta$, then there exists a star subgraph with $\delta$ edges centered at that vertex. Consequently the number of star subgraphs are important in the analysis of $N_{\VVV_\delta}^{(\Psima)}$ and $N_{\VV_\delta}^{(\Psima)}$, hence our interest in the number of star subgraphs.   For $2\leq \delta\leq p-1$, the number of subgraphs in $\Gc_\rho(\Psima)$ that are isomorphic to $\Gamma_\delta$ is denoted by
 $\Nedged^{(\Psima)}$, where $\Gamma_\delta$ denotes a star shaped graph with $\delta$ edges. In the case when $\delta = 1$, we define $N_{E_1}^{(\Psima)}$ to be twice the number of edges in $\Gc_\rho(\Psima)$. $\Nedged^{(\Psima)}$ is referred as the \emph{star subgraph counts}. 
 
\begin{exa}
    \begin{figure}[ht!]
        \labellist
        \small \hair 2pt
        \pinlabel {$2$} at 52 417 
        \pinlabel $1$ at 525 51
        \pinlabel $5$ at 968 239
        \pinlabel $3$ at 579 548
        \pinlabel $4$ at 1049 826
        \pinlabel \rotatebox{-35}{$e_1$} at 241 219
        \pinlabel \rotatebox{20}{$e_2$} at 780 121
        \pinlabel \rotatebox{80}{$e_3$} at 507 334 
        \pinlabel \rotatebox{-30}{$e_4$} at 810 437
        \pinlabel \rotatebox{25}{$e_5$} at 801 725
        \endlabellist
        \centering
        \includegraphics[width=0.3\linewidth]{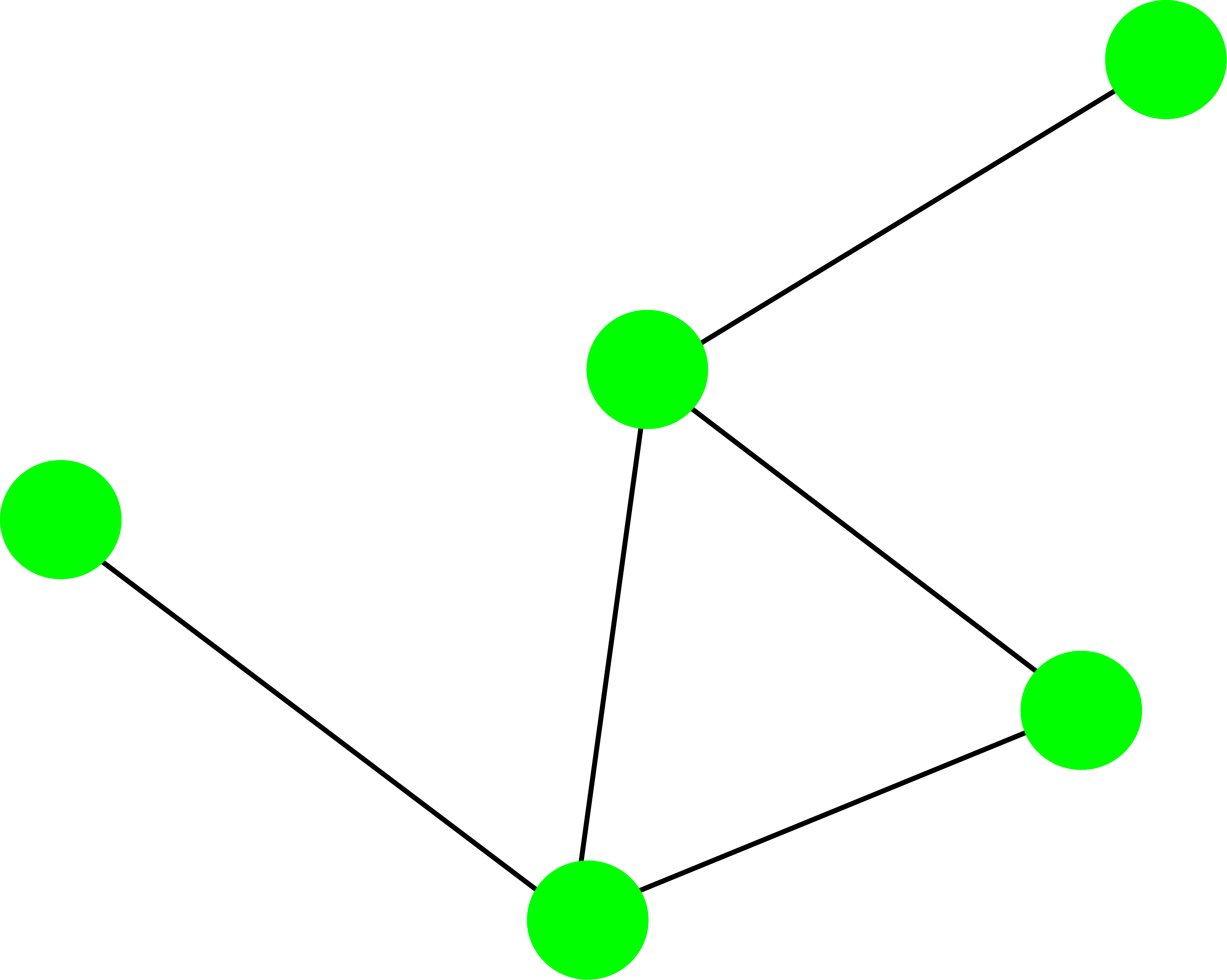}
        \caption{A graph with $5$ vertices and $5$ edges.}
        \label{fig:exaofdef}
    \end{figure}
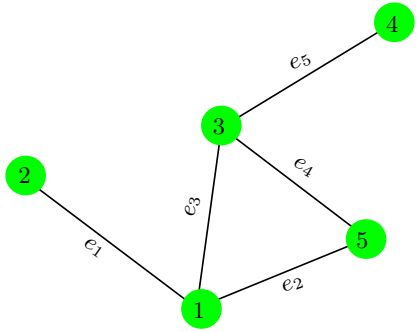

Figure \ref{fig:exaofdef} represents an empirical partial correlation graph. For this graph the  number of vertices of degree $2$ is  $N_{\VVV_2}^{(\Pma)}=1$ and the number of vertices of degree at least $2$ is $N_{\VV_2}^{(\Pma)}=3$. The number of subgraphs isomorphic to $\Gamma_3$ is $N_{E_3}^{(\Pma)} = 2 $. The number of connected vertices is $N_{\VV_1}^{(\Pma)}=5$, and $N_{E_1}^{(\Pma)} = 10 $ as there are $5$ edges.
\end{exa}

 Consider now the case where the sample size $n$ is fixed and there exists a sequence of data matrices $\Xma \in \R ^{n\times p}$ with increasing dimension $p$. Following the procedure described in the paragraph after \eqref{Pdef}, we obtain a sequence of random graphs $G_{\rho}(\Psima)$ with increasing number of vertices. This paper derives finite sample compound Poisson characterizations of the distributions of the $6$ random quantities $\{N_{E_\delta}^{(\Psima)},N_{\breve{V}_\delta}^{(\Psima)},N_{V_\delta}^{(\Psima)}: \Psima \in \{\Rma,\Pma\}  \}$, for $p$ finite and as $p\to\infty$, for suitably chosen $\rho$, under a sparsity assumption on the dispersion parameter $\Sigmama$. \ONE{Such characterizations can be used to test the sparsity structure of the dispersion parameter $\Sigmama$ or to guide the choice of the threshold $\rho$ \cite{hero2011large, hero2012hub}.}
 Throughout the rest of the paper we use $\bar{N}_\delta$ to denote a generic random variable equal to one of the $6$ quantities $\{N_{E_\delta}^{(\Psima)},N_{\breve{V}_\delta}^{(\Psima)},N_{V_\delta}^{(\Psima)}: \Psima \in \{\Rma,\Pma\} \}$. By abuse of terminology, we refer to $\bar{N}_\delta$ generically as vertex counts. We reiterate that in this paper the number $n$ of samples is fixed and the number $p$ of variables  could either be finite or tend to infinity.

\subsection{A unified theorem}
\label{sec:unifiedtheorem}

In this subsection we present a unified theorem that establishes that $\bar{N}_\delta$ converges in distribution to a compound Poisson distribution  when $p\to \infty$. We begin by defining necessary quantities and then we state our main theorem.

For any positive number $\lambda$ and a probability distribution $\bm{\zeta}$ supported on positive integers, let $\CP(\lambda,\bm{\zeta})$ denote the corresponding compound Poisson distribution. Specifically, $\CP(\lambda,\bm{\zeta})$ is the distribution of $ Z= \sum_{i=1}^N Z_i $, where $N$ is distributed as a Poisson random variable with mean $\lambda$, $Z_i \overset{\text{i.i.d.}}{\sim} \bm{\zeta} $ and $N$ is independent of each $Z_i$. Here the random variable $N$ is the number of occurrences of increments and $\bm{\zeta}$ is the distribution of each increment. The parameter $\lambda$ and $\bm{\zeta}$ are often referred to as the arrival rate and the increment distribution, respectively.

As the parameters of the compound Poisson distribution in the next theorem involve a random geometric graph, we define relevant notation. Given a set of points $\{ \vve_i  \}_{i=1}^\delta$ in $\R^{n-2}$, denote by 
 $\textbf{Ge}\left(\{\vve_i  \}_{i=1}^\delta, r   \right)$
 the geometric graph with radius $r$, defined as follows. The vertex set of the graph is $\{ \vve_i  \}_{i=1}^\delta$, and there is an edge between $\vve_i$ and $\vve_j$ if $\|\vve_i-\vve_j\|_2 \leq r$. The graph is called a random geometric graph when the vertices of the geometric graph are random.  A \emph{universal vertex} is a vertex of an undirected graph that is adjacent to all other vertices of the graph \cite{larrion2004clique}. Denote by $\nmd\left(\{\vve_i  \}_{i=1}^\delta, r   \right)$
the number of universal vertices in $\textbf{Ge}\left(\{\vve_i  \}_{i=1}^\delta, r   \right)$. Denote $B^{n-2}$ the unit sphere in $\R^{n-2}$ and denote $\unif(B^{n-2})$ the uniform distribution on $B^{n-2}$. Let $\{\tilde{\uve}_i\}_{i=1}^\delta$ be i.i.d. from $\unif(B^{n-2})$. For the random geometric graph $\textbf{Ge}\left(\{\tilde{\uve}_i\}_{i=1}^\delta, 1   \right)$, we denote the probability that there are exactly $\ell-1$ universal vertices by 
\begin{equation}
\alpha_\ell := \P\left( \nmd\left( \{\tilde{\uve}_i  \}_{i=1}^\delta, 1   \right)= \ell -1\right), \quad \forall \ell\in [\delta+1], \label{eqn:alphaelldef}
\end{equation}
and define a probability distribution $\bm{\zeta}_{n,\delta}$ on $[\delta+1]$: 
\begin{equation}
\bm{\zeta}_{n,\delta}(\ell):=(\alpha_\ell/\ell)/\left(\sum_{s=1}^{\delta+1} (\alpha_s/s) \right),  \quad \forall \ell\in [\delta+1].  \label{eqn:incrementsizelimitdef}
\end{equation}
As will be shown in the next theorem, $\bm{\zeta}_{n,\delta}$ is the increment distribution of a compound Poisson approximation to $\bar{N}_\delta$ when $p\to\infty$. 

We also introduce the following sparsity conditions:
a matrix is said to be row-$\kappa$ sparse if every row has at most $\kappa$ nonzero elements. This is a weaker sparsity condition than the block sparsity condition of \cite{hero2012hub} (see also Subsection \ref{sec:taukappaspa}). The next definition is a stronger sparsity condition than row-$\kappa$ sparsity but remains weaker than block sparsity.

\begin{definition}[$(\tau,\kappa)$ sparsity]
A $p$ by $p$ dimensional symmetric matrix is said to be $(\tau,\kappa)$ sparse if it is  row-$\kappa$ sparse and its lower $p-\tau$ by $p-\tau$ block is diagonal. 
\end{definition}


Another relevant quantity is the normalized determinant defined as follows:
\begin{definition}[Normalized determinant]\label{def:nordet}
    For any symmetric, positive definite matrix $\Ama \in \R^{p\times p}$, its normalized determinant $\mu(\Ama)$ is defined by $$\mu(\Ama):= \prod\limits_{i=1}^p \frac{\lambda_i(\Ama)}{\lambda_p(\Ama)} = \frac{\text{det}(\Ama)}{\left(\lambda_p(\Ama)\right)^p},$$
    where $\lambda_1(\Ama)\leq \lambda_2(\Ama)\leq \cdots \leq \lambda_p(\Ama)$ are the eigenvalues of $\Ama$.
\end{definition}
For $\mathcal{I}\subset [p]$ denote by $\Ama_{\mathcal{I}}$ the set of all $|\mathcal{I}|\times |\mathcal{I}|$ submatrices of $\Ama\in \R^{p\times p}$ obtained by extracting the corresponding rows and columns indexed by $\mathcal{I}$. The set $\Ama_{\mathcal{I}}$ contains $|\mathcal{I}|!$ matrices that are all equivalent up to a permutation applied simultaneously to both rows and columns. Define the \textit{local normalized determinant of degree $m$} of a matrix $\Ama\in \R^{p\times p}$ to be $\mu_m(\Ama) = \min\{\mu(\Ama_{\mathcal{I}}) : \mathcal{I}\subset [p], |\mathcal{I}|=m \} $. Note that $\mu(\Ama_{\mathcal{I}})$ is well defined since $\mu(\cdot)$ is invariant to simultaneous application of a permutation to both rows and columns. For $\Ama\in \R^{p\times p}$ further define the \emph{inverse local normalized determinant}
\begin{equation}
\mu_{n,m} (\Ama) : =\begin{cases} [\mu_{m}(\Ama)]^{-\frac{n-1}{2}}, &  \Ama \text{ symmetric positive definite but not diagonal, }  \\
1, & \Ama \text{ symmetric positive definite and diagonal. }
 \end{cases} \label{eqn:munmA} \end{equation}
By definition $\mu(A)\in (0,1]$ and $\mu_{n,m} (\Ama) \in [1,\infty)$.

Denote $\Gamma(x)$ the gamma function and let $a_n := \frac{\Gamma((\nn-1)/2)}{(n-2)\sqrt{\pi}\Gamma((\nn-2)/2)}$. With the above definitions in place, we now state our main theorem: when $p 
\to \infty$, if the threshold $\rho$ 
approaches $1$ at a particular rate, then the sequence of vertex counts $\bar{N}_\delta$ converges in distribution to a compound Poisson distribution. 

\begin{thm}
    [Compound Poisson limit] \label{cor:Poissonlimit}
    Let $ n \geq 4$ and $\delta$ be fixed positive integers.  
    Let $\Xma\sim \mathcal{VE}(\muve,\Sigmama,g)$.  
    Assume that the threshold $\rho$ is a function of $p$ that satisfies $a_n 2^{\frac{n}{2}}p^{1+\frac{1}{\delta}}(1-\rho)^{\frac{n-2}{2}}\to e_{n,\delta}$ as $p\to \infty$, where $e_{n,\delta}$ is some positive finite constant that possibly depends on $n$ and $\delta$. Denote $\lambda_{\nn,\delta}(e_{n,\delta}) = \frac{1}{\delta !} \left( e_{n,\delta}\right)^\delta \sum_{\ell=1}^{\delta+1}\frac{\alpha_\ell}{\ell}$. Suppose $\Sigmama$, after some row-column permutation, is $(\tau_p,\kkkk_p)$ sparse with $\lim\limits_{p\to \infty} \frac{\tau_p}{p}+\mu_{n,2\delta+2}\left(\Sigmama\right)\frac{\kkkk_p}{p}  \to 0$. Then $\bar{N}_\delta$,  a generic random variable in the set $\{N_{E_\delta}^{(\Psima)},N_{\breve{V}_\delta}^{(\Psima)},N_{V_\delta}^{(\Psima)}: \Psima \in \{\Rma,\Pma\}  \}$, satisfies: 
    \begin{equation}
    \bar{N}_\delta \overset{\Dc}{\to} \CP(\lambda_{n,\delta}(e_{n,\delta}),\bm{\zeta}_{n,\delta}) \text{  as } p\to \infty. \label{eqn:compoilimdis2} 
    \end{equation}
    \end{thm}
 
If only the vertex counts in the empirical correlation graph is of interest, then the $(\tau,\kappa)$ sparsity assumption can be relaxed to row-$\kappa$ sparsity.

\begin{lem}[Compound Poisson limit in empirical correlation graph]
\label{item:Poissonlimita}
 Let $ n \geq 4$ and $\delta$ be fixed positive integers.  
 Let $\Xma\sim \mathcal{VE}(\muve,\Sigmama,g)$.
 Assume that the threshold $\rho$ is a function of $p$ that satisfies $a_n 2^{\frac{n}{2}}p^{1+\frac{1}{\delta}}(1-\rho)^{\frac{n-2}{2}}\to e_{n,\delta}$ as $p\to \infty$, where $e_{n,\delta}$ is some positive finite constant that possibly depends on $n$ and $\delta$. Denote $\lambda_{\nn,\delta}(e_{n,\delta}) = \frac{1}{\delta !} \left( e_{n,\delta}\right)^\delta \sum_{\ell=1}^{\delta+1}\frac{\alpha_\ell}{\ell}$. Suppose $\Sigmama$ is row-$\kkkk_p$ sparse with $\lim\limits_{p\to \infty} \mu_{n,2\delta+2}\left(\Sigmama\right)\frac{\kkkk_p}{p} \to 0$. Then $\tilde{N}_\delta$, a generic random variable in the set $\{N_{E_\delta}^{(\Rma)},N_{\breve{V}_\delta}^{(\Rma)},N_{V_\delta}^{(\Rma)}  \}$, satisfies:
    \begin{equation}
    \tilde{N}_\delta \overset{\Dc}{\to} \CP(\lambda_{n,\delta}(e_{n,\delta}),\bm{\zeta}_{n,\delta}) \text{ as } p\to \infty. \label{eqn:compoilimdis}
    \end{equation}
\end{lem}

\begin{rem} \label{rem:remarkofthm1}
    
    
    The condition $a_n2^{\frac{n}{2}}p^{1+\frac{1}{\delta}}(1-\rho)^{\frac{n-2}{2}}\to e_{n,\delta}>0$ is equivalent to  
    $$
    p^{\frac{2}{n-2}\left(1+\frac{1}{\delta}\right)}(1-\rho)\to \left(\frac{e_{n,\delta}}{a_n2^{\frac{n}{2}}}\right)^{\frac{2}{n-2}}= \frac{1}{2}\left(\frac{e_{n,\delta}}{2a_n}\right)^{\frac{2}{n-2}} ,
    $$ 
    which indicates that $\rho\to 1$ at rate $p^{-\frac{2}{n-2}\left(1+\frac{1}{\delta}\right)}$.  
       { As will be discussed in more detail in Remark \ref{rem:rhorate}, this rate is in fact both necessary and sufficient for the expected counts $\E \tilde{N}_\delta$ to converge to a non-trivial limit. If $\rho$ does not converge to $1$, or converges to  $1$ at a slower rate, then $\E\tilde{N}_\delta$ diverges  to $\infty$, while if $\rho$ converges to $1$ at a faster rate then $\E \tilde{N}_\delta$ converges to $0$. 
       }
\ONE{This particular rate on $\rho$ is consistent with the rate 
    of the existing Poisson approximation results in 
    random geometric graphs \cite{penrose2003random} as will be 
    discussed in Section \ref{sec:rangeo}.} A sequence of correlation thresholds $\rho=\rho_p$ that satisfies this condition is
    \begin{equation}
    \rho_p=1-\frac{1}{2}\left(\frac{e_{n,\delta}}{2a_n p^{1+\frac{1}{\delta}}}\right)^{\frac{2}{n-2}}. \label{eqn:rhopformula}
    \end{equation} 
    Observe that Theorem \ref{cor:Poissonlimit} and Lemma \ref{item:Poissonlimita} hold for any mean $\muve$ and any shaping function $g$ when $\Xma\sim \mathcal{VE}(\muve,\Sigmama,g)$. We will provide intuition for this invariance property in Remark \ref{rem:Udisconsequence}.
\end{rem}

The proofs of Theorem \ref{cor:Poissonlimit} and Lemma \ref{item:Poissonlimita} will be presented in Subsection \ref{sec:compoihighdim}. 
We will call $\CP(\lambda_{n,\delta}(e_{n,\delta}),\bm{\zeta}_{n,\delta})$ the \emph{limiting compound Poisson distribution, approximation or characterization.} Since $\bar{N}_\delta$ and $\CP(\lambda_{n,\delta}(e_{n,\delta}),\bm{\zeta}_{n,\delta})$ are discrete, \eqref{eqn:compoilimdis2} is equivalent to:
    $$
    d_{\TV}\left(\mathscr{L}\left(\bar{N}_\delta\right),\CP(\lambda_{n,\delta}(e_{n,\delta}),\bm{\zeta}_{n,\delta})\right) \to 0 \text{ as } p\to \infty,
    $$
    where $\mathscr{L}\left(\cdot\right)$ represents the probability distribution of the argument,  
    and $d_{\TV}\left(\cdot,\cdot\right)$ is the total variation distance between two probability distributions.
    A variant of Theorem \ref{cor:Poissonlimit} for finite $p$ (Theorem \ref{thm:Poissonultrahigh}), which establishes an upper bound on the total variation distance between $\mathscr{L}\left(\bar{N}_\delta\right)$ and a compound Poisson distribution, will be presented in Subsection \ref{sec:compoihighdim}. \ONE{Theorem \ref{cor:Poissonlimit} and Theorem \ref{thm:Poissonultrahigh} specify asymptotic compound Poisson limits and non-asymptotic 
    bounds on the full distribution of vertex counts. These limits correct and extend the Poisson limits that were falsely claimed to hold for all finite $n,\delta$, 
    although  we note that the compound Poisson limit can be well approximated by the Poisson limit  in the case of moderately large $n$ or large $\delta$ (see Sec. \ref{sec:limitingcompoundPoisson}).}   


In Theorem \ref{cor:Poissonlimit} the $(\tau,\kappa)$ sparsity condition and the condition $\mu_{n,m}(\Sigmama)<\infty$ are assumed. We elaborate on these two conditions in the next two subsections. 


\subsection{$(\qq,\kappa)$ sparsity}
\label{sec:taukappaspa}


The matrix \eqref{examat} below  is an example of a $(\tau,\kappa)$ sparse matrix 
with $\tau=2, \kappa=3$. This $5\times5$ symmetric matrix is $(2,3)$ sparse since each of the  first $2$ rows has at most $3$ nonzero elements and the lower $3\times3$ block is diagonal. 
\begin{equation}
\begin{pmatrix}
5 &  0  & 2  &  0  & 1 \\
0 & 8 & 3 & 0 & 0 \\
2 & 3 & 6 & 0 & 0 \\
0 & 0 & 0 & 7 & 0 \\
1 & 0 & 0 & 0 & 8
\end{pmatrix}   \label{examat}
\end{equation}

If the adjacency matrix of a graph $(\Vc,\Ec)$ is $(\tau,\kappa)$ sparse, then the vertices $\Vc$ can be partitioned into two disjoint subsets $\Vc_1$ and $\Vc_2$ with the following properties: 1) $|\Vc_1|\leq \tau$; 2) there is no edge between any two vertices in $\Vc_2$; 3) the degree of any vertex in $\Vc_1$ is no more than $\kappa-1$; 4) edges connecting vertex in $\Vc_1$ and $\Vc_2$ may exist.

When the dispersion parameter $\Sigmama$ is row-$\kappa$ sparse, the authors of \cite{hero2011large} studied the mean of the quantities $N^{(\Rma)}_{E_1}$ and $N^{(\Rma)}_{V_1}$ and they obtained limits of $\P(N_{E_1}^{(\Rma)}>0)$ or $\P(N_{V_1}^{(\Rma)}>0)$ when $p\to\infty$ while fixing $n$. In \cite{hero2012hub} these results were extended to empirical partial correlation graphs when the dispersion parameter $\Sigmama$ is assumed to be block-$\qq$ sparse up to a row-column permutation, i.e., there exists a permutation matrix $\Tma$ such that
\begin{equation}
\Tma\Sigmama\Tma^\top = \begin{pmatrix}
\Sigmama_{11} &  \Sigmama_{12} \\
\Sigmama_{21} & \Dma_{p-\qq}
\end{pmatrix} \label{eqn:permutationresult}
\end{equation}
where $\Sigmama_{12} = \Sigmama_{21}^\top = \bm{0}\in \R^{\qq\times (p-\qq)}$ and $\Dma_{p-\qq}\in \R^{(p-\qq)\times (p-\qq)}$  is some diagonal matrix. 
In Theorem \ref{cor:Poissonlimit} of this paper $\Sigmama$ is assumed to be $(\tau,\kappa)$ sparse after row-column permutation, i.e. there exists a permutation matrix $\Tma$ such that \eqref{eqn:permutationresult} holds with $\Dma_{p-\qq}\in \R^{(p-\qq)\times (p-\qq)}$  some diagonal matrix and with the first $\tau$ rows $(\Sigmama_{11}\ \Sigmama_{12})$ being row-$\kappa$ sparse.

\begin{figure}[ht]
\begin{subfigure}{.49\textwidth}
  \centering
  \includegraphics[width=.8\linewidth]{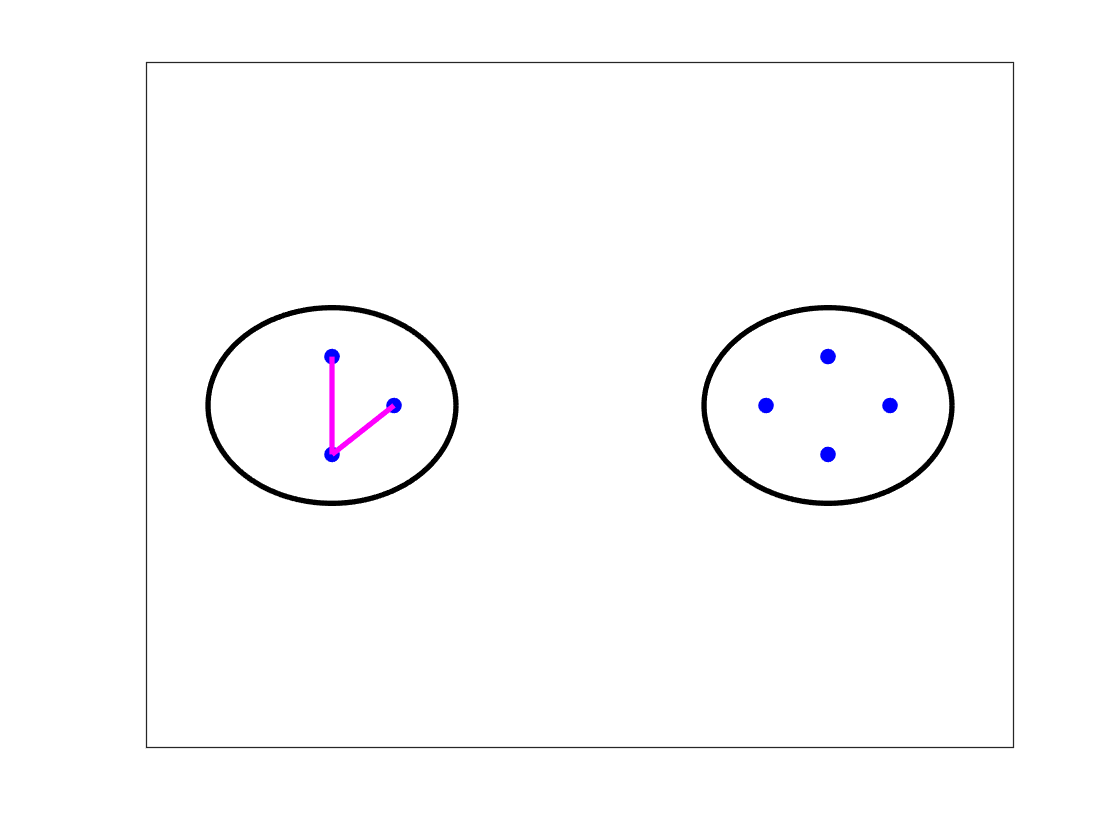}  
  \caption{block sparsity}
\end{subfigure}
\begin{subfigure}{.49\textwidth}
  \centering
  \includegraphics[width=.8\linewidth]{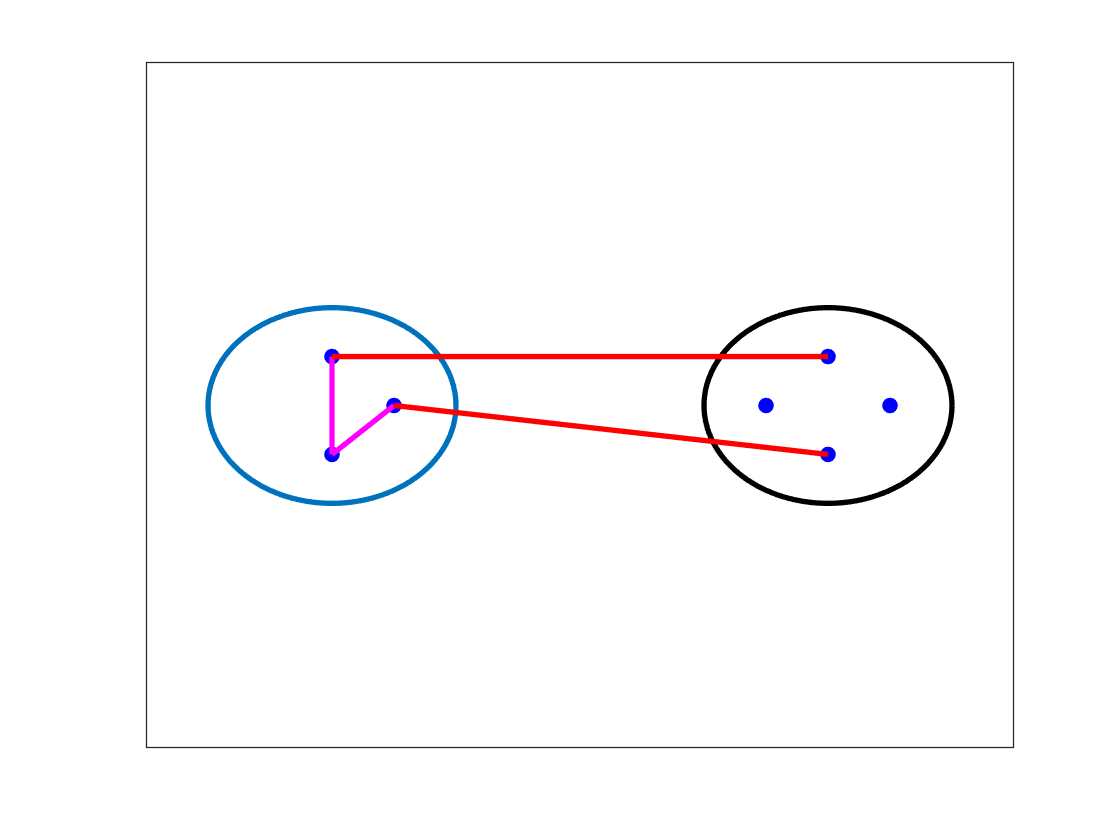}  
  \caption{$(\tau,\kappa)$ sparsity}
\end{subfigure}
\caption{Diagram of the correlation graph $\Gc_0(\tilde{\Sigmama})$ for a $p=7$ dimensional distribution with two different $7\times 7$ correlation matrices. The left panel is associated with a block-3 sparse assumption on $\tilde{\Sigmama}$. Only the $\qq = 3$ variables in the group inside the left circle are correlated: there are no correlations (edges) between the remaining $4$ variables in the right circle and there are no correlations across the two sets of variables in the different circles. The right panel is associated with $(\qq,\kappa)=(3,3)$ sparsity on $\tilde{\Sigmama}$, where two additional edges, representing correlations between variables, exist across the two groups.}
    \label{fig:covasscom}
\end{figure}

It is clear that the $(\tau,\kappa)$ sparsity condition is more general than the block sparsity condition as that there is no restriction on $\Sigmama_{12}=0$. 
Indeed, every block-$\qq$ sparse matrix is $(\tau,\kappa)$ sparse with $\kappa = \tau$. Nevertheless, the $(\qq,\kappa)$ sparsity with $\kappa=\qq$ allows non-zeros in the top-right submatrix, which permits more possible correlations between the variables relative to the block-$\qq$ sparsity in correlation graphical models. To see this, consider the associated graphical model $\Gc_0(\tilde{\Sigmama})$ for a correlation matrix $\tilde{\Sigmama}$. Recall in Section \ref{Hubfra} we define $\Gc_\rho(\cdot)$ as the graph with adjacency matrix obtained by thresholding a matrix with $\rho$. 
In Figure \ref{fig:covasscom}, nodes represent the variables and edges represent the correlation between variables. The left panel is a graphical model associated with the block-$3$ sparse assumption, while the right panel satisfies $(\qq,\kappa)$ sparsity with $(\qq,\kappa)=(3,3)$. The later has more correlations (the red edges) across the two sets of variables in the $2$ circles. The $(\qq,\kappa)$ sparsity condition with $\kappa>\qq$ allows additional correlations between variables. 
{ Notably, $(\tau,\kappa)$ sparsity allows the underlying graphical model to be connected as shown in the following example.} 

\begin{exa}[A connected correlation graph that is $(\tau,\kappa)$ sparse]
\label{exa:connectedtaukappa}
    Start with a set of $p$ vertices $[p]$ with no edges. Denote the neighborhood of vertex $i$ to be $\operatorname{NB}(i)$. We construct a connected $(\tau,\kappa)$ sparse graph on $[p]$ by creating neighborhoods of the first $\tau$ vertices according to the following rule: 
    \begin{align*}
    \operatorname{NB}(i)= & \{ i+1,  \tau+ (i-1) (\kappa-2) +1, \tau+ (i-1) (\kappa-2) +2 ,\ldots, \tau+ i (\kappa-2)   \}, \quad \text{ for } i\in [\tau-1], \\
    \operatorname{NB}(\tau)  = &  \{\tau+ (\tau-1) (\kappa-2) +1, \tau+ (\tau-1) (\kappa-2) +2,\ldots, \tau+ \tau (\kappa-2) + 1 \}.
    \end{align*}
    With this construction, each vertex among the first $\tau$ vertices is connected to $\kappa-1$ other vertices and the remaining $p-\tau$ vertices $[p]\setminus [\tau]$ do not connect to each other. As a result the associated adjacency matrix is $(\tau,\kappa)$ sparse. Moreover, as long as 
    \begin{equation}
    \tau+ \tau (\kappa-2) + 1 = \tau (\kappa-1) + 1 \geq p,  \label{eqn:connectedtaukappa}
    \end{equation}
    the graph is connected. Now, let $\Sigmama$ be a matrix with diagonal $1$ and the locations of the off-diagonal nonzero elements specified by the graph constructed. For simplicity, choose all the off-diagonal nonzero elements to be  $\xi$. If $\xi$ is small enough, then $\Sigmama$ is close to the identity matrix, hence positive definite and is therefore a {\em bone fide} correlation matrix. We have thus constructed a $(\tau,\kappa)$ sparse correlation matrix whose associated graphical model is connected provided \eqref{eqn:connectedtaukappa} is true. Now observe that by specifying $\tau=p^{\alpha}$ and $\kappa=p^{\beta}$ with $\alpha+\beta>1$, \eqref{eqn:connectedtaukappa} is indeed satisfied for large $p$. Moreover, if $\alpha,\beta\in (0,1)$ then $\tau,\kappa=o(p)$. By choosing $\xi$ as a decreasing function of $p$,  the condition number of $\Sigmama$ will not increase too fast w.r.t. $p$ so that the condition that $\lim\limits_{p\to \infty} \frac{\tau}{p}+\mu_{n,2\delta+2}\left(\Sigmama\right)\frac{\kkkk}{p}  \to 0$ in Theorem \ref{cor:Poissonlimit} is satisfied.  
\end{exa}

On the other hand, $(\tau,\kappa)$ sparsity is a stronger assumption than row-$\kappa$ sparsity
, since every $(\tau,\kappa)$ sparse matrix is row-$\kappa$ sparse. $(\tau,\kappa)$ sparsity is thus an intermediate level of sparsity lying between block sparsity and row sparsity.

\begin{rem} \label{rem:firsttaurworemark}
    In Theorem \ref{cor:Poissonlimit} we supposed that the dispersion parameter $\Sigmama$, after some row-column permutation, is $(\qq,\kappa)$ sparse. In this paper we are interested in $\bar{N}_\delta$, which is invariant under permutation of the $p$  variables. Since permutation of the variables is equivalent to row-column permutation of the dispersion parameter $\Sigmama$, without loss of generality, we can assume that the variables have been permuted such that $\Sigmama$ is $(\qq,\kappa)$ sparse. 
\end{rem}

\subsection{Local normalized determinant}
\label{sec:upperboundjointdensity}

In this subsection, we elaborate on the normalized determinant $\mu(\cdot)$ and the inverse local normalized determinant defined in Subsection \ref{sec:unifiedtheorem}. 
\ONE{The normalized determinant measures the closeness between the identity matrix and the dispersion parameter $\Sigmama$. This quantity plays an important role in Theorem \ref{cor:Poissonlimit}.}

 Observe that $\mu(\Ama)\in (0,1]$ and that, $\mu(\Ama)=1$  if and only if $\Ama$ is a multiple of $\Ima_p$. Moreover $\mu(\Ama)$ is close to $1$, and hence bounded away from $0$, as long as all eigenvalues concentrate around a positive number. Below is an example of a sequence of symmetric positive definite matrices with well-concentrated eigenvalues such that their normalized determinants are uniformly bounded away from $0$. 

\begin{exa} \label{exa:normaliseddeterminant}
    Let $\{\alpha_i\}_{i=1}^p$ be a positive real sequence. Let $\{\beta_i\}_{i=1}^\infty $ be a positive, decreasing sequence such that $\sum\limits_{i=1}^\infty \beta_i < \infty$. Consider $\Ama \in \R^{p\times p}$ a symmetric positive definite matrix with eigenvalues $\lambda_i = \alpha_p \exp(-\beta_i)$ for $1\leq i\leq p-1$ and $\lambda_p = \alpha_p$. Then 
    $$
    \mu(\Ama) =  \exp\left(-\sum\limits_{i=1}^{p-1} \beta_i\right).
    $$
    Consider now the case that $p$ is increasing and consider a sequence of matrices $\Ama$ of increasing dimension with the above properties. For this sequence $\mu(\Ama) \geq \exp\left(-\sum\limits_{i=1}^{\infty} \beta_i\right)>0$, i.e. $\mu(\Ama)$ is bounded uniformly away from $0$. 
\end{exa}

It follows by the interlacing property (see Theorem 8.1.7 in \cite{golub2012matrix}) that $\mu_m(\Ama)$ is decreasing with respect to $m\in [p]$ for any symmetric positive definite matrix $\Ama\in \R^{p\times p}$. Thus the inverse local normalized determinant $\mu_{n,m}(\Ama)$ is increasing with respect to $m \in [p]$.

It turns out that the inverse local normalized determinant of the dispersion matrix $\Sigmama$ will play an important role in our study of the distribution of $\bar{N}_\delta$. Indeed, when $\delta\geq 2$, $\Nedged^{(\Rma)}$ can be represented as a sum of indicator functions of the event that a subgraph of $\delta+1$ vertices is isomorphic to $\Gamma_\delta$ (see \eqref{eqn:Nedgedef1} for the precise formula). Each term in the summation involves only $\delta+1$ variables, and thus each pair of two such terms involves at most $2(\delta+1)$ variables. 
Hence $\mu_{2\delta+2}(\Sigmama)$ or $\mu_{n,2\delta+2}(\Sigmama)$ controls the correlation between two indicator terms in the summation of $\Nedged^{(\Rma)}$, which determines the convergence of $\Nedged^{(\Rma)}$ to a compound Poisson distribution. 

In Theorem \ref{cor:Poissonlimit} we assume that $\mu_{n,2\delta+2}\left(\Sigmama\right)\frac{\kkkk_p}{p} \to 0$, which holds when $\mu_{n,2\delta+2}\left(\Sigmama\right)$ is either  bounded or increasing with rate $o(\frac{p}{\kappa_p})$. In the Section \ref{sec:localdeterminantvseigenvalues} of the Appendix,  Lemma \ref{lem:sufass} provides a bound on $\mu_{n,2\delta+2}\left(\Sigmama\right)$ in terms of the condition number and the eigenvalues of $\Sigmama$. 


\section{Compound Poisson characterizations for finite $p$}
\label{mainres}

In this section we establish a compound Poisson approximation for $\mathscr{L}(\bar{N}_\delta)$ for finite $p$, which is then used to prove the asymptotic Theorem \ref{cor:Poissonlimit}. In Subsection \ref{sec:scorerep} Z-score type representations are introduced, which are used in the subsequent development in Subsection \ref{sec:rangeo}. Subsection \ref{sec:rangeo} provides an equivalent formulation of the empirical correlation and partial correlation graphs in terms of random geometric graphs. 
Against the backdrop of the first two subsections, Subsection \ref{sec:closenessnumedge} presents the compound Poisson approximation for star subgraph counts $N_{E_\delta}^{(\Rma)}$. Subsection \ref{sec:allclose} demonstrates that all $6$ quantities in $\{N_{E_\delta}^{(\Psima)},N_{\breve{V}_\delta}^{(\Psima)},N_{V_\delta}^{(\Psima)}: k\in \{\Rma,\Pma\}  \}$ are close in $L^1$ distance. Combining results in Subsection \ref{sec:closenessnumedge} and Subsection \ref{sec:allclose}, a compound Poisson characterization for $\bar{N}_\delta$ for finite $p$ is obtained in Subsection \ref{sec:compoihighdim}, which is then used to deduce Theorem \ref{cor:Poissonlimit} in Subsection \ref{sec:proofoftheorem1}.





\subsection{Z-score type representations of sample correlation and partial correlation}
\label{sec:scorerep}

In this subsection, we define the U-scores and the Y-scores for the sample correlation and the sample partial correlation, respectively. These Z-type scores will serve as the vertex set on which random geometric graphs are constructed in Subsection \ref{sec:rangeo}. 

The matrix of Z-scores $\Zma=[\zve_1,\ldots,\zve_p]\in \R^{n\times p}$ associated with the data matrix $\Xma$ is defined by 
$$
\zve_i := \frac{\xve_i-\bm{1}\bar{x}_i}{\sqrt{S_{ii}(n-1)}},
$$
where $\bar{x}_i$ is the $i$-th coordinate of the sample mean $\bar{\xve}$. The Z-scores specify the sample correlation matrix by the relation $\Rma= \Zma^\top \Zma$. Since $\bm{1}^\top \zve_i =0$ for every $i\in [p]$, we consider an equivalent but lower-dimensional type of Z-scores called the U-scores \cite{hero2011large, hero2012hub}, obtained by projecting each $\zve_i$ on the subspace $\{\wwve\in \R^n: \bm{1}^\top \wwve =0 \}$. Specifically, consider
an orthogonal $\nn \times \nn$ matrix $\Hma = [\nn^{-\frac{1}{2}}\bm{1}, \Hma_{2:\nn}]^\top$. The matrix $\Hma_{2:\nn}\in \R^{n\times(n-1)}$ can be obtained by
Gramm-Schmidt orthogonalization and satisfies the properties
$$
\bm{1}^\top\Hma_{2:\nn} := \bm{0}, \quad \Hma_{2:\nn}^\top\Hma_{2:\nn} = \Ima_{\nn - 1}.
$$
Define the matrix of U-scores as $\Uma=[\uve_1,\ldots,\uve_p]$ with $\uve_i=\Hma_{2:\nn}^\top \zve_i\in \R^{n-1}$. It is clear that 
\begin{equation}
\label{eqn:RUrelation}
\Rma=\Uma^\top \Uma.
\end{equation} 
The normalized outer product of $\Uma$, defined by 
\begin{equation}
\Bma :=\frac{\nn-1}{p}\Uma\Uma^\top\in \R^{(\nn-1)\times (\nn-1)} \label{eqn:Bdef}
\end{equation}
 will play an important role in the analysis of the empirical partial correlation graph. 
\begin{lem}
\label{lem:Binv}
Let $\Xma \sim \mathcal{VE}(\muve,\Sigmama,g)$ and $p\geq n$. Then $\Bma$ is invertible a.s.
\end{lem}

Proofs of all the lemmas and propositions in this paper can be found in the Appendix. By Lemma 1 in \cite{hero2012hub}, provided $\Bma$ is invertible, the Moore-Penrose pseudo-inverse of $\Rma$ can be expressed as follows:
\begin{equation}
\label{Rdagdef}
\Rma^\dagger = \Uma^\top[\Uma\Uma^\top]^{-2}\Uma=\left(\frac{p}{\nn-1}\right)^2\Uma^\top\Bma^{-2}\Uma.
\end{equation}
 It follows from Lemma \ref{lem:Binv} that equation \eqref{Rdagdef} holds $a.s.$ Define 
$\bar{\Yma}=\Bma^{-1}\Uma$ and observe that
\begin{equation}
\label{RdagYbar}
\Rma^\dagger = \left(\frac{p}{\nn-1}\right)^2 \bar{\Yma}^\top\bar{\Yma} \quad a.s. 
\end{equation}
Further define $\yve_i := \bar{\yve}_i/\|\bar{\yve}_i\|_2$ for $i\in [p]$ and $\Yma  := [\yve_1,\yve_2,\ldots,\yve_p]\in \R^{(\nn-1)\times p}$. Then
\begin{equation}
\Pma = \Yma^\top\Yma \quad a.s.  \label{Pfac}
\end{equation}
by equation \eqref{Pdef} and \eqref{RdagYbar}. $\Yma$ in \eqref{Pfac} is referred to as the Y-scores representation of the sample partial correlation matrix, as in \cite{hero2012hub}. 

By \eqref{eqn:RUrelation} and \eqref{Pfac}, the set of U-scores and Y-scores summarize all the information about the sample correlations and partial correlations, and hence they capture all the information about the vertex counts $\bar{N}_\delta$. We will further elaborate on this point in Subsection \ref{sec:rangeo}. Observing that $\Yma$ is a function of $\Uma$, the distribution of $\bar{N}_\delta$ is uniquely determined by $\Uma$. The remainder of this subsection focuses on the distribution of the U-scores.

Denote $S^{n-2}$ the unit sphere in $\R^{n-1}$ and denote $\unif(S^{n-2})$ the uniform distribution on $S^{n-2}$. Let $\{\tilde{\xve}^{(i)}\}_{i=1}^{n-1}\overset{\text{i.i.d.}}{\sim} \Nc(\bm{0},\Sigmama)$ and 
$$
\tilde{\Xma}=[\tilde{\xve}^{(1)},\ldots,\tilde{\xve}^{(n-1)}]^\top =[\tilde{\xve}_1,\ldots,\tilde{\xve}_p] \in \R^{(n-1)\times p}.
$$
In the next lemma we characterize the distribution of $\Uma$ in terms of the distribution of $\tilde{\Xma}$. 

\begin{lem} \label{lem:Uscoresdistribution}
Assume $\Xma \sim \mathcal{VE}(\bm{\mu},\Sigmama,g)$. 
    \begin{enumerate}[label=(\alph*)]
        \item \label{item:Uscoresdistributiona} The matrix $\Uma$ of U-scores has the same distribution as
        $
        \left[ \frac{\tilde{\xve}_1}{\|\tilde{\xve}_1\|_2}, \ldots, \frac{\tilde{\xve}_p}{\|\tilde{\xve}_p\|_2} \right] $. 
        
        \item \label{item:Uscoresdistributionb}
        For each $i\in [p]$, $\uve_i$ is distributed as $\unif(S^{\nn-2})$. Moreover, if $\Sigma_{ij}=0$, then $\uve_i$ and $\uve_j$ are independent.
    \end{enumerate}
\end{lem}

\begin{rem}  \label{rem:Udisconsequence} 
When $\Sigma_{ij}=0$, the $i$-th and the $j$-th columns of the data matrix are uncorrelated but are not necessarily independent. However Lemma \ref{lem:Uscoresdistribution} \ref{item:Uscoresdistributionb} establishes that in fact the corresponding U-scores $\uve_i$ and $\uve_j$ are independent. This is a consequence of the fact that the entire data matrix follows the vector elliptically contoured distribution $\Xma \sim \mathcal{VE}(\bm{\mu},\Sigmama,g)$. Another important implication of Lemma \ref{lem:Uscoresdistribution} is that the distribution of $\Uma$, hence the distribution of $\bar{N}_\delta$, is invariant to the mean $\muve$ and the shaping function $g$. 
\end{rem}


Let $\sigma^{n-2}$ be the spherical measure on $S^{n-2}$, that is $\sigma^{n-2}$ is the uniform probability measure corresponding to $\unif(S^{n-2})$.  
 Denote by $f_{\uve_{j_1},\uve_{j_2},\cdots,\uve_{j_m}}$ the joint density of the $j_1$-th, $j_2$-th, $\ldots$, $j_m$-th column of $\Uma$ with respect to the product measure $\otimes^m \sigma^{n-2} := \underbrace{\sigma^{n-2}\otimes \sigma^{n-2} \otimes \cdots \otimes \sigma^{n-2}}_{m}$. The next lemma establishes that $f_{\uve_{j_1},\uve_{j_2},\cdots,\uve_{j_m}}$ is bounded by the inverse local normalized determinant $\mu_{n,m}(\Sigmama)$. 

\begin{lem} \label{lem:densitybou}
    Assume $\Xma\sim \mathcal{VE}(\bm{\mu}, \Sigmama, g)$. Consider a set of distinct indexes $\mathcal{J}=\{j_i: 1\leq i \leq m\}\subset [p]$.
    \begin{enumerate}[label=(\alph*)]        
        \item \label{item:densityboua}
        The inverse local normalized determinant satisfies
        $$
        \mu_{n,m}(\Sigmama_{\mathcal{J}}) \leq \mu_{n,m} (\Sigmama).
        $$
        
        \item \label{item:densityboub}  
        The joint density of any subset of $m$ columns of U-scores w.r.t. $\otimes^m \sigma^{n-1}$ is upper bounded by $\mu_{n,m} (\Sigmama)$:
        $$f_{\uve_{j_1},\uve_{j_2},\cdots,\uve_{j_m}}(\vve_1,\vve_2,\cdots,\vve_m)  \leq \mu_{n,m}(\Sigmama_{\mathcal{J}}) \leq \mu_{n,m} (\Sigmama),  \quad \forall \  \vve_i\in S^{n-2}, \forall\ i\in [m]. $$
        \item \label{item:densitybouc}
        Let $h:\left(S^{n-2}\right)^m\to \R$ be a nonnegative Borel measurable function. Then 
        \begin{align*}
        \E h(\uve_{j_1},\uve_{j_2},\cdots,\uve_{j_m})
        \leq & \mu_{n,m} (\Sigmama_{\mathcal{J}}) \E h(\uve_{j_1}',\uve_{j_2}',\cdots,\uve_{j_m}') \\
        \leq & \mu_{n,m} (\Sigmama) \E h(\uve_{j_1}',\uve_{j_2}',\cdots,\uve_{j_m}'),    
        \end{align*}
        where $\{\uve_{j_\ell}'\}_{\ell=1}^m$ are i.i.d. uniformly distributed on $S^{n-2}$.
    \end{enumerate}
\end{lem}

According to Lemma \ref{lem:densitybou} \ref{item:densitybouc}, when calculating the expectation of nonnegative function of any $m$ columns of $\Uma$, one may always assume that the associated columns $\{\uve_j\}$ are i.i.d. $\unif(S^{n-1})$ up to an additional multiplicative factor $\mu_{n,m}(\Sigmama)$.

\subsection{Random pseudo geometric graphs}
\label{sec:rangeo}

In this subsection we define random pseudo geometric graphs and provide an equivalent formulation of the empirical correlation and partial correlation graphs in terms of these graphs, for which the vertex sets are, respectively, the U and Y scores. We also define the increment distribution of the compound Poisson that approximates $\bar{N}_\delta$ when $p$ is finite.

By the fact that $\Rma = \Uma^\top \Uma$ and the fact that columns of $\Uma$ have Euclidean norm $1$,
\begin{equation}
R_{ij}=\uve_i^\top\uve_j=1- \frac{\|\uve_i-\uve_j\|_2^2}{2} =  \frac{\|\uve_i+\uve_j\|_2^2}{2} -1. \label{eqn:cortoycol}
\end{equation}
For a threshold $\rho\in [0,1)$, define $r_\rho:=\sqrt{2(1-\rho)} \in (0,\sqrt{2}]$. 
As shown in \cite{hero2011large,hero2012hub}, by \eqref{eqn:cortoycol}, 
\begin{equation}
 \{ |R_{ij}| \geq \rho \} = \{\|\uve_i+\uve_j\|_2 \leq r_\rho\}\cup \{\|\uve_i-\uve_j\|_2 \leq r_\rho\}.   
\label{eqn:corphiequzer}
\end{equation}
An analogous argument yields the following for the empirical partial correlation graph,
\begin{equation}
 \{ |P_{ij}| \geq \rho \} = \{\|\yve_i+\yve_j\|_2 \leq r_\rho\}\cup \{\|\yve_i-\yve_j\|_2 \leq r_\rho\}.   \label{phiequzer}
\end{equation}
Based on \eqref{eqn:corphiequzer}, we now introduce novel geometric connections between empirical correlation  graphs and random geometric graphs. Recall $\{|R_{ij}|\geq \rho\}$ is the event that the magnitude of the sample correlation between the $i$-th and $j$-th variables exceeds the threshold $\rho$, or equivalently, the event that there exists a edge connecting the $i$-th and $j$-th vertices in the empirical correlation graph $\Gc_{\rho}(\Rma)$. Equation \eqref{eqn:corphiequzer} indicates that $\{|R_{ij}|\geq \rho\}$ is the same as the event that the associated U-scores for the $i$-th and $j$-th variables, $\uve_i$ and $\uve_j$, lie in some geometric set on $S^{\nn-2}\times S^{\nn-2}$. This insight provides an equivalent way to construct $\Gc_{\rho}(\Rma)$ through the U-scores. The empirical partial correlation graph $\Gc_{\rho}(\Pma)$ may similarly be constructed through the Y-scores based on \eqref{phiequzer}. These remarks are formalized in the next few paragraphs.

\begin{definition}[Pseudo geometric graph\index{pseudo geometric graph}]
    \label{def:pge}
    Given $m\geq 2$ and a set of points $\{ \vve_i  \}_{i=1}^m$ in $\R^{\mathscr{N}}$, denote by $\textbf{PGe}\left(\{\vve_i  \}_{i=1}^m, r  ;  \mathscr{N} \right)$ the pseudo geometric graph with radius $r$, defined as follows. The vertex set of the graph is $\{ \vve_i  \}_{i=1}^m$, and there is an edge between $\vve_i$ and $\vve_j$ if $\textbf{dist}(\vve_i,\vve_j):= \min\left\{\|\vve_i-\vve_j\|_2,\|\vve_i+\vve_j\|_2 \right \}  \leq r$. The graph is called a random pseudo geometric graph when the vertices are random.
\end{definition}
It is easy to verify that $\textbf{dist}(\cdot,\cdot)$ has the following properties: for $\forall \vve_1,\vve_2,\vve_3\in \R^m$,
\begin{enumerate}
\item     $\textbf{dist}(\vve_1,\vve_2)\geq 0$;
\item   $\textbf{dist}(\vve_1,\vve_2) = 0$ if only if $\vve_1 = \vve_2$ or $\vve_1 = - \vve_2$;
\item   $\textbf{dist}(\vve_1,\vve_2) = \textbf{dist}(\vve_2,\vve_1)$ and $\textbf{dist}(\vve_1,\vve_2) = \textbf{dist}(\vve_1,-\vve_2)$
\item $\textbf{dist}(\vve_1,\vve_2)\leq \textbf{dist}(\vve_1,\vve_3) + \textbf{dist}(\vve_3,\vve_2)$.
\end{enumerate}
That is, $\textbf{dist}(\cdot,\cdot)$ is a pseudo metric on $\R^\mathscr{N}$, hence the name pseudo geometric graph in Definition \ref{def:pge}. 


 With the above definitions, and by the discussions preceding Definition \ref{def:pge}, the empirical correlation graph $\Gc_{\rho}(\Rma)$ is isomorphic to $\textbf{PGe}\left(\{\uve_i  \}_{i=1}^p, r_\rho   \right)$, the random pseudo geometric graph generated by U-scores. Consequently $\{N_{V_\delta}^{(\Rma)},N_{\breve{V}_\delta}^{(\Rma)},N_{E_\delta}^{(\Rma)}\}$ are the numbers of vertices or subgraphs in $\textbf{PGe}\left(\{\uve_i  \}_{i=1}^p, r_\rho   \right)$. 
 An analogous analysis applies to the empirical partial correlation graph and  $\textbf{PGe}\left(\{\yve_i  \}_{i=1}^p, r_\rho   \right)$. This equivalent construction indicates that the distribution of each of the 3 quantities $\{N_{E_\delta}^{(\Psima)},N_{\breve{V}_\delta}^{(\Psima)},N_{V_\delta}^{(\Psima)}  \}$ with $\Psima=\Rma$ ($\Psima=\Pma$) only depends on the pairwise pseudo distances $\textbf{dist}(\cdot,\cdot)$ between columns of $\Uma$ ($\Yma$). 


The random pseudo geometric graph in Definition \ref{def:pge} \ONE{is similar} to the random geometric graph introduced in \cite{penrose2003random}. 
In particular, studied in the monograph \cite{penrose2003random} is the number of \emph{induced} subgraphs isomorphic to a given graph, typical vertex degrees, and other graphical quantities of random geometric graphs. 
\ONE{As a specific example, the rate $a_n2^{n/2}p^{1+\frac{1}{\delta}} (1-\rho)^{\frac{n-2}{2}}\to e(n, \delta)$ as $p\to \infty$ in Theorem \ref{cor:Poissonlimit} is equivalent to $2a_np^{1+\frac{1}{\delta}} r_\rho^{{n-2}}\to e(n, \delta)$, which is consistent with existing Poisson approximation for random geometric graph \cite[ Theorem 3.4]{penrose2003random}.
}
The differences between our random pseudo geometric graph in Definition \ref{def:pge} and the random geometric graphs defined in \cite{penrose2003random} are: 1) our graph has vertices $\uve_i$ lying on the unit sphere instead of on the entire Euclidean space; 2) our graph is induced by distance $\textbf{dist}(\cdot,\cdot)$ instead of the Euclidean distance. Another key difference is that vertices $\uve_i$ are not necessarily independent in Definition \ref{def:pge}. Indeed in our model, the correlations between vertices $\uve_i$ are encoded by a sparse matrix $\Sigmama$ (cf. Lemma \ref{lem:Uscoresdistribution} \ref{item:Uscoresdistributiona}), whereas in \cite{penrose2003random} the vertices associated with the random geometric graph are assumed to be i.i.d. In \cite{penrose2003random} it was \ONE{stated (Example after \cite[Corollary 3.6]{penrose2003random}) without proof} that the number of vertices with degree at least $3$ was approximately compound Poisson. A similar compound Poisson limit is established in Theorem \ref{cor:Poissonlimit}. \ONE{There is recent work on testing whether a given graph is a realization of 
an Erd\H{o}s–R\'enyi random graph or a realization of 
a random geometric graph with vertices i.i.d. uniformly distributed on the sphere 
\cite{bubeck2016testing}. 
{ Therein, the authors study the asymptotics of a count statistic, the (signed) triangles, on a random geometric graph induced by thresholding the distances between $p$ i.i.d. uniformly distributed points on the $n$ dimensional sphere as $n$ goes to $\infty$. In contrast, this paper studies the large $p$ asymptotics of vertex degree counts for dependent non-uniform points on a sphere of fixed dimension $n$.}}

Recall that $\nmd\left(\{\vve_i  \}_{i=1}^m, r  ;  \mathscr{N} \right)$ denotes the number of universal vertices (vertices adjacent to all vertices except itself) in the geometric graph $\textbf{Ge}\left(\{\vve_i  \}_{i=1}^m, r  ;  \mathscr{N} \right)$. 
Analogously, denote by $\pnmd\left(\{\vve_i  \}_{i=1}^m, r  ;  \mathscr{N} \right)$  the number of universal vertices in the pseudo geometric graph $\textbf{PGe}\left(\{\vve_i  \}_{i=1}^m, r  ;  \mathscr{N} \right)$.

Denote by $\deg(\cdot)$ the degree of a given vertex in the graph. Consider $\{\uve'_i\}_{i=1}^{\delta+1}\ \overset{\text{i.i.d.}}{\sim} \ \unif(S^{n-2})$. Denote the conditional probability that there are $\ell$ universal vertices in $\textbf{PGe}\left(\{\uve'_i  \}_{i=1}^{\delta+1}, r_\rho   \right)$ by
\begin{equation}
\alpha_{n,\delta}(\ell,r_\rho) := \P\left(\pnmd\left(\{\uve'_i  \}_{i=1}^{\delta+1}, r_\rho  \right)=\ell|\deg(\uve'_1)=\delta\right),\quad \forall \ell\in [\delta+1]. \label{eqn:alphandeltarho}
\end{equation} 
The conditional probability $\alpha_{n,\delta}(\ell,r_\rho)$ depends on $n$, $\delta$ and the threshold $\rho$ and is abbreviated as $\alpha(\ell,r_\rho)$ when there is no risk of confusion.  
Define a probability  distribution $\bm{\zeta}_{n,\delta,\rho}$ supported on $[\delta+1]$ with
\begin{equation}
\bm{\zeta}_{n,\delta,\rho}(\ell) : = \frac{\alpha(\ell,r_\rho)/\ell}{ \sum_{\ell=1}^{\delta+1} \left(\alpha(\ell,r_\rho)/\ell \right) }.  \label{eqn:zetandelrhodef}
\end{equation}
It will be shown in the next three subsections that the probability distribution $\bm{\zeta}_{n,\delta,\rho}$ is the increment distribution of the compound Poisson approximation to $\bar{N}_\delta$ for finite $p$. 

\subsection{Closeness of the star subgraph counts
to a compound Poisson}
\label{sec:closenessnumedge}
\ONE{In this subsection we present a proposition that establishes an upper bound on the total variation between $\mathscr{L}\left(N_{E_\delta}^{(\Rma)}\right)$ and \ONE{the} compound Poisson distribution.

We first introduce some notation.} Let $\SC(r,\wwve)$ be the spherical cap with radius $r$ at the center $\wwve\in S^{\nn-2}$. Formally,
\begin{equation}
\SC(r,\wwve) = \{\xve\in S^{\nn-2}: \|\xve-\wwve\|_2 \leq r\}. \label{eqn:SCdef}
\end{equation}
Define $P_\nn(r):=\frac{{\text{Area}}(\SC(r,\wwve))}{\text{Area}(S^{\nn-2})}$, where $\text{Area}(\cdot)$ is the area of a subset of $S^{\nn-2}$. $P_\nn(r)$ is the normalized area of the spherical cap with radius $r$. As is shown in (2.6) in \cite{hero2011large}, 
\begin{equation}
P_\nn(r) = \frac{b_\nn}{2} \int_{1-\frac{r^2}{2}}^1 (1-u^2)^{\frac{\nn-4}{2}}du,\quad \text{when }r\in [0,\sqrt{2}], \label{spharefor}
\end{equation}
where $b_\nn = \frac{2\Gamma((\nn-1)/2)}{\sqrt{\pi}\Gamma((\nn-2)/2)}$. It follows by simple calculation that  
$$
P_{\nn}(r)=1-\frac{{\text{Area}}(\SC(\sqrt{4-r^2},\wwve))}{\text{Area}(S^{\nn-2})}=1-P_\nn(\sqrt{4-r^2}) \quad \text{when} \sqrt{2}<r\leq 2,
$$ 
and $P_\nn(r)=1$ when $r>2$. Further properties of $P_n(r)$ are summarized in Lemma \ref{pderivative} in Section \ref{sec:auxlem} of the Appendix. 
\ONE{Define the rate parameter
\begin{align}
   \lambda_{p,n,\delta,\rho}:= \binom{p}{1}\binom{p-1}{\delta}(2P_n(r_\rho))^\delta\sum_{\ell=1}^{\delta+1} \frac{\alpha(\ell,r_\rho)}{\ell}.  \label{eqn:lambdadeffinitep}
\end{align}
The next proposition \ONE{establishes that the distribution of $\Nedged^{(\Rma)}$ is approximated by a compound Poisson distribution $\CP(\lambda_{p,n,\delta,\rho},\bm{\zeta}_{n,\delta,\rho})$ with rate parameter $\lambda_{p,n,\delta,\rho}$ and dispersion parameter $\bm{\zeta}_{n,\delta,\rho}$.}} Here $C$ and $c$, and their subscripted versions, denote generic positive constants that only depend on their subscripts. 


\begin{prop}[Compound Poisson approximation for subgraph counts]
\label{prop:edgecor}
    Let $p\geq n \geq 4$, $\delta\in [p-1]$ and $\gamma>0$ be given. Suppose $\Xma\sim \mathcal{VE}(\muve,\Sigmama,g)$. Suppose $2p^{1+\frac{1}{\delta}}P_n(r_\rho)\leq \gamma$, and $\Sigmama$ is row-$\kkkk$ sparse. Then
    \begin{align*}
    &d_{\TV}\left(\mathscr{L}\left(N_{E_\delta}^{(\Rma)}\right), \CP(\lambda_{p,n,\delta,\rho},\bm{\zeta}_{n,\delta,\rho})\right) \\
    \leq & C_{n,\delta,\gamma} \left(C'_{\delta,\gamma}\right)^{\mu_{n,\delta+1}\left(\Sigmama\right)\frac{\kappa-1}{p}}\left( \mu_{n,2\delta+2}\left(\Sigmama\right)\frac{\kappa}{p} \left(1 + \mu_{n,2\delta+2}\left(\Sigmama\right) \left(\frac{\kappa}{p}\right)^2 \right)  + p^{-\frac{1}{\delta}} \right).
    \end{align*}
    where $C_{n,\delta,\gamma}$ and $C'_{\delta,\gamma}$ are  two constants. 
\end{prop}
\begin{rem} \label{rem:edgecor}
    The condition $2p^{1+\frac{1}{\delta}}P_\nn(r_\rho)\leq \gamma$ specifies an implicit lower bound on the threshold $\rho$. To obtain an explicit lower bound, observe that the condition $2a_np^{1+\frac{1}{\delta}} \left(\sqrt{2(1-\rho)}\right)^{n-2} \leq \gamma $ is sufficient for  $2 p^{1+\frac{1}{\delta}}P_\nn(r_\rho)\leq \gamma$, by Lemma \ref {pderivative} \ref{item:pderivativea} in Section \ref{sec:auxlem} of the Appendix. Solving for $\rho$, we then obtain 
    \begin{equation}
    \label{eqn:rhononasymptotic}
    \rho \geq 1 - \frac{1}{2} \left(\frac{\gamma}{2a_n p^{1+\frac{1}{\delta}}}\right)^{\frac{2}{n-2}}.
    \end{equation}
    The condition \eqref{eqn:rhononasymptotic} is a non-asymptotic version of \eqref{eqn:rhopformula}. \ONE{But the requirement on $\rho$ from $2p^{1+\frac{1}{\delta}}P_\nn(r_\rho)\leq \gamma$ is weaker than \eqref{eqn:rhononasymptotic} since \ONE{the bound} in Lemma \ref {pderivative} \ref{item:pderivativea} is not tight when $p$ is finite.}
     For the upper bound in Proposition \ref{prop:edgecor} to be small, $p$ should be relatively large and, $\Sigmama$ should have relatively small $\mu_{n,2\delta+2}\left(\Sigmama\right)$, and $\Sigmama$ should be row-$\kkkk$ sparse with relative small sparsity level $\kkkk/p$. 
     In Lemma \ref{lem:sufass} in Section \ref{sec:localdeterminantvseigenvalues} of the Appendix, Lemma \ref{lem:sufass}) we bound $\mu_{n,2\delta+2}\left(\Sigmama\right)$ in terms of the condition number and eigenvalues of $\Sigmama$. In the special case when $\Sigmama$ is diagonal, $\mu_{n,2\delta+2}\left(\Sigmama\right) \kkkk/p =1/p$. 
    %
    In the case that
    $\mu_{n,2\delta+2}\left(\Sigmama\right) \kkkk/p$ is small, say $\mu_{n,2\delta+2}\left(\Sigmama\right) \kkkk/p<1$, in the upper bound the term $\mu_{n,2\delta+2}\left(\Sigmama\right) \left(\kkkk/p\right)^2$ can be dropped, resulting in an additional constant factor, since $\mu_{n,2\delta+2}\left(\Sigmama\right)\left(\kkkk/p\right)^2<1$. In other words the effective upper bound is $\mu_{n,2\delta+2}\left(\Sigmama\right)\kkkk/p+p^{-\frac{1}{\delta}}$, neglecting the coefficients depending on $n$, $\delta$ and $\gamma$. 
    Respective expressions for $C'_{\delta,\gamma}$ and $C_{n,\delta,\rho}$ in Proposition \ref{prop:edgecor} are presented in equations \eqref{eqn:C'delgamdef} and \eqref{eqn:Cndelgamdef} in the Appendix.
\end{rem}

Proposition \ref{prop:edgecor} states that for given $n$, $p$, $\delta$ and $\gamma$, if the threshold $\rho$ is properly chosen, and $\Sigmama$ is row-$\kappa$ sparse and has small $\mu_{n,2\delta+2}(\Sigmama)$, then the distribution of $N_{E_\delta}^{(\Rma)}$ can be well approximated by the compound Poisson $\CP(\lambda_{p,n,\delta,\rho},\bm{\zeta}_{n,\delta,\rho})$. We will call $\CP(\lambda_{p,n,\delta,\rho},\bm{\zeta}_{n,\delta,\rho})$ the \emph{non-asymptotic compound Poisson distribution.   } 
\ONE{\ONE{In the Appendix we provide an informal argument (Section \ref{sec:derivationprop:edgecor}) to motivate Proposition \ref{prop:edgecor}   in addition to a complete proof (Section \ref{sec:proofprop:edgecor})}. }

\subsection{A portmanteau result and bounds on pairwise total variations}
\label{sec:allclose}

In this subsection upper bounds for pairwise total variation distances and $L^1$ distances among  $\{N_{E_\delta}^{(\Psima)},N_{\breve{V}_\delta}^{(\Psima)},N_{V_\delta}^{(\Psima)}: k\in \{\Rma,\Pma\}  \}$ are obtained. 
\ONE{Some intuition is presented before stating the portmanteau result Proposition \ref{thm:allclose}.}

\begin{lem} 
    \label{lem:6quantitiesine}
Consider $\delta \in [p-2]$.
$$    
N_{E_\delta}^{(\Rma)} - (\delta+1)N_{E_{\delta+1}}^{(\Rma)}   \leq N_{\breve{V}_\delta}^{(\Rma)} \leq N_{V_\delta}^{(\Rma)} \leq N_{E_\delta}^{(\Rma)}, 
$$
$$    
N_{E_\delta}^{(\Pma)} - (\delta+1)N_{E_{\delta+1}}^{(\Pma)}   \leq N_{\breve{V}_\delta}^{(\Pma)} \leq N_{V_\delta}^{(\Pma)} \leq N_{E_\delta}^{(\Pma)}. 
$$
\end{lem}

It follows directly from Lemma \ref{lem:6quantitiesine} that for $\tilde{N}_\delta \in \left\{ N_{\breve{V}_\delta}^{(\Rma)}, N_{V_\delta}^{(\Rma)} \right\} $, 
\begin{equation}
\E \left|\tilde{N}_\delta - N_{E_\delta}^{(\Rma)} \right| \leq (\delta+1) \E N_{E_{\delta+1}}^{(\Rma)}. \label{eqn:corexaatupp}
\end{equation}
As a result, if $\E N_{E_{\delta+1}}^{(\Rma)}$ is small, then $N_{\breve{V}_\delta}^{(\Rma)}$ and $N_{V_\delta}^{(\Rma)}$ are close to $N_{E_\delta}^{(\Rma)}$ in $L^1$ norm.

To intuitively see why the vertex counts in the empirical partial correlation graph are close to those in the empirical correlation graph, consider large $p$ and suppose that $\Sigmama$ is diagonal. Then by Lemma \ref{lem:Uscoresdistribution} \ref{item:Uscoresdistributionb}, $\{\uve_i\}_{i=1}^p$ are i.i.d. $\unif(S^{n-2})$. According to the law of large numbers, for large $p$,
\begin{equation}
\Bma = \frac{n-1}{p}\sum_{i=1}^p\uve_i\uve_i^\top  \approx (\nn-1)\E\uve_i\uve_i^\top = \Ima_{n-1}, \label{eqn:Bapprox}
\end{equation}
which implies that
\(
\bar{\Yma} = \Bma^{-1} \Uma \approx \Uma
\),
so that $\|\bar{\yve}_i\|_2\approx \|\uve_i\|_2=1$, and thus 
\begin{equation}
\Yma = \left[\frac{\bar{\yve}_1}{\|\bar{\yve}_1\|_2},\ldots,\frac{\bar{\yve}_p}{\|\bar{\yve}_p\|_2}\right] \approx \bar{\Yma} \approx \Uma. \label{eqn:UYclose}
\end{equation}
Recall in Subsection \ref{sec:rangeo} that $N_{E_\delta}^{(\Rma)}$ and $N_{E_\delta}^{(\Pma)}$ are the numbers of subgraphs isomorphic to $\Gamma_\delta$ associated with the random pseudo  geometric graphs   induced respectively by $\Uma$ and $\Yma$. Hence by \eqref{eqn:UYclose}, $N_{E_\delta}^{(\Rma)} \approx N_{E_\delta}^{(\Pma)}$. 
Following the same reasoning, $N_{V_\delta}^{(\Rma)}\approx N_{V_\delta}^{(\Pma)}$ and \(N_{\breve{V}_\delta}^{(\Rma)}\approx N_{\breve{V}_\delta}^{(\Pma)}\). These informal arguments are formalized by Proposition \ref{thm:allclose}, which establishes that all $6$ quantities $\{N_{E_\delta}^{(\Psima)},N_{\breve{V}_\delta}^{(\Psima)},N_{V_\delta}^{(\Psima)}: k\in \{\Rma,\Pma\}  \}$ are mutually close in $L^1$ norm and  their distributions are mutually close in total variation. Proposition \ref{thm:allclose} is therefore called a ``Portmanteau result''.

\begin{prop}[Portmanteau result] \label{thm:allclose}
    { 
        Let $p\geq n \geq 4$ and $\Xma\sim \mathcal{VE}(\bm{\mu},\Sigmama,g)$. Let $\delta \in [p-1]$. Suppose $2p^{1+\frac{1}{\delta}}P_\nn(r_\rho)\leq \gamma$. 
    \begin{enumerate}[label=(\alph*)]
        \item  \label{item:allclosea}
        Suppose $\Sigmama$ is row-$\kappa$ sparse. Then for $\tilde{N}_\delta \in \left\{ N_{\breve{V}_\delta}^{(\Rma)}, N_{V_\delta}^{(\Rma)} \right\} $,
        $$
        d_{\TV}\left(\mathscr{L}(\tilde{N}_\delta ), \mathscr{L}\left(\Nedged^{(\Rma)}\right)  \right) 
        \leq 
        \E \left| \tilde{N}_\delta -    \Nedged^{(\Rma)}  \right| \leq \frac{(\delta+1)^2}{\delta!}\gamma^{\delta+1}\left(1+\mu_{n,\delta+2}(\Sigmama)\frac{\kappa-1}{p}\right)p^{-\frac{1}{\delta}}.$$
        \item \label{item:allclosec}
        Suppose $\Sigmama$, after some row-column permutation, is $(\tau,\kappa)$ sparse with $\tau\leq \frac{p}{2}$. Suppose the condition $\left(\sqrt{\frac{\nn-1}{p}}+\sqrt{\frac{\delta\ln p}{p}}\right) \leq c$ is satisfied for some positive and \ONE{sufficiently} small  constant $c$. Then
        \begin{align*}
        d_{\TV}\left(\mathscr{L}\left(\Nedged^{(\Pma)}\right),\mathscr{L}\left(\Nedged^{(\Rma)}\right)  \right) \leq & \E\left|\Nedged^{(\Pma)}-\Nedged^{(\Rma)}\right|\\
        \leq & C_{E_\delta}^{(\Pma)} \left(1+\frac{\kappa-1}{p}\mu_{n,\delta+1}(\Sigmama)\right)\left(\sqrt{\frac{\ln p}{p}}+\frac{\tau}{p}\right), \numberthis \label{eqn:allclosec}
        \end{align*}
        where $C^{(\Pma)}_{E_\delta} $ is a constant depending on only $n,\delta$ and $\gamma$.
        \item  \label{item:allclosed}
        Suppose the same conditions as in part \ref{item:allclosec} hold. Then
        \begin{align*}
        d_{\TV}\left(\mathscr{L}\left(N_{\VVV_\delta}^{(\Pma)}\right),\mathscr{L}\left(N_{\VVV_\delta}^{(\Rma)}\right)  \right)\leq & \E\left|N_{\VVV_\delta}^{(\Pma)}-N_{\VVV_\delta}^{(\Rma)}\right|\\
        \leq &  C_{\breve{V}_\delta}^{(\Pma)} \left(1+\frac{\kappa-1}{p}\mu_{n,\delta+2}(\Sigmama)\right)\left(\sqrt{\frac{\ln p }{p}}+\frac{\tau}{p}+p^{-\frac{1}{\delta}} \right), \numberthis \label{eqn:allclosed}
         \end{align*}
        where $C_{\breve{V}_\delta}^{(\Pma)}$ is a constant depending  on only $n,\delta$ and $\gamma$.
    \end{enumerate} 
}
\end{prop}
\begin{rem}
\label{rem:allclose}
    In Proposition \ref{thm:allclose} \ref{item:allclosec} and \ref{item:allclosed}, the condition 
    \begin{equation}
    \left(\sqrt{\frac{\nn-1}{p}}+\sqrt{\frac{\delta\ln p}{p}}\right) \leq c \label{eqn:plargendelta}
    \end{equation}
 explicitly specifies a lower bound on $p$. Specifically, observe that the left side of \eqref{eqn:plargendelta} is a decreasing function of $p$, and its limit is $0$ when $p\to\infty$. Thus the smallest positive integer value of $p$ satisfying \eqref{eqn:plargendelta} exists, denoted by $p_0$, which is a function \ONE{of $n, \delta$ and $c$}. \ONE{Therefore} \eqref{eqn:plargendelta} is equivalent to requiring $p\geq p_0$. 
\ONE{The inequalities (\ref{eqn:allclosec}) and (\ref{eqn:allclosed}) are valid as long as $c$ is less than finite constant, given respectively in \eqref{eqn:cdef} and \eqref{eqn:cuppbou2} in the Appendix.  }
 %
   %
    The row-$\kappa$ sparsity condition on $\Sigmama$ guarantees that the vertex counts associated with the empirical correlation graph are close in $L^1$ norm to $N_{E_\delta}^{(\Rma)}$ by Proposition \ref{thm:allclose} \ref{item:allclosea}. In Proposition \ref{thm:allclose} \ref{item:allclosec} and \ref{item:allclosed} the stronger condition of $(\tau,\kappa)$ sparsity suffices to establish that the vertex counts of $G_\rho(\Pma)$ are close in $L^1$ norm to those of $G_\rho(\Rma)$. 
    All $3$ upper bounds in the proposition are bounded by  $ \left(1+\frac{\kappa-1}{p}\mu_{n,\delta+2}(\Sigmama)\right)\left(\sqrt{\frac{\ln p }{p}}+\frac{\tau}{p}+p^{-\frac{1}{\delta}} \right)$, up to a multiplicative constant depending on only $n$, $\delta$ and $\gamma$. The interpretation of this upper bound is similar to the second paragraph in Remark \ref{rem:edgecor}.
    Respective expressions for $C^{(\Pma)}_{E_\delta}$ and $C_{\breve{V}_\delta}^{(\Pma)}$ are presented in equations \eqref{eqn:CEatleadel} and \eqref{eqn:CVexcdel} in Appendix. 
\end{rem}

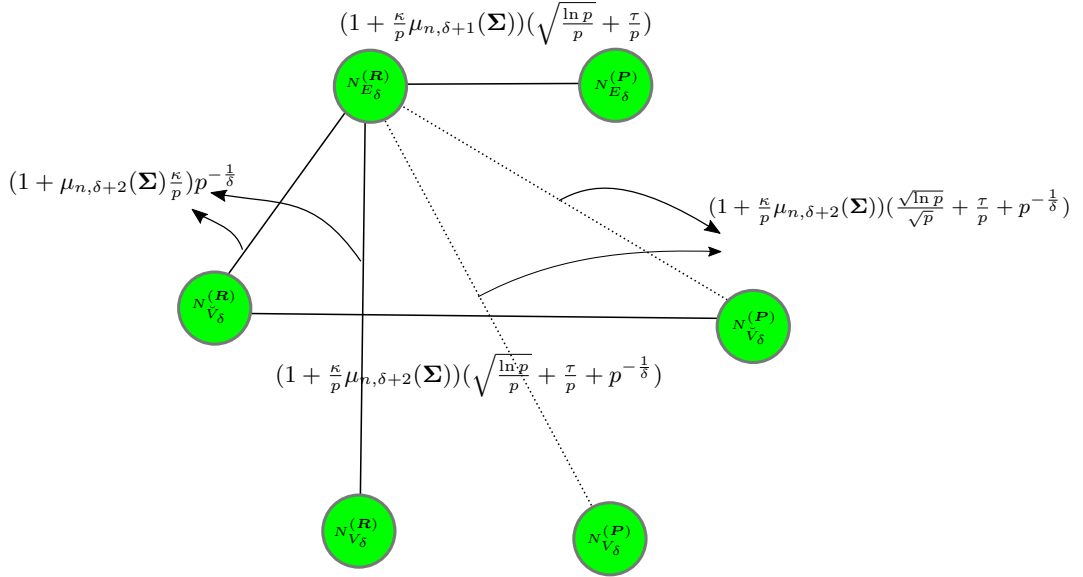
\begin{figure}%
    \labellist
    \small \hair 2pt
    \pinlabel {\tiny $N_{\breve{V}_\delta}^{(\Rma)}$} at 83 602
    \pinlabel {\tiny $N_{V_\delta}^{(\Rma)}$} at 400 98
    \pinlabel {\tiny $N_{E_\delta}^{(\Rma)}$} at 427 1098
    \pinlabel {\tiny $N_{E_\delta}^{(\Pma)}$} at 980 1094
    \pinlabel {\tiny $N_{\breve{V}_\delta}^{(\Pma)}$} at 1287 555
    \pinlabel {\tiny $N_{V_\delta}^{(\Pma)}$} at 965 81
    \pinlabel \rotatebox{0}{$(1+\frac{\kappa}{p}\mu_{n,\delta+1}(\Sigmama))(\sqrt{\frac{\ln p}{p}}+\frac{\tau}{p})$} at 711 1235
    \pinlabel \rotatebox{0 }{$(1+\frac{\kappa}{p}\mu_{n,\delta+2}(\Sigmama))(\sqrt{\frac{\ln p }{p}}+\frac{\tau}{p}+p^{-\frac{1}{\delta}} )$} at 651 450
    \pinlabel \rotatebox{0}{ {\footnotesize $(1+\frac{\kappa}{p}\mu_{n,\delta+2}(\Sigmama))(\frac{\sqrt{\ln p} }{\sqrt{p}}+\frac{\tau}{p}+p^{-\frac{1}{\delta}} )$ } } [l] at 1154 822
    \pinlabel \rotatebox{0}{$(1+\mu_{n,\delta+2}(\Sigmama)\frac{\kappa}{p})p^{-\frac{1}{\delta}}$} [r] at 150 880
    \endlabellist
    \centering
    \includegraphics[width=0.45\linewidth]{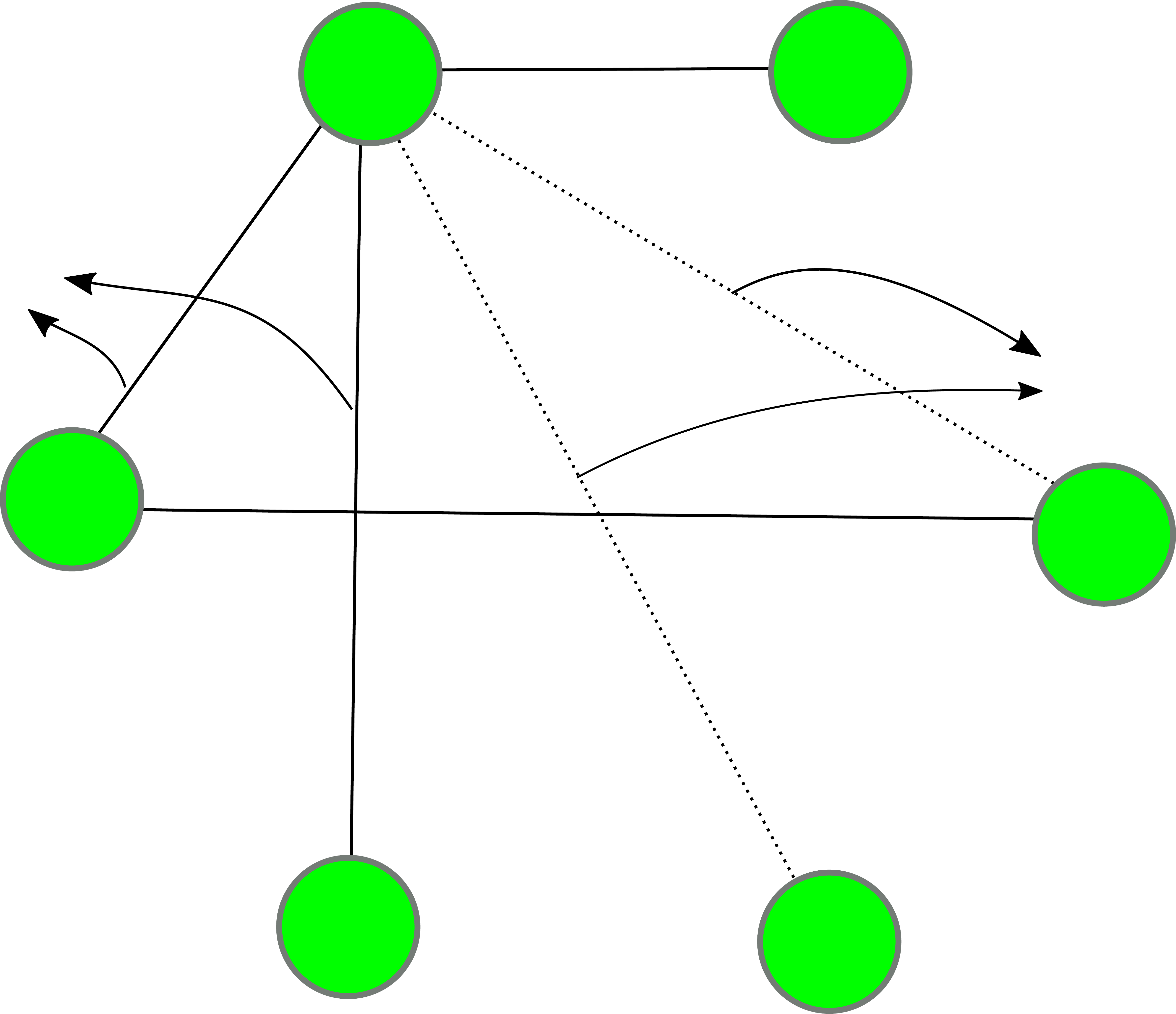}%
    \caption{
    Illustration of the pairwise total variation distances between the vertex counts associated with the empirical correlation and partial correlation graph established by the portmanteau result Proposition \ref{thm:allclose}. The 4 solid edges correspond to existence of direct upper bounds on the total variation distances between two types of vertex counts, where the weights correspond respectively to the $4$ upper bounds (neglecting constant coefficients) in Proposition \ref{thm:allclose}. Dashed edges correspond to indirect upper bounds of the associated total variation distances obtained by applying the triangular inequalities, with weights computed from the solid path connecting the two vertices. 
}
    \label{fig:relations}
\end{figure}

Figure \ref{fig:relations} illustrates the relations established by the portmanteau result Proposition \ref{thm:allclose}.
Dashed edges correspond to indirect upper bounds of the total variation distances between vertices, with weights computed from solid path connecting the two vertices. For instance the weight of dash edge between $N_{E_\delta}^{(\Rma)}$ and $N_{V_\delta}^{(\Pma)}$ is computed from    
\begin{align*}
&d_{\TV}\left(\mathscr{L}\left(N_{\VV_\delta}^{(\Pma)}\right),\mathscr{L}\left(N_{E_\delta}^{(\Rma)}\right)  \right) \\
  \leq & \E \left|N_{E_\delta}^{(\Rma)}-N_{V_\delta}^{(\Pma)}\right| \\
 \leq &  \E \left|N_{E_\delta}^{(\Rma)}-N_{\breve{V}_\delta}^{(\Pma)}\right|+\E \left|N_{E_\delta}^{(\Rma)}-N_{E_\delta}^{(\Pma)}\right| \\ 
 \leq & \E \left|N_{E_\delta}^{(\Rma)}-N_{\breve{V}_\delta}^{(\Rma)}\right|+ \E \left|N_{\breve{V}_\delta}^{(\Rma)}-N_{\breve{V}_\delta}^{(\Pma)}\right| +\E \left|N_{E_\delta}^{(\Rma)}-N_{E_\delta}^{(\Pma)}\right|\\
 \leq & \left(C_{E_\delta}^{(\Pma)} + C_{\breve{V}_\delta}^{(\Pma)}+\frac{(\delta+1)^2}{\delta!} \gamma^{\delta+1}\right)   \left(1+\frac{\kappa-1}{p}\mu_{n,\delta+2}(\Sigmama)\right)\left(\sqrt{\frac{\ln p }{p}}+\frac{\tau}{p}+p^{-\frac{1}{\delta}} \right),
\end{align*}
where the first step follows from Lemma \ref{prop:tvtomean} in Section \ref{sec:auxlem} in the Appendix, the second step follows from Lemma \ref{lem:6quantitiesine}, and the last step follows from Proposition \ref{thm:allclose}. 

Figure \ref{fig:relations} graphically illustrates the implication of Proposition \ref{thm:allclose} that all $6$ quantities $\{N_{E_\delta}^{(\Psima)},N_{\breve{V}_\delta}^{(\Psima)},N_{V_\delta}^{(\Psima)}: k\in \{\Rma,\Pma\}  \}$ are mutually close in total variation. 
As a result, the closeness of one quantity among the 6 to some distribution in total variation implies the closeness of all 6 quantities to that same distribution. By Proposition \ref{prop:edgecor},  $\mathscr{L}\left(N_{E_\delta}^{(\Rma)}\right)$ is close to the non-asymptotic compound Poisson distribution $\CP(\lambda_{p,n,\delta,\rho},\bm{\zeta}_{n,\delta,\rho})$ in total variation, which implies all $6$ quantities $\{N_{E_\delta}^{(\Psima)},N_{\breve{V}_\delta}^{(\Psima)},N_{V_\delta}^{(\Psima)}: k\in \{\Rma,\Pma\}  \}$ are close in total variation to $\CP(\lambda_{p,n,\delta,\rho},\bm{\zeta}_{n,\delta,\rho})$. We state this formally as Theorem \ref{thm:Poissonultrahigh} in the next subsection. 

\subsection{A unified theorem for finite $p$} 
\label{sec:compoihighdim}

The following theorem is a variant of Theorem \ref{cor:Poissonlimit} for finite $p$. It states that if the threshold $\rho$ is properly chosen, and $\Sigmama$ satisfies the $(\tau,\kappa)$ sparsity condition, then $\mathscr{L}(\bar{N}_\delta)$ can be approximated by a compound Poisson distribution.

\begin{thm}[Compound Poisson approximation for finite $p$] \label{thm:Poissonultrahigh}
    Let $ n \geq 4$, $\delta\in [p-1]$, and $\gamma>0$ be given. Consider $\Xma\sim \mathcal{VE}(\bm{\mu},\Sigmama,g)$.  
    Suppose $2p^{1+\frac{1}{\delta}}P_n(r_\rho)\leq \gamma$. 
    Suppose $\Sigmama$, after some row-column permutation, is $(\tau,\kkkk)$ sparse with $\tau \leq \frac{p}{2}$ and $\mu_{n,2\delta+2}(\Sigmama)\frac{\kappa}{p}<1$. Suppose the condition $\left(\sqrt{\frac{\nn-1}{p}}+\sqrt{\frac{\delta\ln p}{p}}\right) \leq c$ is satisfied for some positive and \ONE{sufficiently} small  constant $c$. 
    Then $\bar{N}_\delta$,  a generic random variable in the set $\{N_{E_\delta}^{(\Psima)},N_{\breve{V}_\delta}^{(\Psima)},N_{V_\delta}^{(\Psima)}: k\in \{\Rma,\Pma\}  \}$, satisfies:
    \begin{equation}
    d_{\TV}\left(\mathscr{L}\left(\bar{N}_\delta\right), \CP(\lambda_{p,n,\delta,\rho},\bm{\zeta}_{n,\delta,\rho}) \right)  \leq C_{n,\delta,\gamma}  \left(\mu_{n,2\delta+2}\left(\Sigmama\right)\frac{\kkkk}{p} +p^{-\frac{1}{\delta}}+E(p,\tau)\right), 
    \label{eqn:mainthm2}
    \end{equation}    
    where
    $$
    E(p,\tau) = \begin{cases}
    0 & \text{if } \bar{N}_\delta = N_{E_\delta}^{(\Rma)},\ N_{\breve{V}_\delta}^{(\Rma)} \text{ or } N_{V_\delta}^{(\Rma)}, \\ 
    \sqrt{\frac{\ln p }{p}}+\frac{\tau}{p} & \text{if } \bar{N}_\delta = N_{E_\delta}^{(\Pma)}, N_{\breve{V}_\delta}^{(\Pma)} \text{ or } N_{V_\delta}^{(\Pma)}.
    \end{cases}
    $$
\end{thm}

If only the vertex counts in the empirical correlation graph is of interest, then the $(\tau,\kappa)$ sparsity assumption can be relaxed to row-$\kappa$ sparsity.

\begin{lem}[Compound Poisson approximation in empirical correlation graph]
\label{item:Poissonultrahigha}
 Let $ n \geq 4$, $\delta\in [p-1]$, and $\gamma>0$ be given. Consider $\Xma\sim \mathcal{VE}(\bm{\mu},\Sigmama,g)$.  
    Suppose $2p^{1+\frac{1}{\delta}}P_n(r_\rho)\leq \gamma$. 
    Suppose $\Sigmama$ is row-$\kkkk$ sparse with $\mu_{n,2\delta+2}(\Sigmama)\frac{\kappa}{p}<1$. Then $\tilde{N}_\delta$, a generic random variable in the set $\{N_{E_\delta}^{(\Rma)},N_{\breve{V}_\delta}^{(\Rma)},N_{V_\delta}^{(\Rma)}  \}$, satisfies
    \begin{equation}
    d_{\TV}\left(\mathscr{L}\left(\tilde{N}_\delta\right), \CP(\lambda_{p,n,\delta,\rho},\bm{\zeta}_{n,\delta,\rho}) \right)  \leq C_{n,\delta,\gamma}  \left(\mu_{n,2\delta+2}\left(\Sigmama\right)\frac{\kkkk}{p}+p^{-\frac{1}{\delta}}\right).
    \label{eqn:mainthm}
    \end{equation}
\end{lem}

\begin{rem}    
\label{rem:nonasythm}
    The assumption $\mu_{n,2\delta+2}(\Sigmama)\frac{\kappa}{p}<1$ is only used to obtain simpler expressions for the upper bounds in \eqref{eqn:mainthm2} and \eqref{eqn:mainthm}. Without this assumption, similar inequalities hold with the upper bounds that are the sum the upper bounds in Proposition \ref{prop:edgecor} and Proposition \ref{thm:allclose}. 
    See the third paragraph in Remark \ref{rem:edgecor} for additional discussions.       
    Observe that 
    Theorem \ref{thm:Poissonultrahigh} and Lemma \ref{item:Poissonultrahigha} hold for any mean $\bm{\mu}$ and any shaping function $g$ when $\Xma\sim \mathcal{VE}(\muve,\Sigmama,g)$. This is a consequence of the invariance property of the U-scores distribution (see Remark \ref{rem:Udisconsequence}). 
    \ONE{In particular, none of the constants in Theorem \ref{thm:Poissonultrahigh}  depend on $\muve$ or $g$.}
\end{rem}

Theorem \ref{thm:Poissonultrahigh} and Lemma \ref{item:Poissonultrahigha} directly follow from Proposition \ref{prop:edgecor} and Proposition \ref{thm:allclose} and hence their proofs are omitted. Theorem \ref{thm:Poissonultrahigh} and Lemma \ref{item:Poissonultrahigha} provide an approximation for the family-wise error rate (FWER) \cite{hero2012hub}
$$
\P(\bar{N}_\delta>0)\approx 1- e^{-\lambda_{p,n,\delta,\rho}}.
$$

We end this subsection by making a comparison between Theorem \ref{thm:Poissonultrahigh} and Theorem \ref{cor:Poissonlimit}. 
Theorem \ref{thm:Poissonultrahigh} and Theorem \ref{cor:Poissonlimit} provide compound Poisson approximations to $\bar{N}_\delta$ respectively when $p$ is finite and when $p\to\infty$. By taking the limit as $p\to\infty$, we obtain simpler formulae for parameters of the approximation. Specifically the  increment distribution $\bm{\zeta}_{n,\delta,\rho}$ of the non-asymptotic compound Poisson distribution in Theorem \ref{thm:Poissonultrahigh} depends on conditional probabilities in the random pseudo geometric graph as in \eqref{eqn:alphandeltarho}. On the other hand, the  increment distribution  $\bm{\zeta}_{n,\delta}$ of the limiting compound Poisson distribution in Theorem \ref{cor:Poissonlimit} depends on probabilities in the random geometric graph as in \eqref{eqn:alphaelldef}, which is relatively simpler. For instance, when $\delta=2$, an analytical formula for $\bm{\zeta}_{n,2}$ can be obtained (see Example \ref{exa:limitcompoidelta2}). Obtaining an analytical formula for  $\bm{\zeta}_{n,2,\rho}$ does not seem straightforward. 

\ONE{Despite the fact that the limiting compound Poisson distribution in Theorem \ref{cor:Poissonlimit} is relatively simpler than Theorem \ref{thm:Poissonultrahigh}, it has the disadvantage that the approximation is not accurate unless $p$ is large. Moreover,} Theorem \ref{thm:Poissonultrahigh} is stronger than Theorem \ref{cor:Poissonlimit} in the sense that Theorem \ref{thm:Poissonultrahigh} provides explicit upper bounds for the approximation errors, while Theorem \ref{cor:Poissonlimit} simply provides a limit but no convergence rates. \ONE{ \ONE{For further discussion see Remark \ref{rem:comparisonsfinitelimit} in the next subsection.}}

\subsection{Proof of Theorem \ref{cor:Poissonlimit}}
\label{sec:proofoftheorem1}
In this subsection we present results on the limit of $\CP(\lambda_{p,n,\delta,\rho},\bm{\zeta}_{n,\delta,\rho})$ as $p\to\infty, \rho\to 1$, followed by a proof of Theorem \ref{cor:Poissonlimit}.

To study the limiting distribution of $\CP(\lambda_{p,n,\delta,\rho},\bm{\zeta}_{n,\delta,\rho})$, which requires the limit of the parameter $\alpha(\ell,r_\rho)$, the next two lemmas are useful.

\begin{lem}\label{prop:rangeo}
    Consider $r < 2/\sqrt{5}$ and $\delta\geq 1$. Suppose $\{\uve'_i\}_{i=1}^{\mgeo}\ \overset{\text{i.i.d.}}{\sim} \ \unif(S^{\Ngeo})$. Then for any $\ell \in [\mgeo]$,
    \begin{align*}
    &\P\left(\nmd\left(\{\uve'_i \}_{i=1}^{\delta+1}, r  \right)=\ell|\deg(\uve'_{\delta+1})=\delta\right) \\
    = & \P\left(\pnmd\left(\{\uve'_i  \}_{i=1}^{\delta+1}, r  \right)=\ell|\deg(\uve'_{\delta+1})=\delta\right), \numberthis \label{eqn:newconditionalprobinvariant}
    \end{align*}
    where $\deg(\uve'_{\delta+1})$ on the left (right) side is the degree of vertex $\uve'_{\delta+1}$ in the corresponding random (pseudo) geometric graph.
\end{lem}
Lemma \ref{prop:rangeo} establishes that the conditional distributions of the number of universal vertices in the random geometric graph and the random pseudo geometric graph are identical. 
By Lemma \ref{prop:rangeo} and \eqref{eqn:alphandeltarho}, 
\begin{equation}
\alpha(\ell,r_\rho) = \P\left(\nmd\left(\{\uve'_i  \}_{i=1}^{\delta+1}, r_\rho  \right)=\ell|\deg(\uve'_{\delta+1})=\delta\right) \label{eqn:alphalrhorgg}
\end{equation}
when $r_\rho<2/\sqrt{5}$ or equivalently $\rho>3/5$. Therefore to study the limit of $\alpha(\ell,r_\rho)$ when $r_\rho\to 0$ or equivalently $\rho\to 1$, the limiting form of the right hand side of \eqref{eqn:alphalrhorgg} is required, given by the following lemma.

\begin{lem}
\label{prop:rangeolimit}
    Let $\mgeominus\geq 1$ and $n\geq 3$. Suppose $\{\uve'_i  \}_{i=1}^{\mgeo}\overset{\text{i.i.d.}}{\sim} \unif(S^{\Ngeo})$ and $\{\tilde{\uve}_i\}_{i=1}^{\mgeominus} \overset{\text{i.i.d.}}{\sim} \unif(B^{\Ngeo})$. Then for any $\ell \in [\mgeo]$,
    \begin{align*}
    &\lim_{r\to 0}\P\left(\nmd\left(\{\uve'_i  \}_{i=1}^{\delta+1}, r \right)=\ell|\deg(\uve'_{\delta+1})=\delta \right)\\
    = & \P\left(\nmd\left(\{\tilde{\uve}_i \}_{i=1}^{\delta}, 1  \right)=\ell-1\right). \numberthis \label{eqn:tobeproof2}
\end{align*}
\end{lem}

One immediate consequence of Lemma \ref{prop:rangeolimit} and \eqref{eqn:alphalrhorgg} is that $\lim_{\rho\to 1}\alpha(\ell,r_\rho)=\alpha_\ell$ for any $\ell\in [\delta+1]$. As $p\to \infty$, the condition $2p^{1+\frac{1}{\delta}}P_n(r_\rho)\leq \gamma$ in Theorem \ref{thm:Poissonultrahigh} entails $\rho\to 1$. The following lemma states if the rate of $\rho\to 1$ is coupled with the rate $p\to \infty$, then the non-asymptotic compound Poisson distribution $\CP(\lambda_{p,n,\delta,\rho},\bm{\zeta}_{n,\delta,\rho})$ converges \text{ in distribution} to the limiting compound Poisson distribution $\CP(\lambda_{n,\delta}(e_{n,\delta}),\bm{\zeta}_{n,\delta})$. 
\begin{lem} 
\label{lem:compoundpoissonlimit}
Suppose as $p\to \infty$, $\rho\to 1$ such that $a_n2^{\frac{n}{2}}p^{1+\frac{1}{\delta}}(1-\rho)^{\frac{n-2}{2}}\to e_{n,\delta}$, where $e_{n,\delta}$ is some positive constant that possibly depends on $n$ and $\delta$. Then
\begin{equation}
\CP(\lambda_{p,n,\delta,\rho},\bm{\zeta}_{n,\delta,\rho}) \overset{\Dc}{\to} \CP(\lambda_{n,\delta}(e_{n,\delta}),\bm{\zeta}_{n,\delta}). \label{eqn:jumpsizedislimit}
\end{equation}
\end{lem}

With the above results established we are in a position to prove Theorem \ref{cor:Poissonlimit} and Lemma \ref{item:Poissonlimita}.

\begin{proof}[Proof of Theorem \ref{cor:Poissonlimit} and Lemma \ref{item:Poissonlimita}]
Theorem \ref{cor:Poissonlimit} directly follows from Theorem \ref{thm:Poissonultrahigh} and Lemma \ref{lem:compoundpoissonlimit}. Lemma \ref{item:Poissonlimita} directly follows from Lemma \ref{item:Poissonultrahigha} and Lemma \ref{lem:compoundpoissonlimit}.
\end{proof}
 
\begin{rem}
\label{rem:comparisonsfinitelimit}
\ONE{
 In \eqref{eqn:slowupperbound} of the proof of Lemma \ref{lem:compoundpoissonlimit}, it is shown that part of the upper bound for its error rate is of the order $1-\rho$. This particular rate, however, decreases as $n$ increases. More concisely, if one chooses $\rho$ according to \eqref{eqn:rhopformula} then $1-\rho$ is of the order $p^{-(1+\frac{1}{\delta})\frac{2}{n-2}}$. Hence the convergence to $\CP(\lambda_{n,\delta}(e_{n,\delta}),\bm{\zeta}_{n,\delta})$ in Theorem \ref{cor:Poissonlimit} is only accurate for large $p$. On the other hand, the upper bound in Theorem \ref{thm:Poissonultrahigh} only depends on $\rho$ through $2p^{1+\frac{1}{\delta}}P_n(r_\rho)\leq \gamma$, which holds for small $p$, again if $\rho$ is chosen according to \eqref{eqn:rhopformula} (see the first paragraph of Remark \ref{rem:edgecor} for related discussion). Hence Theorem \ref{thm:Poissonultrahigh} provides an accurate approximation to $\bar{N}_{\delta}$ even for small $p$. The accuracy of these approximations for various values of $\rho$ and $n$ is numerically illustrated in Figure \ref{fig:finiteplimitcomparison} in Section \ref{sec:simulation} of the Appendix.

}
\end{rem}

\section{Convergence of moments}
\label{sec:convergenceofmoments}

 Moment expressions are useful for a number of reasons, including characterizing the behavior of phase transition thresholds and the expected number of false discoveries \cite{hero2011large}. Theorems \ref{cor:Poissonlimit} and \ref{thm:Poissonultrahigh} only prescribe the distribution of $\bar{N}_\delta$ but not the moments. In this subsection, we present approximations to the first moment and second moment of $\bar{N}_\delta$ for finite $p$.

Let $Z\sim \operatorname{CP}(\lambda_{p,n,\delta,\rho},\bm{\zeta}_{n,\delta,\rho})$. Then we can represent $ Z= \sum_{i=1}^N Z_i $, where $N$ is distributed as a Poisson with mean $\lambda_{p,n,\delta,\rho}$, $Z_i \overset{\text{i.i.d.}}{\sim} \bm{\zeta}_{n,\delta,\rho} $ and $N$ is independent of each $Z_i$. The first two moments of $Z$ are: 
\begin{align*}
    \E Z = & \E N \E Z_1 =  \binom{p}{1}\binom{p-1}{\delta}(2P_n(r_\rho))^{\delta},  \numberthis \label{eqn:firstmomentofCP}\\
    \E Z^2 = & \E N \E Z_1^2 + (\E N \E Z_1)^2  \\
    = &\binom{p}{1}\binom{p-1}{\delta}(2P_n(r_\rho))^{\delta} \sum_{\ell =1}^{\delta+1} \ell \alpha(\ell,r_\rho)  +  \left( \binom{p}{1}\binom{p-1}{\delta}(2P_n(r_\rho))^{\delta} \right)^2. \numberthis \label{eqn:secondmomentofCP}
\end{align*}

The next lemma provides an upper bound for the difference between 
the first moment of $N_{E_{\delta}}^{(\Rma)}$ and the first moment of the compound Poisson specified in \eqref{eqn:firstmomentofCP}.

\begin{lem} \label{lem:meanedge}
Let $p\geq n \geq 4$, $\delta\in [p-1]$ and $\gamma>0$ be given. Suppose $\Xma\sim \mathcal{VE}(\muve,\Sigmama,g)$. Suppose $2p^{1+\frac{1}{\delta}}P_n(r_\rho)\leq \gamma$, and $\Sigmama$ is row-$\kkkk$ sparse. Then
$$
\left| \E N_{E_{\delta}}^{(\Rma)}-\E Z \right|\leq  \frac{(\delta+1)}{2((\delta-1)!)}\gamma^{\delta} \mu_{n,\delta+1}\left(\Sigmama\right) \frac{\kappa-1}{p}.
$$
\end{lem}

\begin{rem}
\label{rem:rhorate}
{ 
Lemma \ref{lem:meanedge} implies that the condition on $\rho$ assumed in Theorem \ref{cor:Poissonlimit} is in fact necessary and sufficient for the mean to converge to a finite and strictly positive limit. For simplicity we specialize to the case that $\Sigmama$ is diagonal, which is row-$\kappa$ sparse with  $\kappa=1$. By Lemma \ref{lem:meanedge} and (\ref{eqn:firstmomentofCP}), 
$$
\E N_{E_{\delta}}^{(\Rma)}=\E Z = \binom{p}{1}\binom{p-1}{\delta}(2P_n(r_\rho))^{\delta},
$$
which via Stirling approximation behaves as
$\frac{1}{\delta!}(2p^{1+\frac{1}{\delta}}P_n(r_\rho))^{\delta}$ for large $p$. 
By Lemma \ref{pderivative} \ref{item:pderivativeb}
we further have, again for large $p$, 
$$
\E N_{E_{\delta}}^{(\Rma)} \approx \frac{1}{\delta!}\left(2p^{1+\frac{1}{\delta}}a_n (1-\rho)^{\frac{n-2}{2}}\right)^{\delta}.
$$
Thus, in order that the mean count have a non-degenerate limit when $p\to \infty$, $2p^{1+\frac{1}{\delta}}a_n (1-\rho)^{\frac{n-2}{2}}$ must converge to some strictly positive and finite constant value. If $\rho$ does not converge to $1$, or if it converges to  $1$ at a slower rate, then $\E\tilde{N}_\delta$ diverges  to $\infty$, while if $\rho$ converges to $1$ at a faster rate then $\E \tilde{N}_\delta$ converges to $0$. This is reflected in the phase transition phenomenon discussed in \cite{hero2011large,hero2012hub}. Thus the rate on $\rho$ in Theorem \ref{cor:Poissonlimit} 
is sharp.}
\end{rem}

By combining the preceding lemma and the portmanteau result of Proposition \ref{thm:allclose} one immediately obtains approximations for the first moment of all $6$ quantities $\{N_{E_\delta}^{(\Psima)},N_{\breve{V}_\delta}^{(\Psima)},N_{V_\delta}^{(\Psima)}: \Psima \in \{\Rma,\Pma\}  \}$ for finite $p$. 


We also characterize the second moment approximation of $N_{E_\delta}^{(\Rma)}$ in the next lemma
under the same conditions as Lemma \ref{lem:meanedge}.

\begin{prop}[Second moment bounds for subgraph counts]
\label{lem:2ndmoment}
Let $p\geq n \geq 4$, $\delta\in [p-1]$ and $\gamma>0$ be given. Suppose $\Xma\sim \mathcal{VE}(\muve,\Sigmama,g)$. Suppose $2p^{1+\frac{1}{\delta}}P_n(r_\rho)\leq \gamma$, and $\Sigmama$ is row-$\kkkk$ sparse. Then
$$
\left| \E \left(N_{E_{\delta}}^{(\Rma)}\right)^2-\E Z^2 \right|\leq  C_{n,\delta,\gamma} \left(\mu_{n,2\delta+2}(\Sigmama)\frac{\kappa}{p} +p^{-1/\delta}   \right).
$$
\end{prop}

We extend the preceding proposition to 
other quantities in $\{N_{E_\delta}^{(\Psima)},N_{\breve{V}_\delta}^{(\Psima)},N_{V_\delta}^{(\Psima)}: \Psima \in \{\Rma,\Pma\}  \}$ by generalizing Proposition \ref{thm:allclose}
to $L^1$ distance between the square of the quantities.

\begin{prop}[Portmanteau result for $L^2$ distance] 
\label{prop:allcloseL2}
    { 
        Let $p\geq n \geq 4$ and $\Xma\sim \mathcal{VE}(\bm{\mu},\Sigmama,g)$. Let $\delta \in [p-1]$. Suppose $2p^{1+\frac{1}{\delta}}P_\nn(r_\rho)\leq \gamma$. 
    \begin{enumerate}[label=(\alph*)]
        \item  
        \label{item:allcloseL2a}
        Suppose $\Sigmama$ is row-$\kappa$ sparse. Then for $\tilde{N}_\delta \in \left\{ N_{\breve{V}_\delta}^{(\Rma)}, N_{V_\delta}^{(\Rma)} \right\} $,
        $$
         \E \left| \left(\tilde{N}_\delta\right)^2 -    \left(\Nedged^{(\Rma)}\right)^2  \right|
\leq C_{n,\delta,\gamma}  \left(1+\mu_{n,2\delta+3}(\Sigmama)\frac{\kappa-1}{p}\right)p^{-1/\delta}. $$

        \item \label{item:allcloseL2b}
        Suppose $\Sigmama$, after some row-column permutation, is $(\tau,\kappa)$ sparse with $\tau\leq \frac{p}{2}$. Moreover,  $\left(\sqrt{\frac{\nn-1}{p}}+\sqrt{\frac{\delta\ln p}{p}}\right) \leq c$ hold for some positive and \ONE{sufficiently} small  constant $c$. Then
        \begin{align*}
         \E\left|\left(\Nedged^{(\Pma)}\right)^2-\left(\Nedged^{(\Rma)}\right)^2\right|
        \leq  C_{n,\delta,\gamma}  \left(1+\mu_{n,2\delta+2}(\Sigmama)\frac{\kappa-1}{p}\right) \left(\frac{\sqrt{\ln p}}{\sqrt{p}} + \frac{\tau}{p}\right). 
        \end{align*}
        
        \item  
        \label{item:allcloseL2c}
        Suppose the same conditions as in part \ref{item:allclosec} hold. Then for $\tilde{N}_\delta \in \left\{ N_{\breve{V}_\delta}^{(\Pma)}, N_{V_\delta}^{(\Pma)} \right\} $
        \begin{align*}
         \E\left|\left(\tilde{N}_{\delta}\right)^2-\left(N_{E_\delta}^{(\Pma)}\right)^2\right|
        \leq   C_{n,\delta,\gamma}  \left(1+\mu_{n,2\delta+3}(\Sigmama)\frac{\kappa-1}{p}\right)p^{-1/\delta}. 
         \end{align*}
    \end{enumerate} 
}
\end{prop}

By applying triangle inequalities to Proposition \ref{prop:allcloseL2} \ref{item:allcloseL2b} and \ref{item:allcloseL2c}, one obtain for $\tilde{N}_{\delta}\in \{N_{V_\delta}^{(\Pma)}, N_{\breve{V}_\delta}^{(\Pma)} \}$
$$
\E\left|\left(\tilde{N}_{\delta}\right)^2-\left(N_{E_\delta}^{(\Rma)}\right)^2\right|\leq  C_{n,\delta,\gamma}  \left(1+\mu_{n,2\delta+3}(\Sigmama)\frac{\kappa-1}{p}\right)
\left(\frac{\sqrt{\ln p}}{\sqrt{p}} + \frac{\tau}{p}+p^{-1/\delta}\right).
$$
Thus we have established the $L^1$ distance between the square of each term in $\{N_{E_\delta}^{(\Psima)},N_{\breve{V}_\delta}^{(\Psima)},N_{V_\delta}^{(\Psima)}: \Psima \in \{\Rma,\Pma\}  \}$ and $\left(N_{E_\delta}^{(\Rma)}\right)^2$. By combining Proposition \ref{lem:2ndmoment} and Proposition \ref{prop:allcloseL2} one immediately obtains approximations to the second moment of $\bar{N}_\delta$ for finite $p$. 

While this section focuses on the approximation to the first and second moments of $\bar{N}_\delta$ for finite $p$, their limits can also be obtained when $p\to\infty$ and $\rho\to 1$ at the rate specified in Theorem \ref{cor:Poissonlimit}.

\section{Analytical expressions for the compound Poisson parameters}
\label{sec:limitingcompoundPoisson}

Recall that the limiting compound Poisson distribution in Theorem \ref{cor:Poissonlimit} is in terms of $\alpha_\ell$, while the non-asymptotic compound Poisson distributions in Theorem \ref{thm:Poissonultrahigh} is in terms of $\alpha(\ell,r_\rho)$. Moreover the second moment approximation established in Proposition \ref{lem:2ndmoment} also involves the term $\alpha(\ell,r_\rho)$ due to \eqref{eqn:secondmomentofCP}. Since $\alpha_\ell$ and $\alpha(\ell,r_\rho)$ are expressed in terms of random (pseudo)  geometric graphs, they are tedious to compute. 
 In subsection \ref{sec:approximationstoalphaell}, we obtain simple analytical approximations to these quantities for moderately large $n$ or $\delta$. In subsection \ref{sec:exacttoalphaell}, exact analytical formulae for these quantities are established for the special case $\delta=1$ and $\delta=2$. 

\subsection{Approximations to $\alpha_\ell$ and $\alpha(\ell,r_\rho)$ when $\delta\geq 2$}
\label{sec:approximationstoalphaell}

In this subsection we show that when $n$ or $\delta$ is moderately large, $\alpha_1\approx 1$ and $\alpha_\ell\approx 0$ for $2\leq \ell\leq \delta+1 $ (parallel results for $\alpha(\ell,r_\rho)$ are $\alpha(1,r_\rho)\approx 1$ and $\alpha(\ell,r_\rho)\approx 0$ for $2\leq \ell\leq \delta+1$). These approximations yield simple formulae for parameters of the limiting and non-asymptotic compound Poisson distributions, 
and for parameters of the second moment \eqref{eqn:secondmomentofCP}. Importantly, the compound Poisson distributions are well approximated by Poisson distributions for $n,\delta$ moderately large. 

To show that $\alpha_\ell$ is small for $\ell\geq 2$ and $\alpha_1\approx 1$, it suffices to establish a vanishing upper bound on $\sum_{\ell=2}^{\delta+1}\alpha_\ell$, since $\sum_{\ell=1}^{\delta+1}\alpha_\ell=1$. By definition of $\alpha_\ell$, $\sum_{\ell=2}^{\delta+1}\alpha_\ell=\P\left(\nmd\left(\{\tilde{\uve}_i  \}_{i=1}^\delta, 1   \right) \geq 1 \right)$, where $\{\tilde{\uve}_i\}_{i=1}^{\mgeominus} \overset{\text{i.i.d.}}{\sim} \unif(B^{\Ngeo})$. The next lemma establishes an upper bound on the geometric quantity $\P\left(\nmd\left(\{\tilde{\uve}_i  \}_{i=1}^\delta, 1   \right) \geq 1 \right)$.

\begin{lem} \label{lem:randgeogra}
    Let $n\geq 4$ and $\mgeominus \geq 2$. 
    \begin{enumerate}[label=(\alph*)]
        \item 
        \label{item:randgeograa}
         Consider $\{\tilde{\uve}_i\}_{i=1}^{\mgeominus} \overset{\text{i.i.d.}}{\sim} \unif(B^{\Ngeo})$. Then 
        \begin{align*}
        \sum_{\ell=2}^{\delta+1}\alpha_\ell = & \P\left(\nmd\left(\{\tilde{\uve}_i  \}_{i=1}^\mgeominus, 1  ; \Ngeo \right) \geq  1 \right) \\
        \leq & \mgeominus (\Ngeo)\int_{0}^1 \left(1-\frac{r^2}{4}\right)^{\frac{(\Ngeo)(\mgeominus-1)}{2}} r^{n-3}dr\\
        = & \mgeominus(\Ngeo) 2^{n-3}B\left(\frac{1}{4};\frac{\Ngeo}{2},\frac{(\Ngeo)(\mgeominus-1)}{2}+1\right), \numberthis \label{eqn:tolvarincrementdeltabetafunction}
        \end{align*}
        where $B\left(\cdot;\cdot,\cdot\right)$ is the incomplete beta function.
    \item  
    \label{item:randgeograb}
    $$\int_{0}^1 \left(1-\frac{r^2}{4}\right)^{\frac{(\Ngeo)(\mgeominus-1)}{2}} r^{n-3}dr \leq \begin{cases} \left(\frac{4}{5}\right)^{\frac{(\Ngeo)\mgeominus -1}{2}}+\left(1-\sqrt{\frac{4}{5}}\right)\left(\frac{3}{4}\right)^{\frac{(\Ngeo)(\mgeominus-1)}{2}}, & \mgeominus=2,3, \\ \exp\left(\frac{1}{4}\right) \left(\frac{\mgeominus-1}{\mgeominus}\right)^{\frac{(\Ngeo)(\mgeominus-1)}{2}} \left(\frac{4}{\mgeominus}\right)^{\frac{n-3}{2}}, & \delta \geq 4. \end{cases} $$
\end{enumerate}
\end{lem}

Lemma \ref{lem:randgeogra} \ref{item:randgeograa} establishes an upper bound for the probability that there is at least one universal vertex in the random geometric graph generated by the uniform distribution in the unit ball. Lemma \ref{lem:randgeogra} \ref{item:randgeograb} provides an explicit upper bound for part \ref{item:randgeograa} 
and this upper bound provides the following insight 
when $n$ or $\delta$ is large. In particular, when  $\mgeominus$ is fixed, $\P\left(\nmd\left(\{\tilde{\uve}_i  \}_{i=1}^\mgeominus, 1  ; \Ngeo \right) \geq  1 \right)$ decays exponentially as $n$ increases. While $n$ is fixed, it decays at rate $\mgeominus^{-{\frac{n-3}{2}}}$ as $\mgeominus$ increases. Such decay rates also apply to $\alpha_\ell$ for $\ell\geq 2$ and $|\alpha_1-1|$.

\begin{figure}[ht]
\begin{subfigure}{.49\textwidth}
  \centering
  \includegraphics[width=.8\linewidth]{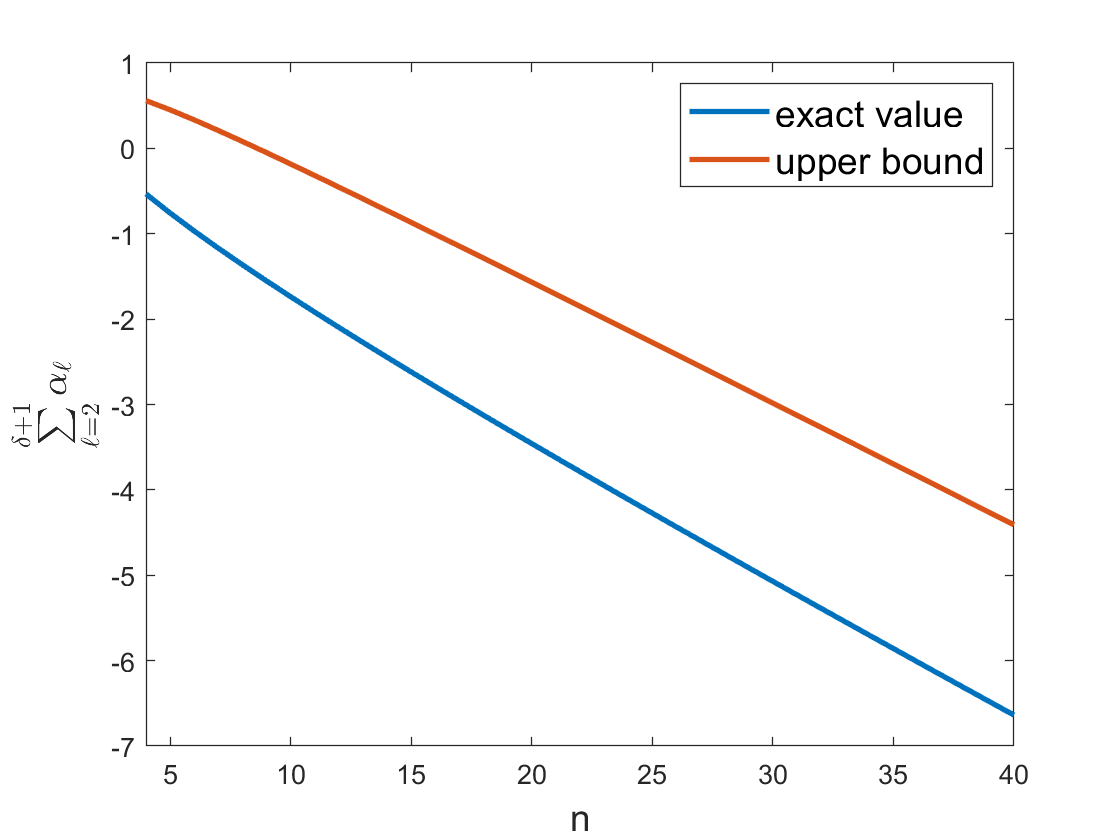}  
  \caption{Log-scale comparison of the decay when $\delta=2$}
\end{subfigure}
\begin{subfigure}{.49\textwidth}
  \centering
  \includegraphics[width=.8\linewidth]{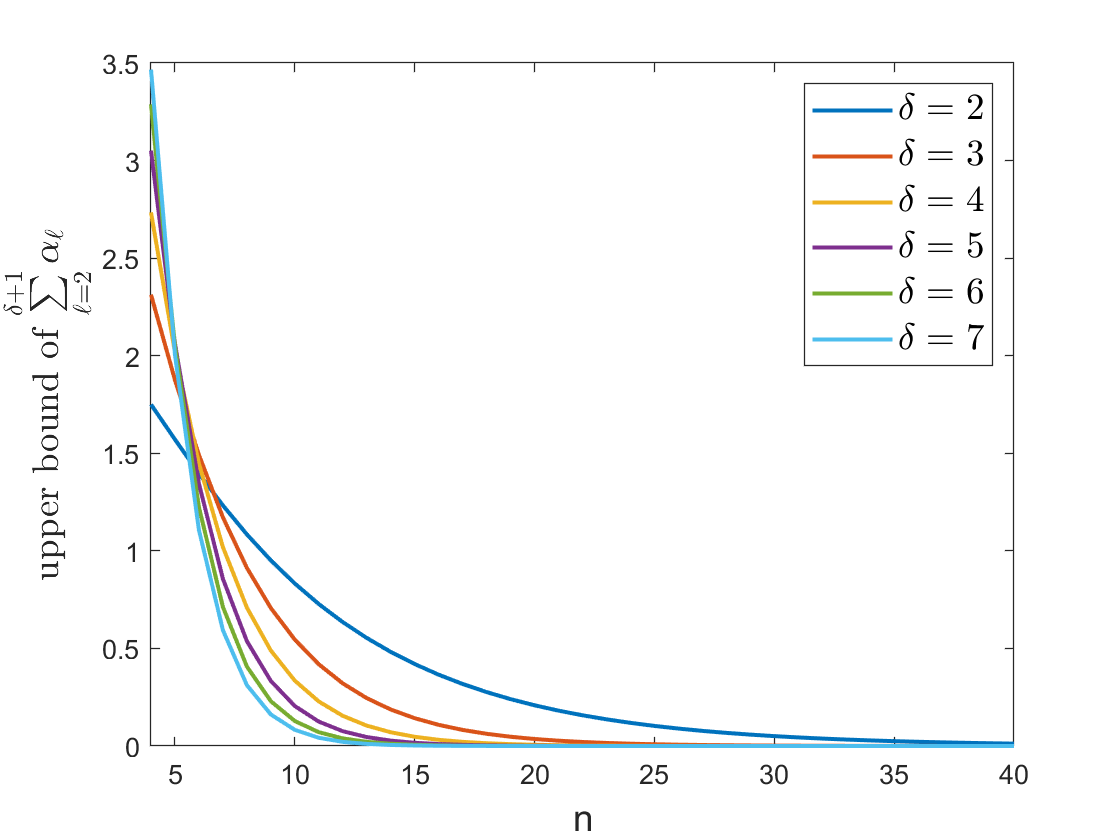}  
  \caption{Family of upper bounds}
\end{subfigure}
\caption{ (a) is a comparison in the log-scale between the upper bound on $\sum_{\ell=2}^{\delta+1}\alpha_\ell$ by \eqref{eqn:tolvarincrementdeltabetafunction} with $\delta=2$ and the exact value of $\sum_{\ell=2}^{\delta+1}\alpha_\ell$ by Example \ref{lem:alphaldelta2}. (b) is the plot of the upper bound on $\sum_{\ell=2}^{\delta+1}\alpha_\ell$ as a function of $n$ for $\delta$ from $2$ to $7$.}
\label{fig:incrementtodeltauppbou}
\end{figure}

As illustrated in Figure \ref{fig:incrementtodeltauppbou} (a), when $\delta=2$, $\sum_{\ell=2}^{\delta+1}\alpha_\ell$ and the upper bound in \eqref{eqn:tolvarincrementdeltabetafunction} both decay to $0$ exponentially fast as $n$ increases. This diagram demonstrates that the upper bound in \eqref{eqn:tolvarincrementdeltabetafunction} captures the decay rate of $\sum_{\ell=2}^{\delta+1}\alpha_\ell$. Figure \ref{fig:incrementtodeltauppbou} (b) plots the upper bounds in \eqref{eqn:tolvarincrementdeltabetafunction} as functions of $n$ for fixed $\delta$. As is clear from the plot, as long as the number of samples $n$ is above $40$, $\sum_{\ell=2}^{\delta+1}\alpha_\ell\approx 0$, which yields $\alpha_\ell\approx0$ for $2\leq \ell \leq \delta+1$ and $\alpha_1\approx 1$.
Moreover, as $\delta$ increases, the number of samples $n$ required for $\sum_{\ell=2}^{\delta+1}\alpha_\ell\approx 0$ decreases. 

Denote the Dirac Distribution at $a$ by $\Dirac(a)$. The next lemma establishes approximations to the limiting compound Poisson distribution in Theorem \ref{cor:Poissonlimit}, with the approximation errors upper bounded in terms of $\sum_{\ell=2}^{\delta+1}\alpha_\ell$. 
\begin{lem}[Approximation to the limiting Compound Poisson distribution]
    \label{prop:jumpsizedecay}
    Consider $n\geq 4$ and $\delta\geq 2$. Let $e_{n,\delta}$ be the same as in Theorem \ref{cor:Poissonlimit}.
    \begin{enumerate}[label=(\alph*)]
    \item\label{item:jumpsizedecaya}
    The increment distribution satisfies 
    $
    d_{\TV}\left(\bm{\zeta}_{n,\delta}, \Dirac(1) \right)\leq \sum_{\ell=2}^{\delta+1}\alpha_\ell.
    $
    \item \label{item:jumpsizedecayb}
     The arrival rate satisfies
    $$
    \left|\lambda_{n,\delta}(e_{n,\delta})-\frac{(e_{n,\delta})^\delta}{\delta!}\right|\leq \frac{3}{2}\frac{(e_{n,\delta})^\delta}{\delta!} \sum_{\ell=2}^{\delta+1}\alpha_\ell.
    $$
    \item \label{item:jumpsizedecayc}
    The limiting compound Poisson distribution satisfies
    $$
    d_{\TV}\left(\CP(\lambda_{n,\delta},\bm{\zeta}_{n,\delta}),\Pois\left(\frac{(e_{n,\delta})^\delta}{\delta!}\right)\right) \leq \frac{5}{2}\frac{(e_{n,\delta})^\delta}{\delta!} \sum_{\ell=2}^{\delta+1}\alpha_\ell.
    $$
    \end{enumerate}
\end{lem} 

The proof of Lemma \ref{prop:jumpsizedecay} \ref{item:jumpsizedecayc} relies on a general upper bound of the total variation distance between two compound Poisson distributions (Lemma \ref{eqn:cptv} in Section \ref{sec:auxlem} of the Appendix). Lemma \ref{eqn:cptv} may be of independent interest in compound Poisson approximation to non-Poissonian distributions (see \cite{barbour2001topics} and the references therein). 
From Lemma \ref{prop:jumpsizedecay} \ref{item:jumpsizedecaya} and Lemma \ref{lem:randgeogra}, the total variation distance between the  increment distribution and Dirac distribution at $1$ decays exponentially as $n$ increases and decays at rate $\delta^{-\frac{n-3}{2}}$ as $\delta$ increases. 
Provided that the threshold $\rho$ is chosen such that  $\frac{(e_{n,\delta})^\delta}{\delta!}$ is not large, the upper bounds in Lemma \ref{prop:jumpsizedecay} \ref{item:jumpsizedecayb} and \ref{item:jumpsizedecayc} have the same interpretation as part \ref{item:jumpsizedecaya}. In that case, the limiting compound Poisson are well approximated by the Poisson distribution $\Pois\left(\frac{(e_{n,\delta})^\delta}{\delta!}\right)$. By combining Lemma \ref{prop:jumpsizedecay} \ref{item:jumpsizedecayc} and Theorem \ref{cor:Poissonlimit}, we immediately  obtain the following result on Poisson approximation to $\mathscr{L}(\bar{N}_\delta)$ as $p\to\infty$.

\begin{prop}[Poisson approximation as $p\to \infty$]
\label{prop:poiapplimit}
Suppose all the conditions in Theorem \ref{cor:Poissonlimit} hold. Then as $p\to \infty$,
$$
d_{\TV}\left(\mathscr{L}(\bar{N}_\delta), \Pois\left(\frac{(e_{n,\delta})^\delta}{\delta!}\right)\right) \leq \frac{5}{2}\frac{(e_{n,\delta})^\delta}{\delta!} \sum_{\ell=2}^{\delta+1}\alpha_\ell.
$$
\end{prop}



We now turn our attention to the parameter $\alpha(\ell,r_\rho)$ of the non-asymptotic compound Poisson distribution. Our results on $\alpha(\ell,r_\rho)$ are similar to those on $\alpha_\ell$ but the proofs are more technical. To establish that $\alpha(1,r_\rho)\approx 1$ and $\alpha(\ell,r_\rho)\approx 0$ for $2\leq \ell\leq \delta+1$, it suffices to obtain a vanishing upper bound on $\sum_{\ell=2}^{\delta+1}\alpha(\ell,r_\rho)$. By \eqref{eqn:alphalrhorgg}, when $\rho>3/5$, $\sum_{\ell=2}^{\delta+1}\alpha(\ell,r_\rho)=\P\left(\nmd\left(\{\uve'_i  \}_{i=1}^\mgeo, r  ; \Ngeo \right) \geq  2 |\deg(\uve'_{\mgeo})=\mgeominus\right)$, where $\{\uve'_i\}_{i=1}^{\mgeominus+1} \overset{\text{i.i.d.}}{\sim} \unif(S^{\Ngeo})$. The next lemma is an analogous result to Lemma \ref{lem:randgeogra}.

\begin{lem} \label{lem:randgeograconduppbou}
    Let $n\geq 4$ and $\delta \geq 2$. 
    \begin{enumerate}[label=(\alph*)]
        \item \label{item:randgeograconduppboua}
         Consider $\{\uve'\}_{i=1}^{\mgeo} \overset{\text{i.i.d.}}{\sim} \unif(S^{\Ngeo})$. Then for $0<r<\sqrt{2}$
        \begin{align*}
       &\P\left(\nmd\left(\{\uve'_i  \}_{i=1}^\mgeo, r  ; \Ngeo \right) \geq  2 |\deg(\uve'_{\mgeo})=\mgeominus\right) \\
        \leq & \bar{h}\left(\frac{1}{\sqrt{1-r^2/4}},n,\delta\right) \delta (\Ngeo)  \int_0^1 \left(1 -\left(\frac{r_1}{2}\right)^2 \right)^{\frac{(n-2){(\mgeominusminus)}}{2}} r_1^{n-3}dr_1\\
        = & \bar{h}\left(\frac{1}{\sqrt{1-r^2/4}},n,\delta\right) \delta (\Ngeo) 2^{n-3} B\left(\frac{1}{4};\frac{n-2}{2},\frac{(n-2)(m-1)}{2}+1\right)         ,
        \end{align*}
        where $\bar{h}(x,n,\delta)=x^{n+\delta-5}x^{(n-2)(\delta-1)}$.
    \item  
    \label{item:randgeograconduppboub}
    When $r\leq \begin{cases} 2\sqrt{1-\sqrt{1-1/5}}, & \delta = 2,3 \\ 2\sqrt{1-\sqrt{1-1/\delta}}, & \delta\geq 4  \end{cases}$, the upper bound in part \ref{item:randgeograconduppboua} is upper bounded by
    \begin{align*}
        \begin{cases} \delta (\Ngeo)\{\left({\frac{\sqrt{5}}{2}}\right)^{\frac{\delta-2}{2}} \left(\frac{2}{\sqrt{5}}\right)^{\frac{(n-2)\delta-1}{2}}+(1-\frac{2}{\sqrt{5}})\left(\frac{\sqrt{5}}{2}\right)^{\frac{n+\delta-5}{2}}\left(\frac{3\sqrt{5}}{8}\right)^{\frac{(\Ngeo)(\delta-1)}{2}}\}, & \delta=2,3, \\ \delta (\Ngeo)\exp\left(\frac{1}{4}\right) \left(\sqrt{\frac{\delta}{\delta-1}}\right)^{\frac{\delta-2}{2}} \left(\sqrt{\frac{\delta-1}{\delta}}\right)^{\frac{(n-2)(\delta-1)}{2}} \left(\frac{4}{\sqrt{\delta(\delta-1)}}\right)^{\frac{n-3}{2}}, & \delta \geq 4. \end{cases}.
\end{align*}

\item 
\label{item:randgeograconduppbouc}
When $\rho> 3/5$, $\sum_{\ell=2}^{\delta+1} \alpha(\ell,r_\rho)$ is upper bounded by the upper bound in part \ref{item:randgeograconduppboua} with $r$ replaced by $r_\rho$. When $\rho \geq \begin{cases} 4/\sqrt{5}-1, & \delta = 2,3 \\ 2\sqrt{1-1/\delta}-1, & \delta\geq 4  \end{cases}$,
$\sum_{\ell=2}^{\delta+1} \alpha(\ell,r_\rho)$ is upper bounded by the upper bound in part \ref{item:randgeograconduppboub} with $r$ replaced by $r_\rho$.
\end{enumerate}
\end{lem}

Lemma \ref{lem:randgeograconduppbou} \ref{item:randgeograconduppboua} establishes an upper bound for the conditional probability that there are at least two universal vertices, conditioned on the existence of one universal vertex, in the random geometric graph over points generated by the uniform distribution on the sphere. 
The interpretations of Lemma \ref{lem:randgeograconduppbou} \ref{item:randgeograconduppboub}, \ref{item:randgeograconduppbouc} are similar to  those of Lemma \ref{lem:randgeogra} \ref{item:randgeograb} and $\sum_{\ell=2}^{\delta+1} \alpha_\ell$, and are therefore omitted.

The next lemma establishes approximations to the non-asymptotic compound Poisson distribution in Theorem \ref{thm:Poissonultrahigh} and to its second moment, with the approximation errors upper bounded in terms of $\sum_{\ell=2}^{\delta+1}\alpha(\ell,r_\rho)$. 
Denote  $\bar{\lambda}_{p,n,\delta,\rho} := \binom{p}{1}\binom{p-1}{\delta}(2P_n(r_\rho))^{\delta}$.

\begin{lem}[Approximation to the non-asymptotic Compound Poisson distribution]
    \label{prop:jumpsizedecaynonasymptotic}
    Consider $n\geq 4$ and $\delta\geq 2$. 
    \begin{enumerate}[label=(\alph*)]
    \item\label{item:jumpsizedecaynonasymptotica}
    The increment distribution satisfies 
    $
    d_{\TV}\left(\bm{\zeta}_{n,\delta,\rho}, \Dirac(1) \right)\leq \sum_{\ell=2}^{\delta+1}\alpha(\ell,r_\rho).
    $
    \item \label{item:jumpsizedecaynonasymptoticb}
     The arrival rate satisfies
    $
    \left|\lambda_{p,n,\delta,\rho} - \bar{\lambda}_{p,n,\delta,\rho}\right|\leq \frac{3}{2}\bar{\lambda}_{p,n,\delta,\rho} \sum_{\ell=2}^{\delta+1}\alpha(\ell,r_\rho).
    $
    \item \label{item:jumpsizedecaynonasymptoticc}
    The non-asymptotic compound Poisson distribution satisfies
    \begin{align*}
    d_{\TV}\left(    \CP(\lambda_{p,n,\delta,\rho},\bm{\zeta}_{n,\delta,\rho}), \Pois\left(\bar{\lambda}_{p,n,\delta,\rho}\right) \right) 
    \leq & \frac{5}{2}\bar{\lambda}_{p,n,\delta,\rho} \sum_{\ell=2}^{\delta+1}\alpha(\ell,r_\rho).
\end{align*}
    \item \label{item:jumpsizedecaynonasymptoticd}
    Let $Z\sim \operatorname{CP}(\lambda_{p,n,\delta,\rho},\bm{\zeta}_{n,\delta,\rho})$. Then the second moment of the non-asymptotic compound Poisson $\E Z^2$ in \eqref{eqn:secondmomentofCP} satisfies
    $$
    \left|\E Z^2 - (\bar{\lambda}_{p,n,\delta,\rho}+(\bar{\lambda}_{p,n,\delta,\rho})^2) \right|\leq \frac{3}{2}\bar{\lambda}_{p,n,\delta,\rho} \sum_{\ell=2}^{\delta+1}\alpha(\ell,r_\rho).
    $$
    \end{enumerate}
\end{lem} 

The proof of Lemma \ref{prop:jumpsizedecaynonasymptotic} and its interpretation are analogous to those of Lemma \ref{prop:jumpsizedecay}, and are therefore omitted. By combining Lemma \ref{prop:jumpsizedecaynonasymptotic} \ref{item:jumpsizedecaynonasymptoticc} and Theorem \ref{thm:Poissonultrahigh}, we immediately  obtain the following result on Poisson approximation to $\mathscr{L}(\bar{N}_\delta)$ for finite $p$. One can derive an equivalent proof of the following proposition by first applying Poisson approximation from Chen-Stein's method (e.g. \cite[Theorem 1]{arratia1990poisson}) to obtain an upper bound $d_{\TV}(\mathscr{L}(N_{E_\delta}^{(\Rma)}), \Pois\left(\bar{\lambda}_{p,n,\delta,\rho}\right))$, which is then combined with the portmanteau result Proposition \ref{thm:allclose}. 

\begin{prop}[Poisson approximation for finite $p$]
\label{prop:poiappfinitep}
Suppose all the conditions in Theorem \ref{thm:Poissonultrahigh} hold. Then $\bar{\lambda}_{p,n,\delta,\rho}\leq \frac{\gamma^\delta}{\delta !}$ and
\begin{align*}
d_{\TV}\left(\mathscr{L}(\bar{N}_\delta), \Pois\left(\bar{\lambda}_{p,n,\delta,\rho}\right) \right) 
\leq &  C_{n,\delta,\gamma}  \left(\mu_{n,2\delta+2}\left(\Sigmama\right)\frac{\kkkk}{p} +p^{-\frac{1}{\delta}}+E(p,\tau)\right) + \frac{5}{2}\frac{\gamma^\delta}{\delta !} \sum_{\ell=2}^{\delta+1}\alpha(\ell,r_\rho),
\end{align*}
where $E(p,\tau)$ is defined in Theorem \ref{thm:Poissonultrahigh}.
\end{prop}

We have thus established that the limiting compound Poisson distribution $\CP(\lambda_{n,\delta}(e_{n,\delta}),  \bm{\zeta}_{n,\delta})$ in Theorem \ref{cor:Poissonlimit} can be approximated 
by $\Pois\left(\frac{(e_{n,\delta})^\delta}{\delta!}\right)$, and that the non-asymptotic compound Poisson distribution $\CP(\lambda_{p,n,\delta,\rho}, \bm{\zeta}_{n,\delta,\rho})$ in Theorem \ref{thm:Poissonultrahigh} can be approximated by $\Pois(\binom{p}{1}\binom{p-1}{\delta}(2P_n(r_\rho))^{\delta})$ for sufficiently large $n$ or $\delta$. By combining these results with Theorem \ref{thm:Poissonultrahigh} and Theorem \ref{cor:Poissonlimit} we then obtain Poisson approximations to $\mathscr{L}(\bar{N}_\delta)$. In Section \ref{sec:simulation} of the Appendix, Figure \ref{fig:Poissonfinitevslimit} provides numerical simulations to demonstrate the accuracy of the Poisson approximation to $\mathscr{L}(N_\delta)$. See Figure \ref{fig:Poissonfinitevslimit} and the ensuing discussions for more details.

\subsection{Exact formulae for $\alpha_\ell$ and $\alpha(\ell,\rho)$ when $\delta=1$ and $\delta=2$}
\label{sec:exacttoalphaell}
 In this subsection we provide analytical expressions for $\alpha_\ell$ and $\alpha(\ell,r_\rho)$ when $\delta=1$ and expressions for  $\alpha_\ell$ when $\delta=2$. 

\begin{exa}[$\alpha_\ell$ and $\alpha(\ell,r_\rho)$ when $\delta=1$]\label{exa:limitcompoidelta1}
  When $\delta=1$, $\alpha_2=1$ since in the random geometric graph $\textbf{Ge}(\{\tilde{\uve}_i\}_{i=1}^1; 1, n-2)$ the number of universal vertices (vertices of degree $0$) is $1$. In this case $\alpha_1=1-\alpha_2=0$. Similarly, $\alpha(2,r_\rho)=1$ when $\delta=1$ since in the random pseudo geometric graph $\textbf{PGe}(\{\uve'_i\}_{i=1}^2; r_\rho, n-1)$ the number of universal vertices is $2$ as long as there exists one universal vertex. In this case $\alpha(1,r_\rho)=1-\alpha(2,r_\rho)=0$.
\end{exa}

\begin{rem}[Compound Poisson approximations when $\delta=1$]\label{rem:limitcompoidelta1}
Using the results for $\alpha_\ell$ in Example \ref{exa:limitcompoidelta1}, we obtained $\bm{\zeta}_{n,1}=\Dirac(2)$ and $\lambda_{\nn,1}(e_{n,1}) = \frac{1}{2}  e_{n,1}$. Then the limiting compound Poisson distribution in Theorem \ref{cor:Poissonlimit} is $\CP(\frac{1}{2}  e_{n,1}, \Dirac(2))$. On the other hand, by the results for $\alpha(\ell,r_\rho)$ in Example \ref{exa:limitcompoidelta1}, $\bm{\zeta}_{n,1,\rho}=\Dirac(2)$ and $\lambda_{p,n,1,\rho}=p(p-1)P_n(r_\rho)$. Therefore the non-asymptotic compound Poisson distribution in Theorem \ref{thm:Poissonultrahigh} is $\CP(p(p-1)P_n(r_\rho), \Dirac(2))$. The intuition for the result that the increment distributions $\bm{\zeta}_{n,1}$ and $\bm{\zeta}_{n,1,\rho}$ are $\Dirac(2)$ is as follows. Recall that $N_{E_1}^{(\Psima)}$ is twice of the number of edges. Thus $N_{E_1}^{(\Psima)}$ always increases by $2$ whenever there is a new edge and $N_{{\VVV}_1}^{(\Psima)}$ always increases by $2$ since the increment is always a new pair of vertices of degree $1$. 
Note that $N_{V_1}^{(\Psima)}$ has increment close to $2$ since $N_{V_1}^{(\Psima)}\approx N_{{\VVV}_1}^{(\Psima)}$ by Lemma  \ref{lem:6quantitiesine}. Given that the increment distributions $\bm{\zeta}_{n,1}$ and $\bm{\zeta}_{n,1,\rho}$ are $\Dirac(2)$, we have the following equivalent Poisson approximations to  $\bar{N}_1/2$: $\mathscr{L}(\bar{N}_1/2)\approx \Pois(p(p-1)P_n(r_\rho))$ and  $\bar{N}_1/2\overset{\Dc}{\to} \Pois(\frac{1}{2}  e_{n,1})$ as $p\to \infty$ and $a_n 2^{\frac{n}{2}}p^{2}(1-\rho)^{\frac{n-2}{2}}\to e_{n,1}$.
As a comparison, Proposition 1 and its proof in \cite{hero2011large} under row-$\kappa$ sparsity condition established that $N_{E_1}^{(\Rma)}/2$ converges to a Poisson distribution and obtained the limits of $\E N_{V_1}^{(\Rma)}$ and $\P(N_{V_1}^{(\Rma)}>0)$. Proposition 1 and Proposition 3 in \cite{hero2012hub} under block sparsity condition extended the results in \cite{hero2011large} to the corresponding versions in the empirical partial correlation graph, i.e. the same conclusions hold with $\Rma$ replaced by $\Pma$. Our results in Theorems \ref{thm:Poissonultrahigh},  \ref{cor:Poissonlimit} (and Lemmas \ref{item:Poissonultrahigha},  \ref{item:Poissonlimita}) with $\delta=1$ characterize the full distributions of the $6$ quantities $\{N_{E_\delta}^{(\Psima)},N_{\breve{V}_\delta}^{(\Psima)},N_{V_\delta}^{(\Psima)}: \Psima \in \{\Rma,\Pma\}  \}$, and our results in Section \ref{sec:convergenceofmoments} characterize their first and second moments,  which together contain and extend the aforementioned previous results. Moreover, our results are established with weaker sparsity assumptions and provide concise formulae for the parameters. 
\end{rem}

We next explicitly characterize $\alpha_\ell$ when $\delta=2$.

\begin{lem}[$\alpha_\ell$ when $\delta=2$]
\label{lem:alphaldelta2}
    When $\delta=2$, $\alpha_2=0$, $\alpha_3 = \frac{3}{2} I_{\frac{3}{4}}(\frac{n-1}{2},\frac{1}{2})$ and $\alpha_1=1-\alpha_3$, where $I_x(a,b)$ is the regularized incomplete Beta function.
\end{lem}

\begin{rem}[Limiting compound Poisson approximation when $\delta=2$]\label{exa:limitcompoidelta2}
 When $\delta=2$,
 by Lemma \ref{lem:alphaldelta2}, $\alpha_2=0$, $\alpha_3 = \frac{3}{2} I_{\frac{3}{4}}(\frac{n-1}{2},\frac{1}{2})$ and $\alpha_1=1-\alpha_3$. Then $\sum_{\ell=1}^3\alpha_\ell/\ell=1- I_{\frac{3}{4}}(\frac{n-1}{2},\frac{1}{2})$.
 Thus, the parameters for $\CP(\lambda_{n,2}(e_{n,2}),\bm{\zeta}_{n,2})$ in Theorem \ref{cor:Poissonlimit} are
\begin{equation}
\bm{\zeta}_{n,2}(1)=\frac{1-\frac{3}{2} I_{\frac{3}{4}}(\frac{n-1}{2},\frac{1}{2})}{1- I_{\frac{3}{4}}(\frac{n-1}{2},\frac{1}{2})}, \quad \bm{\zeta}_{n,2}(2)=0,  \quad \bm{\zeta}_{n,2}(3)=\frac{\frac{1}{2} I_{\frac{3}{4}}(\frac{n-1}{2},\frac{1}{2})}{1- I_{\frac{3}{4}}(\frac{n-1}{2},\frac{1}{2})}  \label{eqn:incrementdelta2}
\end{equation}
and 
$$
\lambda_{n,2}(e_{n,2}) = \frac{1}{2}(e_{n,2})^2\left(1- I_{\frac{3}{4}}\left(\frac{n-1}{2},\frac{1}{2}\right)\right).
$$
Note that Proposition 1 in \cite{hero2012hub} states that for any fixed $n$ and $\delta$, $\P(N_{V_{\delta}}^{(\Rma)}>0)$ converges to $\P(N_\delta^*>0)$ when $p\to\infty$, where $N_\delta^*$ is the Poisson random variable specified in the proposition. Since the increment distribution $\bm{\zeta_{n,2}} \neq \Dirac(1)$ by \eqref{eqn:incrementdelta2}, the limit of $\P(N_{V_{\delta}}^{(\Rma)}>0)$ is $\P(\CP(\lambda_{n,2}(e_{n,2}),\bm{\zeta}_{n,2})>0)$ and hence 
\cite[Proposition 1]{hero2012hub} is incorrect. 
{ However, as shown in Subsection \ref{sec:approximationstoalphaell} 
the result \cite[Proposition 1]{hero2012hub} is still useful when $n$ or $\delta$ is large.}   
Specifically, the distribution of 
each of the count variables $\{N_{E_\delta}^{(\Psima)},N_{\breve{V}_\delta}^{(\Psima)},N_{V_\delta}^{(\Psima)}:i\in\{\Rma,\Pma\}, \Psima \in\{\Rma,\Pma\} \}$ can be accurately approximated by a Poisson distribution but with an additional error term depending on $n$ and $\delta$. This error term is small 
for large $n$ or large $\delta$ (see Proposition \ref{prop:poiapplimit} and Proposition \ref{prop:poiappfinitep}).
\end{rem}

\section{Conclusions and discussions}
\label{sec:conclusions}

In this paper, we studied the number of highly connected vertices in both the empirical correlation graph and the empirical partial correlation graph by adopting unified framework. More specifically, we showed that the distributions of the number of hubs $N_{V_\delta}^{(\Psima)}$ or $N_{\breve{V}_\delta}^{(\Psima)}$ and the star subgraph counts $N_{E_\delta}^{(\Psima)}$ are close to common compound Poisson, when $p$ is finite and when $p$ approaches infinity. We also establish that their first and second moments are close to that of the compound Poisson distribution. The parameters in the compound Poisson distributions are characterized in terms of random geometric graphs and random pseudo geometric graphs. The parameters are also approximated by simple formulae, which implies that the approximating compound Poisson distributions can be further approximated by Poisson distributions for reasonably large sample size $n$ or a reasonably hub degree $\delta$. 

There are multiple avenues for future research. Numerical experiments suggest that the results in this paper hold beyond the $(\tau,\kappa)$ sparsity condition. A future line of work is to characterize the compound Poisson approximations for weaker sparsity conditions. 
 Generalization of the moment convergence results to beyond the second moments are also of interest.



\bibliographystyle{IEEEtran}
\bibliography{screening}

\newpage

\begin{appendix}


\subsection{List of symbols in the main text}
\label{sec:symbols}

\begin{itemize}
\setlength\itemsep{0mm}
 \item 
 $\|\cdot\|_2$ :  Euclidean norm. 
 
 \item
 $|\cdot|$ : the cardinality of a set.
 
 \item 
 $C$ and $c$ : denote generic positive universal constants that often differ from line to line.
 	$C$ and $c$ with subscripts are positive constants depending on the parameter in their subscripts and may differ from line to line.

 \item 
 $\Xma$: $n\times p$ matrix of observations. It has density $f_{\Xma}(\Xma)$ given in \eqref{eqn:ellden}.
 
 \item 
 $\Rma$ : $p\times p$ sample correlation matrix. Also denoted $\Psima^{(\mathbf R)}$.
 
 \item 
 $\Pma$ : $p\times p$ sample partial correlation matrix. Also denoted $\Psima^{(\mathbf P)}$.
 
 \item 
 $\rho \in [0, 1]$ : screening threshold applied to elements of matrices $\Rma$ or $\Pma$.

\item 
$\Phima^{(\Psima)}$ : the binary indicator matrix of the non-zero entries of the thresholded sample correlation or sample partial correlation matrix $\Psima$ with entries $\Phi_{ij}^{(\Psima)}= \Phi_{ij}^{(\Psima)}(\rho)=1(|\Psi_{ij}|\geq \rho)$, where $\Psima \in\{\Rma,\Pma\}$. 

\item
$\Gc_{\rho}(\Psima)$ : the graph associated with adjacency matrix $\Phima^{(\Psima)}$. Specifically, it is the empirical correlation graph when 
 $\Psima =\Rma$ or, equivalently, $\Psima=\Rma$ and it is the empirical partial correlation graph when  $\Psima =\Pma$ or, equivalently, $\Psima=\Pma$ (see Subsection \ref{Hubfra}).

   \item
  $N_{V_\delta}^{(\Psima)}$ : number of vertices of degree at least $\delta$ in $\Gc_{\rho}(\Psima)$.
 
  \item
  $N_{\VVV_\delta}^{(\Psima)}$ : number of vertices of degree exactly $\delta$ in $\Gc_{\rho}(\Psima)$. The quantity $\{N_{\VVV_\delta}^{(\Psima)}\}_{\Psima =0}^{n-1}$ is the empirical vertex degree distribution associated with the sample correlation graph ($\Psima =\Rma$) or partial correlation graph  ($\Psima =\Pma$). 

  \item
  $\Gamma_{\delta}$ : a star graph with $\delta$ edges.
 
  \item
  $N_{E_\delta}^{(\Psima)}$ : the number of subgraphs isomorphic to $\Gamma_{\delta}$ in $\Gc_{\rho}(\Psima)$ when $\delta \geq 2$.  $N_{E_\delta}^{(\Psima)}$ is defined as twice the number of edges in $\Gc_{\rho}(\Psima)$ when $\delta =1$.
 
  \item
  $\bar{N}_{\delta}$ :  a generic random variable denoting any of the quantities $\{N_{E_\delta}^{(\Psima)},N_{\breve{V}_\delta}^{(\Psima)},N_{V_\delta}^{(\Psima)}: \Psima \in \{\Rma,\Pma\}  \}$.
  \item
  $\Uma=[\uve_1,\ldots,\uve_p]$ : $(n-1) \times p$ matrix of correlation U-scores such that $\Uma^\top \Uma = \Rma $.  
  \item
  $\Yma=[\yve_1,\ldots,\yve_p]$ : $(n-1) \times p$ matrix of partial correlation Y-scores such that $\Yma^\top \Yma =\Pma$. 
  \item
  $r_{\rho}=\sqrt{2(1-\rho)}$ : spherical cap radius parameter in terms of $\rho$.
  \item
  $S^{n-2}$ : unit sphere in $\R^{n-1}$, {i.e. $S^{n-2}=\{\wwve\in \R^{n-1}: \|\wwve\|_2=1$\}. Its area is denoted by $\textbf{Area}(S^{n-2})$ or $|S^{n-2}|$.
  }
  \item
  $B^{n-2}$ : unit ball in $\R^{n-2}$, {i.e. $B^{n-2}=\{\wwve\in \R^{n-2}: \|\wwve\|_2\leq 1$\}. Its volume is denoted by $|B^{n-2}|$.}
  
  \item
  $\unif(S^{n-2})$ : uniform distribution on $S^{n-2}$.
  
 \item
  $\unif(B^{n-2})$ : uniform distribution in $B^{n-2}$.

  \item
  $P_n(r)$ : normalised area of a spherical cap on $S^{n-2}$ with radius $r$ (see \eqref{spharefor}). 
 
  \item
  $\textbf{Ge}\left(\{\vve_i  \}_{i=1}^\delta, r   \right)$ :  the geometric graph generated by a set of points $\{\vve_i  \}_{i=1}^\delta\subset \R^{n-2}$ with ball radius parameter $r$ (see Subsection \ref{Hubfra}).
  
  \item Universal vertex : a vertex of an undirected graph that is adjacent to all other vertices
  
  \item
  $\nmd\left(\{\vve_i  \}_{i=1}^\delta, r   \right)$ :
the number of universal vertices in $\textbf{Ge}\left(\{\vve_i  \}_{i=1}^\delta, r   \right)$. 
\item
$\textbf{PGe}\left(\{\vve_i  \}_{i=1}^m, r  ;  \mathscr{N} \right)$ : the pseudo geometric graph generated by $\{\vve_i  \}_{i=1}^m\subset \R^{\mathscr{N}}$ with radius paraneter $r$ (see Definition \ref{def:pge}).

\item
$\pnmd\left(\{\vve_i  \}_{i=1}^m, r  ;  \mathscr{N} \right)$ : the number of universal vertices in $\textbf{PGe}\left(\{\vve_i  \}_{i=1}^m, r  ;  \mathscr{N} \right)$.

\item
$\{\uve'_i\}_i$ : a sequence of i.i.d. random vectors drawn from $\unif(S^{n-2})$, the uniform distribution on $S^{n-2}$.

\item 
$\alpha(\ell,r_\rho)=\alpha_{n,\delta}(\ell,r_\rho)$: the conditional probability that there are $\ell$ universal vertices in $\textbf{PGe} \left(\{\uve'_i  \}_{i=1}^{\delta+1}, r_\rho   \right)$ (see \eqref{eqn:alphandeltarho}). 


\item
$\bm{\zeta}_{n,\delta,\rho}$ : the discrete distribution supported on $[\delta+1]$ with $\bm{\zeta}_{n,\delta,\rho}(\ell)\propto \alpha(\ell,r_\rho)/\ell$ (see \eqref{eqn:zetandelrhodef}).

\item 
$\CP(\lambda,\bm{\zeta})$ : 
the compound Poisson distribution, where the arrival rate $\lambda$ is the rate for the underlying Poisson random variable and the increment distribution $\bm{\zeta}$ is the distribution of each increment.

\item 
$\CP(\lambda_{p, n,\delta,\rho},\bm{\zeta}_{n,\delta,\rho})$ : 
the non-asymptotic compound Poisson distribution in Theorem \ref{thm:Poissonultrahigh}.\\

\item  

  $\overset{\Dc}{\to}$ : denotes convergence in distribution.

\item  
$\lambda_{p, n,\delta,\rho}$ : a quantity defined as $\binom{p}{1}\binom{p-1}{\delta}(2P_n(r_\rho))^\delta\sum_{\ell=1}^{\delta+1} \frac{\alpha(\ell,r_\rho)}{\ell}$. It is the arrival rate of $\CP(\lambda_{p, n,\delta,\rho},\bm{\zeta}_{n,\delta,\rho})$. 

\item
$\{\tilde{\uve}_i\}_i$ : a sequence of i.i.d. random vectors drawn from $\unif(B^{n-2})$, the uniform distribution in $B^{n-2}$.

\item 
$\alpha_\ell$ : the probability that there are exactly $\ell-1$ universal vertices in $\textbf{Ge}\left(\{\tilde{\uve_i}\}_{i=1}^\delta, 1   \right)$ (see \eqref{eqn:alphaelldef}). $\alpha_\ell=\lim_{\rho\to 1}\alpha(\ell,r_\rho)$ (see Lemma \ref{prop:rangeolimit} and the paragraph beneath it).


\item
$\bm{\zeta}_{n,\delta}$ : the discrete distribution supported on $[\delta+1]$ with $\bm{\zeta}_{n,\delta}(\ell)\propto \alpha_\ell/\ell$ (see \eqref{eqn:incrementsizelimitdef}). $\bm{\zeta}_{n,\delta}(\ell)=\lim_{\rho\to 1} \bm{\zeta}_{n,\delta,\rho}(\ell)$ (see Lemma \ref{lem:compoundpoissonlimit}).

\item 
$\lambda_{n,\delta}(e_{n,\delta})$ : a quantity defined as $\frac{1}{\delta !} \left( e_{n,\delta}\right)^\delta \sum_{\ell=1}^{\delta+1}\frac{\alpha_\ell}{\ell}$.  $\lambda_{n,\delta}(e_{n,\delta})=\lim_{p\rightarrow \infty} \lambda_{p,n,\delta,\rho}$ when  $\rho$ satisfies $a_n2^{\frac{n}{2}}p^{1+\frac{1}{\delta}}(1-\rho)^{\frac{n-2}{2}}\to e_{n,\delta}$, where $a_n=\frac{\Gamma((\nn-1)/2)}{(n-2)\sqrt{\pi}\Gamma((\nn-2)/2)}$ and $e_{n,\delta}$ is some positive constant (see Lemma \ref{lem:compoundpoissonlimit}).  

\item 
$\CP(\lambda_{n,\delta},\bm{\zeta}_{n,\delta})$ : 
the limiting compound Poisson distribution in Theorem \ref{cor:Poissonlimit}. Moreover, $\CP(\lambda_{p,n,\delta,\rho},\bm{\zeta}_{n,\delta,\rho}) \overset{\Dc}{\to} \CP(\lambda_{n,\delta}(e_{n,\delta}),\bm{\zeta}_{n,\delta})$, when $p\to \infty$ and $\rho$ satisfies  $a_n2^{\frac{n}{2}}p^{1+\frac{1}{\delta}}(1-\rho)^{\frac{n-2}{2}}\to e_{n,\delta}$. (see Lemma \ref{lem:compoundpoissonlimit})

\item
$\mu(\Ama)$ : normalized determinant of square matrix $\Ama$ (see Definition \ref{def:nordet}).

\item
$\mu_m(\Ama)$ : local normalized determinant of degree $m$ (see  the paragraph after Definition \ref{def:nordet}).

\item
$\mu_{n,m} (\Ama)$ : inverse local normalized determinant. It is powers of $\mu_m(\Ama)$ (see \eqref{eqn:munmA}).

\item 

$C_\delta^<$ : the index set defined as $\{\vec{i}=(i_0,i_1,\cdots,i_\delta)\in [p]^{\delta+1}: i_1<i_2<\cdots<i_\delta, \text{ and } i_\ell\not= i_0,\ \forall\ \ell\in [\delta] \}$. Each $\vec{i}\in C_\delta^<$ corresponds to a group of $\delta$ vertices indexed by $\{i_j\}_{j=1}^{\delta}$ and a center indexed by $i_0$. 


\item 

$\Phi_{\vec{i}}^{(\Rma)}$ : defined as $\prod_{j=1}^\delta \Phi_{i_0i_j}^{(\Rma)}$. 
It is the indicator random variable that vertex ${i_0}$ is adjacent to each vertex ${i_j}$ for $j\in [\delta]$ in $\Gc_\rho(\Rma)$. Equivalently, it is the indicator function of the event that there exist a star subgraph in $\Gc_{\rho}(\Rma)$ with center at $i_0$ and leaves $\{i_j\}_{j=1}^{\delta}$.

\item 

$S_{\vec{i}}$ : defined as $\left\{\vec{j}\in C_\delta^< \backslash \{\vec{i}\}: \bigcup_{\ell=0}^\delta \{j_\ell\}  = \bigcup_{\ell=0}^\delta \{i_\ell\}  \right\}$ for $\vec{i}\in C_\delta^<$. It is the set of indices in $C_\delta^<$ that share the same vertices with $\vec{i}$ but with center different from $i_0$.  


\item
$U_{\vec{i}}$ : defined as $\sum_{\vec{j}\in S_{\vec{i}}} \Phi_{\vec{j}}^{(\Rma)}$. It is the sum of highly dependent terms of $\Phi_{\vec{i}}^{(\Rma)}$. (see the discussions preceding and following \eqref{eqn:Uidef}).
\end{itemize}

\subsection{Controlling local normalized determinant by extreme eigenvalues}
\label{sec:localdeterminantvseigenvalues}
\begin{lem}     \label{lem:sufass}
Let $\Ama$ be a symmetric positive definite matrix. Consider a sequence of symmetric positive definite matrices $\Sigmama\in \R^{p\times p}$ with increasing dimension $p$.
    \begin{enumerate}[label=(\alph*)]
        \item \label{item:sufassnewa}
        $\mu_{n,m}(\Ama)$ is bounded by powers of the largest local condition number:
        $$
        \mu_{n,m}(\Ama) \leq  \begin{cases} \max\limits_{{\mathcal{I}}\subset [p]}\left(\frac{\lambda_{\text{max}}(\Ama_{\mathcal{I}})}{\lambda_{\text{min}}(\Ama_{\mathcal{I}})}\right)^{\frac{m(n-1)}{2}}, & \Ama   \text{ not diagonal, }  \\
        1, & \Ama \text{ diagonal. }
        \end{cases}
        $$
        \item \label{item:sufassa}
        $\mu_{n,m}(\Ama)$ is bounded by powers of the condition number:
        $$
        \mu_{n,m}(\Ama) \leq  \begin{cases} \left(\frac{\lambda_{\text{max}}(\Ama)}{\lambda_{\text{min}}(\Ama)}\right)^{\frac{m(n-1)}{2}}, &  \Ama \text{ not diagonal, }  \\
        1, & \Ama \text{ diagonal. }
        \end{cases}
        $$
        \item  If $\lambdamin\left(\Sigmama\right) \geq \underline{\lambda}$ and $\lambdamax\left(\Sigmama\right) \leq \overline{\lambda}$ for all $p$, then
        $$
        \mu_{n,m}(\Sigmama) \leq  \begin{cases} \left(\frac{\overline{\lambda}}{\underline{\lambda}}\right)^{\frac{m(n-1)}{2}}, &\Sigmama \text{ not diagonal, }  \\
        1, & \Sigmama \text{  diagonal. }
        \end{cases}
        $$
        \item  Let $M>0$  be a constant. Suppose 
        for all $p$, $\sup_{1\leq i\leq p} \Sigma_{ii}\leq M $. Moreover suppose $\lambda_{\text{min}}\left(\Sigmama\right) \geq \underline{\lambda}$ for all $p$. Then
        $$
        \mu_{n,m}(\Sigmama) \leq  \begin{cases} \left(\frac{Mm}{\underline{\lambda}}\right)^{\frac{m(n-1)}{2}}, &\Sigmama \text{ not diagonal, }  \\
        1, & \Sigmama \text{  diagonal. }
        \end{cases}
        $$
    \end{enumerate}
\end{lem}
\begin{proof}
(a)
It follows directly by definition of $\mu_{n,m}(A)$ and $\mu(A_{\mathcal{I}})\geq \left(\frac{\lambda_{\text{min}}(A_{\mathcal{I}})}{\lambda_{\text{max}}(A_{\mathcal{I}})}\right)^{m}$.\\
    
    \noindent (b)
    Since $\mu_{n,m}(A)$ is increasing in $m$ as discussed after Example \ref{exa:normaliseddeterminant}, $\mu_{n,m}(A)\leq \mu_{n,p}(A)$. The proof is then complete by applying \ref{item:sufassnewa} to $\mu_{n,p}(A)$.\\
        
        \noindent (c)
        It follows directly from \ref{item:sufassa}.\\

        \noindent (d)
        Let $\mathcal{I}\subset [p]$ with $|\mathcal{I}|=m$. Since $\Sigmama_{\mathcal{I}}$ is symmetric positive definite, $|(\Sigmama_{\mathcal{I}})_{ij}|\leq \sqrt{(\Sigmama_{\mathcal{I}})_{ii}(\Sigmama_{\mathcal{I}})_{jj} }\leq M$ and thus $\lambda_{\text{max}}(\Sigmama_{\mathcal{I}}) = \|\Sigmama_{\mathcal{I}}\|_2 \leq \|\Sigmama_{\mathcal{I}}\|_F \leq Mm  $. By the interlacing property (see Theorem 8.1.7 in \cite{golub2012matrix}), $\lambdamin(\Sigmama_{\mathcal{I}}))\geq \lambdamin(\Sigmama).$
        The proof is then completed by applying part \ref{item:sufassa}.
\end{proof}

\subsection{Proofs in Subsection \ref{sec:scorerep}}

\subsubsection{Proofs of Lemma \ref{lem:Binv} and Lemma \ref{lem:Uscoresdistribution}}
\label{sec:Uscoresproof}

\begin{proof}[Proof of Lemma \ref{lem:Binv}] 
It suffices to prove that $\Uma$ has rank $n-1$ $a.s.$ By Lemma \ref{lem:Uscoresdistribution} \ref{item:Uscoresdistributiona}, it suffices to prove that $\left[\frac{\tilde{\xve}_1}{\|\tilde{\xve}_1\|_2},\ldots,\frac{\tilde{\xve}_p}{\|\tilde{\xve}_p\|_2}\right]$ has rank $n-1$ $a.s.$ Note that $\tilde{\Xma}$ is of rank $n-1$  $a.s.$, since it has a density with respect to Lebesgue measure on $\R^{(n-1)p}$ and $p\geq n$. Then $\left[\frac{\tilde{\xve}_1}{\|\tilde{\xve}_1\|_2},\ldots,\frac{\tilde{\xve}_p}{\|\tilde{\xve}_p\|_2}\right]$ have rank $\nn-1$ $a.s.$ as well since it is obtained by normalizing the columns of $\tilde{\Xma}$. 
\end{proof}

\begin{proof}[Proof of Lemma \ref{lem:Uscoresdistribution}]
(a)
Consider $\hat{\xve}^{(i)}\overset{\text{i.i.d.}}{\sim}\Nc(\bm{0},\Sigmama)$. Then $\hat{\Xma} := [\hat{\xve}^{(1)},\ldots,\hat{\xve}^{(n)}]^\top = [\hat{\xve}_1,\ldots,\hat{\xve}_p] \in \R^{n\times p}$ has the distribution $\mathcal{VE}(\bm{0},\Sigmama,g_0)$ where $g_0(w)=(2\pi)^{-\frac{\nn p}{2}}\exp(-\frac{1}{2}w)$. 
The sample mean $\bar{\hat{\xve}}$ of $\hat{\Xma}$ is given as a row vector
$$
\bar{\hat{\xve}}=\frac{1}{\nn}\Sum_{i=1}^{\nn}\hat{\xve}^{(i)} = \frac{1}{\nn} \hat{\Xma}^\top\bm{1}=(\bar{\hat{x}}_1,\ldots,\bar{\hat{x}}_p),
$$
where $\bar{\hat{x}}_i:=n^{-1}\sum_{j=1}^p {\hat{X}_{ij}}$. The sample covariance matrix $\hat{\Sma}$ of $\hat{\Xma}$ is
 \begin{equation*}
 \hat{\Sma}=\frac{1}{\nn-1}\sum_{i=1}^{\nn} (\hat{\xve}^{(i)}-\bar{\hat{\xve}})(\hat{\xve}^{(i)}-\bar{\hat{\xve}})^\top
 =\frac{1}{\nn-1}\hat{\Xma}^\top \left(\Ima-\frac{1}{n}\bm{1}\bm{1}^\top\right) \hat{\Xma}.
 \end{equation*}
The Z-scores of $\hat{\Xma}$ is
$$
\hat{\zve}_i := \frac{\hat{\xve}_i-\bm{1}\bar{\hat{x}}_i}{\sqrt{\hat{S}_{ii}(n-1)}}= \frac{\left(\Ima-\frac{1}{n}\bm{1}\bm{1}^\top\right) \hat{\xve}_i}{\sqrt{\hat{S}_{ii}(n-1)}} .
$$
Recall $\Hma_{2:n}\in \R^{n\times (n-1)}$ is a matrix such that
 $\nn \times (\nn-1)$ matrix $\Hma = [\nn^{-\frac{1}{2}}\bm{1}, \Hma_{2:\nn}]^\top$ is orthogonal, which satisfies $\Hma_{2:n} \Hma_{2:n}^\top = \Ima -\frac{1}{n}\bm{1}\bm{1}^T$. The U-scores of $\hat{\Xma}$ is denoted $\hat{\Uma}=[\hat{\uve}_1,\ldots,\hat{\uve}_p]$ with  $\hat{\uve}_i=\Hma_{2:n}^\top z_i$ for $i\in [p]$.

By the theorem in Section 6 of \cite{anderson1992nonnormal}, the distribution of $\Uma$ is invariant to $g$ and $\muve$. In particular, $\Uma$ has the same distribution as  $\hat{\Uma}$. It suffices to study the distribution of $\hat{\Uma}$.

 The proof of Theorem 3.3.2 in \cite{anderson2003introduction} establishes that 
\begin{equation*}
\check{\Xma} = \Hma_{2:\nn}^\top \hat{\Xma} \in \R ^{(\nn-1) \times p}    
\end{equation*}
has i.i.d. rows $\{\check{\xve}^{(i)}\}_{i=1}^{\nn-1}$ distributed as $\Nc(\bm{0},\Sigmama)$, and
\begin{equation*}
\hat{\Sma} = \frac{1}{\nn-1}\hat{\Xma}^\top \left(\Ima-\frac{1}{n}\bm{1}\bm{1}^\top\right) \hat{\Xma} =\frac{1}{\nn-1} \left(\Hma_{2:n}^\top \hat{\Xma}\right)^\top \Hma_{2:n}^\top \hat{\Xma}= \frac{1}{\nn-1} \check{\Xma}^\top \check{\Xma}. 
\end{equation*}
Then $\hat{S}_{ii} = \|\check{\xve}_i\|^2/(n-1) $ where $\check{\xve}_i$ is the $i$-th column of $\check{\Xma}$. Thus
$$
\hat{\uve}_i = \Hma_{2:n}^\top \frac{\left(\Ima-\frac{1}{n}\bm{1}\bm{1}^\top\right) \hat{\xve}_i}{\sqrt{\hat{S}_{ii}(n-1)}} = \frac{\Hma_{2:n}^\top \Hma_{2:n} \Hma_{2:n}^\top \hat{\xve}_i }{\sqrt{\|\check{\xve}_i\|^2}}=\frac{\check{\xve}_i}{\|\check{\xve}_i\|_2}.
$$
The proof is then complete by observing that $\check{\Xma}$ and $\tilde{\Xma}$ has the same distribution.\\

\noindent
(b)
Part \ref{item:Uscoresdistributionb} directly follows from part \ref{item:Uscoresdistributiona}.
\end{proof}

\subsubsection{Proof of Lemma \ref{lem:densitybou}} \label{sec:proofdensitybou}


Part \ref{item:densityboua} of Lemma \ref{lem:densitybou} follows from $\mu_m(\Sigmama_{\mathcal{J}})=\mu(\Sigmama_{\mathcal{J}})\geq \mu_m(\Sigmama)$. Part \ref{item:densitybouc} immediately follows from part \ref{item:densityboub} and in the remainder of this subsection we give the proof of part \ref{item:densityboub}.

It suffices to prove the following statement:
$$
f_{\uve_{j_1},\uve_{j_2},\cdots,\uve_{j_m}}(\vve_1,\vve_2,\cdots,\vve_m) \leq   \mu_{n,m}(\Sigmama_{\mathcal{J}}), \quad \forall \  \vve_i\in S^{n-2}, \forall\ i\in [m].
$$
For notation convenience, we only present the proof for $m=p$ and $\mathcal{J}=[p]$ since the proof of the general $m$ and $\mathcal{J}$ follows the same proof procedure. When $m=p$ and $\mathcal{J}=[p]$, the statement of Lemma \ref{lem:densitybou} \ref{item:densityboub} becomes:

The joint density of columns of U-scores w.r.t. $\otimes^p \sigma^{n-1}$ is upper bounded  by $\mu_{n,p}(\Sigmama)$: \begin{equation}
f_{\uve_1,\uve_2,\cdots,\uve_p}(\vve_1,\vve_2,\cdots,\vve_p)\leq \mu_{n,p}(\Sigmama),  \quad \forall \  \vve_i\in S^{n-2}, \forall\ i\in [p]. \label{eqn:jointuuppbou}
\end{equation}

\begin{proof}[Proof of \eqref{eqn:jointuuppbou}]
		Recall $\{\tilde{\xve}^{(i)}\}_{i=1}^{\nn-1}\subset \R^p$, the rows of $\tilde{\Xma}$, are i.i.d. copy of $\Nc(\bm{0},\Sigmama)$; $\{\tilde{\xve}_i\}_{i=1}^{p}$ are the columns of $\tilde{\Xma}$; $\uve_i\vcentcolon=\frac{\tilde{\xve}_i}{\|\tilde{\xve}_i\|_2}\in \R^{\nn-1}$ has distribution $\unif(S^{\nn-2})$ for $i\in [p]$. 
		
		When $\Sigmama$ is symmetric positive definite and diagonal, $\{\tilde{\xve}_i\}_{i=1}^{p}$ are independent, which implies that $\{\uve_i\}_{i=1}^{p}$ are independent. Thus in this case, the joint density of columns of U-scores w.r.t. $\otimes^p \sigma^{n-1}$ is $1$.
		
		Consider general symmetric positive definite $\Sigmama$. The probability density of $\tilde{\Xma}$ w.r.t. the Lebesgue measure on $\R^{(n-1)p}$ is 
		$$
		f_{\tilde{\Xma}}(\tilde{\Xma}) = A\exp\left( -\frac{1}{2}\sum_{j=1}^{n-1}
		\left(\tilde{\xve}^{(j)}\right)^\top\Sigmama^{-1}\tilde{\xve}^{(j)}  \right),
		$$ 
		where the constant $A=\text{det}(\Sigmama)^{-\frac{n-1}{2}}(2\pi)^{-\frac{(n-1)p}{2}}$.
		Use the spherical transform for each column $\tilde{\xve}_i = \left(\tilde{\Xma}_{ji}:1\leq j \leq n-1\right)^\top$:
		\begin{equation*}
		\begin{cases}
		\tilde{\Xma}_{1i}=R_i\cos(\theta_{1i}),\\
		\tilde{\Xma}_{2i}=R_i\sin(\theta_{1i})\cos(\theta_{2i}),\\ \vdots & \text{\ \ for } 1\leq i \leq p,\\
		\tilde{\Xma}_{(n-2)i}=R_i\sin(\theta_{1i})\sin(\theta_{2i})\cdots\sin(\theta_{(n-3)i})\cos(\theta_{(n-2)i}),\\ \tilde{\Xma}_{(n-1)i}=R_i\sin(\theta_{1i})\sin(\theta_{2i})\cdots\sin(\theta_{(n-3)i})\sin(\theta_{(n-2)i}), 
		\end{cases}
		\end{equation*}
		where for each $i\in [p]$: $R_i\geq 0, \theta_{ji}\in[0,\pi]$ for $1\leq j\leq n-3$ and $\theta_{(n-2)i}\in[0,2\pi)$.
		
		Denote $\bm{R}=(R_i:1\leq i\leq p)$ and $\bm{\Theta}=(\theta_{ji}: 1\leq i\leq p, 1\leq j\leq (n-2))$.
		Then the joint density of $(\bm{R},\bm{\Theta})$ is:
		\begin{align*}
		&f_{\bm{R},\bm{\Theta}}(\bm{R},\bm{\Theta})\\
		= &
		A\exp\left( -\frac{1}{2}\sum_{j=1}^{n-1}
		\left(\bm{h}^{(j)}\right)^\top\Sigmama^{-1}\bm{h}^{(j)}  \right) \prod_{i=1}^{p} \left(R_i^{n-2}\prod_{j=1}^{n-2}\sin^{n-2-j}(\theta_{ji})\right),
		\end{align*}
		where 
		$$
		\bm{h}^{(j)}=\left(R_i\cos(\theta_{ji})\prod\limits_{q=1}^{j-1}\sin(\theta_{qi}):1\leq i\leq p\right)^\top\in \R^{p} \text{ for } 1\leq j\leq n-2
		$$ 
		and  
		$$
		\bm{h}^{(n-1)}=\left(R_i\prod\limits_{q=1}^{n-2}\sin(\theta_{qi}):1\leq i\leq p\right)^\top\in \R^p. 
		$$
		
		Then the density of $\bm{\Theta}$ is:
		\begin{align*}
		&f_{\bm{\Theta}}(\bm{\Theta}) \\
		= &
		A \left(\prod_{i=1}^{p} \prod_{j=1}^{n-2}\sin^{n-2-j}(\theta_{ji}) \right) \int_{[0,\infty)^p}\exp\left(-\frac{1}{2}\sum_{j=1}^{n-1}
		\left(\bm{h}^{(j)}\right)^\top\Sigmama^{-1}\bm{h}^{(j)} \right) \prod_{i=1}^{p} R_i^{n-2} dR_i\\
		\leq & 
		A \left(\prod_{i=1}^{p} \prod_{j=1}^{n-2}\sin^{n-2-j}(\theta_{ji}) \right) \int_{[0,\infty)^p} \exp\left( -\frac{1}{2}\lambda_{\text{min}}(\Sigmama^{-1})\sum_{j=1}^{n-1}
		\|\bm{h}^{(j)}\|_2^2  \right) \prod_{i=1}^{p} R_i^{n-2}dR_i\\
		=& 
		A \left(\prod_{i=1}^{p} \prod_{j=1}^{n-2}\sin^{n-2-j}(\theta_{ji}) \right) \int_{[0,\infty)^p} \exp\left( -\frac{1}{2}[\lambda_{\text{max}}(\Sigmama)]^{-1}\sum_{i=1}^{p}
		R_i^2 \right) \prod_{i=1}^{p} R_i^{n-2} d R_i \\
		=& 
		A \left(\prod_{i=1}^{p} \prod_{j=1}^{n-2}\sin^{n-2-j}(\theta_{ji}) \right)  \left(\int_{[0,\infty)}\exp\left( -\frac{1}{2}[\lambda_{\text{max}}(\Sigmama)]^{-1}
		R_1^2  \right)  R_1^{n-2} dR_1\right)^p \\
		\overset{(m)}{=}& 
		A \left(\prod_{i=1}^{p} \prod_{j=1}^{n-2}\sin^{n-2-j}(\theta_{ji}) \right) \left( \left[\lambda_{\text{max}}(\Sigmama)\right]^{\frac{n-1}{2}} \Gamma\left(\frac{n-1}{2}\right)2^{\frac{n-3}{2}}
		\right)^p\\
		\overset{(mm)}{=}&\left[\frac{\left(\lambda_{\text{max}}(\Sigmama)\right)^p}{\text{det}(\Sigmama)}\right]^\frac{n-1}{2} \frac{1}{(\text{Area}(S^{n-2}))^p} \prod_{i=1}^{p} \prod_{j=1}^{n-2}\sin^{n-2-j}(\theta_{ji})
		\end{align*}
		where equality $(m)$ follows from the integration of Chi distribution with degree $n-1$, and equality $(mm)$ follows from $|S^{(n-2)}|=2\pi^{\frac{n-1}{2}}/\Gamma((n-1)/2)$. The proof is complete by noticing $f_{\bm{\Theta}}(\bm{\Theta})$ is joint density of columns of U-scores expressed in spherical coordinate
		and $\frac{1}{(\text{Area}(S^{n-2}))^p} \prod_{i=1}^{p} \prod_{j=1}^{n-2}\sin^{n-2-j}(\theta_{ji})$ is the joint distribution of $p$ independent $\unif(S^{n-2})$ expressed as spherical coordinates.		 
\end{proof}


\subsection{\ONE{An informal derivation of Proposition \ref{prop:edgecor}}}
\label{sec:derivationprop:edgecor}

In this subsection we provide motivations to the following two questions: 1) why is the distribution of the star subgraph counts $N_{E_\delta}^{(\Rma)}$ approximately compound Poisson? 2) what are the associated parameters of the compound Poisson approximation to $N_{E_\delta}^{(\Rma)}$? 

We first introduce some notation to answer the first question. 
Denote 
\begin{equation}
C_\delta^{<}:=\{\vec{i}=(i_0,i_1,\cdots,i_\delta)\in [p]^{\delta+1}: i_1<i_2<\cdots<i_\delta, \text{ and } i_\ell\not= i_0,\ \forall\ \ell\in [\delta] \}. \label{eqn:Cdeltadef}
\end{equation}
Each $\vec{i}\in C_\delta^<$ corresponds to a group of $\delta$ vertices indexed by $\{i_j\}_{j=1}^{\delta}$ and a center indexed by $i_0$. For $\vec{i}\in C_\delta^<$, denote by
\begin{equation}
\Phi_{\vec{i}}^{(\Rma)}= \prod_{j=1}^\delta \Phi_{i_0i_j}^{(\Rma)} = 1\left(\bigcap_{j=1}^{\delta}\{\textbf{dist}(\uve_{i_0},\uve_{i_j})\leq r_\rho\}\right) \label{eqn:NEdeltaRdef}
\end{equation}
the indicator function of the event that vertex $i_0$ is connected to each vertex $i_j$ for $j\in [\delta]$ in the empirical correlation graph. Equivalently, when $\delta\geq 2$, $\Phi_{\vec{i}}^{(\Rma)}$ is the indicator function of the event that there exists a star subgraph in $\Gc_{\rho}(\Rma)$ with center at $i_0$ and with leaves $\{i_j\}_{j=1}^{\delta}$.
Then by definition
\begin{equation}
\Nedged^{(\Rma)}  =\sum_{\vec{i}\in C_\delta^<}\ \  \Phi_{\vec{i}}^{(\Rma)}. \label{eqn:Nedgedef1}
\end{equation}
\index{$\Nedged^{(\Rma)}$}

If $\Phi_{\vec{i}}^{(\Rma)}$ for different $\vec{i}\in C^<_{\delta}$ are independent or weakly dependent, then the distribution of $\Nedged^{(\Rma)}$, as a sum of independent or weakly dependent indicator random variables, might be expected to be approximately Poisson. This, however, is not the case due to high dependency among many terms in the summation. Specifically, for any $\vec{i}\in C_\delta^<$, $\Phi_{\vec{i}}^{(\Rma)}$ is highly dependent on $\Phi_{\vec{j}}^{(\Rma)}$ for any $\vec{j} \in S_{\vec{i}}$ where 
\begin{align}
S_{\vec{i}} & :=\left\{\vec{j}\in C_\delta^< \backslash \{\vec{i}\}: \bigcup_{\ell=0}^\delta \{j_\ell\}  = \bigcup_{\ell=0}^\delta \{i_\ell\}  \right\}.
\label{eqn:Sidef}
\end{align}
$S_{\vec{i}}$ is the set of indexes sharing the same vertices with $\vec{i}$ but with center different from $i_0$ and thus $|S_{\vec{i}}| = \delta$.
Indeed, provided $\Phi_{\vec{i}}^{(\Rma)}=1$, which is equivalent to $\textbf{dist}(\uve_{i_0},\uve_{i_j})\leq r_\rho$ for $\forall j\in [\delta]$,
$$
\textbf{dist}(\uve_{i_1},\uve_{i_j})\leq \textbf{dist}(\uve_{i_1},\uve_{i_0})+\textbf{dist}(\uve_{i_0},\uve_{i_j})\leq 2r_\rho, \text{ for } \forall 2\leq j \leq \delta.
$$ 
That is, $\{\uve_{i_j}\}_{j=2}^{\delta+1}$ are all close to $\uve_{i_1}$ and hence it is likely there are edges connecting them. In other words it is likely
 for  $\vec{i'}=(i_1,i_0,\cdots,i_\delta)$, $\Phi_{\vec{i'}}^{(\Rma)}=1$, where we without loss of generality assume $i_0<i_j$ for $2\leq j \leq \delta$. 
 Let 
\begin{equation}
U_{\vec{i}} = \sum_{\vec{j}\in S_{\vec{i}}} \Phi_{\vec{j}}^{(\Rma)} \label{eqn:Uidef}
\end{equation}
be the sum of highly dependent terms of $\Phi_{\vec{i}}^{(\Rma)}$. In summary, if there is an increment for $\Nedged^{(\Rma)}$, say $\Phi_{\vec{i}}^{(\Rma)}=1$, there is a certain probability that $U_{\vec{i}}$ is greater than $0$ due to the high dependence, causing each increment of $\Nedged^{(\Rma)}$ to be greater than $1$ with a certain probability. This is typical behavior for a compound Poisson random variable.

To answer the second question, we next informally derive the parameters of the compound Poisson approximation to $N_{E_\delta}^{(\Rma)}$ in the special case when $\Sigmama$ is diagonal. We require some additional definitions. Let $\left[\vec{i}\right]=\{i_0,i_1,\cdots,i_\delta\}$ be the unordered set of indexes, which is referred as \emph{index group}, of any $\vec{i}\in C_\delta^<$ and define $\left[C_\delta^{<}\right]:=\left\{\left[\vec{i}\right]: \vec{i}\in C_\delta^< \right\}$. It follows that $\left|\left[C_\delta^{<}\right]\right| =\binom{p}{\delta+1}$, 
where $|\cdot|$ is the cardinality of a set. For a given group of $\delta+1$ indexes $\left[\vec{i}\right]$, $\Phi_{\vec{i}}^{(\Rma)} + U_{\vec{i}}$ is the increment associated to this group and its value is between $0$ and $\delta+1$. 
In the following informal argument we assume that the event that two different index groups both have non-zero increments has probability zero. This assumption will be verified in the proof of Proposition \ref{prop:edgecor}. Consequently the probability of increment size $\ell$ for $\ell\geq 1$ is proportional to the expectation of the fraction of the number of groups with increment $\ell$:
\begin{align*}
&\E \frac{1}{\left|\left[C_\delta^{<}\right]\right|} \sum_{\left[\overset{}{\vec{i}}\right]\in \left[C_\delta^{<}\right] } 1\left( \Phi_{\vec{i}}^{(\Rma)} + U_{\vec{i}} =  \ell \right) \\
=& \frac{1}{\left|\left[C_\delta^{<}\right]\right|} \frac{1}{\ell} \E \sum_{\left[\overset{}{\vec{i}}\right]\in \left[C_\delta^{<}\right] } \left(\Phi_{\vec{i}}^{(\Rma)} + U_{\vec{i}}\right) 1\left( \Phi_{\vec{i}}^{(\Rma)} + U_{\vec{i}} = \ell \right) \\
=& \frac{1}{\left|\left[C_\delta^{<}\right]\right|} \frac{1}{\ell} \E \sum_{\vec{i}\in C_\delta^{<} } \Phi_{\vec{i}}^{(\Rma)} 1\left( \Phi_{\vec{i}}^{(\Rma)} + U_{\vec{i}} = \ell \right) \\
=& \frac{1}{\left|\left[C_\delta^{<}\right]\right|}  \frac{1}{\ell}  \sum_{\vec{i}\in C_\delta^{<} } \P\left(\Phi_{\vec{i}}^{(\Rma)}=1\right) \P \left( \Phi_{\vec{i}}^{(\Rma)} + U_{\vec{i}} = \ell | \Phi_{\vec{i}}^{(\Rma)}=1 \right). \numberthis \label{eqn:jumprateint} 
\end{align*}
Since $\Sigmama$ is assumed to be diagonal, $\{\uve_i\}_{i=1}^p$ are i.i.d. $\unif(S^{n-2})$ by Lemma \ref{lem:Uscoresdistribution} \ref{item:Uscoresdistributionb} and hence \eqref{eqn:jumprateint} becomes
\begin{align*}
\E \frac{1}{\left|\left[C_\delta^{<}\right]\right|} \sum_{\left[\vec{i}\right]\in \left[C_\delta^{<}\right] } 1\left( \Phi_{\vec{i}}^{(\Rma)} + U_{\vec{i}} = \ell \right) = & \frac{1}{\left|\left[C_\delta^{<}\right]\right|}
\frac{1}{\ell}  \left|C_\delta^<\right| (2P_n(r_\rho))^{\delta} \alpha(\ell,r_\rho)  \\
= & 
\frac{\delta+1}{\ell}   (2P_n(r_\rho))^{\delta} \alpha(\ell,r_\rho), 
\end{align*}
where $\alpha(\ell,r_\rho)$ is defined in \eqref{eqn:alphandeltarho}
and $\P\left(\Phi_{\vec{i}}^{(\Rma)}=1\right)=(2P_n(r_\rho))^{\delta}$ by conditioning on $\uve_{i_0}$. As a consequence, the probability of increment size $\ell$ for $\ell\geq 1$ is:  
\begin{align*}
\frac{\delta+1}{\ell}   (2P_n(r_\rho))^{\delta} \alpha(\ell,r_\rho)/ \sum_{\ell =1}^{\delta+1}\left(\frac{\delta+1}{\ell}   (2P_n(r_\rho))^{\delta} \alpha(\ell,r_\rho)\right) 
= &\bm{\zeta}_{n,\delta,\rho}(\ell),
\end{align*}
where the equation follows from 
\eqref{eqn:zetandelrhodef}. This argument implies that $\bm{\zeta}_{n,\delta,\rho}$ is the increment distribution of the compound Poisson approximation. 

Since the mean of a compound Poisson distribution is the product of the arrival rate and the mean of the increment distribution, the arrival rate $\lambda$ of the compound Poisson approximation satisfies the following mean constraint:
$$
  \lambda \E \bm{\zeta}_{n,\delta,\rho} = \E \Nedged^{(\Rma)},
$$
where $\E \bm{\zeta}_{n,\delta,\rho}$ is the mean of $\bm{\zeta}_{n,\delta,\rho}$. One can easily verify $\E \bm{\zeta}_{n,\delta,\rho}= 1/\sum_{\ell=1}^{\delta+1} \left(\alpha(\ell,r_\rho)/\ell \right) $ and $\E \Nedged^{(\Rma)} = \binom{p}{1}\binom{p-1}{\delta}\E\Phi_{\vec{i}}^{(\Rma)} = \binom{p}{1}\binom{p-1}{\delta}(2P_n(r_\rho))^\delta $ since $\Sigmama$ is diagonal. Hence the arrival rate for the compound Poisson is
\begin{align}
\lambda = \binom{p}{1}\binom{p-1}{\delta}(2P_n(r_\rho))^\delta\sum_{\ell=1}^{\delta+1} \frac{\alpha(\ell,r_\rho)}{\ell}= \lambda_{p,n,\delta,\rho}.  
\end{align}
This informal argument motivates the compound Poisson approximation $\CP(\lambda_{p,n,\delta,\rho},\bm{\zeta}_{n,\delta,\rho})$ to $\Nedged^{(\Rma)}$ when $\Sigmama$ is diagonal.

\subsection{Proof of Proposition \ref{prop:edgecor}}
\label{sec:proofprop:edgecor}

\subsubsection{Auxiliary lemmas for Proposition \ref{prop:edgecor}}
Recall for any $\delta\geq 1$, $C_\delta^<$ is defined in \eqref{eqn:Cdeltadef}. For $\vec{i}\in C_\ell^{<}$, define a symmetric positive definite matrix $\Sigmama_{\vec{i}}\in \R^{(\ell+1)\times (\ell+1)}$ to be the submatrix of $\Sigmama$, consisting of rows and columns $\Sigmama$ indexed by the ordered components $(i_0,i_1,\ldots,i_{\ell})$ of $\vec{i}$. Let $\left[\vec{i}\right]=\{i_0,i_1,\cdots,i_\ell\}$ be the unordered set of indices of any $\vec{i}\in C_\ell^<$. Then $\Sigmama_{\vec{i}}\in \Sigmama_{[\vec{i}]}$
and $\mu_{n,\ell+1}\left(\Sigmama_{\vec{i}}\right)=\mu_{n,\ell+1}\left(\Sigmama_{[\vec{i}]}\right)$, where $\Sigmama_{[\vec{i}]}$ and $\mu_{n,\ell+1}\left(\Sigmama_{[\vec{i}]}\right)$ are defined in the paragraph after Definition \ref{def:nordet}.

\begin{lem} \label{item:prod}
Suppose $\Xma \sim \mathcal{VE}(\bm{\mu},\Sigmama,g)$. Let $\ell\in[p-1]$. Consider $\vec{i}=(i_0,i_1,\cdots,i_\ell)\in C^<_\ell$. 
		\begin{equation*}
		\E \prod_{q=1}^\ell\Phi_{i_0i_q}^{(\Rma)} \leq \mu_{n,\ell+1}\left(\Sigmama_{\vec{i}}\right) \left( 2P_\nn(r_\rho)\right)^\ell. \end{equation*}
Moreover, when $\Sigmama_{\vec{i}}$ is diagonal, in the last expression the equality holds and $\mu_{n,\ell+1}\left(\Sigmama_{\vec{i}}\right)=1$. 
\end{lem}
\begin{proof}
		$ \prod\limits_{q=1}^\ell\Phi_{i_0i_q}^{(\Rma)}$ is a nonnegative Borel Measurable function of $\uve_j$ for $j\in \left[\vec{i}\right]$. By Lemma \ref{lem:densitybou} \ref{item:densitybouc}, it suffices to show 
		$$
		\E \prod_{q=1}^\ell\Phi_{i_0i_q}^{(\Rma)} \leq \left( 2P_\nn(r)\right)^\ell
		$$
		for the case $\uve_j$ for $j\in \left[\vec{i}\right]$ are $\ell+1$ independent $\unif(S^{n-2})$. The last inequality holds with equality, which follows from that the terms in the product on the left hand side are independent conditioned on $\uve_i$. 
\end{proof}

 Lemma \ref{item:prod} suggests differentiating whether $\Sigmama_{\vec{i}}$ is diagonal or not since $\mu_{n,\ell+1}\left(\Sigmama_{\vec{i}}\right)=1$ when $\Sigmama_{\vec{i}}$ is diagonal. The next lemma establishes an upper bound on the number of $\vec{i}\in C_\delta^<$ such that $\Sigmama_{\vec{i}}$ is not diagonal.

\begin{lem} 
	\label{lem:numofnondia}
	Let $\Sigmama$ be row-$\kappa$ sparse. Let $\delta\in [p-1]$. Then
	$$
	\sum_{\substack{\vec{i}\in C_\delta^{<}\\ \Sigmama_{\vec{i}} \text{ not diagonal}}}\ 1 \leq 	\frac{\delta(\delta+1)}{2}(\kappa-1)\binom{p}{\delta}\leq 
	\frac{(\delta+1)}{2((\delta-1)!)}   p^{\delta} (\kappa-1).
	$$	
\end{lem}
\begin{proof}
Note that 
$$
\sum_{\substack{\vec{i}\in C_\delta^{<}\\ \Sigmama_{\vec{i}} \text{ diagonal}}}\ 1 \geq \frac{1}{\delta!}p(p-\kappa)\ldots(p-\delta\kappa),
$$
where the $\frac{1}{\delta!}$ is due to in our definition $\vec{i}$ the index $i_1<\ldots <i_{\delta}$ are sorted. Then 
\begin{align*}
\sum_{\substack{\vec{i}\in C_\delta^{<}\\ \Sigmama_{\vec{i}} \text{ not diagonal}}}\ 1 \leq  \binom{p}{1}\binom{p-1}{\delta} -
\frac{1}{\delta!}p\prod_{\ell=1}^{\delta}(p-\ell\kappa)\leq \frac{\delta(\delta+1)}{2}(\kappa-1)\binom{p}{\delta},
\end{align*}
where the last inequality follows from Lemma \ref{basine} \ref{item:basinei}.
\end{proof}

Note $\kappa=1$,  the Lemma \ref{lem:numofnondia} shows $\sum\limits_{\substack{\vec{i}\in C_\delta^{<}\\ \Sigmama_{\vec{i}} \text{ not diagonal}}}\ 1 = 0$, which means $\Sigmama$ is diagonal matrix. Next we present a lemma to bound $
\sum_{\vec{i}\in C_\ell^<}\mu_{n,\ell+1}\left(\Sigmama_{\vec{i}}\right)
$.

\begin{lem} \label{lem:summuupp}
	$$
	\sum_{\vec{i}\in C_\ell^<}\mu_{n,\ell+1}\left(\Sigmama_{\vec{i}}\right) \leq \frac{p^{\ell+1}}{\ell!}\left(1+\ell^2\mu_{n,\ell+1}(\Sigmama)\frac{\kappa-1}{p}\right).
	$$
\end{lem}
\begin{proof}
		\begin{align*}
	\sum_{\vec{i}\in C_\ell^<}\mu_{n,\ell+1}(\Sigmama_{\vec{i}}) =& \sum_{\substack{\vec{i}\in C_\ell^<\\ \Sigmama_{\vec{i} }\text{ diagonal }}}1+ \sum_{\substack{\vec{i}\in C_\ell^<\\ \Sigmama_{\vec{i} } \text{ not diagonal }}} \mu_{n,\ell+1}(\Sigmama_{\vec{i}})\\
	\leq & \binom{p}{1}\binom{p-1}{\ell} + \mu_{n,\ell+1}(\Sigmama) \sum_{\substack{\vec{i}\in C_\ell^<\\ \Sigmama_{\vec{i} } \text{ not diagonal }}}1 \\
	\leq & \frac{p^{\ell+1}}{\ell!}\left(1+\ell^2\mu_{n,\ell+1}(\Sigmama)\frac{\kappa-1}{p}\right),
	\end{align*}	
	where the first inequality follows from the Lemma \ref{lem:densitybou} \ref{item:densityboua}, and the second inequality follows from Lemma \ref{lem:numofnondia}.
\end{proof}

\begin{lem}\label{lem:eleine}
	Let $\Xma\sim \mathcal{VE}(\bm{\mu}, \Sigmama, g)$. Let $\{i_q\}_{q=0}^\alpha, \{j_q\}_{q=0}^\beta\subset [p]$ be respectively a sequence of $\alpha+1$ and $\beta+1$ distinct integers. Let $m\in [\min\{\alpha,\beta\}]$. Suppose $i_q=j_q$ for $q\in [m]$ and $i_q \not= j_{q'}$ for $q,q'\not\in [m]$. Denote $\mathcal{I} = \bigcup_{q=0}^\alpha \{i_q\} \bigcup \left(\bigcup_{q'=0}^{\beta} \{j_q'\}\right)$ and then $|\mathcal{I}|=\alpha+\beta-m+2$.
	\begin{enumerate}[label=(\alph*)]
		\item \label{item:eleineb}
		Then
		\begin{equation}
		\E \left(\prod\limits_{q=1}^\alpha\Phi_{i_0i_q}^{(\Rma)}\right) \left(\prod\limits_{q'=1}^\beta\Phi_{j_0j_{q'}}^{(\Rma)}\right) \leq 
		\mu_{n,|\mathcal{I}|}\left(\Sigmama_{\mathcal{I}}\right) (2P_\nn(r_\rho))^{\alpha+\beta-m} \left(2P_\nn(2r_\rho)\right). 
		\label{eqn:prodmupp2}
		\end{equation}

		\item \label{item:eleinec}
		Then
		\begin{align*}
		\E\Phi_{i_0j_0}^{(\Rma)} \left(\prod\limits_{q=1}^\alpha\Phi_{i_0i_q}^{(\Rma)}\right) \left(\prod\limits_{q'=1}^\beta\Phi_{j_0j_{q'}}^{(\Rma)}\right) 
		\leq &
 	 \mu_{n,|\mathcal{I}|}\left(\Sigmama_{\mathcal{I}}\right) 	\left(2P_n(r_\rho)\right)^{\alpha+\beta-m+1} \sum_{\ell=2}^{m+2} \alpha_{n,m+1}(\ell,r_\rho)\\
		\leq &
		\mu_{n,|\mathcal{I}|}\left(\Sigmama_{\mathcal{I}}\right) (2P_\nn(r_\rho))^{\alpha+\beta-m+1}. \numberthis \label{eqn:prodmupp3}
		\end{align*}
		The inequality \eqref{eqn:prodmupp3} also holds with $m=0$. 		
	\end{enumerate}
\end{lem}
\begin{proof}
\noindent	(a)	
	By Lemma \ref{lem:densitybou} \ref{item:densitybouc}, it suffices to prove
		\eqref{eqn:prodmupp2} without $\mu_{n,|\mathcal{I}|}\left(\Sigmama_{\mathcal{I}}\right)$ for the case that $\{\uve_j\}$ for $j\in \mathcal{I}$ are independent $\unif(S^{n-2})$. Conditioned on $\uve_{i_0}$ and $\uve_{j_0}$, $\{\Phi_{i_0i_q}^{(\Rma)}\Phi_{j_0i_q}^{(\Rma)}\}_{q=1}^m$ are i.i.d., $\{\Phi_{i_0i_q}^{(\Rma)}\}_{q=m+1}^\alpha\bigcup\{\Phi_{j_0j_{q'}}^{(\Rma)}\}_{q'=m+1}^\beta$ are i.i.d. and moreover, every term in $\{\Phi_{i_0i_q}^{(\Rma)}\Phi_{j_0i_q}^{(\Rma)}\}_{q=1}^m$ is independent of every term in $\{\Phi_{i_0i_q}^{(\Rma)}\}_{q=m+1}^\alpha\bigcup\{\Phi_{j_0j_{q'}}^{(\Rma)}\}_{q'=m+1}^\beta$. Thus
		\begin{align*}
		&\E\left[ \left.\left(\prod\limits_{q=1}^\alpha\Phi_{i_0i_q}^{(\Rma)}\right) \left(\prod\limits_{q'=1}^\beta\Phi_{j_0j_{q'}}^{(\Rma)}\right) \right|\uve_{i_0},\uve_{j_0}\right]\\
		=&
		\E\left[\left.\prod_{q=1}^m \Phi_{i_0i_q}^{(\Rma)}\Phi_{j_0i_q}^{(\Rma)}\right|\uve_{i_0},\uve_{j_0}\right] \left(\E\left[ \left.\Phi_{i_0i_\alpha}^{(\Rma)}{1}(\alpha>m)+\Phi_{j_0j_\beta}^{(\Rma)}{1}(\alpha=m,\beta>m) \right|\uve_{i_0},\uve_{j_0}\right]\right)^{\alpha+\beta-2m}\\
		=& 
		\E\left[\left. \prod_{q=1}^m \Phi_{i_0i_q}^{(\Rma)}\Phi_{j_0i_q}^{(\Rma)}\right|\uve_{i_0},\uve_{j_0}\right] \left(2P_n(r_\rho)\right)^{\alpha+\beta-2m} \numberthis \label{eqn:conine2} \\
		=& 
		\left(\E\left[\left.  \Phi_{i_0i_q}^{(\Rma)}\Phi_{j_0i_q}^{(\Rma)}\right|\uve_{i_0},\uve_{j_0}\right]\right)^m \left(2P_n(r_\rho)\right)^{\alpha+\beta-2m} \numberthis \label{eqn:conine}
		\end{align*}
		where for the first equality the convention $0^0=1$ is used if $\alpha=\beta=m$. Notice \eqref{eqn:conine} also holds for $m=0$.
		
		Denote $ \overline{\SC}(r,\wwve) =  \SC(r,\wwve) \cup \SC(r,-\wwve)$. Then conditioned on $\uve_{i_0}$ and $\uve_{j_0}$,
		\begin{align*}
		&\Phi_{i_0i_1}^{(\Rma)}\Phi_{j_0i_1}^{(\Rma)}\\
		=& 1\left(\uve_{i_1}\in \overline{\SC}(r,\uve_{i_0}) \cap \overline{\SC}(r,\uve_{j_0})\right)\\ =&1(\|\uve_{i_0}-\uve_{j_0}\|_2\leq 2r_\rho\text{ or } \|\uve_{i_0}+\uve_{j_0}\|_2\leq 2r_\rho) 1\left(\uve_{i_1}\in \overline{\SC}(r,\uve_{i_0}) \cap \overline{\SC}(r,\uve_{j_0})\right),
		\end{align*}
		where the last equality follows by noticing that $\overline{\SC}(r_\rho,\uve_{i_0}) \cap \overline{\SC}(r_\rho,\uve_{j_0})$ is non-empty only when $\|\uve_{i_0}-\uve_{j_0}\|_2\leq 2r_\rho$ or $\|\uve_{i_0}+\uve_{j_0}\|_2\leq 2r_\rho$. Plugging the above inequality into \eqref{eqn:conine}, we obtain
		\begin{align*}
		&\E\left[ \left.\left(\prod\limits_{q=1}^\alpha\Phi_{i_0i_q}^{(\Rma)}\right) \left(\prod\limits_{q'=1}^\beta\Phi_{j_0j_{q'}}^{(\Rma)}\right) \right|\uve_{i_0},\uve_{j_0}\right]\\
		= & 1(\|\uve_{i_0}-\uve_{j_0}\|_2\leq 2r_\rho\text{ or } \|\uve_{i_0}+\uve_{j_0}\|_2\leq 2r_\rho)\left(\E\left[\left.\Phi_{i_0i_1}^{(\Rma)}\Phi_{j_0i_1}^{(\Rma)}\right|\uve_{i_0},\uve_{j_0}\right]\right)^m \left(2P_n(r_\rho)\right)^{\alpha+\beta-2m}\\
		\leq &  1(\|\uve_{i_0}-\uve_{j_0}\|_2\leq 2r_\rho\text{ or } \|\uve_{i_0}+\uve_{j_0}\|_2\leq 2r_\rho)\left(2P_n(r_\rho)\right)^m \left(2P_n(r_\rho)\right)^{\alpha+\beta-2m}.
		\end{align*}
		The result then follows by taking expectation w.r.t. $\uve_{i_0}$ and $\uve_{j_0}$.\\
		
		\noindent (b)
		Similar to the proof of \ref{item:eleineb}, it suffices to prove
		\eqref{eqn:prodmupp3} without $\mu_{n,|\mathcal{I}|}\left(\Sigmama_{\mathcal{I}}\right)$ for the case $\uve_j$ for $j\in \mathcal{I}$ are independent $\unif(S^{n-2})$. Conditioned on $\uve_{i_0}$ and $\uve_{j_0}$,
		\begin{align*}
		&\E \left[\Phi_{i_0j_0}^{(\Rma)} \left.\left(\prod\limits_{q=1}^\alpha\Phi_{i_0i_q}^{(\Rma)}\right) \left(\prod\limits_{q'=1}^\beta\Phi_{j_0j_{q'}}^{(\Rma)}\right) \right|\uve_{i_0},\uve_{j_0}\right]\\
		=&
		\Phi_{i_0j_0}^{(\Rma)} \E\left[ \prod_{q=1}^m \left.\Phi_{i_0i_q}^{(\Rma)}\Phi_{j_0i_q}^{(\Rma)}\right|\uve_{i_0},\uve_{j_0}\right] \left(2P_n(r_\rho)\right)^{\alpha+\beta-2m}.
		\end{align*}
		where the equality follows from \eqref{eqn:conine2}. Then
		\begin{align*}
		&\E \left[\Phi_{i_0j_0}^{(\Rma)} \left(\prod\limits_{q=1}^\alpha\Phi_{i_0i_q}^{(\Rma)}\right) \left(\prod\limits_{q'=1}^\beta\Phi_{j_0j_{q'}}^{(\Rma)}\right)\right]\\
		= & 
		\E\left[ \prod_{q=0}^m \Phi_{i_0i_q}^{(\Rma)} \prod_{q'=0}^{m} \Phi_{j_0 i_{q'}}^{(\Rma)} \right] \left(2P_n(r_\rho)\right)^{\alpha+\beta-2m}\\
		= & 
		\P\left( \prod_{q=0}^m \Phi_{i_0i_q}^{(\Rma)} =1 \right) \P\left( \left. \prod_{q'=0}^{m} \Phi_{j_0 i_{q'}}^{(\Rma)} = 1\right|\prod_{q=0}^m \Phi_{i_0i_q}^{(\Rma)} =1 \right) \left(2P_n(r_\rho)\right)^{\alpha+\beta-2m} \\
		\overset{(*)}{\leq} & 
		\left(2P_n(r_\rho)\right)^{\alpha+\beta-m+1} \sum_{\ell=2}^{m+2} \alpha_{n,m+1}(\ell,r_\rho)\\
		\leq &
		\left(2P_n(r_\rho)\right)^{\alpha+\beta-m+1}, 
		\end{align*}
		where step $(*)$ follows from the definition of $\alpha_{n,\delta}(\ell,r_\rho)$ in \eqref{eqn:alphandeltarho}. Notice \eqref{eqn:conine} also holds for $m=0$.
		
\end{proof}

\subsubsection{Lemmas on double summations}
\label{sec:quadraticterms}

Denote $\vec{i}\cup \vec{j} =  \left[\vec{i}\right]\bigcup\left[\vec{j}\right] $ for any $\vec{i}\in C_q^<$ and any $\vec{j}\in C_\delta^<$. Consider any $\theta_{\vec{i},\vec{j}}$ that is a non-negative function of $\uve_{\ell}$ for $\ell\in \vec{i}\cup\vec{j}$ defined for $\vec{i}\in C_q^<$ and $\vec{j}\in C_\delta^<$ with $1\leq  \delta\leq q\leq p-1$. In this section  an upper bound on $\E \sum_{\vec{i}\in C_q^< }\sum_{\vec{j}\in C_\delta^<} \theta_{\vec{i},\vec{j}}$ is presented. The results in this subsection will be used in the proofs of Proposition \ref{prop:edgecor} and Proposition \ref{prop:allcloseL2}.

For $i\in [p]$, let 
\begin{equation}
\mathcal{NZ}(i):=\{m\in [p]: \Sigmama_{im}\not =0\} \label{eqn:NZdef}
\end{equation}
denote the index of the variables that has non zero correlation with the $i$-th variable. For $\vec{i}\in C_q^<$, define $\mathcal{NZ}\left(\vec{i}\right) := \bigcup\limits_{\ell=0}^q\mathcal{NZ}(i_{\ell}) $. Since $\Sigmama$ is row-$\kappa$ sparse, for any $\vec{i}\in C_q^<$, $\left|\mathcal{NZ}\left(\vec{i}\right)\right|\leq (q+1)\kkkk$, and
	\begin{equation}
	p_{\vec{i}}:=\left|[p]\backslash \mathcal{NZ}\left(\vec{i}\right)\right|\geq p-(q+1)\kkkk. \label{eqn:p'upperbounew}
	\end{equation}
Note that $p_{\vec{i}}$ denotes the number of variables that are independent of variables with index in $[\vec{i}]$. 
	\\
	
	For $\vec{i}\in C_q^<$, define 
	\begin{align}
	 J_{\vec{i}} & := \left\{\vec{j}\in C_\delta^<:  \bigcup_{\ell=0}^\delta \{j_\ell\}  \subset \bigcup_{\ell=0}^q \{i_\ell\} \right\}\label{eqn:Sidefgeneral},\\
	 T_{\vec{i}} & :=\left\{\vec{j}\in C_\delta^<: \left(\bigcup_{\ell=0}^\delta \{j_\ell\} \right)\bigcap \left( \mathcal{NZ}\left(\vec{i}\right) \right) = \emptyset  \right\}, \label{eqn:Tidef} \\
	 N_{\vec{i}} & := C_q^< \backslash J_{\vec{i}}  \backslash T_{\vec{i}}. \label{eqn:Nidef}
	 \end{align}
	 Here $J_{\vec{i}}$ is the set of indices in $C_\delta^<$ consisting of coordinates as subsets of $\left[\vec{i}\right]$; $T_{\vec{i}}$ is the set of indices in $C_\delta^<$ consisting of coordinates outside the neighborhood of $\vec{i}$; $N_{\vec{i}}$ is the set of "correlated but not highly correlated" indices in $C_\delta^<$, i.e. the set of indices of which at least one coordinate is in the neighborhood of $\vec{i}$, but excluding those sets of indices of which the set of coordinates are subsets as that of $\vec{i}$.
	
	The strategy is to decompose 
	$$
	\E \sum_{\vec{i}\in C_q^< }\sum_{\vec{j}\in C_\delta^<} \theta_{\vec{i},\vec{j}} =  \E \sum_{\vec{i}\in C_q^< } \sum_{\vec{j}\in  J_{\vec{i}}}\theta_{\vec{i},\vec{j}} + \E \sum_{\vec{i}\in C_q^< } \sum_{\vec{j}\in  T_{\vec{i}}}  \theta_{\vec{i},\vec{j}} +\E \sum_{\vec{i}\in C_q^< } \sum_{\vec{j}\in N_{\vec{i}}}\theta_{\vec{i},\vec{j}}
	$$
	and bound each of the three terms.

	The next result is an upper bound on the first two terms.

\begin{lem}
	\label{lem:quadraticsum}
	Let $p\geq n \geq 4$ and $\Xma\sim \mathcal{VE}(\bm{\mu},\Sigmama,g)$. Suppose $\Sigmama$ is row-$\kappa$ sparse.
		Consider any $\theta_{\vec{i},\vec{j}}$ that is a non-negative function of $\uve_{\ell}$ for $\ell\in \vec{i}\cup\vec{j}$ defined for $\vec{i}\in C_q^<$ and $\vec{j}\in C_\delta^<$ with $1\leq  \delta\leq q\leq p-1$.  
		\begin{enumerate}[label=(\alph*)]
		\item \label{item:quadraticsuma}
		Suppose there exist positive constants $a,z$ such that $\E \theta_{\vec{i},\vec{j}} \leq \mu_{n,q+1}(\Sigmama_{\vec{i}}) a z^q$ for any $\vec{j}\in J_{\vec{i}}$. Then
		\begin{align*}
		\sum_{\vec{i}\in C_q^<} \sum_{\vec{j}\in J_{\vec{i}}} \E \theta_{\vec{i},\vec{j}} \leq & ap (pz)^q \frac{q+1}{\delta!(q-\delta)!} \left(1+q^2\mu_{n,q+1}(\Sigmama)\frac{\kappa-1}{p}\right)
		\end{align*}
		\item \label{item:quadraticsumb}
		Suppose there exist positive constants $a,z$ such that $\E \theta_{\vec{i},\vec{j}} \leq \mu_{n,q+\delta+2}(\Sigmama_{\vec{i}\cup\vec{j}}) a z^{q+\delta}$ for any  $\vec{j}\in T_{\vec{i}}$.
		Then
		\begin{align*}
		\sum_{\vec{i}\in C_q^<}  \sum_{\vec{j}\in T_{\vec{i}}} \E \theta_{\vec{i},\vec{j}} \leq& a p^2 (pz)^{q+\delta} \frac{3}{\delta! (q-1)!} \left(1+ \mu_{n,q+\delta+2}(\Sigmama) \frac{\kappa-1}{p} \right) .
		\end{align*}
		\end{enumerate}
	\end{lem}
	
\begin{proof} (a)
     Since $|J_{\vec{i}}|=\binom{q+1}{1}\binom{q}{\delta}$
     \begin{align*}
         \sum_{\vec{i}\in C_q^<} \sum_{\vec{j}\in J_{\vec{i}}} \E \theta_{\vec{i},\vec{j}} 
         \leq & a z^q \binom{q+1}{1}\binom{q}{\delta} \sum_{\vec{i}\in C_q^<}\mu_{n,q+1}(\Sigmama_{\vec{i}})\\
         \leq & a p (pz)^q \frac{q+1}{\delta!(q-\delta)!} \left(1+q^2\mu_{n,q+1}(\Sigmama)\frac{\kappa-1}{p}\right),
     \end{align*}
     where the last step follows from Lemma \ref{lem:summuupp}.\\

\noindent
(b)     
\begin{align*}
&\sum_{\vec{i}\in C_q^<}  \sum_{\vec{j}\in T_{\vec{i}}} \E \theta_{\vec{i},\vec{j}} \\
\leq & 
az^{q+\delta} \sum_{\vec{i}\in C_q^<}  \sum_{\vec{j}\in T_{\vec{i}}} \mu_{n,q+\delta+2}(\Sigmama_{\vec{i}\cup \vec{j}})  \\
\leq & 
az^{q+\delta}
\left(\sum_{\substack{\vec{i}\in C_q^<\\ \Sigmama_{\vec{i}} \text{ diagonal }}}  \sum_{\substack{\vec{j}\in T_{\vec{i}}\\ \Sigmama_{\vec{j}} \text{ diagonal }  }} 1
+ \mu_{n,q+\delta+2}(\Sigmama) \sum_{\substack{\vec{i}\in C_q^<\\ \Sigmama_{\vec{i}} \text{ diagonal }}}  \sum_{\substack{\vec{j}\in T_{\vec{i}}\\ \Sigmama_{\vec{j}} \text{ not diagonal }  }}1+\right.\\
&\quad \left.\mu_{n,q+\delta+2}(\Sigmama) \sum_{\substack{\vec{i}\in C_q^<\\ \Sigmama_{\vec{i}} \text{ not diagonal }}}  \sum_{\vec{j}\in T_{\vec{i}}} 1\right)
\\
\leq & az^{q+\delta}
\left(\binom{p}{1}\binom{p-1}{q}  \binom{p}{1}\binom{p-1}{\delta} 
+ \mu_{n,q+\delta+2}(\Sigmama) \binom{p}{1}\binom{p-1}{q}  \sum_{\substack{\vec{j}\in C_\delta^<\\ \Sigmama_{\vec{j}} \text{ not diagonal }  }}1+ \right. \\
&\quad \left.\mu_{n,q+\delta+2}(\Sigmama) \sum_{\substack{\vec{i}\in C_q^<\\ \Sigmama_{\vec{i}} \text{ not diagonal }}}  \binom{p}{1}\binom{p-1}{\delta} \right)\\
\leq & 
a p^2 (pz)^{q+\delta} \frac{3}{\delta! (q-1)!} \left(1+ \mu_{n,q+\delta+2}(\Sigmama) \frac{\kappa-1}{p} \right) ,
\end{align*}
where the second inequality follows from that for $\vec{j}\in T_{\vec{i}}$, $\Sigmama_{\vec{i}\cup \vec{j} }$ is diagonal if and only if $\Sigmama_{\vec{i}}$ and $\Sigmama_{\vec{j}}$ are both diagonal; and the last step follows from Lemma \ref{lem:numofnondia}.\\
\end{proof}
	
	To control $\E \sum_{\vec{i}\in C_q^< } \sum_{\vec{j}\in N_{\vec{i}}}\theta_{\vec{i},\vec{j}}$, we further partition $N_{\vec{i}}$ into $6$ subsets as follows. For $\vec{i}\in C_q^<$ with $q\geq \delta$, define
	\begin{align*}
	\mathcal{K}_1\left(\vec{i}\right) & := \left\{\vec{j}\in N_{\vec{i}}:j_0=i_0 \right\},\\
	\mathcal{K}_2\left(\vec{i}\right) & := \left\{\vec{j}\in N_{\vec{i}}:j_0\not=i_0, j_0\in \bigcup_{\ell=1}^q \{i_\ell\}, i_0\in \bigcup_{\ell=1}^\delta\{j_\ell\}\right\}, \\
	\mathcal{K}_3\left(\vec{i}\right) & := \left\{\vec{j}\in N_{\vec{i}}:j_0\not=i_0, j_0\not\in \bigcup_{\ell=1}^q\{i_\ell\}, i_0\in \bigcup_{\ell=1}^\delta\{j_\ell\}\right\},\\
	\mathcal{K}_4\left(\vec{i}\right) & := \left\{\vec{j}\in N_{\vec{i}}:j_0\not=i_0, j_0\in \bigcup_{\ell=1}^q \{i_\ell\}, i_0\not\in \bigcup_{\ell=1}^\delta\{j_\ell\} \right\},\\
	\mathcal{K}_5\left(\vec{i}\right) & := \left\{\vec{j}\in N_{\vec{i}}:j_0\not=i_0, j_0\not\in \bigcup_{\ell=1}^q\{i_\ell\}, i_0\not\in \bigcup_{\ell=1}^\delta\{j_\ell\},   \left|\left(\bigcup_{\ell=1}^q\{i_\ell\}\right)\bigcap
	\left(\bigcup_{\ell=1}^\delta \{j_\ell\}\right)  \right| \geq 1      \right\},\\
	\mathcal{K}_6\left(\vec{i}\right) & := \left\{\vec{j}\in N_{\vec{i}}:j_0\not=i_0, j_0\not\in \bigcup_{\ell=1}^q\{i_\ell\}, i_0\not\in \bigcup_{\ell=1}^\delta\{j_\ell\},   \left(\bigcup_{\ell=1}^q\{i_\ell\}\right) \bigcap \left(\bigcup_{\ell=1}^\delta \{j_\ell\}\right)  = \emptyset      \right\}.
	\end{align*}
	Then $N_{\vec{i}}= \cup_{w=1}^6\Kc_w\left(\vec{i}\right)$. Let $D_{\vec{i}}^m=\{\vec{j}\in N_{\vec{i}}: \left| \left(\cup_{\ell=1}^q\{i_\ell\}\right) \bigcap \left(\cup_{\ell=1}^\delta\{j_\ell\}\right)  \right|=m \}.$
	 We are now in a good position to present 
	a lemma on $\E \sum_{\vec{i}\in C_q^< } \sum_{\vec{j}\in N_{\vec{i}}} \theta_{\vec{i},\vec{j}}$. 
	
	\begin{lem}
	\label{lem:b2}
	Let $p\geq n \geq 4$ and $\Xma\sim \mathcal{VE}(\bm{\mu},\Sigmama,g)$. Suppose $\Sigmama$ is row-$\kappa$ sparse.
		Consider any $\theta_{\vec{i},\vec{j}}$ that is a non-negative function of $\uve_{\ell}$ for $\ell\in \vec{i}\cup\vec{j}$ defined for $\vec{i}\in C_q^<$ and $\vec{j}\in C_\delta^<$ with $1\leq  \delta\leq q\leq p-1$. Suppose there exist positive constants $a,b,z$ such that $\theta_{\vec{i},\vec{j}}$ satisfies: 
		\begin{align*}
		\E \theta_{\vec{i},\vec{j}} \leq \mu_{n,|\vec{i}\cup \vec{j}|}(\Sigmama_{\vec{i}\cup \vec{j}}) az^{q+\delta-m}, &\quad \forall\ \vec{j}\in \mathcal{K}_w\left(\vec{i}\right) \cap D_{\vec{i}}^m,\ \forall\ 0\leq m\leq \delta-1,\ \forall w\in\{1,3,4\}; \\
		\E \theta_{\vec{i},\vec{j}} \leq \mu_{n,|\vec{i}\cup \vec{j}|}(\Sigmama_{\vec{i}\cup \vec{j}}) az^{q+\delta-m-1}, &\quad \forall\ \vec{j}\in \mathcal{K}_2\left(\vec{i}\right) \cap D_{\vec{i}}^m,\ \forall\ 0\leq m\leq \delta-2; \\
		\E \theta_{\vec{i},\vec{j}} \leq \mu_{n,|\vec{i}\cup \vec{j}|}(\Sigmama_{\vec{i}\cup \vec{j}}) abz^{q+\delta-m}, &\quad \forall\ \vec{j}\in \mathcal{K}_5\left(\vec{i}\right) \cap D_{\vec{i}}^m,\ \forall\ 1\leq m\leq \delta;\\
		\E \theta_{\vec{i},\vec{j}} \leq \mu_{n,|\vec{i}\cup \vec{j}|}(\Sigmama_{\vec{i}\cup \vec{j}}) az^{q+\delta}, &\quad \forall\ \vec{j}\in \Kc_6\left(\vec{i}\right).
		\end{align*}
		Then
		\begin{align*}
		\sum_{\vec{i}\in C_q^<} \sum_{\vec{j}\in N_{\vec{i}}} \E \theta_{\vec{i},\vec{j}} \leq& ap(pz)^{q+1}\left(1+ \mu_{n,q+\delta+1}\left(\Sigmama\right) (3q^2)\frac{\kappa-1}{p}\right)\left(1+pz\right)^{\delta-1} \delta\frac{4+b/z}{(\delta-1)!}+\\
		&\quad\quad\quad ap^2 (pz)^{q+\delta} \frac{(\delta+1)(q+1)}{\delta! q!}  \mu_{n,q+\delta+2}\left(\Sigmama\right)\frac{\kappa-1}{p}.
		\end{align*}
	\end{lem}
	
	\begin{proof}
	Since
	\begin{align*}
	N_{\vec{i}} = \bigcup_{w=1}^6\mathcal{K}_w\left(\vec{i}\right),
	\end{align*}
    \begin{equation}
    \sum_{\vec{i}\in C_q^<} \sum_{\vec{j}\in N_{\vec{i}}} \E \theta_{\vec{i},\vec{j}} =  \sum_{w=1}^6 \sum_{\vec{i}\in C_q^<}I_w\left(\vec{i}\right), \label{eqn:b2dec}
    \end{equation}
    with 
    $$
    I_w\left(\vec{i}\right) :=  \sum_{\vec{j}\in \mathcal{K}_w\left(\vec{i}\right)} \E\theta_{\vec{i},\vec{j}}.
    $$	
    
	\noindent \textbf{Case 1}: $p \geq q+\delta +2 $\\	
	 Obviously $\mathcal{K}_1\left(\vec{i}\right)=\bigcup_{m=0}^{\delta-1} \left(\mathcal{K}_1\left(\vec{i}\right)\cap D_{\vec{i}}^m \right)$. Then for any $\vec{i}\in C_q^<$ satisfying $\Sigmama_{\vec{i}}$ diagonal,  
	\begin{align*}
	&\left|\vec{j}\in \mathcal{K}_1\left(\vec{i}\right) \bigcap D_{\vec{i}}^m : \Sigmama_{\vec{i}\bigcup\vec{j}} \text{ diagonal } \right| \\
	= & 
	\binom{q}{m} \frac{1}{(\delta-m)!} \sum_{  j_1\in [p]\backslash \mathcal{NZ}\left(\vec{i}\right)  }\  \sum_{ \substack{ j_2\in [p]\backslash \mathcal{NZ}\left(\vec{i}\right) \\ j_2\not \in \mathcal{NZ}(j_1) } } \ \ \sum_{ \substack{ j_3\in [p]\backslash \mathcal{NZ}\left(\vec{i}\right) \\ j_3\not \in \cup_{\ell=1}^2 \mathcal{NZ}(j_\ell) } } \cdots \sum_{ \substack{ j_{\delta-m}\in [p]\backslash \mathcal{NZ}\left(\vec{i}\right) \\ j_{\delta-m}\not \in \cup_{\ell=1}^{\delta-m-1} \mathcal{NZ}(j_\ell) } } 1 \\
	\geq &
	\binom{q}{m} \frac{1}{(\delta-m)!}\prod_{\ell=0}^{\delta-m-1}\left(p_{\vec{i}}-\ell \kappa \right), \numberthis \label{eqn:K1mdia}
	\end{align*}
	where in the first inequality we assume without loss of generality that the components of $\vec{j}$ distinct from $\vec{i}$ are $j_1,j_2,\cdots,j_{\delta-m}$.
	Then,
	\begin{align*}
	&\left|\vec{j}\in \mathcal{K}_1\left(\vec{i}\right) \bigcap D_{\vec{i}}^m : \Sigmama_{\vec{i}\bigcup\vec{j}} \text{ not diagonal } \right|\\
	\leq &\binom{q}{m} \binom{p-1-q}{\delta-m} -\binom{q}{m} \frac{1}{(\delta-m)!}\prod_{\ell=0}^{\delta-m-1}\left(p_{\vec{i}}-\ell \kappa \right)\\
	= &  
	\frac{\binom{q}{m}}{(\delta-m)!}\left(\prod_{\ell=0}^{\delta-m-1}\left(p-1-q-\ell  \right) -\prod_{\ell=0}^{\delta-m-1}\left(p_{\vec{i}}-\ell  \right)+\prod_{\ell=0}^{\delta-m-1}\left(p_{\vec{i}}-\ell  \right) -\prod_{\ell=0}^{\delta-m-1}\left(p_{\vec{i}}-\ell \kappa \right) \right) \\
	\leq &  
	\binom{q}{m} \frac{1}{(\delta-m)!}(\delta-m)p^{\delta-m-1}(q+\delta)(\kappa-1), \numberthis \label{eqn:K1mnotdia}
	\end{align*}
	where the first inequality follows from \eqref{eqn:K1mdia},
	and the second inequality follows from Lemma  \ref{basine} \ref{item:basineh}, \ref{item:basinei} and \eqref{eqn:p'upperbou}. 
	
	Then
	\begin{align*}
	&\sum_{\vec{i}\in C_q^<}I_1 \left(\vec{i}\right) \\
	=& 
	\sum_{m=0}^{\delta-1}\ \left( \sum_{ \substack{\vec{i}\in C_q^<\\ \Sigmama_{\vec{i}} \text{ diagonal}}} \  \sum_{ \substack{\vec{j}\in \mathcal{K}_1\left(\vec{i}\right)\bigcap D_{\vec{i}}^m\\ \Sigmama_{\vec{i}\cup\vec{j} \text{ diagonal}}}} 
	+
	\sum_{ \substack{\vec{i}\in C_q^<\\ \Sigmama_{\vec{i}} \text{ diagonal}}} \  \sum_{ \substack{\vec{j}\in \mathcal{K}_1\left(\vec{i}\right)\bigcap D_{\vec{i}}^m\\ \Sigmama_{\vec{i}\cup\vec{j} \text{ not diagonal}}}} 
	+
	\sum_{ \substack{\vec{i}\in C_q^<\\ \Sigmama_{\vec{i}} \text{ not diagonal}}} \  \sum_{ \vec{j}\in \mathcal{K}_1\left(\vec{i}\right)\bigcap D_{\vec{i}}^m} \right) \E\theta_{\vec{i},\vec{j}}\\
	\leq & 
	\sum_{m=0}^{\delta-1}\ \left( \sum_{ \substack{\vec{i}\in C_q^<\\ \Sigmama_{\vec{i}} \text{ diagonal}}} \  \sum_{ \substack{\vec{j}\in \mathcal{K}_1\left(\vec{i}\right)\bigcap D_{\vec{i}}^m\\ \Sigmama_{\vec{i}\cup\vec{j} \text{ diagonal}}}} 
	+\right.\\
	&
	\left.\mu_{n,q+\delta+1}\left(\Sigmama\right)\left(\sum_{ \substack{\vec{i}\in C_q^<\\ \Sigmama_{\vec{i}} \text{ diagonal}}} \  \sum_{ \substack{\vec{j}\in \mathcal{K}_1\left(\vec{i}\right)\bigcap D_{\vec{i}}^m\\ \Sigmama_{\vec{i}\cup\vec{j} \text{ not diagonal}}}} +\sum_{ \substack{\vec{i}\in C_q^<\\ \Sigmama_{\vec{i}} \text{ not diagonal}}} \  \sum_{ \vec{j}\in \mathcal{K}_1\left(\vec{i}\right)\bigcap D_{\vec{i}}^m}\right) \right) a z^{q+\delta-m}\\
	\leq & 
	\sum_{m=0}^{\delta-1}\ \left( \binom{p}{1}\binom{p-1}{q}  \binom{q}{m}\binom{p-1-q}{\delta-m} \right.\\ 
	&\left. \ \ \ \ \ +\mu_{n,q+\delta+1}\left(\Sigmama\right)\binom{p}{1}\binom{p-1}{q} \binom{q}{m} \frac{\delta-m}{(\delta-m)!}p^{\delta-m-1} (q+\delta)(\kappa-1) \right. \\
	&\quad \quad \quad \left. +\mu_{n,q+\delta+1}\left(\Sigmama\right)\frac{(q+1)}{2((q-1)!)}   p^{q} (\kappa-1)  \binom{q}{m}\binom{p-1-q}{\delta-m}  \right) a z^{q+\delta-m}\\
	\leq & 
	\sum_{m=0}^{\delta-1}\ \left(   \frac{1}{m!(\delta-m)!} +\mu_{n,q+\delta+1}\left(\Sigmama\right)\frac{\kappa-1}{p} \left( \frac{(q+\delta)(\delta-m)}{m!(\delta-m)!}    +      \frac{q(q+1)}{m!(\delta-m)!2}\right)  \right) ap\left( pz\right)^{q+\delta-m}\\
	\leq & 
	ap\left(pz\right)^{q+1} \left(   1   +\mu_{n,q+\delta+1}\left(\Sigmama\right)(3q^2)    \frac{\kappa-1}{p}    \right) \frac{1}{(\delta-1)!} \sum_{m=0}^{\delta-1} \frac{(\delta-1)!}{m!(\delta-1-m)!}   \left( pz\right)^{\delta-1-m} \\
	= &  
	ap\left(pz\right)^{q+1} \left(   1   +\mu_{n,q+\delta+1}\left(\Sigmama\right)(3q^2)    \frac{\kappa-1}{p}    \right) \frac{1}{(\delta-1)!} \left(1+pz\right)^{\delta-1}, \numberthis \label{eqn:I1temp}
	\end{align*}
	where the first inequality follows from $\mu_{n,q+\delta-m+1}\left(\Sigmama\right)\leq \mu_{n,q+\delta+1}\left(\Sigmama\right)$, and the second inequality follows from Lemma \ref{lem:numofnondia} and \eqref{eqn:K1mnotdia}.
	
	Obviously, 
	\begin{align*}
	\mathcal{K}_2\left(\vec{i}\right)= \bigcup_{m=0}^{\delta-2}\left(\mathcal{K}_2\left(\vec{i}\right) \bigcap D_{\vec{i}}^m \right) , & \quad \mathcal{K}_3\left(\vec{i}\right)= \bigcup_{m=0}^{\delta-1}\left(\mathcal{K}_3\left(\vec{i}\right) \bigcap D_{\vec{i}}^m \right) ,\\ \mathcal{K}_4\left(\vec{i}\right)= \bigcup_{m=0}^{\delta-1}\left(\mathcal{K}_4\left(\vec{i}\right) \bigcap D_{\vec{i}}^m \right), & \quad \mathcal{K}_5\left(\vec{i}\right) = \bigcup_{m=1}^\delta \left(\mathcal{K}_5\left(\vec{i}\right) \bigcap D_{\vec{i}}^m  \right) .
	\end{align*}
	Then, following a similar analysis to $\mathcal{K}_1(\vec{i})$, additionally with Lemma \ref{lem:eleine}, we obtain
	\begin{align*}
	\sum_{\vec{i}\in C_q^<}I_2 \left(\vec{i}\right) \leq &  ap\left(pz\right)^{q+1} \left(   1   +\mu_{n,q+\delta}\left(\Sigmama\right)(3q^2)    \frac{\kappa-1}{p}    \right) \frac{1}{(\delta-2)!} \left(1+pz\right)^{\delta-2}1(\delta \geq 2), \numberthis \label{eqn:I2temp} \\   
	\sum_{\vec{i}\in C_q^<}I_3 \left(\vec{i}\right) \leq &  ap(pz)^{q+1} \left(   1   +\mu_{n,q+\delta+1}\left(\Sigmama\right)(3q^2)    \frac{\kappa-1}{p}    \right) \frac{1}{(\delta-1)!} \left(1+pz\right)^{\delta-1}, \numberthis \label{eqn:I3temp} \\
	\sum_{\vec{i}\in C_q^<}I_4 \left(\vec{i}\right) \leq &  ap\left(pz\right)^{q+1} \left(   1   +\mu_{n,2\delta+1}\left(\Sigmama\right)(3q^2)    \frac{\kappa-1}{p}    \right) \frac{1}{(\delta-1)!} \left(1+pz\right)^{\delta-1}, \numberthis \label{eqn:I4temp} \\
	\sum_{\vec{i}\in C_q^<}I_5 \left(\vec{i}\right) 
	\leq &  a\left(\frac{b}{z}\right) p\left(pz\right)^{q+1}  \left(   1   +\mu_{n,q+\delta+1}\left(\Sigmama\right)(3\delta^2)    \frac{\kappa-1}{p}    \right) \frac{1}{(\delta-1)!} \left(1+pz\right)^{\delta-1}. \numberthis  \label{eqn:I5temp}
	\end{align*} 
	 The detailed derivation of the above inequalities are omitted for the sake of brevity.
	 
	 Observe that 
	 $$
	 \mathcal{K}_6\left(\vec{i}\right) =
	 \left\{\vec{j}\in C_q^<: \left(\bigcup_{\ell=0}^q\{i_\ell\}\right) \bigcap \left(\bigcup_{\ell=0}^\delta \{j_\ell\}\right)  = \emptyset , \  \exists \ell\in [\delta]\cup\{0\} \text{ such that } j_\ell \in \mathcal{NZ}\left(\vec{i}\right)   \right\} .
	 $$ 
	 Then,
	 \begin{align*}
	 \left|\mathcal{K}_6\left(\vec{i}\right)\right|=& \binom{p-1-q}{1}\binom{p-2-q}{\delta}-\binom{p_{\vec{i}}}{1}\binom{p_{\vec{i}}-1}{\delta}\\
	 \leq & \frac{1}{\delta!}(\delta+1)p^{\delta}(q+1)(\kappa-1), \numberthis \label{eqn:k6bou}
	 \end{align*}
	 where the last inequality follows from Lemma  \ref{basine} \ref{item:basineh} and \eqref{eqn:p'upperbou}. 
	 Thus,
	 \begin{align*}
	 \sum_{\vec{i}\in C_q^<}I_6\left(\vec{i}\right) 
	 \leq & \binom{p}{1}\binom{p-1}{q} \frac{1}{\delta!}(\delta+1)p^{\delta}(q+1)(\kappa-1) \mu_{n,q+\delta+2}\left(\Sigmama\right) az^{q+\delta},\\
	 \leq & ap^2 (pz)^{q+\delta} \frac{(\delta+1)(q+1)}{\delta! q!}  \mu_{n,q+\delta+2}\left(\Sigmama\right)\frac{\kappa-1}{p}, \numberthis \label{eqn:I6temp}
	 \end{align*}
	 where the first inequality follows from \eqref{eqn:k6bou} and Lemma \ref{lem:densitybou} \ref{item:densityboua}.
	 
	 \noindent \textbf{Case 2}: $p < q+\delta +2 $ \\
	 We have imposed the condition $p \geq 2\delta +2 $ to derive \eqref{eqn:I1temp}, \eqref{eqn:I2temp}, \eqref{eqn:I3temp}, \eqref{eqn:I4temp}, \eqref{eqn:I5temp} and \eqref{eqn:I6temp}. However, 
	 one can verify directly that these inequalities also hold when $p < q+\delta +2 $. We omit these tedious verifications here and take it for granted \eqref{eqn:I1temp}, \eqref{eqn:I2temp}, \eqref{eqn:I3temp}, \eqref{eqn:I4temp}, \eqref{eqn:I5temp} and \eqref{eqn:I6temp} holds for all $1\leq \delta\leq q\leq p-1$.
	 
	 Thus combining \eqref{eqn:b2dec}, \eqref{eqn:I1temp}, \eqref{eqn:I2temp}, \eqref{eqn:I3temp}, \eqref{eqn:I4temp}, \eqref{eqn:I5temp} and \eqref{eqn:I6temp}, yield
	 \begin{align*}
	 \sum_{\vec{i}\in C_q^<} \sum_{\vec{j}\in N_{\vec{i}}} \E \theta_{\vec{i},\vec{j}} & \leq 
	 ap(pz)^{q+1}\left(1+ \mu_{n,q+\delta+1}\left(\Sigmama\right) (3q^2)\frac{\kappa-1}{p}\right)\left(1+pz\right)^{\delta-1} \delta\frac{4+b/z}{(\delta-1)!}+\\
	 &\quad\quad\quad ap^2 (pz)^{q+\delta} \frac{(\delta+1)(q+1)}{\delta! q!}  \mu_{n,q+\delta+2}\left(\Sigmama\right)\frac{\kappa-1}{p}. 
	 \end{align*}
	\end{proof}

Combining Lemma \ref{lem:quadraticsum} and Lemma \ref{lem:b2} immediately yields the following lemma.
\begin{lem}
\label{lem:b2wholesum}
Let $p\geq n \geq 4$ and $\Xma\sim \mathcal{VE}(\bm{\mu},\Sigmama,g)$. Suppose $\Sigmama$ is row-$\kappa$ sparse. Consider any $\theta_{\vec{i},\vec{j}}$ that is a non-negative function of $\uve_{\ell}$ for $\ell\in \vec{i}\cup\vec{j}$ defined for $\vec{i}\in C_q^<$ and $\vec{j}\in C_\delta^<$ with $1\leq  \delta\leq q\leq p-1$. Suppose there exist $a,z,b$ such that all conditions in Lemma \ref{lem:quadraticsum} and in Lemma \ref{lem:b2} hold. Moreover suppose $b/z\leq c_{n,\delta,q}$ for some positive constant $c_{n,\delta,q}$ that only depends  on $n$, $q$ and $\delta$. Then 
\begin{align*}
& \sum_{\vec{i}\in C_q^<} \sum_{\vec{j}\in C_\delta^<} \E \theta_{\vec{i},\vec{j}} \\
\leq  &
C_{n,q,\delta} \left(p^{1+\frac{1}{\delta}}z\right)^q \left(1+(p^{1+\frac{1}{\delta}}z)^{\delta}\right)\left(1+pz\right)^{\delta-1} \left(1+\mu_{n,q+\delta+2}(\Sigmama)\frac{\kappa-1}{p}\right) ap^{1-\frac{q}{\delta}}
\end{align*}
\end{lem}

\subsubsection{Proof of Proposition \ref{prop:edgecor}}

In what follows in this subsection, for the sake of brevity, we write $C^<$ for $C_\delta^<$, and write $\Phi_{\vec{i}}$ for $\Phi_{\vec{i}}^{(\Rma)}$, for any $\vec{i}\in C^<$. 

\begin{proof}[Proof of Proposition \ref{prop:edgecor}]
	 Recall 
	 $$
	 N_{E_{\delta}}^{(\Rma)}=\sum_{\vec{i}\in C_\delta^<}\ \  \prod_{j=1}^\delta \Phi_{i_0i_j}^{(\Rma)}=\sum_{\vec{i}\in C_\delta^<}\ \  \Phi_{\vec{i}}.
	 $$

	 To apply a Compound Poisson Approximation result, some additional notation need to be introduced. For $\vec{i}\in C^<$, let $S_{\vec{i}}$ be defined in \eqref{eqn:Sidef}, and let $T_{\vec{i}}$, $N_{\vec{i}}$ be defined respectively as in \eqref{eqn:Tidef}, \eqref{eqn:Nidef} with $q=\delta$. Here $T_{\vec{i}}$ is the set of indices consisting of coordinates outside the neighborhood of $\vec{i}$; $N_{\vec{i}}$ is the set of "correlated but not highly correlated" indices, i.e. the set of indices of which at least one component but not every component is in the neighborhood of $\vec{i}$. 
	 Denote
	 \begin{equation}
	 W_{\vec{i}} = \sum_{\vec{j}\in T_{\vec{i}}} \Phi_{\vec{j}}, \quad Z_{\vec{i}} = \sum_{\vec{j}\in N_{\vec{i}}} \Phi_{\vec{j}}, \label{eqn:WiZidef}
	 \end{equation}
	 and recall $U_{\vec{i}} = \sum_{\vec{j}\in S_{\vec{i}}}\Phi_{\vec{j}}^{(\Rma)}$ is defined in \eqref{eqn:Uidef}. Then $W_{\vec{i}}$ is independent of $U_{\vec{i}}$ and $\Phi_{\vec{i}}$. Further denote 
	 \begin{align*}
	 \lambda_0 & = \sum_{\vec{i}\in C^<} \E \left( \frac{\Phi_{\vec{i}}}{\Phi_{\vec{i}}+U_{\vec{i}}} 1\left(\Phi_{\vec{i}}+U_{\vec{i}} \geq 1\right)\right),\\
	 \zeta_{0\ell} & = \frac{1}{\lambda_0 \ell}\sum_{\vec{i}\in C^<} \E\left( \Phi_{\vec{i}} 1\left(\Phi_{\vec{i}}+ U_{\vec{i}} = \ell \right) \right), \quad \forall \ell \geq 1 \numberthis \label{eqn:defofzeta0l}
	 \end{align*}
	 and a probability distribution $\bm{\zeta}_0$ on positive integers  with $\bm{\zeta}_0(\ell) = \zeta_{0\ell}$. The mean of $\bm{\zeta}_0$ is $\E \bm{\zeta}_0 = \sum_{\ell \geq 1} \ell \zeta_{0\ell}$. Moreover, let $b_1 =\sum_{\vec{i}\in C^<} \E \Phi_{\vec{i}} \E(\Phi_{\vec{i}} \ + U_{\vec{i}}+Z_{\vec{i}}) $ and 
	 \begin{equation}
	 b_2 = \sum_{\vec{i}\in C^<}\E \left( \Phi_{\vec{i}} Z_{\vec{i}} \right). \label{eqn:b2def}
	 \end{equation}

	 In this proof we write $\lambda$ and $\bm{\zeta}$ for $\lambda_{p,n,\delta,\rho}$ and $\bm{\zeta}_{n,\delta,\rho}\left( \ell \right)$ respectively when there is no confusion. By the compound Poisson Stein's approximation, i.e. (5.19) and (5.16) in \cite{barbour2001topics}, 
	 \begin{equation}
	 d_{\TV}\left(\mathscr{L}\left(N_{E_\delta}^{(\Rma)}\right), \CP(\lambda,\bm{\zeta})\right) \leq e^{\lambda_0}  \left(b_1+b_2+\lambda_0  d_{W}(\bm{\zeta}_0',\bm{\zeta}') \E \bm{\zeta}_0+ |\lambda_0\E\bm{\zeta_0} - \lambda\E\bm{\zeta} |\right), \label{eqn: compoibou}
	 \end{equation}
	 where $\bm{\zeta}_{0}'(\ell) = \ell \zeta_{0\ell}/\E \bm{\zeta}_{0} $ and $\bm{\zeta}'(\ell) = \ell \bm{\zeta}(\ell)/\E \bm{\zeta} $ for $\ell \in \mathbb{Z}_+$, the set of all positive integers. In \eqref{eqn: compoibou}, the distance $d_W$ is the Wasserstein $L_1$ metric on probability measures over the set of positive integers $\mathbb{Z}_+$:
	 $$
	 d_W(P,Q) = \sup_{f\in \text{Lip}_1} \left| \int f dP - \int f dQ \right|
	 $$
	 where $\text{Lip}_1=\{f:|f(r)-f(s)|\leq |r-s|, r,s \in \mathbb{Z}_+ \}$.

	By Lemma \ref{prop:was1dis}, 
	\begin{align*}
	\lambda_0 d_{W}(\bm{\zeta}_0',\bm{\zeta}') \E \bm{\zeta}_0 \leq & \lambda_0 \E \bm{\zeta}_0 \frac{\delta}{2}\sum_{\ell=1}^{\delta+1}\ell\left|\frac{\lambda_0\zeta_{0\ell}}{\lambda_0\E\bm{\zeta}_0} - \frac{\lambda\bm{\zeta}(\ell)}{\lambda\E\bm{\zeta}} \right|\\
	 =  & \frac{\delta}{2}\sum_{\ell=1}^{\delta+1}\ell\left|\left(\lambda_0\zeta_{0\ell} - \lambda \bm{\zeta}(\ell)\right) + \left(\lambda\E\bm{\zeta}-\lambda_0 \E \bm{\zeta}_0 \right)\frac{\lambda\bm{\zeta}(\ell)}{\lambda\E\bm{\zeta}} \right|\\
	 \leq & \frac{\delta}{2}\sum_{\ell=1}^{\delta+1}\ell\left|\lambda_0\zeta_{0\ell} - \lambda \bm{\zeta}(\ell)\right| + \frac{\delta}{2}\left| \lambda\E\bm{\zeta}-\lambda_0 \E \bm{\zeta}_0  \right|\sum_{\ell=1}^{\delta+1}\ell\frac{\lambda\bm{\zeta}(\ell)}{\lambda\E\bm{\zeta}} \\
	 \leq & \delta\sum_{\ell=1}^{\delta+1}\ell\left|\lambda_0\zeta_{0\ell} - \lambda \bm{\zeta}(\ell)\right|.
	\end{align*}
	Plugging the above inequalities into \eqref{eqn: compoibou},

\begin{equation}
d_{\TV}\left(\mathscr{L}\left(N_{E_\delta}^{(\Rma)}\right), \CP(\lambda,\bm{\zeta})\right) \leq e^{\lambda_0}  \left(b_1+b_2+ \left(\delta+1\right)\sum_{\ell=1}^{\delta+1}\ell\left|\lambda_0\zeta_{0\ell} - \lambda \bm{\zeta}(\ell)\right|\right). \label{eqn:compoibou1}
\end{equation}
	It remains to estimate the quantities in the right hand side of \eqref{eqn:compoibou1}. \\

    \noindent{\textbf{Part I. Upper bound for} $\lambda_0$ \textbf{and} $\sum_{\ell=1}^{\delta+1}\ell\left|\lambda_0\zeta_{0\ell} - \lambda \bm{\zeta}(\ell)\right|$ } \\
    For $\ell \in [\delta+1]$,
    
    \begin{align*}
    &\left|\lambda_0\zeta_{0\ell} - \lambda\bm{\zeta}\left( \ell \right) \right|\\
    = &\left|\lambda_0\zeta_{0\ell} - \frac{1}{\ell}\binom{p}{1}\binom{p-1}{\delta}(2P_n(r_\rho))^\delta\alpha(\ell,r_\rho) \right| \\
    \leq &\frac{1}{ \ell}\sum_{\vec{i}\in C^<} \left|\E\left( \Phi_{\vec{i}} 1\left(\Phi_{\vec{i}}+ U_{\vec{i}} = \ell \right) \right) - (2P_n(r_\rho))^\delta\alpha(\ell,r_\rho) \right|\\
    = & \frac{1}{ \ell}\sum_{ \substack{\vec{i}\in C^< \\ \Sigmama_{\vec{i}} \text{ not diagonal} }} \ \  \left|\E\left( \Phi_{\vec{i}} 1\left(\Phi_{\vec{i}}+ U_{\vec{i}} = \ell \right) \right) - (2P_n(r_\rho))^\delta\alpha(\ell,r_\rho) \right| \\
    \leq & \frac{1}{ \ell} (\mu_{n,\delta+1}\left(\Sigmama\right)+1) (2P_n(r_\rho))^\delta\alpha(\ell,r_\rho)\sum_{ \substack{\vec{i}\in C^< \\ \Sigmama_{\vec{i}} \text{ not diagonal} }} \ \  1 \\
    \leq &  \frac{\alpha(\ell,r_\rho)}{ \ell} \mu_{n,\delta+1}\left(\Sigmama\right) \gamma^\delta \frac{(\delta+1)}{(\delta-1)!}    \frac{\kappa-1}{p}, \numberthis \label{eqn:lambdazetabou}
     \end{align*}
    where the first inequality follows from the definition of $\zeta_{0\ell}$ in \eqref{eqn:defofzeta0l}, the second inequality follows from Lemma \ref{item:prod} and Lemma \ref{lem:densitybou} \ref{item:densityboua}, and the last inequality follows from Lemma \ref{lem:numofnondia} and $\mu_{n,\delta+1}\left(\Sigmama\right)\geq 1$. 
    
    Then,
    \begin{equation}
    |\lambda_0\E\bm{\zeta_0} - \lambda\E\bm{\zeta} | \leq \sum_{\ell=1}^{\delta+1} \ell\left|\lambda_0\zeta_{0\ell} - \lambda\bm{\zeta(\ell)} \right|\leq  \mu_{n,\delta+1}\left(\Sigmama\right) \gamma^\delta \frac{(\delta+1)}{(\delta-1)!}    \frac{\kappa-1}{p}, \label{eqn:lambdameandif}
    \end{equation}
    where the last inequality follows from  \eqref{eqn:lambdazetabou}. As an immediate consequences,
    \begin{align*}
    \lambda_0 
    &\leq \lambda_0\E\bm{\zeta_0}  \\
    & \leq |\lambda_0\E\bm{\zeta_0} - \lambda\E\bm{\zeta} |+ \lambda\E\bm{\zeta} \\
    & \leq   \mu_{n,\delta+1}\left(\Sigmama\right) \gamma^\delta \frac{(\delta+1)}{(\delta-1)!}    \frac{\kappa-1}{p} + \binom{p}{1}\binom{p-1}{\delta}(2P_n(r_\rho))^\delta \\
    & \leq \gamma^\delta \frac{(\delta+1)}{(\delta-1)!}\left( \mu_{n,\delta+1}\left(\Sigmama\right)     \frac{\kappa-1}{p} +1 \right), \numberthis \label{eqn:lambda0uppbou}
    \end{align*}
    where the third inequality follows from \eqref{eqn:lambdameandif}.\\

	\noindent {\textbf{Part II. Upper bound for } $b_1$}\\
	
	Since $N_{\vec{i}}\cup S_{\vec{i}}\cup \{\vec{i}\}=C^<\backslash T_{\vec{i}}$,
	\begin{equation}
	b_1 = \sum_{\vec{i}\in C^<} \ \  \sum_{\vec{j}\in C^< \backslash T_{\vec{i}} } \E\Phi_{\vec{i}} \E \Phi_{\vec{j}}.  \label{eqn:b1uppbounew}
	\end{equation}
	
	Given $\vec{i}\in C^<$,	by \eqref{eqn:p'upperbounew} with $q=\delta$, 
	\begin{equation}
	p_{\vec{i}}:=\left|[p]\backslash \mathcal{NZ}\left(\vec{i}\right)\right|\geq p-(\delta+1)\kkkk. \label{eqn:p'upperbou}
	\end{equation}
	Since $\left|T_{\vec{i}}\right| = p_{\vec{i}}\binom{p_{\vec{i}}-1}{\delta}$, 
	\begin{equation}
	\left| C^<\backslash T_{\vec{i}}\right|=p\binom{p-1}{\delta}-p_{\vec{i}}\binom{p_{\vec{i}}-1}{\delta}
	\leq \frac{1}{\delta !} (\delta+1)\left(\prod_{\alpha=0}^{\delta-1} (p-\alpha) \right)   (p-p_{\vec{i}})\leq \frac{(\delta+1)^2}{\delta !}  p^{\delta}\kkkk, \label{eqn:neiborsize}
	\end{equation}
	where the first inequality follows from Lemma \ref{basine} \ref{item:basineh}.  

	One straightforward upper bound is
	\begin{align*}
	b_1 & = \sum_{\vec{i}\in C^{<}} \  \sum_{\vec{j}\in C^<\backslash T_{\vec{i}}} \E\left( \prod_{l=1}^\delta \Phi_{i_0i_l}^{(\Rma)} \right) \E\left( \prod_{l'=1}^\delta \Phi_{j_0j_{l'}}^{(\Rma)} \right)  \\
	&\leq p\binom{p-1}{\delta}\frac{(\delta+1)^2}{\delta !}  p^{\delta}\kkkk\left(\mu_{n,\delta+1}\left(\Sigmama\right)\right)^2(2P_n(r_\rho))^{2\delta} \\
	& \leq \left(\mu_{n,\delta+1}\left(\Sigmama\right)\right)^2\frac{(\delta+1)^2}{(\delta !)^2}\left(2p^{1+\frac{1}{\delta}}P_n(r_\rho)\right)^{2\delta}\frac{\kkkk}{p},
	\end{align*}
	where the first inequality follows from Lemma \ref{item:prod}, Lemma \ref{lem:densitybou} \ref{item:densityboua} and \eqref{eqn:neiborsize}. The $\left(\mu_{n,\delta+1}\left(\Sigmama\right)\right)^2$ in the above upper bound is not very satisfactory, and can be improved by a more involved analysis. 
	
	Observe for given $\vec{i}\in C^<$,
	\begin{align*}
	&\left|\{\vec{j}\in C^<\backslash T_{\vec{i}}: \Sigmama_{\vec{j}} \text{ not diagonal} \}  \right|\\
	=& |C^<\backslash T_{\vec{i}}|   - \left( \left|\{\vec{j}\in C^<: \Sigmama_{\vec{j}} \text{ diagonal} \}  \right| -\left|\{\vec{j}\in T_{\vec{i}}: \Sigmama_{\vec{j}} \text{ diagonal} \}  \right| \right)\\
	\leq &  \frac{1}{\delta !} (\delta+1)\left(\prod_{\alpha=0}^{\delta-1} (p-\alpha) \right)   (p-p_{\vec{i}}) - \left( \left|\{\vec{j}\in C^<: \Sigmama_{\vec{j}} \text{ diagonal} \}  \right| -\left|\{\vec{j}\in T_{\vec{i}}: \Sigmama_{\vec{j}} \text{ diagonal} \}  \right| \right), \numberthis \label{eqn:neinotdiagonalsize}
	\end{align*}
	where the inequality follows from \eqref{eqn:neiborsize}. 
	Then
	\begin{align*}
	&  \left|\{\vec{j}\in C^<: \Sigmama_{\vec{j}} \text{ diagonal} \}  \right| -\left|\{\vec{j}\in T_{\vec{i}}: \Sigmama_{\vec{j}} \text{ diagonal} \}  \right|  \\
	= & 
	\frac{1}{\delta !}\sum_{j_0=1}^p\ \sum_{j_1\in [p]\backslash \mathcal{NZ}(j_0)}\ \cdots \sum_{j_\delta\in [p]\backslash \bigcup\limits_{l=0}^{\delta-1}\mathcal{NZ}(j_l)}1- \frac{1}{\delta !}\sum_{\substack{j_0\in [p] \\ j_0 \not \in \mathcal{NZ}(\vec{i})} }\ \sum_{\substack{j_1\in [p]\backslash \mathcal{NZ}(j_0)\\ j_1 \not \in \mathcal{NZ}(\vec{i}) }}\  \cdots \sum_{\substack{j_\delta\in [p]\backslash \bigcup\limits_{l=0}^{\delta-1}\mathcal{NZ}(j_l)\\ j_\delta \not \in \mathcal{NZ}(\vec{i})}}1  \numberthis \label{eqn:substract} \\
	=& 
	\frac{1}{\delta !} \sum_{m=0}^\delta \left(\sum_{\substack{j_0\in [p] \\ j_0 \not \in \mathcal{NZ}(\vec{i})} }\ 
	\sum_{\substack{j_1\in [p]\backslash \mathcal{NZ}(j_0)\\ j_1 \not \in \mathcal{NZ}(\vec{i}) }}\  \cdots\right.\\ &\left.\quad \sum_{\substack{j_{m-1}\in [p]\backslash \bigcup\limits_{l=0}^{m-2}\mathcal{NZ}(j_l)\\ j_{m-1} \not \in \mathcal{NZ}(\vec{i})}}\  \sum_{\substack{j_m\in [p]\backslash \bigcup\limits_{l=0}^{m-1}\mathcal{NZ}(j_l)\\ j_m \in \mathcal{NZ}(\vec{i})}} \ \ \sum_{\substack{j_{m+1}\in [p]\backslash \bigcup\limits_{l=0}^{m}\mathcal{NZ}(j_l)}}  \cdots \sum_{\substack{j_\delta\in [p]\backslash \bigcup\limits_{l=0}^{\delta-1}\mathcal{NZ}(j_l)}}1 \right)\\
	= & \frac{1}{\delta !} \sum_{\substack{j_0\in [p] \\ j_0 \in \mathcal{NZ}(\vec{i})} }\ \sum_{\substack{j_1\in [p]\backslash \mathcal{NZ}(j_0) }}\   \cdots \sum_{\substack{j_\delta\in [p]\backslash \bigcup\limits_{l=0}^{\delta-1}\mathcal{NZ}(j_l)}}1+
	\frac{1}{\delta !} \sum_{m=1}^\delta \left( \right.\\
	&\left.\sum_{\substack{ j_m \in \mathcal{NZ}(\vec{i})}}  \sum_{\substack{j_0\in [p] \\ j_0 \not \in \mathcal{NZ}(\vec{i})\\ j_0\not \in \mathcal{NZ}(j_m) } }\ \sum_{\substack{j_1\in [p]\backslash \mathcal{NZ}(j_0)\\ j_1 \not \in \mathcal{NZ}(\vec{i}) \\ j_1 \not \in \mathcal{NZ}(j_m) }}\  \cdots \sum_{\substack{j_{m-1}\in [p]\backslash \bigcup\limits_{l=0}^{m-2}\mathcal{NZ}(j_l)\\ j_{m-1} \not \in \mathcal{NZ}(\vec{i})\\ j_{m-1}\not \in \mathcal{NZ}(j_m)}}\   \ \ \sum_{\substack{j_{m+1}\in [p]\backslash \bigcup\limits_{l=0}^{m}\mathcal{NZ}(j_l)}}  \cdots \sum_{\substack{j_\delta\in [p]\backslash \bigcup\limits_{l=0}^{\delta-1}\mathcal{NZ}(j_l)}}1 \right) \\
	\geq & \frac{1}{\delta !} (p-p_{\vec{i}}) \prod_{\beta=1}^\delta (p-\beta \kkkk)+
	\frac{1}{\delta !} \sum_{m=1}^\delta (p-p_{\vec{i}}) \left(\prod_{\alpha=1}^m (p_{\vec{i}}-\alpha \kkkk)\right) \left(\prod_{\beta=m+1}^\delta (p-\beta \kkkk)\right) \\
	\geq & 
	\frac{(\delta+1)}{\delta !}  (p-p_{\vec{i}}) \prod_{\alpha=1}^\delta (p_{\vec{i}}-\alpha \kkkk), \numberthis \label{eqn:neidiagonal}
	\end{align*}
	where the second equality follows by writing \eqref{eqn:substract} as a telescoping sum with the convention that the summation over $j_{-1}$ for $m=0$ and the summation over $j_{\delta+1}$ for $m=\delta$ vanish, and the third equality follows from changing the order of the summation for $m\geq 1$. Plugging \eqref{eqn:neidiagonal} into \eqref{eqn:neinotdiagonalsize},
	\begin{align*}
	&\left|\{\vec{j}\in C^<\backslash T_{\vec{i}}: \Sigmama_{\vec{j}} \text{ not diagonal} \}  \right|\\
	\leq &  \frac{1}{\delta !} (\delta+1)   (p-p_{\vec{i}}) \left(\prod_{\alpha=0}^{\delta-1} (p-\alpha) - \prod_{\alpha=0}^{\delta-1} (p_{\vec{i}}-1-\alpha) + \prod_{\alpha=1}^{\delta} (p_{\vec{i}}-\alpha)-  \prod_{\alpha=1}^\delta (p_{\vec{i}}-\alpha \kkkk) \right)\\
	\leq & \frac{1}{\delta !} (\delta+1)   (p-p_{\vec{i}}) \left(\delta p^{\delta-1} (p-p_{\vec{i}}+1) + \frac{\delta(\delta+1)}{2} p^{\delta-1} (\kkkk-1) \right) \\  
	\leq & \frac{3\delta(\delta+1)^3}{\delta !}     p^{\delta-1} \kkkk^2, \numberthis \label{eqn:neinotdiasize2}
	\end{align*}
	where the second inequality follows from Lemma \ref{basine} \ref{item:basineh} and Lemma \ref{basine} \ref{item:basinei}, and the last inequality follows from \eqref{eqn:p'upperbou}. 

	Then for any $\vec{i}\in C^<$,
	\begin{align*}
	 \sum_{\vec{j}\in C^<\backslash T_{\vec{i}}}\mu_{n,\delta+1}\left(\Sigmama_{\vec{j}}\right) \leq & \sum_{\substack{\vec{j}\in C^<\backslash T_{\vec{i}}\\ \Sigmama_{\vec{j}} \text{ diagonal}}}1 + \mu_{n,\delta+1}\left(\Sigmama\right) \sum_{\substack{\vec{j}\in C^<\backslash T_{\vec{i}}\\ \Sigmama_{\vec{j}} \text{ not diagonal}}}1\\
	 \leq &  \frac{(\delta+1)^2}{\delta !}  p^{\delta}\kkkk	+ \mu_{n,\delta+1}\left(\Sigmama\right) \frac{3\delta(\delta+1)^3}{\delta !}     p^{\delta-1} \kkkk^2, \\
	 \leq &   \frac{3\delta(\delta+1)^3}{\delta !}     p^{\delta} \kkkk\left(1+ \mu_{n,\delta+1}\left(\Sigmama\right) \frac{\kappa}{p} \right), \numberthis \label{eqn:sumti}
	\end{align*}
	where the first inequality follows from Lemma \ref{lem:densitybou} \ref{item:densityboua}, and the second inequality follows from \eqref{eqn:neiborsize}, \eqref{eqn:neinotdiasize2}.
	
	Then following \eqref{eqn:b1uppbounew}, 	
	\begin{align*}
	b_1 & \leq \sum_{\vec{i}\in C^{<}} \  \sum_{\vec{j}\in C^<\backslash T_{\vec{i}}} \mu_{n,\delta+1}\left(\Sigmama_{\vec{i}}\right)\mu_{n,\delta+1}\left(\Sigmama_{\vec{j}}\right) (2P_\nn(r_\rho))^{2\delta}\\
	&\leq \frac{p^{\delta+1}}{\delta !}\left(1+\delta^2\mu_{n,\delta+1}(\Sigmama)\frac{\kappa-1}{p}\right)
	  \frac{3\delta(\delta+1)^3}{\delta !}     p^{\delta} \kkkk\left(1+ \mu_{n,\delta+1}\left(\Sigmama\right) \frac{\kappa}{p} \right) (2P_\nn(r_\rho))^{2\delta} \\ 
	& \leq \left(3\frac{\delta^3(\delta+1)^3}{(\delta!)^2}(2p^{1+\frac{1}{\delta}}P_\nn(r_\rho))^{2\delta}\right)\  \frac{\kkkk}{p} \left(1+\mu_{n,\delta+1}(\Sigmama)\frac{\kappa}{p}\right)^2, \numberthis \label{eqn:b1upperbou3}
	\end{align*}
	where the first inequality follows from Lemma \ref{item:prod}, the second inequality follows from Lemma \ref{lem:summuupp} and \eqref{eqn:sumti}. 
	
	\vspace{2mm}

	\noindent{\textbf{Part III. Upper bound for } $b_2$}\\
	Let $\Kc_w(\vec{i})$ and $D_{\vec{i}}^m$ be the same as in Subsection \ref{sec:quadraticterms} with $q=\delta$. By Lemma \ref{lem:densitybou} \ref{item:densitybouc}, Lemma  \ref{item:prod} and Lemma \ref{lem:eleine},
	it follows that the conditions in Lemma \ref{lem:b2} with $q=\delta$ and $\theta_{\vec{i},\vec{j}}=\Phi_{\vec{i}}\Phi_{\vec{j}}$ are satisfied with $a=1$, $b=2P_n(2r_\rho)1(\delta\geq 2)+2P_n(r_\rho)1(\delta= 1)$ and $z=2P_n(r_\rho)$. Moreover, $b/z\leq 2^{n-2}1(\delta\geq 2)+1$ by Lemma \ref{pderivative} \ref{item:pderivatived}. Thus by Lemma \ref{lem:b2} with $q=\delta$ and $\theta_{\vec{i},\vec{j}}=\Phi_{\vec{i}}\Phi_{\vec{j}}$, 
	\begin{align*}
	b_2 & \leq p(2pP_n(r_\rho))^{\delta+1}\left(1+ \mu_{n,2\delta+2}\left(\Sigmama\right) (3\delta^2)\frac{\kappa-1}{p}\right)\left(1+2pP_n(r_\rho)\right)^{\delta-1} \delta\frac{5+2^{n-1}1(\delta \geq 2)}{(\delta-1)!}+\\
	&\quad\quad\quad \frac{(\delta+1)^2}{(\delta!)^2}p^2(2pP_n(r_\rho))^{2\delta}\mu_{n,2\delta+2}\left(\Sigmama\right)\frac{\kappa-1}{p}. \numberthis \label{eqn:b2upperbou2new} \\
	\leq & \gamma^{\delta+1}p^{-\frac{1}{\delta}}\left(1+ \mu_{n,2\delta+2}\left(\Sigmama\right) (3\delta^2)\frac{\kappa-1}{p}\right)\left(1+\gamma p^{-\frac{1}{\delta}}\right)^{\delta-1} \delta\frac{5+2^{n-1}1(\delta \geq 2)}{(\delta-1)!}+\\
	&\quad\quad\quad \frac{(\delta+1)^2}{(\delta!)^2}\gamma^{2\delta}\mu_{n,2\delta+2}\left(\Sigmama\right)\frac{\kappa-1}{p}, \numberthis \label{eqn:b2upperbou2}
	\end{align*} 
    where the last step follows from $2p^{1+\frac{1}{\delta}}P_n(r_\rho)\leq \gamma$.

	By combining \eqref{eqn:compoibou1}, \eqref{eqn:lambdameandif}, \eqref{eqn:lambda0uppbou}, \eqref{eqn:b1upperbou3}, \eqref{eqn:b2upperbou2}, together with the assumption $2p^{1+\frac{1}{\delta}}P_\nn(r_\rho)\leq \gamma$,
	\begin{align*}
	&d_{\TV}\left(\mathscr{L}\left(N_{E_\delta}^{(\Rma)}\right), \CP(\lambda,\bm{\zeta})\right) \\
	\leq & C_{n,\delta,\gamma} \left(C'_{\delta,\gamma}\right)^{\mu_{n,\delta+1}\left(\Sigmama\right)(\kappa-1)/p}\left( \mu_{n,2\delta+2}\left(\Sigmama\right)\kappa/p \left(1 + \mu_{n,2\delta+2}\left(\Sigmama\right) \left(\kappa/p\right)^2 \right)  + p^{-\frac{1}{\delta}} \right) ,
	\end{align*}
	where 
	\begin{equation}
	C'_{\delta,\gamma} = \exp\left(\gamma^{\delta}\frac{\delta+1}{(\delta-1)!}\right), \label{eqn:C'delgamdef}
	\end{equation}
	and
	\begin{equation}
	C_{n,\delta,\gamma} = C\frac{\delta^6+\delta^22^{n-1}1(\delta \geq 2)}{\delta!} \gamma^{\delta+1}(1+\gamma)^\delta C'_{\delta,\gamma}. \label{eqn:Cndelgamdef}
	\end{equation}
\end{proof}

\subsection{Proofs in Subsection \ref{sec:allclose}  } \label{sec:proofofallclose}

\subsubsection{Proof of Lemma \ref{lem:6quantitiesine}}
\begin{proof}[Proof of Lemma \ref{lem:6quantitiesine}]
    $N_{\breve{V}_\delta}^{(\Rma)} \leq N_{V_\delta}^{(\Rma)} \leq N_{E_\delta}^{(\Rma)}$ follows trivially from their definitions. It remains to show
    \begin{equation}
    N_{E_\delta}^{(\Rma)}\leq  (\delta+1)N_{E_{\delta+1}}^{(\Rma)} + N_{\breve{V}_\delta}^{(\Rma)}. \label{eqn:tobeproved}
    \end{equation}
     To see this, consider $\delta\geq 2$ and any vertex $i$ and denote its degree by $m$. If $m<\delta$, then it contributes zero to both sides of \eqref{eqn:tobeproved}. If $m=\delta$, then it contributes $1$ to both sides of \eqref{eqn:tobeproved}. If $m>\delta$, it contributes $\binom{m}{\delta}$ to left hand side of \eqref{eqn:tobeproved}, while contributes $(\delta+1)\binom{m}{\delta+1}=(m-\delta)\binom{m}{\delta}$. The above observation proves \eqref{eqn:tobeproved}. The case $\delta=1$ is similar and is omitted.
    
    The above proof indeed applies to any graph and, in particular, the empirical partial correlation graph. So the second equation in the statement of the lemma holds. 
\end{proof}

\subsubsection{Proof of Proposition \ref{thm:allclose} \ref{item:allclosea}}

By \eqref{eqn:corexaatupp}, it suffices to establish an upper bound on $\E N_{E_{\delta+1}}^{(\Rma)}$.

\begin{lem} \label{lem:NEbou}
{Let $\Xma\sim \mathcal{VE}(\bm{\mu},\Sigmama,g)$. Suppose $\Sigmama$ is row-$\kappa$ sparse. Let $\ell\in [p-1]$. Then
$$
\E N_{E_{\ell}}^{(\Rma)} \leq \frac{1}{\ell!}\left(1+\ell^2\mu_{n,\ell+1}(\Sigmama)\frac{\kappa-1}{p}\right)p\left( 2pP_\nn(r_\rho)\right)^\ell.
$$}
\end{lem}
\begin{proof}
	\begin{align*}
	\E N_{E_{\ell}}^{(\Rma)} = & \sum_{\vec{i}\in C_{\ell}^< }\E  \prod_{j=1}^{\ell} \Phi_{i_0i_j}^{(\Rma)}  \\
	\leq & \sum_{\vec{i}\in C_\ell^<}\mu_{n,\ell+1}\left(\Sigmama_{\vec{i}}\right) \left( 2P_\nn(r_\rho)\right)^\ell \\
	\leq &   \frac{1}{\ell!}\left(1+\ell^2\mu_{n,\ell+1}(\Sigmama)\frac{\kappa-1}{p}\right)p\left( 2pP_\nn(r_\rho)\right)^\ell,
	\end{align*}
	where the first inequality follows from Lemma \ref{item:prod}, and the second inequality follows from Lemma \ref{lem:summuupp}.
\end{proof}

\begin{proof}[Proof of Proposition \ref{thm:allclose} \ref{item:allclosea}] It follows from \eqref{eqn:corexaatupp}, Lemma \ref{lem:NEbou} and Lemma \ref{prop:tvtomean}.
\end{proof}

\subsubsection{Proof of Proposition \ref{thm:allclose} \ref{item:allclosec}}
\label{sec:proofthmallclosec}

Similar to \eqref{eqn:NEdeltaRdef} and \eqref{eqn:Nedgedef1}, denote
\begin{equation*}
\Phi_{\vec{i}}^{(\Rma)}= \prod_{j=1}^\delta \Phi_{i_0i_j}^{(\Rma)} = 1\left(\bigcap_{j=1}^{\delta}\{\textbf{dist}(\uve_{i_0},\uve_{i_j})\leq r_\rho\}\right). \label{eqn:NEdeltaPdef}
\end{equation*}
Then by definition
\begin{equation}
\Nedged^{(\Pma)}  =\sum_{\vec{i}\in C_\delta^<}\ \  \Phi_{\vec{i}}^{(\Pma)}. \label{eqn:NEdeltaPdef1}
\end{equation}
\index{$\Nedged^{(\Rma)}$}

By \eqref{eqn:Nedgedef1} and \eqref{eqn:NEdeltaPdef1},
$$
\left|\Nedged^{(\Pma)}-\Nedged^{(\Rma)}\right| \leq \sum_{\vec{i}\in C_\delta^<}|\Phi_{\vec{i}}^{(\Pma)} - \Phi_{\vec{i}}^{(\Rma)}|.
$$
The next three lemmas establish upper bound on $|\Phi_{\vec{i}}^{(\Pma)} - \Phi_{\vec{i}}^{(\Rma)}|$. 

We may suppose $\Sigmama$ is $(\tau,\kappa)$ sparse throughout this proof and the proof of Proposition \ref{thm:allclose} \ref{item:allclosed} since the conclusion is invariant to permutation of the variables by Remark \ref{rem:firsttaurworemark}. As a result, the U-scores may be partitioned into $\hat{\Uma}\in \R^{(n-1)\times \tau}$ consisting of the first $\tau$ columns and $\check{\Uma}\in \R^{(n-1)\times (p-\tau)}$ consisting the remaining $p-\tau$ columns.

Denote $[\tau]=\{1,2,\cdots,\tau \}$.  
Define a matrix $\check{\Bma}$ by
\begin{equation}
\check{\Bma} = \frac{\nn-1}{p-\tau}\check{\Uma}[\check{\Uma}]^\top
=
\frac{\nn-1}{p-\tau}\Sum_{i\in [p]\backslash [\tau] }^{p}\uve_i\uve_i^\top.  \label{Bp'def}
\end{equation}
Denote $\check{\Qma}=\sqrt{\nn-1}\check{\Uma}$. Observe that $\check{\Qma}$ has exactly $p-\tau$ independent columns and each column $\sqrt{\nn-1}\uve_i\sim \unif(\sqrt{\nn-1}S^{\nn-2})$. These observations immediately give us part \ref{item:bpproa} of the following.

\begin{lem}
	\label{lem:bppro}
	Let $\{\uve_\alpha\}_{\alpha =1 }^ p$ be columns of $\Uma$ defined in Section \ref{sec:scorerep}. Let $\check{\Bma}$ be defined as in equation \eqref{Bp'def}. 
	\begin{enumerate}[label=(\alph*)]
		\item \label{item:bpproa} Suppose $\Sigmama$ is $(\tau,\kappa)$ sparse. $\check{\Bma}=\frac{1}{p-\tau}\Qma\Qma^\top$, where $\Qma\in R^{(\nn-1)\times (p-\tau)}$ has independent columns with each column distributed as $\unif(\sqrt{\nn-1}S^{\nn-2})$.
		
		\item \label{item:bpprob} $|\lambdamax\left(\frac{p}{p-\tau}\Bma\right)-\lambdamax(\check{\Bma})|\leq \frac{\nn-1}{p-\tau}\qq$, and $\lambdamin\left(\frac{p}{p-\tau}\Bma\right)\geq \lambdamin(\check{\Bma})$.
	\end{enumerate}
\end{lem}
\begin{proof}
	(b) Recall $\Bma=\frac{\nn-1}{p}\sum_{i=1}^p \uve_i\uve_i^\top$. Then,
	$$
	\frac{p}{p-\tau}\Bma-\check{\Bma} = \frac{\nn-1}{p-\tau}\sum_{i\in [\tau] }\uve_i\uve_i^\top.
	$$
	By Lemma \ref{perthe} \ref{item:perthea}, we have:
	\begin{align*}
	|\lambdamax\left(\frac{p}{p-\tau}\Bma\right)-\lambdamax\left(\check{\Bma}\right)| &\leq \left\|\frac{\nn-1}{p-\tau}\sum_{i\in [\tau] }\uve_i\uve_i^\top\right\|_2\\
	&\leq \frac{\nn-1}{p-\tau}\sum_{i\in [\tau] }\|\uve_i\uve_i^\top\|_2 \\
	&\leq \frac{\nn-1}{p-\tau}\qq,
	\end{align*}
	where for the last inequality, we use the fact that $\uve_i \in S^{\nn-2}$. Moreover, by Lemma \ref{perthe} \ref{item:perthec}, we get $\lambdamin\left(\frac{p}{p-\tau}\Bma\right)\geq \lambdamin(\check{\Bma})$.
\end{proof}

Denote $h_0\left(\check{\Bma}\right)= \frac{\lambdamax\left(\check{\Bma}\right)+\frac{n-1}{p-\tau}\tau}{\lambdamin\left(\check{\Bma}\right)}=\frac{S_\text{max}\left(\check{\Bma}\right)+\frac{n-1}{p-\tau}\tau}{S_\text{min}\left(\check{\Bma}\right)}$ to be the perturbational condition number of $\check{\Bma}$, where $\lambdamax\left(\check{\Bma}\right)$, $\lambdamin\left(\check{\Bma}\right)$, $S_{\text{max}}\left(\check{\Bma}\right)$, and $S_{\text{min}}\left(\check{\Bma}\right)$ are respectively the largest eigenvalue, smallest eigenvalue, largest singular value and smallest singular value of $\check{\Bma}$. 

\begin{lem}
	\label{uyrelation}
	Suppose $p\geq n$. Let $\{\uve_\alpha\}_{\alpha =1 }^ p$ and $\{\yve_\alpha\}_{\alpha =1 }^p$ be defined as in Section \ref{sec:scorerep}.
	Consider distinct $i,j$ satisfying $1\leq i, j \leq p$. Then with probability $1$,
	$$
	\frac{1}{h_0\left(\check{\Bma}\right)}\|\uve_i-\uve_j\|_2 \leq \|\yve_i-\yve_j\|_2 \leq h_0\left(\check{\Bma}\right)\|\uve_i-\uve_j\|_2,
	$$
	and
	$$
	\frac{1}{h_0\left(\check{\Bma}\right)}\|\uve_i+\uve_j\|_2 \leq \|\yve_i+\yve_j\|_2 \leq h_0\left(\check{\Bma}\right)\|\uve_i+\uve_j\|_2.
	$$
\end{lem}

\begin{proof}
	Recall $\yve_\alpha=\bar{\yve}_\alpha/\|\bar{\yve}_\alpha\|_2$ and $\bar{\yve}_\alpha = \Bma^{-1}\uve_\alpha$ $a.s.$, for $\alpha=i,j$. Apply the upper bound in Lemma \ref{isoine},
	\begin{align*}
	\|\yve_i-\yve_j\|_2  
	&\leq \frac{\lambdamax(\Bma^{-1})}{\lambdamin(\Bma^{-1})}\|\uve_i-\uve_j\|_2 \quad a.s. \\
	& = 
	\frac{\lambdamax\left(\frac{p}{p-\tau}\Bma\right)}{\lambdamin\left(\frac{p}{p-\tau}\Bma\right)}\|\uve_i-\uve_j\|_2  \\ 
	& \leq \frac{\lambdamax\left(\check{\Bma}\right)+\frac{n-1}{p-\tau}\tau}{\lambdamin\left(\check{\Bma}\right)}\|\uve_i-\uve_j\|_2 ,
	\end{align*}
	where  the last inequality follows from Lemma \ref{lem:bppro} \ref{item:bpprob}. The lower bound of the first desired display follows similarly, by the lower bound in Lemma \ref{isoine}.
	
	The second desired expression follows analogously.
\end{proof}

For $\{i,j\}\in [p]$ with $i\not =j$,  $q\in\{-1,+1\}$, define
\begin{align*}
S_{ij}^{(q)}(r_\rho)& =\{\|\yve_i-q\yve_j\|_2 \leq r_\rho \}, \quad F_{ij}^{(q)}(r_\rho)=\{\|\uve_i-q\uve_j\|_2 \leq r_\rho \},  \numberthis \label{eqn:Fijdef} \\
G_{ij}^{(q)}(r_\rho)& =\left\{\|\uve_i-q\uve_j\|_2 \leq \frac{1}{h_0\left(\check{\Bma}\right)}r_\rho \right \}, \quad H_{ij}^{(q)}(r_\rho)=\left \{\|\uve_i-q\uve_j\|_2 \leq h_0\left(\check{\Bma}\right)r_\rho \right\}.
\end{align*}
Define $F_{ij}(r_\rho)=F_{ij}^{(-1)}(r_\rho)\cup F_{ij}^{(+1)}(r_\rho)$. $G_{ij}(r_\rho)$, $H_{ij}(r_\rho)$, $S_{ij}(r_\rho)$ are defined similarly. Using these notation, then $\Phi_{ij}^{(\Pma)}(\rho)=1\left(S_{ij}(r_\rho)\right)$, and $\Phi_{ij}^{(\Rma)}(\rho)=1\left(F_{ij}(r_\rho)\right)$. For $\vec{i}\in C_\delta^<$, denote 
$$
H_{\vec{i}}(r_\rho)=\bigcap_{\ell=1}^{\delta} H_{i_0i_\ell}(r_\rho),\quad H_{\vec{i},-m}(r_\rho)=\bigcap_{\substack{\ell=1\\ \ell\neq m}}^{\delta} H_{i_0i_\ell}(r_\rho).
$$

When it is clear from the context, the dependence of $r_\rho$ for the above quantities will be suppressed. By Lemma \ref{uyrelation}, with probability $1$,
\begin{equation}
G_{ij}^{(q)} \subset S_{ij}^{(q)}\subset H_{ij}^{(q)}, \quad G_{ij}^{(q)} \subset F_{ij}^{(q)}\subset H_{ij}^{(q)}. \label{eqn:FGJSsetinclusion}
\end{equation}

\begin{lem}
\label{lem:inddifbou}
	Suppose $p\geq n$. Consider $\delta\in [p-1]$. For any $\vec{i}\in C_\delta^<$, with probability $1$,
	$$
	\left| \Phi_{\vec{i}}^{(\Pma)}-\Phi_{\vec{i}}^{(\Rma)} \right|\leq \xi_{\vec{i}},
	$$
	where 
	$$
\xi_{\vec{i}} := 1\left( \bigcup\limits_{m=1}^\mylll
\left ( \left( H_{i_0i_m} \big \backslash G_{i_0i_m}   \right) \bigcap  H_{\vec{i},-m} \right )\right). 
$$
\end{lem}
\begin{proof}
	Notice 
	$\Phi_{\vec{i}}^{(\Rma)}=1\left(\bigcap\limits_{m=1}^\mylll F_{i_0i_m}\right)$ and $\Phi_{\vec{i}}^{(\Pma)}=1\left(\bigcap\limits_{m=1}^\mylll S_{i_0i_m}\right)$. Let $\bigtriangleup$ denote the symmetrization difference of two sets. Then
	\begin{align*}
	\left| \Phi_{\vec{i}}^{(\Pma)}-\Phi_{\vec{i}}^{(\Rma)} \right|=  
	1\left(\left(\bigcap\limits_{m=1}^\mylll F_{i_0i_m}\right)\bigtriangleup \left(\bigcap\limits_{m=1}^\mylll S_{i_0i_m}\right) \right) 
	\leq \xi_{\vec{i}},
	\end{align*}
	where the inequality follows from \eqref{eqn:FGJSsetinclusion} and Lemma \ref{setrelation} \ref{prop:setinclusiona}.	
\end{proof}

To obtain an upper bound on the expectation of $\xi_{\vec{i}}$, we first bound the expectation on a high-probability set. Define the set $\Ec(t)$, with $t$ being a parameter to be determined, by
\begin{align*}
\Ec(t) &= \left\{ \left[1-C_1\left(\sqrt{\frac{\nn-1}{p-\tau}}+\frac{t}{\sqrt{p-\tau}}\right)\right]^2 \leq \lambdamin(\check{\Bma})\right\}\bigcap \\
&\bigcap \left\{\lambdamax(\check{\Bma})\leq  \left[1+C_1\left(\sqrt{\frac{\nn-1}{p-\tau}}+\frac{t}{\sqrt{p-\tau}}\right)\right]^2 \right\}, \numberthis \label{Epsdef}
\end{align*}
to be the set such that \eqref{eqn:eigupplowbou} in Lemma \ref{proeig} holds, i.e. the constant $C_1$ in $\Ec(t)$ is the same constant as $C$ in \eqref{eqn:eigupplowbou}. By Lemma \ref{lem:bppro} \ref{item:bpproa} and Lemma \ref{proeig},
\begin{equation}
\P(\Ec^c(t))\leq 2\exp(-c_1t^2). \label{eqn:proepst}
\end{equation}
{ Since $\tau\leq\frac{p}{2}$},
\begin{equation}
\frac{n-1}{p-\tau}\tau \leq 2(n-1)\frac{\tau}{p}, \label{eqn:ntaupupp}
\end{equation}
and
\begin{equation}
C_1\left(\sqrt{\frac{\nn-1}{p-\tau}}+\frac{t}{\sqrt{p-\tau}}\right)\leq \sqrt{2}C_1\left(\sqrt{\frac{\nn-1}{p}}+\frac{t}{\sqrt{p}}\right). \label{eqn:c1nptauupp}
\end{equation}
Moreover, on $\Ec(t)$, and assuming 
\begin{equation}
\sqrt{2}C_1\left(\sqrt{\frac{\nn-1}{p}}+\frac{t}{\sqrt{p}}\right) \leq \frac{1}{2}, \label{eqn:postempass}
\end{equation}
one has
\begin{multline*}
h_0(\Bma)  \leq \frac{\left(1+C_1\left(\sqrt{\frac{\nn-1}{p-\tau}}+\frac{t}{\sqrt{p-\tau}}\right)\right)^2+\frac{n-1}{p-\tau}\tau}{\left(1-C_1\left(\sqrt{\frac{\nn-1}{p-\tau}}+\frac{t}{\sqrt{p-\tau}}\right)\right)^2}\\
\leq 1+ 16 \sqrt{2}C_1\left(\sqrt{\frac{\nn-1}{p}}+\frac{t}{\sqrt{p}}\right) + 8(n-1)\frac{\tau}{p}
\vcentcolon= \theta_1(t), \numberthis \label{ctempdef}
\end{multline*}
where the second inequality follows from \eqref{eqn:ntaupupp}, \eqref{eqn:c1nptauupp} and Lemma \ref{basine} \ref{item:basinej}.

For $\vec{i}\in C_\delta^<$, denote 
$$
F_{\vec{i}}(r_\rho)=\bigcap_{\ell=1}^{\delta} F_{i_0i_\ell}(r_\rho),\quad F_{\vec{i},-m}(r_\rho)=\bigcap_{\substack{\ell=1\\ \ell\neq m}}^{\delta} F_{i_0i_\ell}(r_\rho).
$$

\begin{lem}
\label{lem:inddifuppexpbou}
Let $p\geq n\geq 4$, $\delta\in [p-1]$ and $\Xma\sim \mathcal{VE}(\bm{\mu},\Sigmama,g)$. Suppose $\Sigmama$, after some row-column permutation, is $(\tau,\kappa)$ sparse with $\tau\leq \frac{p}{2}$. Let $t$ be any positive number, and suppose \eqref{eqn:postempass} holds. Then for any $\vec{i}\in C_\delta^<$, with probability $1$,
$$
\xi_{\vec{i}} 1\left(\Ec(t)\right) \leq \eta_{\vec{i}}(t),
$$
where 
\begin{align}
    \eta_{\vec{i}}(t) := &1 \left(  \bigcup_{m=1}^\mylll\left( \left(F_{i_0i_m}(\theta_1(t)r_\rho)\big \backslash F_{i_0i_m}\left(\frac{r_\rho}{\theta_1(t)}\right) \right)  \bigcap  F_{\vec{i},-m}(\theta_1(t)r_\rho) \right)\right). \label{eqn:etavecidef}
\end{align}
Moreover,
\begin{align*}
&\E 1  \left( \left(F_{i_0i_m}(\theta_1(t)r_\rho)\big \backslash F_{i_0i_m}\left(\frac{r_\rho}{\theta_1(t)}\right) \right)  \bigcap  F_{\vec{i},-m}(\theta_1(t)r_\rho) \right)\\
\leq & \mu_{n,\delta+1}(\Sigmama_{\vec{i}}) 2\left(P_\nn(r_\rho \theta_1(t))-P_\nn\left(\frac{r_\rho}{\theta_1(t)} \right)\right) \left(2P_\nn\left(\theta_1(t)r_\rho\right)\right)^{\mylll-1}, \numberthis \label{eqn:etasingletermbound}
\end{align*}
and
	\begin{align*}
	 \E \eta_{\vec{i}}(t) 
	\leq &  \delta \mu_{n,\delta+1}(\Sigmama_{\vec{i}}) 2\left(P_\nn(r_\rho \theta_1(t))-P_\nn\left(\frac{r_\rho}{\theta_1(t)} \right)\right) \left(2P_\nn\left(\theta_1(t)r_\rho\right)\right)^{\mylll-1} \\ 
	\leq &  \mu_{n,\delta+1}(\Sigmama_{\vec{i}}) \mylll \nn  (\theta_1(t))^{\nn\delta}\left(\theta_1(t)-\frac{1}{\theta_1(t)}\right) \left( 2P_\nn\left(r_\rho\right) \right)^{\delta}.
	\end{align*}
\end{lem}

\begin{proof} 
By \eqref{ctempdef}, $H_{ij}(r_\rho) \cap \Ec(t)\subset F_{ij}(\theta_1(t)r_\rho)$ and $G_{ij}(r_\rho) \cap \Ec(t)\supset F_{ij}\left(\frac{r_\rho}{\theta_1(t)}\right)$. Then
	\begin{align*}
	\xi_{\vec{i}} 1\left(\Ec(t)\right)
	\leq  \eta_{\vec{i}}(t). \numberthis \label{eqn:A1upp1}
	\end{align*}

	\begin{align*}
	&\E 1  \left( \left(F_{i_0i_m}(\theta_1(t)r_\rho)\big \backslash F_{i_0i_m}\left(\frac{r_\rho}{\theta_1(t)}\right) \right)  \bigcap  F_{\vec{i},-m}(\theta_1(t)r_\rho) \right) \\
	 \leq & \mu_{n,\delta+1}(\Sigmama_{\vec{i}})    \P \left( \left(\bigcup\limits_{q\in\{-1,+1\}}\left\{ \frac{r_\rho}{\theta_1(t)} <\|\uve'_{i_0}-q\uve'_{i_m}\|_2 \leq \theta_1(t)r_\rho \right\}\right)  \bigcap \right.\\
	& \ \ \ \ \ \ \ \left. \left(\bigcap_{\substack{\alpha=1 \\ \alpha\not=m}}^\mylll \left(\bigcup\limits_{q\in\{-1,+1\}} \{ \|\uve'_{i_0}-q\uve'_{i_\alpha}\|_2\leq \theta_1(t)r_\rho \} \right)\right)\right) , \numberthis \label{eqn:expindupp}
	\end{align*}
	where the last inequality follows from Lemma \ref{lem:densitybou} \ref{item:densitybouc} with 
	$$
	\uve'_{i_0},\uve'_{i_1},\cdots,\uve'_{i_\delta}\overset{\text{i.i.d.}}{\sim} \unif(S^{n-2}).
	$$
	For any $\bm{w}\in S^{n-2}$, define $\Omega_{\wve}^{(q)} :=\{\vve \in S^{\nn-2}:\frac{1}{\theta_1(t)} r_\rho < \|\vve-q\wve\|_2 \leq r_\rho \theta_1(t)\}$. Then 
	$$\P\left(\uve'_{im}\in  \bigcup_{q\in\{-1,+1\}} \Omega_{\wve}^{(q)} \right)=2\left(P_\nn(r_\rho \theta_1(t))-P_\nn\left(\frac{1}{\theta_1(t)}r_\rho \right)\right).$$ 
	By conditioning on $\uve'_{i_0}$, the term in right hand side of \eqref{eqn:expindupp} equals to
	\begin{align*}
	\mu_{n,\delta+1}(\Sigmama_{\vec{i}}) 2\left(P_\nn(r_\rho \theta_1(t))-P_\nn\left(\frac{1}{\theta_1(t)}r_\rho \right)\right) \left(2P_\nn\left(\theta_1(t)r_\rho\right)\right)^{\mylll-1},
	\end{align*}
	which then proves \eqref{eqn:etasingletermbound}.
	
	By the union bound,
	\begin{align*}
	    &\E \eta_{\vec{i}}(t) \\
	    \leq &
	\sum\limits_{m=1}^\mylll \E 1 \left( \left(F_{i_0i_m}(\theta_1(t)r_\rho)\big \backslash F_{i_0i_m}\left(\frac{r_\rho}{\theta_1(t)}\right) \right)  \bigcap  F_{\vec{i},-m}(\theta_1(t)r_\rho) \right) \\
	\overset{(*)}{\leq} & 
	\delta \mu_{n,\delta+1}(\Sigmama_{\vec{i}}) 2\left(P_\nn(r_\rho \theta_1(t))-P_\nn\left(\frac{1}{\theta_1(t)}r_\rho \right)\right) \left(2P_\nn\left(\theta_1(t)r_\rho\right)\right)^{\mylll-1}\\
	\overset{(**)}{\leq} & \delta\mu_{n,\delta+1}(\Sigmama_{\vec{i}})
	2(n-2)P_\nn(r_\rho)\left(\theta_1(t)\right)^{n-3}  \left(\theta_1(t)-\frac{1}{\theta_1(t)}\right)\left(2P_\nn(r_\rho \theta_1(t))\right)^{\delta-1}\\
	\overset{(***)}{\leq} & \delta\mu_{n,\delta+1}(\Sigmama_{\vec{i}})
	2(n-2)P_\nn(r_\rho) \left(\theta_1(t)\right)^{n-3}\left(\theta_1(t)-\frac{1}{\theta_1(t)}\right) \left(\left(\theta_1(t)\right)^{n-2} 2P_\nn\left(r_\rho\right) \right)^{(\delta-1)},  
	\end{align*}
	where $(*)$ follows from \eqref{eqn:etasingletermbound}, $(**)$ follows from Lemma \ref{pderivative} \ref{item:pderivativec},
	and $(***)$ follows from Lemma \ref{pderivative}  \ref{item:pderivatived}.

\end{proof}

\begin{lem}
\label{lem:ideltaedgeindicatorbou} 
	Let $p\geq n\geq 4$, $\delta\in [p-1]$ and $\Xma\sim \mathcal{VE}(\bm{\mu},\Sigmama,g)$.
	Let $t$ be any positive number, and suppose \eqref{eqn:postempass} holds. Suppose $\Sigmama$, after some row-column permutation, is $(\tau,\kappa)$ sparse with $\tau\leq \frac{p}{2}$. Then 
	$$
	\left|N_{E_{\delta}}^{(\Pma)}-N_{E_{\delta}}^{(\Rma)} \right| 1\left(\Ec(t) \right) \leq  \sum_{\vec{i}\in C_\delta^<} \eta_{\vec{i}}(t)
	$$
	and
	$$
	\E \sum_{\vec{i}\in C_\delta^<} \eta_{\vec{i}}(t)
	\leq \frac{Cn^2}{(\delta-1)!}\left(1+\delta^2\frac{\kappa-1}{p}\mu_{n,\delta+1}(\Sigmama)\right)  (\theta_1(t))^{\nn\delta}\left(\sqrt{\frac{1}{p}}+\frac{t}{\sqrt{p}} + \frac{\tau}{p}\right) p\left( 2pP_\nn\left(r_\rho\right) \right)^{\delta},
	$$
	where $C$ is an universal constant.		
\end{lem}

\begin{proof}
\begin{align*}
     \left|N_{E_{\delta}}^{(\Pma)}-N_{E_{\delta}}^{(\Rma)} \right| 1\left(\Ec(t) \right) \leq & 
	\sum_{\vec{i}\in C_\delta^<}\ \ \E\left|   \Phi_{\vec{i}}^{(\Pma)}
	 - \Phi_{\vec{i}}^{(\Rma)}\right| 1\left(\Ec(t) \right) \\
	 \leq & \sum_{\vec{i}\in C_\delta^<} \eta_{\vec{i}}(t),
\end{align*}
where the last inequality follows from Lemma \ref{lem:inddifbou} and Lemma \ref{lem:inddifuppexpbou}. 

By Lemma \ref{lem:inddifuppexpbou},
	\begin{align*}
	&  \sum_{\vec{i}\in C_\delta^<} \E \eta_{\vec{i}} \\
	 \leq & 
	 \sum_{\vec{i}\in C_\delta^<}\mu_{n,\delta+1}(\Sigmama_{\vec{i}}) \mylll \nn  (\theta_1(t))^{\nn\delta}\left(\theta_1(t)-\frac{1}{\theta_1(t)}\right) \left( 2P_\nn\left(r_\rho\right) \right)^{\delta}\\
	 \leq & \frac{p^{\delta+1}}{\delta!}\left(1+\delta^2\frac{\kappa-1}{p}\mu_{n,\delta+1}(\Sigmama)\right)\mylll \nn  (\theta_1(t))^{\nn\delta}\left(\theta_1(t)-\frac{1}{\theta_1(t)}\right) \left( 2P_\nn\left(r_\rho\right) \right)^{\delta}\\
	 \leq & \frac{C\nn}{(\delta-1)!}\left(1+\delta^2\frac{\kappa-1}{p}\mu_{n,\delta+1}(\Sigmama)\right)   (\theta_1(t))^{\nn\delta}\left(\sqrt{\frac{\nn}{p}}+\frac{t}{\sqrt{p}} + n\frac{\tau}{p}\right) \left( 2p^{1+\frac{1}{\delta}}P_\nn\left(r_\rho\right) \right)^{\delta} , \numberthis \label{eqn:temptemptemp1}
	\end{align*}
	where the third inequality follows from Lemma \ref{lem:summuupp} and the last inequality follows from Lemma \ref{basine} \ref{item:basinek} and \eqref{ctempdef}. 
\end{proof}

\begin{lem} \label{lem:Nedgedelp}
	Let $p\geq n \geq 4$, $\delta\in [p-1]$ and $\Xma\sim \mathcal{VE}(\bm{\mu},\Sigmama,g)$. Suppose $\Sigmama$, after some row-column permutation, is $(\tau,\kappa)$ sparse with $\tau\leq \frac{p}{2}$.  Suppose  $2p^{1+\frac{1}{\delta}}P_\nn(r_\rho)\leq \gamma$ and $\left(\sqrt{\frac{\nn-1}{p}}+\sqrt{\frac{\delta\ln p}{p}}\right) \leq c$ hold for some positive and \ONE{sufficiently} small universal constant $c$. Then
	\begin{align*}
	\E\left|\Nedged^{(\Pma)}-\Nedged^{(\Rma)}\right|\leq C_{E_\delta}^{(\Pma)} \left(1+\frac{\kappa-1}{p}\mu_{n,\delta+1}(\Sigmama)\right)\left(\sqrt{\frac{\ln p}{p}}+\frac{\tau}{p}\right),
	\end{align*}
	where $C^{(\Pma)}_{E_\delta} $ is  defined in \eqref{eqn:CEatleadel}.
\end{lem}
\begin{proof}
	\begin{align*}
	\E\left|\Nedged^{(\Pma)}-\Nedged^{(\Rma)}\right|
	&\leq  \E \left|\Nedged^{(\Pma)}-\Nedged^{(\Rma)}\right|1\left(\Ec(t)\right) + \binom{p}{1}\binom{p-1}{\delta}\P(\Ec^c(t)), \\
	&\leq \E \left|\Nedged^{(\Pma)}-\Nedged^{(\Rma)}\right|1\left(\Ec(t)\right) + \frac{p^{\delta+1}}{\delta!}2\exp(-c_1t^2), \numberthis \label{eqn:Nedgedelpr}
	\end{align*}
	where the first inequality follows from $0\leq \Nedged^{(\Psima)} \leq \binom{p}{1}\binom{p-1}{\delta}$ for both $\Psima =\Rma$ and $\Psima = \Pma$, and the second inequality follows from \eqref{eqn:proepst}. 
	
	Choose $t=c_\delta\sqrt{\ln p}$ with  $c_\delta=\sqrt{\frac{5\delta}{2c_1}} \geq \sqrt{\left(\frac{3}{2}+\delta\right)/c_1}$ such that 
	$$
	2\exp(-c_1t^2)\leq 2\exp\left(-\left(\frac{3}{2}+\delta\right)\ln p\right)=\frac{2}{p^{\frac{3}{2}+\delta}}.
	$$
	Moreover, for any 
	\begin{equation}
	    c<\frac{1}{2\max\left \{\sqrt{\frac{5}{2c_1}},1\right\}\sqrt{2}C_1}, \label{eqn:cdef}
	\end{equation} 
	the inequality
	$$
	\left(\sqrt{\frac{\nn-1}{p}}+\sqrt{\frac{\delta \ln p}{p}}\right) \leq c
	$$
	implies
	\begin{equation}
	\sqrt{2}C_1\left(\sqrt{\frac{\nn-1}{p}}+c_\delta\sqrt{\frac{\ln p}{p}}\right) \leq \frac{1}{2}, \label{eqn:CEatleadelpre}
	\end{equation}
	which is \eqref{eqn:postempass} with $t= c_\delta \sqrt{\ln p}$. Then apply Lemma \ref{lem:ideltaedgeindicatorbou} with $t=c_\delta \sqrt{\ln p}$ to \eqref{eqn:Nedgedelpr},
	\begin{align*}
	  &\E\left|\Nedged^{(\Pma)}-\Nedged^{(\Rma)}\right|\\
	 \leq & \frac{Cn^2}{(\delta-1)!}\left(1+\delta^2\frac{\kappa-1}{p}\mu_{n,\delta+1}(\Sigmama)\right)  \left(\theta_1\left(c_\delta\sqrt{\ln p}\right)\right)^{\nn\delta}\left(\sqrt{\frac{1}{p}}+\frac{\sqrt{\delta \ln p}}{\sqrt{p}} + \frac{\tau}{p}\right) \gamma^{\delta} + \frac{2}{\delta!\sqrt{p}}\\
	 \leq & \frac{Cn^2\sqrt{\delta}}{(\delta-1)!}\left(1+\delta^2\frac{\kappa-1}{p}\mu_{n,\delta+1}(\Sigmama)\right)  (\theta_1(c_\delta \sqrt{\ln p}))^{\nn\delta}\left(\frac{\sqrt{\ln p}}{\sqrt{p}} + \frac{\tau}{p}\right) \gamma^{\delta} + \frac{2}{\delta!\sqrt{p}} \\
	 \leq & C_{E_\delta}^{(\Pma)} \left(1+\frac{\kappa-1}{p}\mu_{n,\delta+1}(\Sigmama)\right)\left(\sqrt{\frac{\ln p}{p}}+\frac{\tau}{p}\right),
	\end{align*}
	where 
	\begin{align*}
	C_{E_\delta}^{(\Pma)} = & \frac{Cn^2\delta^{\frac{5}{2}}}{(\delta-1)!} \left(\theta_1\left(c_\delta \sqrt{\ln p}\right)\right)^{\nn\delta} \gamma^{\delta}+\frac{2}{\delta!\sqrt{\ln p}}\\
	\leq & \frac{Cn^2\delta^{\frac{5}{2}}}{(\delta-1)!} \left(4n+5\right)^{\nn\delta} \gamma^{\delta}+\frac{2}{\delta!}, \numberthis \label{eqn:CEatleadel}
	\end{align*}
	where the last step follows from $\theta_1\left(c_\delta \sqrt{\ln p}\right)\leq 9+ 4(n-1)=4n+5$ by \eqref{eqn:CEatleadelpre} and $\tau\leq p/2$.
\end{proof}

\begin{proof}[Proof of Proposition \ref{thm:allclose} \ref{item:allclosec}]
	 It follows directly from Lemma \ref{lem:Nedgedelp} and Lemma \ref{prop:tvtomean}.
\end{proof}

\subsubsection{Proof of Proposition \ref{thm:allclose} \ref{item:allclosed} }
By Lemma \ref{lem:6quantitiesine}, 
\[
N_{E_{\delta}}^{(\Pma)}-(\delta+1)N_{E_{\delta+1}}^{(\Pma)} - N_{E_{\delta}}^{(\Rma)} \leq N_{\VVV_\delta}^{(\Pma)}-N_{\VVV_\delta}^{(\Rma)} \leq N_{E_{\delta}}^{(\Pma)}-N_{E_{\delta}}^{(\Rma)} +(\delta+1) N_{E_{\delta+1}}^{(\Rma)},
\]
which implies
\begin{align*}
\left|N_{\VVV_\delta}^{(\Pma)}-N_{\VVV_\delta}^{(\Rma)}\right|\leq  \left|N_{E_{\delta}}^{(\Pma)} - N_{E_{\delta}}^{(\Rma)}\right|+ (\delta+1)\left|N_{E_{\delta+1}}^{(\Pma)}-N_{E_{\delta+1}}^{(\Rma)}\right|+(\delta+1)N_{E_{\delta+1}}^{(\Rma)}. \numberthis \label{eqn:Nvexacdifbou}
\end{align*}

\begin{lem} \label{lem:Nvexadelpr}
		Let $p\geq n \geq 4$, $\delta\in [p-1]$ and $\Xma\sim \mathcal{VE}(\bm{\mu},\Sigmama,g)$. Suppose $\Sigmama$, after some row-column permutation, is $(\tau,\kappa)$ sparse with $\tau\leq \frac{p}{2}$.  Suppose  $2p^{1+\frac{1}{\delta}}P_\nn(r_\rho)\leq \gamma$ and $\left(\sqrt{\frac{\nn-1}{p}}+\sqrt{\frac{\ln p}{p}}\right) \leq c$ hold for some positive and \ONE{sufficiently} small  constant $c$. Then
	$$\E\left|N_{\VVV_\delta}^{(\Pma)}-N_{\VVV_\delta}^{(\Rma)}\right|\leq  C_{\breve{V}_\delta}^{(\Pma)} \left(1+\frac{\kappa-1}{p}\mu_{n,\delta+2}(\Sigmama)\right)\left(\sqrt{\frac{\ln p }{p}}+\frac{\tau}{p}+p^{-\frac{1}{\delta}} \right) $$
	where $C_{\breve{V}_\delta}^{(\Pma)}$ is  defined in \eqref{eqn:CVexcdel}.
\end{lem}
\begin{proof}
	Let $\Ec(t)$ be the same as in \eqref{Epsdef} with $t$ to be determined. Consider $\delta\in [p-2]$.
	\begin{align*}
	&	\E \left|N_{\VVV_\delta}^{(\Pma)}-N_{\VVV_\delta}^{(\Rma)}\right| \\
	\leq & \E \left|N_{\VVV_\delta}^{(\Pma)}-N_{\VVV_\delta}^{(\Rma)}\right| 1\left(\Ec(t)\right) + p \P (\Ec^c(t))   \\
	\leq &  \E |N_{E_{\delta}}^{(\Pma)} - N_{E_{\delta}}^{(\Rma)}|1\left(\Ec(t)\right)+ (\delta+1)\E|N_{E_{\delta+1}}^{(\Pma)}-N_{E_{\delta+1}}^{(\Rma)}|1\left(\Ec(t)\right)+(\delta+1)\E N_{E_{\delta+1}}^{(\Rma)} +2p e^{-c_1t^2}, \numberthis  \label{eqn:Nvexadelpr}
	\end{align*}
	where the first inequality follows from $0\leq N_{\VVV_\delta}^{(\Psima)}\leq p$ for $\Psima =\Rma$ and $\Pma$, the second inequality follows from \eqref{eqn:Nvexacdifbou} and \eqref{eqn:proepst}. If $\delta=p-1$, then
	\begin{align*}
	\E \left|N_{\VVV_{p-1}}^{(\Pma)}-N_{\VVV_{p-1}}^{(\Rma)}\right| 
	\leq & \E \left|N_{\VVV_{p-1}}^{(\Pma)}-N_{\VVV_{p-1}}^{(\Rma)}\right| 1\left(\Ec(t)\right) + p \P (\Ec^c(t))   \\
	\leq &  \E \left|N_{E_{p-1}}^{(\Pma)} - N_{E_{p-1}}^{(\Rma)}\right|1\left(\Ec(t)\right) +2p \exp(-c_1t^2),
	\end{align*}
    which shows that \eqref{eqn:Nvexadelpr} also holds for $\delta=p-1$ with the convention $N_{E_{p}}^{(\Pma)}=N_{E_{p}}^{(\Rma)}=0$.
	
	Choose $t=\sqrt{\frac{3}{c_1}\ln p}:=c_2 \sqrt{\ln p}$, such that $p\exp(-c_1t^2)=\frac{1}{p^2}$. Moreover, for any 
	\begin{equation}
	c<\frac{1}{2\max\{c_2,1\}\sqrt{2}C_1}, \label{eqn:cuppbou2}
	\end{equation}
	the inequality
	$$
	\left(\sqrt{\frac{\nn-1}{p}}+\sqrt{\frac{\ln p}{p}}\right) \leq c
	$$
	implies
	\begin{equation}
	\sqrt{2}C_1\left(\sqrt{\frac{\nn-1}{p}}+c_2\sqrt{\frac{\ln p}{p}}\right) \leq \frac{1}{2}, \label{eqn:CEatleadelprenewnew}
	\end{equation}
	which is \eqref{eqn:postempass} with $t= c_2 \sqrt{\ln p}$. With $t= c_2 \sqrt{\ln p}$ Lemma \ref{lem:ideltaedgeindicatorbou} becomes:
	\begin{align*}
	& \E \left|N_{V_{\delta}}^{(\Pma)}-N_{V_{\delta}}^{(\Rma)} \right| 1\left(\Ec(t) \right)\\
	\leq & \frac{Cn^2}{(\delta-1)!}\left(1+\delta^2\frac{\kappa-1}{p}\mu_{n,\delta+1}(\Sigmama)\right)  (\theta_1(c_2\sqrt{\ln p}))^{\nn\delta}\left(\frac{\sqrt{\ln p}}{\sqrt{p}} + \frac{\tau}{p}\right) p\left( 2pP_\nn\left(r_\rho\right) \right)^{\delta}, \numberthis \label{eqn:Etemp1}
	\end{align*}
	
	Then for $\delta\in [p-1]$ applying \eqref{eqn:Etemp1} with $\delta$, $\delta+1$ and Lemma \ref{lem:NEbou} to \eqref{eqn:Nvexadelpr}, together with $\mu_{n,\delta+1}(\Sigmama)\leq \mu_{n,\delta+2}(\Sigmama)$,
	\begin{align*}
	&\E \left|N_{\VVV_\delta}^{(\Pma)}-N_{\VVV_\delta}^{(\Rma)}\right|\\
	\leq & \frac{Cn^2}{(\delta-1)!}\left(1+\delta^2\frac{\kappa-1}{p}\mu_{n,\delta+2}(\Sigmama)\right)  \left(\theta_1\left(c_2\sqrt{\ln p}\right)\right)^{\nn\delta}\left(\frac{\sqrt{\ln p}}{\sqrt{p}} + \frac{\tau}{p}\right) p\left( 2pP_\nn\left(r_\rho\right) \right)^{\delta}+ \\
	&
	{ (\delta+1)\frac{Cn^2}{\delta!}\left(1+(\delta+1)^2\frac{\kappa-1}{p}\mu_{n,\delta+2}(\Sigmama)\right)  \left(\theta_1\left(c_2\sqrt{\ln p}\right)\right)^{\nn(\delta+1)}\left(\frac{\sqrt{\ln p}}{\sqrt{p}} + \frac{\tau}{p}\right) p\left( 2pP_\nn\left(r_\rho\right) \right)^{\delta+1} } \\
	& \ \ \ \ +(\delta+1)\frac{1}{(\delta+1)!}\left(1+(\delta+1)^2\mu_{n,\delta+2}(\Sigmama)\frac{\kappa-1}{p}\right)p\left( 2pP_\nn(r_\rho)\right)^{\delta+1}
	  + \frac{2}{p^2} 	\\
	\leq & \frac{Cn^2}{(\delta-1)!}\left(1+(\delta+1)^2\frac{\kappa-1}{p}\mu_{n,\delta+2}(\Sigmama)\right)  \left(\theta_1\left(c_2\sqrt{\ln p}\right)\right)^{\nn(\delta+1)}\left(\frac{\sqrt{\ln p}}{\sqrt{p}} + \frac{\tau}{p}\right) \gamma^{\delta}\left(1+ \frac{\delta+1}{\delta}\gamma p^{-\frac{1}{\delta}}\right) \\
	& \ \ \ \ +\frac{1}{\delta!}\left(1+(\delta+1)^2\mu_{n,\delta+2}(\Sigmama)\frac{\kappa-1}{p}\right)\gamma^{\delta+1}p^{-\frac{1}{\delta}}
	+ \frac{2}{p^2} 	\\
	\leq & C_{\breve{V}_\delta}^{(\Pma)} \left(1+\frac{\kappa-1}{p}\mu_{n,\delta+2}(\Sigmama)\right)\left(\sqrt{\frac{\ln p }{p}}+\frac{\tau}{p}+p^{-\frac{1}{\delta}} \right),
	\end{align*}
	where in the last inequality
	\begin{align*}
	C_{\breve{V}_\delta}^{(\Pma)}=&\frac{Cn^2(\delta+1)^2}{(\delta-1)!}  \left(\theta_1\left(c_2\sqrt{\ln p}\right)\right)^{\nn(\delta+1)} \gamma^{\delta}(1+\gamma)\left(1+ \frac{\delta+1}{\delta}\gamma p^{-\frac{1}{\delta}}\right) 
	+ \frac{2}{p^{2-\frac{1}{\delta}}}\\
	\leq & \frac{Cn^2(\delta+1)^2}{(\delta-1)!}  \left(4n+5\right)^{\nn(\delta+1)} \gamma^{\delta}(1+\gamma)\left(1+ \frac{\delta+1}{\delta}\gamma p^{-\frac{1}{\delta}}\right) 
	+ 2. \numberthis \label{eqn:CVexcdel}
	\end{align*}
	where the last step follows from $\theta_1\left(c_2 \sqrt{\ln p}\right)\leq 9+ 4(n-1)=4n+5$ by \eqref{eqn:CEatleadelprenewnew} and $\tau\leq p/2$.
\end{proof}

\begin{proof}[Proof of Proposition \ref{thm:allclose} \ref{item:allclosed}]
Lemmas \ref{lem:Nvexadelpr} and \ref{prop:tvtomean} complete the proof of Proposition \ref{thm:allclose} \ref{item:allclosed}.
\end{proof}

\subsection{Proofs in Subsection \ref{sec:compoihighdim}}

\subsubsection{{Proof of Lemma \ref{prop:rangeo}}}
\label{sec:proofofprop:rangeo}

To utilize the notation we have defined in this paper, we make the following adjustments to the notation throughout this subsection. In this proof it suffices to prove the conclusion for any $\delta+1$ i.i.d. random points from $\unif(S^{n-2})$. Without loss of generality  assume in this subsection that the first $\delta+1$ U-scores $\{\uve_i\}_{i=1}^{\delta+1}$ are independent. Another adjustment is to replace $r$ by $r_\rho$. With these adjustments Lemma \ref{prop:rangeo} is equivalent to the following: when $r_\rho <2/\sqrt{5}$, $\delta\geq 1$, for any $\ell\in [\delta+1]$,
\begin{align*}
	&\P\left(\nmd\left(\{\uve_i  \}_{i=1}^{\delta+1}, r_\rho  \right)=\ell|\deg(\uve_{\delta+1})=\delta\right) \\
	= & \P\left(\pnmd\left(\{\uve_i  \}_{i=1}^{\delta+1}, r_\rho  \right)=\ell|\deg(\uve_{\delta+1})=\delta\right). \numberthis \label{eqn:proofequilem1}
	\end{align*}

Take $\vec{i}=(\delta+1,1,\cdots,\delta)$. Recall the notation $\vec{i},\Phi_{\vec{i}}=\Phi^{(\Rma)}_{\vec{i}},U_{\vec{i}}$ are defined in Subsection \ref{sec:closenessnumedge}, where the dependence of $\Rma$ in $\Phi_{\vec{i}}^{(\Rma)}$ is suppressed throughout this subsection for the sake of clearer exposition. Then $\pnmd\left(\{\uve_i  \}_{i=1}^{\delta+1}, r_\rho  \right)= U_{\vec{i}}+ \Phi_{\vec{i}}$. Moreover, the
event $\{\deg(\uve_{\delta+1})=\delta\}$ in $\textbf{PGe}\left(\{\uve_i  \}_{i=1}^{\delta+1}, r_\rho   \right)$ is the same as $\{\Phi_{\vec{i}} =1 \}$. Define $F_{ij}^{(q)}=\{\|\uve_i-q\uve_j\|_2 \leq r_\rho \}$. Then $1\left(\bigcap\limits_{j=1}^\delta F_{j(\delta+1)}^{(+1)}\right)$ is the indicator function of the event that the degree of vertex $\uve_{\delta+1}$ in $\textbf{Ge}\left(\{\uve_i  \}_{i=1}^{\delta+1}, r_\rho   \right)$ is $\delta$. Hence Lemma \ref{prop:rangeo} is also equivalent to the following: when $r_\rho <2/\sqrt{5}$, $\delta\geq 1$, for any $\ell\in [\delta+1]$,
\begin{equation}
\P\left(U_{\vec{i}}+ \Phi_{\vec{i}}=\ell|\Phi_{\vec{i}}=1\right) = \P\left(\nmd\left(\{\uve_i  \}_{i=1}^{\delta+1}, r_\rho   \right) =\ell \left| 1\left(\bigcap_{j=1}^\delta F_{j(\delta+1)}^{(+1)}\right) =1 \right. \right). \label{eqn:pseudogeogeoequiv}
\end{equation}

\begin{proof}[Proof of \eqref{eqn:pseudogeogeoequiv}]
	For $\vec{q} =(q_1,q_2,\cdots,q_\delta) \in \{-1,+1\}^{\delta}$, denote $F_{\delta+1}^{(\vec{q})} = \bigcap_{j=1}^\delta F_{j(\delta+1)}^{(q_j)}$.
	Observe that
	$$
	\{\Phi_{\vec{i}} =1 \} = \bigcap_{j=1}^{\delta} \ \ \bigcup_{q_j \in \{+1,-1\} } F_{j(\delta+1)}^{(q_j)} =   
	= \bigcup_{\vec{q}\in \{-1,+1\}^{\delta} } \ \ \ F_{\delta+1}^{(\vec{q})}.
	$$
	Since $r_\rho < 2/\sqrt{5}< \sqrt{2} $, $F_{j(\delta+1)}^{(-1)}$ and $F_{j(\delta+1)}^{(+1)}$ are disjoint for every $j\in [\delta]$, which implies $F_{\delta+1}^{(\vec{q})}$ for different $\vec{q} \in \{-1,+1\}^{\delta}  $ are disjoint. Hence, 
	\begin{equation}
	\P (\Phi_{\vec{i}} =1, U_{\vec{i}}=\ell -1) = \sum_{\vec{q} \in \{-1,+1\}^{\delta}} \P(F_{\delta+1}^{(\vec{q})}, U_{\vec{i}} =\ell -1). \label{eqn:temp1}
	\end{equation}
	Next observe that $1(F_{\delta+1}^{(\vec{q})})$ is a function of $\uve_1,\cdots,\uve_{\delta+1}$, and hence it has the same distribution when replacing $\uve_i$ by $-\uve_i$ for any $i\in [\delta]$. Moreover, replacing $\uve_i$ by $-\uve_i$ for any $i\in [\delta]$ wouldn't change $U_{\vec{i}}$. As a result, \eqref{eqn:temp1} implies  
	\begin{equation}
	\P (\Phi_{\vec{i}} =1, U_{\vec{i}}=\ell -1) = 2^{\delta} \P\left(F_{\delta+1}^{(\vec{q}_0)}, U_{\vec{i}} =\ell -1\right), \label{eqn:temp2}
	\end{equation}
	where $\vec{q}_0 =(+1,+1,\cdots,+1)$ is the vector in $\R ^{\delta}$ with all its components $+1$.

	Consider $\omega\in F_{\delta+1}^{(\vec{q}_0)}$. Then $\Phi_{\vec{i}}(\omega)=1$ or equivalently, $\Phi_{i(\delta+1)}^{(\Rma)}(\omega)=1$ for any $i \in [\delta]$. Then 
	\begin{equation}
	U_{\vec{i}}(\omega) = \sum_{\vec{j}\in S_{\vec{i}}} \Phi_{\vec{j}} (\omega) = \sum_{i=1}^\delta \prod_{\substack{j=1\\ j\not = i} }^{\delta+1} \Phi^{(\Rma)}_{ij} (\omega)=\sum_{i=1}^\delta \prod_{\substack{j=1\\ j\not = i} }^{\delta} \Phi^{(\Rma)}_{ij} (\omega). \label{eqn:Uitemptemp2}
	\end{equation}
	Since for any distinct $i,j\in [\delta]$, $\|\uve_i(\omega) - \uve_j (\omega) \|_2 \leq \| \uve_i -\uve_{\delta+1}(\omega)\|_2 + \| \uve_i (\omega) -\uve_{\delta+1}(\omega)\|_2 \leq 2r_\rho < 4/\sqrt{5} $, $\|\uve_i(\omega) + \uve_j (\omega) \|_2 =\sqrt{4- \|\uve_i(\omega) - \uve_j (\omega) \|_2^2 } > 2/\sqrt{5} > r_\rho  $. Thus  $\Phi^{(\Rma)}_{ij}(\omega)  =  1_{F_{ij}^{(+1)}}(\omega) $. That is, in the set $F_{\delta+1}^{(\vec{q}_0)}$, \eqref{eqn:Uitemptemp2} becomes
	\begin{equation}
	U_{\vec{i}} = \sum_{i=1}^\delta \prod_{\substack{j=1\\ j\not = i} }^\delta 1(F_{ij}^{(+1)}) = \nmd\left(\{\uve_i  \}_{i=1}^{\delta}, r_\rho   \right),  \label{eqn:Uidef2}
	\end{equation}
	which implies
	\begin{align*}
	\left(\Phi_{\vec{i}} + U_{\vec{i}}\right) 1\left(F_{\delta+1}^{(\vec{q}_0)}\right)  =& \left(1+ \nmd\left(\{\uve_i  \}_{i=1}^{\delta}, r_\rho   \right) \right) 1\left(F_{\delta+1}^{(\vec{q}_0)}\right) \\
	= & \nmd\left(\{\uve_i  \}_{i=1}^{\delta+1}, r_\rho   \right) 1\left(F_{\delta+1}^{(\vec{q}_0)}\right) . \numberthis \label{eqn:temptemp1}
	\end{align*}
	Thus,
	\begin{align*}
	\P\left(U_{\vec{i}}+ \Phi_{\vec{i}}=\ell|\Phi_{\vec{i}}=1\right) =& 
	\frac{2^{\delta} \P\left(F_{\delta+1}^{(\vec{q}_0)}, U_{\vec{i}} =\ell -1\right)}{\P \left(\Phi_{\vec{i}}=1\right)}\\
	=& 
	\frac{2^{\delta} \P\left(F_{\delta+1}^{(\vec{q}_0)}, U_{\vec{i}} =\ell -1\right)}{\left(2P_n(r_\rho)\right)^\delta} \\
	=& 
	\frac{ \P\left(F_{\delta+1}^{(\vec{q}_0)}, \nmd\left(\{\uve_i  \}_{i=1}^{\delta+1}, r_\rho   \right) =\ell \right)}{\left(P_n(r_\rho)\right)^\delta} \\
	=& 
	\P\left(\nmd\left(\{\uve_i  \}_{i=1}^{\delta+1}, r_\rho   \right) =\ell \left| 1\left(\bigcap_{j=1}^\delta F_{j(\delta+1)}^{(+1)}\right) =1 \right. \right),
	\end{align*}
	where the first equality follows from \eqref{eqn:temp2}, the second equality follows from Lemma \ref{item:prod}, and the third equality follows from \eqref{eqn:temptemp1}.

\end{proof}

\subsubsection{Proofs of Lemma \ref{prop:rangeolimit} and Lemma \ref{lem:compoundpoissonlimit}}
\label{sec:proofofprop:rangeolimit}

\begin{proof}[Proof of Lemma \ref{prop:rangeolimit}]
In the set $\{\deg(\uve'_{\delta+1})=\delta\}$, it follows that
$$
\nmd\left(\{\uve'_i  \}_{i=1}^{\delta+1}, r   \right)  = \nmd\left( \{\uve'_i  \}_{i=1}^\delta, r   \right)+1.
$$
Thus, 
\begin{align*}
&\P\left(\deg(\uve'_{\delta+1})=\delta, \nmd\left(\{\uve'_i  \}_{i=1}^{\delta+1}, r   \right) =\ell\right) \\
 =&  \E 1\left(\{\deg(\uve'_{\delta+1})=\delta\} \bigcap \{ \nmd\left( \{\uve'_i  \}_{i=1}^\delta, r   \right)= \ell -1   \} \right) \\
 =&\E\left( \E \left( 1\left(\{\deg(\uve'_{\delta+1})=\delta\} \bigcap \{ \nmd\left( \{\uve'_i  \}_{i=1}^\delta, r \right)= \ell -1   \} \right) \big| \uve'_{\delta+1} \right) \right) \\
 =& \E \left( 1\left(\{\deg(\uve'_{\delta+1})=\delta\} \bigcap \{ \nmd\left( \{\uve'_i  \}_{i=1}^\delta, r  \right)= \ell -1   \} \right) \big| \uve'_{\delta+1} = \vve_0 \right), \numberthis \label{eqn:tempcond}
\end{align*}
where the last equality follows from that
$$
\E \left( 1\left(\{\deg(\uve'_{\delta+1})=\delta\} \bigcap \{ \nmd\left( \{\uve'_i  \}_{i=1}^{\mgeominus}, r  \right)= \ell -1   \} \right) \big| \uve'_{\delta+1} \right)
$$ 
as a function of the random variable $\uve'_{\delta+1}$, due to the rotation invariance property of the distribution $\unif(S^{n-2})$, equals to
 $$
 \E \left( 1\left(\{\deg(\uve'_{\delta+1})=\delta\} \bigcap \{ \nmd\left( \{\uve'_i  \}_{i=1}^\delta, r  \right)= \ell -1   \} \right) \big| \uve'_{\delta+1} = \vve_0 \right)
 $$ 
 a.s.  with $\vve_0 = (1,0,0,\cdots,0)\in S^{n-2}$.

Under the condition $\uve'_{\delta+1} = \vve_0$, $1\left(\{\deg(\uve'_{\delta+1})=\delta\}\right) = \prod_{i=1}^{\delta} 1\left( \uve'_i \in \SC(r, \vve_0) \right)$, where $\SC(r, \vve_0)$ is defined in \eqref{eqn:SCdef}. Consider the coordinate system for $\uve'_i = \left(u'_{ji}:1\leq j \leq n-1\right)^\top$ in the region $\SC(r, \vve_0)$:
\begin{equation*}
\begin{cases}
u'_{1i}=1-\frac{r^2 r_i^2}{2},\\
u'_{2i}=r r_i\sqrt{1-\frac{r^2 r_i^2}{4}} \cos(\theta_{2i}),\\
\vdots & \text{\ \ for } 1\leq i \leq \delta,\\
u'_{ji}=r r_i\sqrt{1-\frac{r^2 r_i^2}{4}} \cos(\theta_{ji})\prod\limits_{m=2}^{j-1}\sin(\theta_{mi}),\\
\vdots \\
u'_{(n-2)i}=r r_i\sqrt{1-\frac{r^2 r_i^2}{4}}\sin(\theta_{2i})\cdots\sin(\theta_{(n-3)i})\cos(\theta_{(n-2)i}),\\
u'_{(n-1)i}=r r_i\sqrt{1-\frac{r^2 r_i^2}{4}}\sin(\theta_{2i})\cdots\sin(\theta_{(n-3)i})\sin(\theta_{(n-2)i}), 
\end{cases}
\end{equation*}
where for each $i\in [\delta]$:
\begin{equation}
 r_i\in [0,1], \theta_{ji}\in[0,\pi] \text{ for } 2\leq j\leq n-3 \text{ and } \theta_{(n-2)i}\in[0,2\pi).  \label{eqn:intregion}
\end{equation}
Then,
\begin{align*}
&r^{-(n-2)\delta}\E \left( 1\left(\{\deg(\uve'_{\delta+1})=\delta\} \bigcap \{ \nmd\left( \{\uve'_i  \}_{i=1}^\delta, r  \right)= \ell -1   \} \right) \big| \uve'_{\delta+1} = \vve_0 \right) \\
= & r^{-(n-2)\delta} \E \prod_{i=1}^{\delta} 1\left( \uve'_i \in \SC(r, \vve_0) \right) 1\left( \nmd\left( \{\uve'_i  \}_{i=1}^\delta, r \right)= \ell -1\right) \\
\overset{(*)}{=}\  & \frac{r^{-(n-2)\delta}}{|S^{n-2}|^\delta} \int\cdots\int_{\Omega_0} 1\left( \nmd\left( \{\uve'_i  \}_{i=1}^\delta, r \right)= \ell -1\right)\\
&\ \ \ \times \prod_{i=1}^{\delta}\left( r^{n-2}r_i^{n-3}\left(1-\frac{r^2 r^2_i}{4}\right)^{\frac{n-4}{2}} dr_i \prod_{j=2}^{n-2} \left(\sin^{n-2-j}(\theta_{ji}) d\theta_{ji} \right)\right)  \\
=& \frac{1}{|S^{n-2}|^\delta} \int\cdots\int_{\Omega_0} 1\left( \nmd\left( \{\uve'_i  \}_{i=1}^\delta, r \right)= \ell -1\right)\\
&\ \ \ \times\prod_{i=1}^{\delta}\left( r_i^{n-3}\left(1-\frac{r^2 r^2_i}{4}\right)^{\frac{n-4}{2}} dr_i \prod_{j=2}^{n-2} \left(\sin^{n-2-j}(\theta_{ji}) d\theta_{ji} \right)\right), \numberthis \label{eqn:temp4}
\end{align*}
where $\Omega_0$ in equality $(*)$ is the region described in \eqref{eqn:intregion}. Denote by $f(r)$ the integrand in \eqref{eqn:temp4}. $f(r)$ is a function of $r_i$, $\theta_{ji}$ for $2\leq j\leq n-2$ and $1\leq i \leq \delta$, of which the dependences are suppressed.

Note  $\nmd\left(\{\vve_i  \}_{i=1}^\mgeominus, r   \right)$ is a function of $(1\left(\|\vve_i-\vve_j\|_2\leq r\right): 1\leq i<j \leq \mgeominus)$ and it does not depend on the specific location of each vertex. Thus in \eqref{eqn:temp4} $ \nmd\left( \{\uve'_i  \}_{i=1}^\delta, r \right)$ is a function of $(1\left(\|\uve'_i-\uve'_j\|_2< r\right): 1\leq i<j \leq \mgeominus)$ since $\|\uve'_i-\uve'_j\|_2 = r$, as a set of Lebesgue measure $0$, contributes nothing to the integral.
We then write 
\begin{align*}
\nmd\left( \{\uve'_i  \}_{i=1}^\delta, r \right) = & \chi(1\left(\|\uve'_i-\uve'_j\|_2< r\right): 1\leq i<j \leq \mgeominus)\\
= & \chi\left(1\left(\frac{1}{r}\|\uve'_i-\uve'_j\|_2< 1\right): 1\leq i<j \leq \mgeominus\right). 
\end{align*}
Intrinsically, $1\left(\|\uve'_i-\uve'_j\|_2< r\right)$ is the indicator random variable about whether there is an edge between vertex $i$ and $j$, and the function $\chi$ is the function that takes all edge information among $\delta$ vertices as input and outputs the number of universal vertices.

Then as $r\to 0^+$, 
\begin{align*}
\lim_{r\to 0^+} f(r)= &  \prod_{i=1}^{\delta}\left( r_i^{n-3}  \prod_{j=2}^{n-2} \left(\sin^{n-2-j}(\theta_{ji})  \right)\right)\\
& \times\lim_{r\to 0^+} 1\left(\chi\left(1\left(\frac{1}{r}\|\uve'_i-\uve'_j\|_2< 1\right): 1\leq i<j \leq \mgeominus\right) = \ell-1 \right). \numberthis \label{eqn:temp3}
\end{align*}
Observe 
\begin{align*}
&\lim_{r \to 0^+}\left(\frac{1}{r}\|\uve'_i-\uve'_j\|_2\right)^2 \\
= & \left( r_i \cos(\theta_{2i}) - r_j \cos(\theta_{2j})  \right)^2+\sum_{q=3}^{n-2} \left( r_i \prod_{m=2}^{q-1} \sin(\theta_{mi})\cos(\theta_{qi}) - r_j \prod_{m=2}^{q-1} \sin(\theta_{mj})\cos(\theta_{qj})  \right)^2 \\
& \quad + \left( r_i \prod_{m=2}^{n-2} \sin(\theta_{mi}) - r_j \prod_{m=2}^{n-2} \sin(\theta_{mj})  \right)^2. \numberthis \label{eqn:temp5}
\end{align*}
On $\Omega_0$, for $1\leq i \leq \delta$, define 
\begin{equation}
\begin{cases}
\tilde{u}'_{1i} = r_i \cos(\theta_{2i}),\\	
\tilde{u}'_{ji}=r_i \cos(\theta_{(j+1)i}) \prod\limits_{m=2}^{j}\sin(\theta_{mi}),\text{ for } 2\leq j \leq n-3 & \\
\tilde{u}'_{(n-2)i}= r_i\prod\limits_{m=2}^{n-2}\sin(\theta_{mi}),\\
\end{cases} \label{eqn:newballparametrization}
\end{equation}
and $\tilde{\uve}'_i =(\tilde{u}'_{ji}:1\leq j \leq n-2)\in B^{n-2}$. Then by \eqref{eqn:temp5}
$$
\lim_{r\to 0^+}\frac{1}{r}\|\uve'_i-\uve'_j\|_2 = \| \tilde{\uve}'_i -\tilde{\uve}'_j\|_2, 
$$
which, together with \eqref{eqn:temp3}, imply 
\begin{align*}
&\lim_{r \to 0^+} f(r) \\
= & \prod_{i=1}^{\delta}\left( r_i^{n-3}  \prod_{j=2}^{n-2} \left(\sin^{n-2-j}(\theta_{ji})  \right)\right)    1\left(\chi\left(1\left(\|\tilde{\uve}'_i-\tilde{\uve}'_j\|_2< 1\right): 1\leq i<j \leq \mgeominus\right) = \ell-1 \right)\\
=&\prod_{i=1}^{\delta}\left( r_i^{n-3}  \prod_{j=2}^{n-2} \left(\sin^{n-2-j}(\theta_{ji})  \right)\right)    1\left(\chi\left(1\left(\|\tilde{\uve}'_i-\tilde{\uve}'_j\|_2\leq 1\right): 1\leq i<j \leq \mgeominus\right) = \ell-1 \right)\\
=& \prod_{i=1}^{\delta}\left( r_i^{n-3}  \prod_{j=2}^{n-2} \left(\sin^{n-2-j}(\theta_{ji})  \right)\right)    1\left(\nmd\left( \{\tilde{\uve}'_i  \}_{i=1}^\delta, 1   \right)= \ell -1\right),
\end{align*}
where the second equality holds $a.s.$ with respect to the Lebesgue measure on $\Omega_0$,

Moreover, $|f(r)|\leq 1$, which is integrable over the bounded set $\Omega_0$. Applying Dominated Convergence Theorem to \eqref{eqn:temp4}, 
\begin{align*}
&\lim_{r\to 0^+}r^{-(n-2)\delta}\E \left( 1\left(\{\deg(\uve'_{\delta+1})=\delta\} \bigcap \{ \nmd\left( \{\uve'_i  \}_{i=1}^\delta, r \right)= \ell -1   \} \right) \big| \uve'_{\delta+1} = \vve_0 \right) \\
=& \frac{1}{|S^{n-2}|^\delta} \int\cdots\int_{\Omega_0} 1\left( \nmd\left( \{\tilde{\uve}'_i  \}_{i=1}^\delta, 1   \right)= \ell -1\right) \prod_{i=1}^{\delta} r_i^{n-3} dr_i \prod_{j=2}^{n-2} \sin^{n-2-j}(\theta_{ji}) d\theta_{ji}  \\
= & \frac{|B^{n-2}|^\delta}{|S^{n-2}|^\delta} \P\left( \nmd\left( \{\tilde{\uve}_i  \}_{i=1}^\delta, 1   \right)= \ell -1\right), \numberthis \label{eqn:templimit}
\end{align*}
where the parametrization \eqref{eqn:newballparametrization} and the region $\Omega_0$ coincide with the spherical coordinates for $B^{\Ngeo}$. 

Thus
\begin{align*}
&\lim_{r \to 0^+} \P\left(\nmd\left(\{\uve'_i  \}_{i=1}^{\delta+1}, r \right) =\ell \left| \deg(\uve'_{\delta+1})=\delta \right. \right) \\
= &  \lim_{r \to 0^+} r^{-(n-2)\delta}\P\left(\deg(\uve'_{\delta+1})=\delta, \nmd\left(\{\uve'_i  \}_{i=1}^{\delta+1}, r \right) =\ell\right) \frac{r^{(n-2)\delta}}{\left(P_n(r)\right)^\delta} \\
 = & \frac{|B^{n-2}|^\delta}{|S^{n-2}|^\delta} \P\left( \nmd\left( \{\tilde{\uve}_i  \}_{i=1}^\delta, 1   \right)= \ell -1\right)  \frac{1}{(a_n)^\delta} \\
 = & \frac{1}{(a_n)^\delta} \frac{|B^{n-2}|^\delta}{|S^{n-2}|^\delta} \P\left( \nmd\left( \{\tilde{\uve}_i  \}_{i=1}^\delta, 1   \right)= \ell -1\right) \\
 = & \P\left( \nmd\left( \{\tilde{\uve}_i  \}_{i=1}^\delta, 1   \right)= \ell -1\right), 	\end{align*}
where the second equality follows from \eqref{eqn:tempcond}, \eqref{eqn:templimit} and Lemma \ref{pderivative} \ref{item:pderivativeb}.
\end{proof}

\begin{proof}[Proof of Lemma \ref{lem:compoundpoissonlimit}] 
By \eqref{eqn:alphalrhorgg}, Lemma \ref{prop:rangeolimit} and \eqref{eqn:alphaelldef},
\begin{equation}
\lim_{\rho\to 1^-}\alpha(\ell,r_\rho) = \lim_{r_\rho\to 0^+}\alpha(\ell,r_\rho) = \alpha_\ell, \quad \forall \ell \in [\delta+1], \label{eqn:alphalrholimit}
\end{equation}
and thus,
\begin{equation}
\lim_{\rho\to 1^-}\bm{\zeta}_{n,\delta,\rho}(\ell) = \bm{\zeta}_{n,\delta}(\ell), \quad \forall \ell \in [\delta+1]. \label{eqn:alphalimit}
\end{equation}

Note that
\begin{align*}
   \left| 2p^{1+\frac{1}{\delta}}P_n(r_\rho) - e_{n,\delta}\right| \leq & \left|2p^{1+\frac{1}{\delta}}P_n(r_\rho)- 2a_n p^{1+\frac{1}{\delta}} r_\rho^{\nn-2} \right|+ \left| 2a_n p^{1+\frac{1}{\delta}} r_\rho^{\nn-2}- e_{n,\delta}\right| \\
   \leq &  2a_n p^{1+\frac{1}{\delta}} r_\rho^{\nn-2} \left(\frac{n-4}{8}{r_\rho^2}\right)  + \left| 2a_n p^{1+\frac{1}{\delta}} r_\rho^{\nn-2}- e_{n,\delta}\right|\\
   =  &  2a_n p^{1+\frac{1}{\delta}} r_\rho^{\nn-2} \left(\frac{n-4}{4}(1-\rho)\right)  + \left| 2a_n p^{1+\frac{1}{\delta}} r_\rho^{\nn-2}- e_{n,\delta}\right|. \numberthis \label{eqn:slowupperbound}
\end{align*}
where the second inequality follows from Lemma \ref{pderivative} \ref{item:pderivativea}.
Then the preceding expression and \eqref{eqn:alphalrholimit} yield
\begin{equation}
\lim_{p\to\infty}\lambda_{p,n,\delta,\rho} =\lim_{p\to\infty} \frac{1}{\delta!}p^{\delta+1}(2P_n(r_\rho))^\delta\sum_{\ell=1}^{\delta+1} \frac{\alpha(\ell,r_\rho)}{\ell}=\lim_{p\to\infty} \frac{1}{\delta!}\left(e_{n,\delta}\right)^\delta \sum_{\ell=1}^{\delta+1} \frac{\alpha_\ell}{\ell}=\lambda_{\nn,\delta}(e_{n,\delta}). \label{eqn:lambdalimit}
\end{equation}
\eqref{eqn:alphalimit} and \eqref{eqn:lambdalimit} immediately yield the conclusion. 
\end{proof}

\subsection{Proofs in Section \ref{sec:convergenceofmoments}}

\subsubsection{Proofs of Lemma \ref{lem:meanedge} and Proposition \ref{lem:2ndmoment}}
\begin{proof}[Proof of Lemma \ref{lem:meanedge}]
	\begin{align*}
	\E N_{E_{\delta}}^{(\Rma)}-\binom{p}{1}\binom{p-1}{\delta}(2P_n(r_\rho))^{\delta} 
	= & \sum_{ \substack{\vec{i}\in C_{\delta}^<\\ \Sigmama_{\vec{i}}\text{ not diagonal}} } \left(\E  \prod_{j=1}^{\delta} \Phi_{i_0i_j}^{(\Rma)}-(2P_n(r_\rho))^{\delta}\right) \numberthis \label{eqn:meantemp1}  \\
	\leq & 
	\left(\mu_{n,\delta+1}\left(\Sigmama\right)-1\right) \left( 2P_\nn(r_\rho)\right)^\delta \sum_{ \substack{\vec{i}\in C_{\delta}^<\\ \Sigmama_{\vec{i}}\text{ not diagonal}} }1,
	\end{align*}
	where the first inequality follows from Lemma \ref{lem:densitybou} \ref{item:densitybouc} and Lemma \ref{item:prod}. By \eqref{eqn:meantemp1},
	$$
	\E N_{E_{\delta}}^{(\Rma)}-\binom{p}{1}\binom{p-1}{\delta}(2P_n(r_\rho))^{\delta}\geq - (2P_n(r_\rho))^{\delta} \sum_{ \substack{\vec{i}\in C_{\delta}^<\\ \Sigmama_{\vec{i}}\text{ not diagonal}} } 1 .
	$$
	Combining the preceding two expressions,
	\begin{align*}
	   \left| \E N_{E_{\delta}}^{(\Rma)}-\binom{p}{1}\binom{p-1}{\delta}(2P_n(r_\rho))^{\delta} \right|
	\leq & \max\{1,\mu_{n,\delta+1}\left(\Sigmama\right)-1)\} \left( 2P_\nn(r_\rho)\right)^\delta \sum_{ \substack{\vec{i}\in C_{\delta}^<\\ \Sigmama_{\vec{i}}\text{ not diagonal}} }1\\
	\leq &
	\mu_{n,\delta+1}\left(\Sigmama\right) (2P_n(r_\rho))^{\delta} \frac{(\delta+1)}{2((\delta-1)!)}   p^{\delta} (\kappa-1)\\
	\leq & 
	 \frac{(\delta+1)}{2((\delta-1)!)}\gamma^{\delta} \mu_{n,\delta+1}\left(\Sigmama\right) \frac{\kappa-1}{p},
	\end{align*}
	where the second inequality follows from Lemma \ref{lem:numofnondia}.
\end{proof}

\begin{proof}[Proof of Proposition \ref{lem:2ndmoment}] 
Recall for $\vec{i}\in C_\delta^<$, $\Phi_{\vec{i}}^{(\Rma)}$ is defined in \eqref{eqn:NEdeltaRdef}, $U_{\vec{i}}$ is defined in \eqref{eqn:Uidef}, and $Z_{\vec{i}}$ and $W_{\vec{i}}$ are defined in \eqref{eqn:WiZidef}. Then, 
$$
N_{E_\delta}^{(\Rma)} = \sum_{\vec{i}\in C^<_{\delta}} \Phi_{\vec{i}}^{(\Rma)} = \Phi_{\vec{i}}^{(\Rma)} + U_{\vec{i}} + Z_{\vec{i}} + W_{\vec{i}}.
$$
Then, 
\begin{equation}
\left(N_{E_\delta}^{(\Rma)}\right)^2 = \sum_{\vec{i}\in C^<_{\delta}} \Phi_{\vec{i}}^{(\Rma)} \left( \Phi_{\vec{i}}^{(\Rma)} + U_{\vec{i}} + Z_{\vec{i}} + W_{\vec{i}} \right). \label{eqn:squaredecom}
\end{equation}
\textbf{Step 1}:\\
Since $\Phi_{\vec{i}}^{(\Rma)} + U_{\vec{i}}$ takes value in $[\delta+1]\bigcup \{0\}$,
\begin{align*}
    \E  \Phi_{\vec{i}}^{(\Rma)} \left( \Phi_{\vec{i}}^{(\Rma)} + U_{\vec{i}}\right) = & \sum_{\ell=1}^{\delta+1} \ell\ \E  \Phi_{\vec{i}}^{(\Rma)} \1 \left( \Phi_{\vec{i}}^{(\Rma)} + U_{\vec{i}} = \ell \right) \\
    =& \sum_{\ell=1}^{\delta+1} \ell\  \P \left( \Phi_{\vec{i}}^{(\Rma)} =1 \right) \P \left( \Phi_{\vec{i}}^{(\Rma)} + U_{\vec{i}} = \ell  | \Phi_{\vec{i}}^{(\Rma)} =1 \right).
\end{align*}
For $\vec{i}\in C^<_\delta$ such that $\Sigmama_{\vec{i}}$ diagonal, we have $\P \left( \Phi_{\vec{i}}^{(\Rma)} =1 \right) = (2P_n(r_\rho))^{\delta}$ by Lemma \ref{item:prod} and moreover $\P \left( \Phi_{\vec{i}}^{(\Rma)} + U_{\vec{i}} = \ell  | \Phi_{\vec{i}}^{(\Rma)} =1 \right) = \alpha(\ell,r_\rho)$ by \eqref{eqn:alphalrhorgg}. Thus in this case, 
$$
\E  \Phi_{\vec{i}}^{(\Rma)} \left( \Phi_{\vec{i}}^{(\Rma)} + U_{\vec{i}}\right) = (2P_n(r_\rho))^{\delta} \sum_{\ell=1}^{\delta+1} \ell  \alpha(\ell,r_\rho).
$$
Moreover, when $\vec{i}\in C^<_\delta$ such that $\Sigmama_{\vec{i}}$ is not diagonal, by Lemma \ref{lem:densitybou} \ref{item:densitybouc} 
$$
\E  \Phi_{\vec{i}}^{(\Rma)} \left( \Phi_{\vec{i}}^{(\Rma)} + U_{\vec{i}}\right) \leq \mu_{n,\delta+1}(\Sigmama) (2P_n(r_\rho))^{\delta} \sum_{\ell=1}^{\delta+1} \ell  \alpha(\ell,r_\rho).
$$
Then by the preceding two expressions,
\begin{align*}
    &\sum_{\vec{i}\in C^<_{\delta}} \E \Phi_{\vec{i}}^{(\Rma)} \left( \Phi_{\vec{i}}^{(\Rma)} + U_{\vec{i}}\right) - \binom{p}{1}\binom{p-1}{\delta} (2P_n(r_\rho))^{\delta} \sum_{\ell=1}^{\delta+1} \ell  \alpha(\ell,r_\rho)\\
    = & 
    \sum_{\substack{\vec{i}\in C^<_{\delta}\\ \Sigmama_{\vec{i}} \text{ not diagonal} }} \left( \E  \Phi_{\vec{i}}^{(\Rma)} \left( \Phi_{\vec{i}}^{(\Rma)} + U_{\vec{i}}\right) -  (2P_n(r_\rho))^{\delta} \sum_{\ell=1}^{\delta+1} \ell  \alpha(\ell,r_\rho) \right) \numberthis \label{eqn:firsttermtemp}  \\
    \leq 
    & \left(\mu_{n,\delta+1}(\Sigmama) -1 \right) (2P_n(r_\rho))^{\delta} \sum_{\ell=1}^{\delta+1} \ell  \alpha(\ell,r_\rho) \sum_{\substack{\vec{i}\in C^<_{\delta}\\ \Sigmama_{\vec{i}} \text{ not diagonal} }} 1.
\end{align*}
By \eqref{eqn:firsttermtemp}, 
\begin{align*}
&\sum_{\vec{i}\in C^<_{\delta}} \E \Phi_{\vec{i}}^{(\Rma)} \left( \Phi_{\vec{i}}^{(\Rma)} + U_{\vec{i}}\right) - \binom{p}{1}\binom{p-1}{\delta} (2P_n(r_\rho))^{\delta} \sum_{\ell=1}^{\delta+1} \ell  \alpha(\ell,r_\rho) \\
\geq &
- (2P_n(r_\rho))^{\delta} \sum_{\ell=1}^{\delta+1} \ell  \alpha(\ell,r_\rho) \sum_{\substack{\vec{i}\in C^<_{\delta}\\ \Sigmama_{\vec{i}} \text{ not diagonal} }} 1.
\end{align*}
By combining the preceding two expressions, 
\begin{align*}
&\left|\sum_{\vec{i}\in C^<_{\delta}} \E \Phi_{\vec{i}}^{(\Rma)} \left( \Phi_{\vec{i}}^{(\Rma)} + U_{\vec{i}}\right) - \binom{p}{1}\binom{p-1}{\delta} (2P_n(r_\rho))^{\delta}\sum_{\ell=1}^{\delta+1} \ell  \alpha(\ell,r_\rho)\right| \\
\leq &
\mu_{n,\delta+1}(\Sigmama) (2P_n(r_\rho))^{\delta} \sum_{\ell=1}^{\delta+1} \ell  \alpha(\ell,r_\rho) \sum_{\substack{\vec{i}\in C^<_{\delta}\\ \Sigmama_{\vec{i}} \text{ not diagonal} }} 1 \\
\leq & \mu_{n,\delta+1}(\Sigmama) (2P_n(r_\rho))^{\delta} \sum_{\ell=1}^{\delta+1} \ell  \alpha(\ell,r_\rho) \frac{(\delta+1)}{2((\delta-1)!)}   p^{\delta} (\kappa-1) \\
\leq & \mu_{n,\delta+1}(\Sigmama) (2p^{1+1/\delta}P_n(r_\rho))^{\delta}  \frac{(\delta+1)^2}{2((\delta-1)!)}    \frac{\kappa-1}{p}
\numberthis \label{eqn:2ndmomentstep1}
\end{align*}
where the second inequality follows from Lemma \ref{lem:numofnondia}.\\

\noindent
\textbf{Step 2}:\\
\begin{align*}
&\sum_{\vec{i}\in C^<_{\delta}} \E \Phi_{\vec{i}}^{(\Rma)} W_{\vec{i}} - \left(\binom{p}{1}\binom{p-1}{\delta}(2P_n(r_\rho))^{\delta}\right)^2 \\
= & 
\sum_{\vec{i}\in C^<_{\delta}} \sum_{\vec{j}\in T_{\vec{i}}} \E \Phi_{\vec{i}}^{(\Rma)} \Phi_{\vec{j}}^{(\Rma)}-\sum_{\vec{i}\in C^<_{\delta}} \sum_{\vec{j}\in C^<_{\delta}} (2P_n(r_\rho))^{2\delta}\\
= & 
\sum_{\vec{i}\in C^<_{\delta}} \sum_{\vec{j}\in T_{\vec{i}}} \E \Phi_{\vec{i}}^{(\Rma)} \E\Phi_{\vec{j}}^{(\Rma)} - \sum_{\vec{i}\in C^<_{\delta}} \sum_{\vec{j}\in T_{\vec{i}}} (2P_n(r_\rho))^{2\delta}  -\sum_{\vec{i}\in C^<_{\delta}} \sum_{\vec{j}\in C^<_{\delta}\backslash T_{\vec{i}}} (2P_n(r_\rho))^{2\delta}\\
= & 
-\sum_{\vec{i}\in C^<_{\delta}} \sum_{\vec{j}\in C^<_{\delta}\backslash T_{\vec{i}}} (2P_n(r_\rho))^{2\delta}
\sum_{\substack{\vec{i}\in C^<_{\delta}\\ \Sigmama_{\vec{i}} \text{ not diagonal} }} \sum_{\vec{j}\in T_{\vec{i}}} \left(\E \Phi_{\vec{i}}^{(\Rma)} \Phi_{\vec{j}}^{(\Rma)} 
-  (2P_n(r_\rho))^{2\delta}\right)  \\
&  +\sum_{\substack{\vec{i}\in C^<_{\delta}\\ \Sigmama_{\vec{i}} \text{ diagonal} }} \  \sum_{ \substack{\vec{j}\in T_{\vec{i}}\\ \Sigmama_{\vec{j}} \text{ not diagonal}   }} \left(\E \Phi_{\vec{i}}^{(\Rma)} \Phi_{\vec{j}}^{(\Rma)} -  (2P_n(r_\rho))^{2\delta}\right), 
\numberthis \label{eqn:3rdtermtemp}
\end{align*}
where the last equality follows from $\E \Phi_{\vec{i}}^{(\Rma)} =(2P_n(r_\rho))^{\delta}$ for $\vec{i}\in \left\{ \vec{j}\in C_\delta^<:  \Sigmama_{\vec{j}} \text{ diagonal} \right\}$  by  Lemma \ref{item:prod}. Then by \eqref{eqn:3rdtermtemp},
\begin{align*}
&\sum_{\vec{i}\in C^<_{\delta}} \E \Phi_{\vec{i}}^{(\Rma)} W_{\vec{i}} - \left(\binom{p}{1}\binom{p-1}{\delta}(2P_n(r_\rho))^{\delta}\right)^2\\ 
\leq & \left(\mu_{n,2\delta+2}(\Sigmama) -  1\right)(2P_n(r_\rho))^{2\delta} 
\left(\sum_{\substack{\vec{i}\in C^<_{\delta}\\ \Sigmama_{\vec{i}} \text{ not diagonal} }} \sum_{\vec{j}\in T_{\vec{i}}} 1 +\sum_{\substack{\vec{i}\in C^<_{\delta}\\ \Sigmama_{\vec{i}} \text{ diagonal} }} \  \sum_{ \substack{\vec{j}\in T_{\vec{i}}\\ \Sigmama_{\vec{j}} \text{ not diagonal}   }} 1\right),
\end{align*}
where the inequality follows from Lemma \ref{lem:densitybou} \ref{item:densitybouc}. On the other hand, by \eqref{eqn:3rdtermtemp}, 
\begin{align*}
   & \sum_{\vec{i}\in C^<_{\delta}} \E \Phi_{\vec{i}}^{(\Rma)} W_{\vec{i}} - \left(\binom{p}{1}\binom{p-1}{\delta}(2P_n(r_\rho))^{\delta}\right)^2 \\ 
   \geq &  -  (2P_n(r_\rho))^{2\delta} \left(\sum_{\substack{\vec{i}\in C^<_{\delta}\\ \Sigmama_{\vec{i}} \text{ not diagonal} }} \sum_{\vec{j}\in T_{\vec{i}}} 1 +\sum_{\substack{\vec{i}\in C^<_{\delta}\\ \Sigmama_{\vec{i}} \text{ diagonal} }} \  \sum_{ \substack{\vec{j}\in T_{\vec{i}}\\ \Sigmama_{\vec{j}} \text{ not diagonal}   }} 1 \right) -\sum_{\vec{i}\in C^<_{\delta}} \sum_{\vec{j}\in C^<_{\delta}\backslash T_{\vec{i}}} (2P_n(r_\rho))^{2\delta}.
\end{align*}
Combining the preceding two displays,
\begin{align*}
& \left|\sum_{\vec{i}\in C^<_{\delta}} \E \Phi_{\vec{i}}^{(\Rma)} W_{\vec{i}} - \left(\binom{p}{1}\binom{p-1}{\delta}(2P_n(r_\rho))^{\delta}\right)^2\right| \\ 
\leq & 
\mu_{n,2\delta+2}(\Sigmama) (2P_n(r_\rho))^{2\delta} 
\left(\sum_{\substack{\vec{i}\in C^<_{\delta}\\ \Sigmama_{\vec{i}} \text{ not diagonal} }} \sum_{\vec{j}\in T_{\vec{i}}} 1 +\sum_{\substack{\vec{i}\in C^<_{\delta}\\ \Sigmama_{\vec{i}} \text{ diagonal} }} \  \sum_{ \substack{\vec{j}\in T_{\vec{i}}\\ \Sigmama_{\vec{j}} \text{ not diagonal}   }} 1\right) \\
&\quad \quad \quad + \sum_{\vec{i}\in C^<_{\delta}} \sum_{\vec{j}\in C^<_\delta\backslash T_{\vec{i}}} (2P_n(r_\rho))^{2\delta} \\
\leq & 
\mu_{n,2\delta+2}(\Sigmama) (2P_n(r_\rho))^{2\delta} \left(2\sum_{\substack{\vec{i}\in C^<_{\delta}\\ \Sigmama_{\vec{i}} \text{ not diagonal} }} \sum_{\vec{j}\in C_\delta^<} 1\right)  + \sum_{\vec{i}\in C^<_{\delta}} \sum_{\vec{j}\in C^<_\delta\backslash T_{\vec{i}}} (2P_n(r_\rho))^{2\delta}\\
\overset{(*)}{\leq} & 
\mu_{n,2\delta+2}(\Sigmama) (2P_n(r_\rho))^{2\delta} 2\binom{p}{1}\binom{p-1}{\delta} \frac{(\delta+1)}{2((\delta-1)!)}   p^{\delta} (\kappa-1) \\
&\quad \quad + \binom{p}{1}\binom{p-1}{\delta} \frac{(\delta+1)^2}{\delta !}  p^{\delta}\kkkk (2P_n(r_\rho))^{2\delta}\\
\leq&  
\frac{2(\delta+1)^2}{(\delta!)^2} (2p^{1+1/\delta}P_n(r_\rho))^{2\delta} \mu_{n,2\delta+2}(\Sigmama)\frac{\kappa}{p}, \numberthis \label{eqn:2ndmomentstep2}
\end{align*}
where step $(*)$ follows from Lemma \ref{lem:numofnondia} and \eqref{eqn:neiborsize}.

\noindent
\textbf{Step 3}:\\
Notice $\sum_{\vec{i}\in C^<_{\delta}} \E \Phi_{\vec{i}}^{(\Rma)} Z_{\vec{i}} = b_2$ as in \eqref{eqn:b2def} and thus satisfies the bound \eqref{eqn:b2upperbou2new}. Then by \eqref{eqn:squaredecom},
\begin{align*}
    &\left| \E \left(N_{E_{\delta}}^{(\Rma)}\right)^2-\binom{p}{1}\binom{p-1}{\delta} (2P_n(r_\rho))^{\delta}\sum_{\ell=1}^{\delta+1} \ell  \alpha(\ell,r_\rho) -\left(\binom{p}{1}\binom{p-1}{\delta}(2P_n(r_\rho))^{\delta}\right)^2 \right| \\
    \leq & \left|\sum_{\vec{i}\in C^<_{\delta}} \E \Phi_{\vec{i}}^{(\Rma)} \left( \Phi_{\vec{i}}^{(\Rma)} + U_{\vec{i}}\right) - \binom{p}{1}\binom{p-1}{\delta} (2P_n(r_\rho))^{\delta} \sum_{\ell=1}^{\delta+1} \ell  \alpha(\ell,r_\rho) \right|\\
    &+ \left|\sum_{\vec{i}\in C^<_{\delta}} \E \Phi_{\vec{i}}^{(\Rma)} W_{\vec{i}} - \left(\binom{p}{1}\binom{p-1}{\delta}(2P_n(r_\rho))^{\delta}\right)^2\right| + b_2\\
    \leq & C_{n,\delta} \left( (2p^{1+1/\delta}P_n(r_\rho))^{\delta} \left(1+2p^{1+1/\delta}P_n(r_\rho)\right)^{\delta} \mu_{n,2\delta+2}(\Sigmama)\frac{\kappa}{p}\right.\\
    &\left. \quad +  p(2pP_n(r_\rho))^{\delta+1} \left(1+2pP_n(r_\rho)\right)^{\delta}\right), \numberthis \label{eqn:2ndmomentdelta}
\end{align*}
where the last inequality follows from \eqref{eqn:2ndmomentstep1}, \eqref{eqn:2ndmomentstep2} and \eqref{eqn:b2upperbou2new}. The proof is then completed by using the inequality $2p^{1+1/\delta}P_n(r_\rho)\leq \gamma$.
\end{proof}

\subsubsection{Proof of Proposition \ref{prop:allcloseL2}}

\begin{proof}[Proof of Proposition \ref{prop:allcloseL2} \ref{item:allcloseL2a}]\ \\
By taking the square of each of the terms in Lemma \ref{lem:6quantitiesine}, 
$$
\left(N_{E_\delta}^{(\Rma)}\right)^2 - 2(\delta+1)N_{E_\delta}^{(\Rma)} N_{E_{\delta+1}}^{(\Rma)}   \leq \left(N_{\breve{V}_\delta}^{(\Rma)}\right)^2 \leq \left(N_{V_\delta}^{(\Rma)}\right)^2 \leq \left(N_{E_\delta}^{(\Rma)}\right)^2,
$$
which then implies for $\bar{N}_\delta \in \left\{ N_{\breve{V}_\delta}^{(\Rma)}, N_{V_\delta}^{(\Rma)} \right\} $
\begin{equation}
\left|\left(\bar{N}_{\delta}\right)^2-\left(N_{E_\delta}^{(\Rma)}\right)^2\right| \leq 2(\delta+1)N_{E_\delta}^{(\Rma)} N_{E_{\delta+1}}^{(\Rma)} = 2(\delta+1)\sum_{\vec{i}\in C_{\delta+1}^<} \sum_{\vec{j}\in C_{\delta}^<} \Phi_{\vec{i}}^{(\bm{R})} \Phi_{\vec{j}}^{(\bm{R})}. \label{eqn:squarecorrelationuppbou}
\end{equation}
It suffices to establish an upper bound on $\E N_{E_{\delta+1}}^{(\Rma)} N_{E_\delta}^{(\Rma)} = \E \sum_{\vec{i}\in C_{\delta+1}^<} \sum_{\vec{j}\in C_{\delta}^<} \Phi_{\vec{i}}^{(\bm{R})} \Phi_{\vec{j}}^{(\bm{R})} $. 

Observe for $\vec{j}\in J_{\vec{i}}$, 
$$
\E \Phi_{\vec{i}}^{(\bm{R})} \Phi_{\vec{j}}^{(\bm{R})} \leq \E \Phi_{\vec{i}}^{(\bm{R})} \leq \mu_{n,\delta+2}(\Sigmama_{\vec{i}}) (2P_n(r_\rho))^{\delta+1}.
$$

For $\vec{j} \in T_{\vec{i}}$, $\left[\vec{j}\right]\cap \left[\vec{i}\right]=\emptyset$. Thus, if $\Sigmama_{\vec{i}\cup \vec{j}}$ is diagonal, $\E \Phi_{\vec{i}}^{(\bm{R})} \Phi_{\vec{j}}^{(\bm{R})} = \E \Phi_{\vec{i}}^{(\bm{R})} \E \Phi_{\vec{j}}^{(\bm{R})} =(2P_n(r_\rho))^{2\delta+1}$ by Lemma \ref{item:prod}. Then, for the general case when $\Sigmama_{\vec{i}\cup \vec{j}}$ is not necessarily diagonal, by Lemma \ref{lem:densitybou} \ref{item:densitybouc},
$$
\E \Phi_{\vec{i}}^{(\bm{R})} \Phi_{\vec{j}}^{(\bm{R})} \leq \mu_{n,2\delta+3}(\Sigmama_{\vec{i}\cup\vec{j}}) (2P_n(r_\rho))^{2\delta+1}
$$

By Lemma \ref{lem:densitybou} \ref{item:densitybouc}, Lemma  \ref{item:prod} and Lemma \ref{lem:eleine}, it is straightforward  that the conditions in Lemma \ref{lem:b2} with $q=\delta+1$ and $\theta_{\vec{i},\vec{j}}=\Phi_{\vec{i}}^{(\Rma)}\Phi_{\vec{j}}^{(\Rma)}$ are satisfied with $a=1$, $b=2P_n(2r_\rho)1(\delta\geq 2)+2P_n(r_\rho)1(\delta= 1)$ and $z=2P_n(r_\rho)$. Moreover, $b/z\leq 2^{n-2}1(\delta\geq 2)+1$ by Lemma \ref{pderivative} \ref{item:pderivatived}. 

Thus, Lemma \ref{lem:b2wholesum} with $q=\delta+1$ and $\theta_{\vec{i},\vec{j}}=\Phi_{\vec{i}}^{(\Rma)}\Phi_{\vec{j}}^{(\Rma)}$, $a=1$, $b=2P_n(2r_\rho)1(\delta\geq 2)+2P_n(r_\rho)1(\delta= 1)$ and $z=2P_n(r_\rho)$, together with the fact that $pz\leq p^{1+\frac{1}{\delta}}z = 2p^{1+\frac{1}{\delta}}P_n(r_\rho) \leq \gamma$,
yield
\begin{align}
\E N_{E_{\delta+1}}^{(\Rma)} N_{E_\delta}^{(\Rma)} = \E \sum_{\vec{i}\in C_{\delta+1}^<} \sum_{\vec{j}\in C_{\delta}^<} \Phi_{\vec{i}}^{(\bm{R})} \Phi_{\vec{j}}^{(\bm{R})}
\leq C_{n,\delta,\gamma}  \left(1+\mu_{n,2\delta+3}(\Sigmama)\frac{\kappa-1}{p}\right)p^{-1/\delta}.
\label{eqn:NEdelta1NEdelta}
\end{align}
The proof is then complete by the preceding expression and \eqref{eqn:squarecorrelationuppbou}.

\end{proof}

We now present a few lemmas that are used in the proof of Proposition \ref{prop:allcloseL2} \ref{item:allcloseL2b} and \ref{item:allcloseL2c}. Recall that $F_{ij}(r_\rho),H_{ij}(r_\rho), G_{ij}(r_\rho), F_{\vec{i}}(r_\rho)$ are defined in Section \ref{sec:proofthmallclosec}.

\begin{lem}
\label{lem:inddifbou2ndmoment}
Suppose $p\geq n$. $1\leq \delta\leq q \leq p-1$. Then for any $\vec{i}\in C_q^<$, $\vec{j}\in C_{\delta}^<$, with probability $1$,
	$$
	\left|\Phi_{\vec{i}}^{(\bm{P})}\Phi_{\vec{j}}^{(\bm{P})}-\Phi_{\vec{i}}^{(\bm{R})}\Phi_{\vec{j}}^{(\bm{R})}\right|\leq \xi_{\vec{i},\vec{j}},
	$$
	where 
	$$
\xi_{\vec{i},\vec{j}} := 1\left( \bigcup\limits_{m=1}^q \left ( \left( H_{i_0i_m} \big \backslash G_{i_0i_m}   \right) \bigcap  H_{\vec{i},-m} \bigcap H_{\vec{j}} \right ) \bigcup \bigcup\limits_{\ell=1}^{\delta} \left ( \left( H_{j_0j_{\ell}} \big \backslash G_{j_0j_{\ell}}   \right) \bigcap  H_{\vec{j},-\ell} \bigcap H_{\vec{i}} \right ) \right). 
$$
\end{lem}
\begin{proof}
\begin{align*}
\left|\Phi_{\vec{i}}^{(\bm{P})}\Phi_{\vec{j}}^{(\bm{P})}-\Phi_{\vec{i}}^{(\bm{R})}\Phi_{\vec{j}}^{(\bm{R})}\right| =& 1\left(\left(\bigcap\limits_{m=1}^q F_{i_0i_m} \bigcap \bigcap\limits_{\ell=1}^{\delta} F_{j_0j_\ell} \right)\bigtriangleup \left(\bigcap\limits_{m=1}^q S_{i_0i_m} \bigcap \bigcap\limits_{\ell=1}^{\delta} S_{j_0j_\ell} \right) \right)\leq \xi_{\vec{i},\vec{j}}
\end{align*}
where the inequality follows from \eqref{eqn:FGJSsetinclusion} and Lemma \ref{setrelation} \ref{prop:setinclusiona}.
\end{proof}

\begin{lem}
\label{lem:xietarelationship}
Let $p\geq n\geq 4$, $1\leq \delta\leq q \leq p-1$ and $\Xma\sim \mathcal{VE}(\bm{\mu},\Sigmama,g)$. Suppose $\Sigmama$, after some row-column permutation, is $(\tau,\kappa)$ sparse with $\tau\leq \frac{p}{2}$. Let $t$ be any positive number, and suppose \eqref{eqn:postempass} holds. Then for any $\vec{i}\in C_q^<$, $\vec{j}\in C_{\delta}^<$, with probability $1$,
$$
\xi_{\vec{i},\vec{j}} 1\left(\Ec(t)\right) \leq \eta_{\vec{i},\vec{j}}(t),
$$
where 
\begin{align*}
    \eta_{\vec{i},\vec{j}}(t) := &
    1 \left(  \bigcup_{m=1}^q \left( \left(F_{i_0i_m}(\theta_1(t)r_\rho)\big \backslash F_{i_0i_m}\left(\frac{r_\rho}{\theta_1(t)}\right) \right)  \bigcap  F_{\vec{i},-m}(\theta_1(t)r_\rho) \bigcap F_{\vec{j}}(\theta_1(t)r_\rho) \right) \bigcup \right.\\
    &\left.  \quad  \bigcup_{\ell=1}^\delta \left( \left(F_{j_0j_\ell}(\theta_1(t)r_\rho)\big \backslash F_{j_0j_\ell}\left(\frac{r_\rho}{\theta_1(t)}\right) \right)  \bigcap  F_{\vec{j},-\ell}(\theta_1(t)r_\rho) \bigcap F_{\vec{i}}(\theta_1(t)r_\rho) \right)
    \right).
\end{align*}
\end{lem}
\begin{proof}
By \eqref{ctempdef}, $H_{ij}(r_\rho) \cap \Ec(t)\subset F_{ij}(\theta_1(t)r_\rho)$ and $G_{ij}(r_\rho) \cap \Ec(t)\supset F_{ij}\left(\frac{r_\rho}{\theta_1(t)}\right)$. Then
	$$
	\xi_{\vec{i},\vec{j}} 1\left(\Ec(t)\right)
	\leq  \eta_{\vec{i},\vec{j}}(t). 
	$$
\end{proof}

\begin{lem}
\label{lem:etadoublesumexpectation}
Let $p\geq n\geq 4$, $1\leq \delta\leq q \leq p-1$ and $\Xma\sim \mathcal{VE}(\bm{\mu},\Sigmama,g)$. Suppose $\Sigmama$, after some row-column permutation, is $(\tau,\kappa)$ sparse with $\tau\leq \frac{p}{2}$. Let $t$ be any positive number, and suppose \eqref{eqn:postempass} holds. Suppose additionally $2p^{1+\frac{1}{\delta}}P_n(r_\rho)\leq \gamma$. Then
\begin{align*}
&\E \sum_{\vec{i}\in C_q^<} \sum_{\vec{j}\in C_\delta^<}  \eta_{\vec{i},\vec{j}}(t) \\
\leq &
C_{n,q,\delta,\gamma} \left(\theta_1(t)\right)^{n(2\delta+q)} \left(1+\mu_{n,q+\delta+2}(\Sigmama)\frac{\kappa-1}{p}\right) \left(\sqrt{\frac{1}{p}}+\frac{t}{\sqrt{p}} + \frac{\tau}{p}\right) p^{1-\frac{q}{\delta}}.
\end{align*}
\end{lem}
\begin{proof}
 Note
\begin{align*}
  \sum_{\vec{i}\in C_q^<}  \sum_{\vec{j}\in C_\delta^<} \eta_{\vec{i},\vec{j}}(t) = \left(\sum_{\vec{i}\in C_q^<}  \sum_{\vec{j}\in J_{\vec{i}}} + \sum_{\vec{i}\in C_q^<}  \sum_{\vec{j}\in T_{\vec{i}}} +\sum_{\vec{i}\in C_q^<}  \sum_{\vec{j}\in N_{\vec{i}}}\right) 
\eta_{\vec{i},\vec{j}}(t).
\end{align*}

\noindent
\textbf{Step 1: $\vec{j}\in T_{\vec{i}}$}\\
By the union bound for indicator functions,
\begin{equation}
\eta_{\vec{i},\vec{j}}(t) \leq  \eta_{\vec{i}}(t) 1\left(F_{\vec{j}}(\theta_1(t)r_\rho)\right) +  \eta_{\vec{j}}(t) 1\left(F_{\vec{i}}(\theta_1(t)r_\rho)\right), \label{eqn:etaunionbound}
\end{equation}
where $\eta_{\vec{i}}(t)$ is defined in \eqref{eqn:etavecidef} with $\delta$ replaced by $q$. Then for $\vec{j}\in T_{\vec{i}}$
\begin{equation}
\E\eta_{\vec{i},\vec{j}}(t) \leq  \E\eta_{\vec{i}}(t) \P\left(F_{\vec{j}}(\theta_1(t)r_\rho)\right) +  \E\eta_{\vec{j}}(t) \P\left(F_{\vec{i}}(\theta_1(t)r_\rho)\right). \label{eqn:etadoubleindexuppbouTi}
\end{equation}
Moreover, for $\vec{j}\in T_{\vec{i}}$, $\Sigmama_{\vec{i}\cup \vec{j} }$ is diagonal if and only if $\Sigmama_{\vec{i}}$ and $\Sigmama_{\vec{j}}$ are both diagonal. 

Now suppose $\Sigmama_{\vec{i}\cup \vec{j} }$ is diagonal. By conditioning on $\uve_{j_0}$ 
$$
\P\left(F_{\vec{j}}(\theta_1(t)r_\rho)\right) = \left(2P_n(\theta_1(t)r_\rho)\right)^\delta, \quad \P\left(F_{\vec{i}}(\theta_1(t)r_\rho)\right) = \left(2P_n(\theta_1(t)r_\rho)\right)^q.
$$
The preceding expression, \eqref{eqn:etadoubleindexuppbouTi}, Lemma \ref{lem:inddifuppexpbou} applied to $\E\eta_{\vec{j}}(t)$, and Lemma \ref{lem:inddifuppexpbou} with $\delta=q$ applied to $\E\eta_{\vec{i}}(t)$ yield
$$
\E \eta_{\vec{i},\vec{j}}(t) \leq (\delta+q)  2\left(P_\nn(r_\rho \theta_1(t))-P_\nn\left(\frac{r_\rho}{\theta_1(t)} \right)\right) \left(2P_\nn\left(\theta_1(t)r_\rho\right)\right)^{q+\mylll-1}.
$$

For the general case that $\Sigmama_{\vec{i}\cup \vec{j} }$ is not necessarily diagonal, by Lemma \ref{lem:densitybou} \ref{item:densitybouc}, for any $\vec{j}\in T_{\vec{i}}$
$$
\E \eta_{\vec{i},\vec{j}}(t) \leq \mu_{n,q+\delta+2}(\Sigmama_{\vec{i}\cup\vec{j}}) (\delta+q)  2\left(P_\nn(r_\rho \theta_1(t))-P_\nn\left(\frac{r_\rho}{\theta_1(t)} \right)\right) \left(2P_\nn\left(\theta_1(t)r_\rho\right)\right)^{q+\mylll-1}.
$$
Then the condition in Lemma \ref{lem:quadraticsum} \ref{item:quadraticsumb} is satisfied with $\theta_{\vec{i},\vec{j}}=\eta_{\vec{i},\vec{j}}$, $z=2P_\nn(r_\rho \theta_1(t))$ and $a = (\delta+q)  \frac{P_\nn(r_\rho \theta_1(t))-P_\nn\left(\frac{r_\rho}{\theta_1(t)} \right)}{ P_\nn(r_\rho \theta_1(t)) } $.\\

\noindent
\textbf{Step 2: $\vec{j}\in J_{\vec{i}}$}\\
\eqref{eqn:etaunionbound} implies
\begin{align}
    \eta_{\vec{i},\vec{j}}(t)\leq & \eta_{\vec{i}}(t)+ \sum_{\ell=1}^{\delta} 1\left(F_{j_0j_\ell}(\theta_1(t)r_\rho)\big \backslash F_{j_0j_\ell}\left(\frac{r_\rho}{\theta_1(t)}\right) \right) 1(F_{\vec{i}}(\theta_1(t)r_\rho)). \label{eqn:step1etadoubleindexuppbou}
\end{align}
 For $\vec{j}\in J_{\vec{i}}$, $j_0,j_\ell\in \left[ \vec{i} \right]$. If $i_0\in \{j_0,j_\ell\}$, without loss of generality, say $i_0=j_\ell$ and $j_0=i_\alpha$ for some $1\leq \alpha\leq q$. Then 
\begin{align*}
&\E 1\left(F_{j_0j_\ell}(\theta_1(t)r_\rho)\big \backslash F_{j_0j_\ell}\left(\frac{r_\rho}{\theta_1(t)}\right) \right) 1(F_{\vec{i}}(\theta_1(t)r_\rho)) \\
= &
\P \left( \left(F_{i_0i_\alpha}(\theta_1(t)r_\rho)\big \backslash F_{i_0i_\alpha}\left(\frac{r_\rho}{\theta_1(t)}\right) \right)  \bigcap  F_{\vec{i},-\alpha}(\theta_1(t)r_\rho) \right)\\
\leq & 
\mu_{n,q+1}(\Sigmama_{\vec{i}}) 2\left(P_\nn(r_\rho \theta_1(t))-P_\nn\left(\frac{r_\rho}{\theta_1(t)} \right)\right) \left(2P_\nn\left(\theta_1(t)r_\rho\right)\right)^{q-1}, \numberthis \label{eqn:etadoubleindexcase1}
\end{align*}
where the last step follows from Lemma \ref{lem:inddifuppexpbou} with $\delta$ replace by $q$.

If $i_0\not \in \{j_0,j_\ell\}$, without loss of generality, let $j_0=i_\alpha$, $j_\ell=i_{\beta}$ for some $1\leq \alpha \neq \beta \leq q$. Suppose for now that $\Sigmama_{\vec{i}}$ is diagonal, and then
\begin{align*}
&\E 1\left(F_{j_0j_\ell}(\theta_1(t)r_\rho)\big \backslash F_{j_0j_\ell}\left(\frac{r_\rho}{\theta_1(t)}\right) \right) 1(F_{\vec{i}}(\theta_1(t)r_\rho)) \\
= & 
\E 1\left(F_{i_\alpha i_\beta}(\theta_1(t)r_\rho)\big \backslash F_{i_\alpha i_\beta}\left(\frac{r_\rho}{\theta_1(t)}\right) \right) 1(F_{\vec{i}}(\theta_1(t)r_\rho))  \\
\overset{(*)}{=} & 
(2P_n(\theta_1(t)r_\rho))^{q-2} \E 1\left(F_{i_\alpha i_\beta}(\theta_1(t)r_\rho)\setminus F_{i_\alpha i_\beta} \left(\frac{r_\rho}{\theta_1(t)}\right) \right) 1(F_{i_0i_\alpha}(\theta_1(t)r_\rho))1(F_{i_0i_\beta}(\theta_1(t)r_\rho))
\\
\overset{(**)}{\leq} &
2\left (P_n(\theta_1(t)r_\rho)-P_n\left(\frac{r_\rho}{\theta_1(t)}\right)\right) (2P_n(\theta_1(t)r_\rho))^{q-1} , 
\end{align*}
where the step $(*)$ follows from conditioning on $\uve_{i_0},\uve_{i_\alpha},\uve_{i_\beta}$, and the step $(**)$ follows from dropping the term $1(F_{i_0i_\beta}(\theta_1(t)r_\rho))$ and then conditioning on $\uve_{i_\alpha}$. For the general case that $\Sigmama_{\vec{i}}$ is not necessarily diagonal, by Lemma \ref{lem:densitybou} \ref{item:densitybouc}, 
\begin{align*}
&\E 1\left(F_{j_0j_\ell}(\theta_1(t)r_\rho)\big \backslash F_{j_0j_\ell}\left(\frac{r_\rho}{\theta_1(t)}\right) \right) 1(F_{\vec{i}}(\theta_1(t)r_\rho)) \\
\leq &
\mu_{n,q+1}(\Sigmama_{\vec{i}}) 2\left (P_n(\theta_1(t)r_\rho)-P_n\left(\frac{r_\rho}{\theta_1(t)}\right)\right) (2P_n(\theta_1(t)r_\rho))^{q-1}. \numberthis \label{eqn:etadoubleindextemp4}
\end{align*}

By combining \eqref{eqn:step1etadoubleindexuppbou}, \eqref{eqn:etadoubleindexcase1}, \eqref{eqn:etadoubleindextemp4} and Lemma \ref{lem:inddifuppexpbou} with $\delta$ replace by $q$,
$$
\E \eta_{\vec{i},\vec{j}}(t) \leq (q+\delta)\mu_{n,q+1}(\Sigmama_{\vec{i}}) 2\left (P_n(\theta_1(t)r_\rho)-P_n\left(\frac{r_\rho}{\theta_1(t)}\right)\right) (2P_n(\theta_1(t)r_\rho))^{q-1}.
$$
Then the conditions in Lemma \ref{lem:quadraticsum} \ref{item:quadraticsuma} with $\theta_{\vec{i},\vec{j}}=\eta_{\vec{i},\vec{j}}$, $z=2P_n(r_\rho\theta_1(t))$ and $a=(q+\delta) \frac{P_n(\theta_1(t)r_\rho)-P_n\left(\frac{r_\rho}{\theta_1(t)}\right)}{P_n(\theta_1(t)r_\rho)}$ is satisfied.\\

\noindent
\textbf{Step 3: $\vec{j}\in N_{\vec{i}}$}\\
It is straightforward by Lemma \ref{lem:densitybou} \ref{item:densitybouc}, Lemma  \ref{item:prod} and Lemma \ref{lem:eleine} that the conditions in Lemma \ref{lem:b2} with  $\theta_{\vec{i},\vec{j}}=\eta_{\vec{i},\vec{j}}$ are satisfied with $a=a_1=(q+\delta) \frac{P_n(\theta_1(t)r_\rho)-P_n\left(\frac{r_\rho}{\theta_1(t)}\right)}{P_n(\theta_1(t)r_\rho)}$, $b=2P_n(2r_\rho\theta_1(t))1(\delta\geq 2)+2P_n(r_\rho\theta_1(t))1(\delta= 1)$ and $z=2P_n(r_\rho\theta_1(t))$. Moreover, $b/z\leq 2^{n-2}1(\delta\geq 2)+1$ by Lemma \ref{pderivative} \ref{item:pderivatived}.

By Lemma \ref{lem:b2wholesum} with $\theta_{\vec{i},\vec{j}}=\eta_{\vec{i},\vec{j}}$, $a=(q+\delta) \frac{P_n(\theta_1(t)r_\rho)-P_n\left(\frac{r_\rho}{\theta_1(t)}\right)}{P_n(\theta_1(t)r_\rho)}$, $b=2P_n(2r_\rho\theta_1(t))1(\delta\geq 2)+2P_n(r_\rho\theta_1(t))1(\delta= 1)$ and $z=2P_n(r_\rho\theta_1(t))$
\begin{align*}
&\sum_{\vec{i}\in C_q^<} \sum_{\vec{j}\in C_\delta^<} \E \eta_{\vec{i},\vec{j}} \\
\leq & 
C_{n,q,\delta} \left(p^{1+\frac{1}{\delta}}z\right)^q \left(1+(p^{1+\frac{1}{\delta}}z)^{\delta}\right)\left(1+pz\right)^{\delta-1} \left(1+\mu_{n,q+\delta+2}(\Sigmama)\frac{\kappa-1}{p}\right) ap^{1-\frac{q}{p}}.  \numberthis \label{eqn:etadoublesumtemp}
\end{align*}

\noindent
\textbf{Step 4}\\
Observe that
\begin{align}
    pz \leq  p^{1+\frac{1}{\delta}} z \leq  \left(\theta_1(t)\right)^{n-2} 2p^{1+\frac{1}{\delta}}P_n(r_\rho) 
       \leq  \left(\theta_1(t)\right)^{n-2} \gamma \label{eqn:pzbou}
\end{align}
where the second inequality follows from Lemma \ref{pderivative} \ref{item:pderivatived}. Moreover, by Lemma \ref{pderivative} \ref{item:pderivativec} and the fact that $\theta_1(t)\geq 1$,
\begin{equation}
a\leq (q+\delta) \frac{P_n(\theta_1(t)r_\rho)-P_n\left(\frac{r_\rho}{\theta_1(t)}\right)}{P_n(r_\rho)} \leq (q+\delta)(n-2)\left(\theta_1(t)\right)^{n-3}\left(\theta_1(t)-\frac{1}{\theta_1(t)}\right). \label{eqn:aupperbound}
\end{equation}

Plugging \eqref{eqn:pzbou} and \eqref{eqn:aupperbound} into \eqref{eqn:etadoublesumtemp} and by the fact that $\theta_1(t)\geq 1$,  
\begin{align*}
    \sum_{\vec{i}\in C_q^<} \sum_{\vec{j}\in C_\delta^<} \eta_{\vec{i},\vec{j}}
    \leq &  
    C_{n,q,\delta,\gamma} \left(\theta_1(t)\right)^{n(2\delta+q)} \left(1+\mu_{n,q+\delta+2}(\Sigmama)\frac{\kappa-1}{p}\right) \left(\theta_1(t)-\frac{1}{\theta_1(t)}\right)p^{1-\frac{q}{p}}\\
    \leq &  C_{n,q,\delta,\gamma} \left(\theta_1(t)\right)^{n(2\delta+q)} \left(1+\mu_{n,q+\delta+2}(\Sigmama)\frac{\kappa-1}{p}\right) \left(\sqrt{\frac{\nn}{p}}+\frac{t}{\sqrt{p}} + n\frac{\tau}{p}\right) p^{1-\frac{q}{p}},  
\end{align*}
where the last inequality follows from Lemma \ref{basine} \ref{item:basinek} and \eqref{ctempdef}. 
\end{proof}

\begin{lem} 
\label{lem:NEqNEdelta}
	Let $p\geq n \geq 4$ and $\Xma\sim \mathcal{VE}(\bm{\mu},\Sigmama,g)$. Suppose $\Sigmama$, after some row-column permutation, is $(\tau,\kappa)$ sparse with $\tau\leq \frac{p}{2}$. Consider $1\leq \delta \leq p-2$ and let $q\in \{\delta,\delta+1\}$. Suppose  $2p^{1+\frac{1}{\delta}}P_\nn(r_\rho)\leq \gamma$ and $\left(\sqrt{\frac{\nn-1}{p}}+\sqrt{\frac{\delta\ln p}{p}}\right) \leq c$ hold for some positive and \ONE{sufficiently} small universal constant $c$. Then
	\begin{align*}
	\E \left|N_{E_q}^{(\Pma)} N_{E_\delta}^{(\Pma)}-N_{E_q}^{(\Rma)} N_{E_\delta}^{(\Rma)}\right|\leq C_{n,\delta,\gamma}  \left(1+\mu_{n,q+\delta+2}(\Sigmama)\frac{\kappa-1}{p}\right) \left(\frac{\sqrt{\ln p}}{\sqrt{p}} + \frac{\tau}{p}\right) p^{1-\frac{q}{\delta}}.
	\end{align*}
\end{lem}

\begin{proof}

For $\Psima \in \{\Rma,\Pma\}$,
$$
N_{E_q}^{(\Psima)} N_{E_\delta}^{(\Psima)} = \sum_{\vec{i}\in C_q^<} \sum_{\vec{j}\in C_\delta^<} \Phi_{\vec{i}}^{(\Psima)} \Phi_{\vec{j}}^{(\Psima)}.
$$
Thus,
\begin{align*}
&\E \left|N_{E_q}^{(\Pma)} N_{E_\delta}^{(\Pma)}-N_{E_q}^{(\Rma)} N_{E_\delta}^{(\Rma)}\right| \\
\leq & \E \sum_{\vec{i}\in C_q^<} \sum_{\vec{j}\in C_\delta^<} \left|\Phi_{\vec{i}}^{(\Pma)} \Phi_{\vec{j}}^{(\Pma)} -\Phi_{\vec{i}}^{(\Rma)} \Phi_{\vec{j}}^{(\Rma)}\right|1\left(\Ec(t)\right) + \binom{p}{1}\binom{p-1}{q}\binom{p}{1}\binom{p-1}{\delta}\P(\Ec^c(t))\\
\leq & \E \sum_{\vec{i}\in C_q^<} \sum_{\vec{j}\in C_\delta^<} \eta_{\vec{i},\vec{j}} + \frac{p^{q+\delta+2}}{\delta!q!} 2\exp(-c_1t^2) \numberthis \label{eqn:Nedgedelprnew}
\end{align*}
where the first inequality follows from $0\leq \Nedged^{(\Psima)} \leq \binom{p}{1}\binom{p-1}{\delta}$ for both $\Psima =\Rma$ and $\Psima = \Pma$, and the second inequality follows from Lemma \ref{lem:inddifbou2ndmoment}, Lemma \ref{lem:xietarelationship} and \eqref{eqn:proepst}.

	Choose $t=s_0\sqrt{\ln p}$ with  $s_0 = \sqrt{\left(\frac{9}{2}+2\delta\right)/c_1}$. Notice that $s_0 \geq \sqrt{\left(\frac{3}{2}+q+\delta + \frac{q}{\delta} \right)/c_1}$ since $q\in \{\delta,\delta+1\}$. Then
	$$
	2\exp(-c_1t^2)\leq 2\exp\left(-\left(\frac{3}{2}+q+\delta+\frac{q}{\delta}\right)\ln p\right)=\frac{2}{p^{\frac{3}{2}+q+\delta+\frac{q}{\delta}}}.
	$$
	Moreover, for any $c<\frac{1}{2\max\left \{\sqrt{\left(\frac{9}{2}+2\delta\right)/c_1},1\right\}\sqrt{2}C_1} $, 
	$$
	\left(\sqrt{\frac{\nn-1}{p}}+\sqrt{\frac{\delta \ln p}{p}}\right) \leq c
	$$
	implies
	\begin{equation}
	\sqrt{2}C_1\left(\sqrt{\frac{\nn-1}{p}}+s_0 \sqrt{\frac{\ln p}{p}}\right) \leq \frac{1}{2}, \label{eqn:CEatleadelprenew}
	\end{equation}
	which is \eqref{eqn:postempass} with $t= s_0 \sqrt{\ln p}$. Then apply Lemma \ref{lem:etadoublesumexpectation} with $t=s_0 \sqrt{\ln p}$ to \eqref{eqn:Nedgedelprnew},
	\begin{align*}
	  &\E \left|N_{E_q}^{(\Pma)} N_{E_\delta}^{(\Pma)}-N_{E_q}^{(\Rma)} N_{E_\delta}^{(\Rma)}\right|\\
	 \leq & C_{n,q,\delta,\gamma} \left(\theta_1(s_0\ln p)\right)^{n(2\delta+q)} \left(1+\mu_{n,q+\delta+2}(\Sigmama)\frac{\kappa-1}{p}\right) \left(\sqrt{\frac{1}{p}}+\frac{s_0\sqrt{\ln p}}{\sqrt{p}} + \frac{\tau}{p}\right) p^{1-\frac{q}{\delta}}\\
	 &\ \ \ \ \ \ \  + \frac{2}{\delta!q!\sqrt{p}} p^{1-\frac{q}{\delta}}\\
	 \leq & C_{n,q,\delta,\gamma}  \left(1+\mu_{n,q+\delta+2}(\Sigmama)\frac{\kappa-1}{p}\right) \left(s_0\frac{\sqrt{\ln p}}{\sqrt{p}} + \frac{\tau}{p}\right) p^{1-\frac{q}{\delta}} + \frac{2}{\delta!q!\sqrt{p}}p^{1-\frac{q}{\delta}} \\
	 \leq & C_{n,q,\delta,\gamma}  \left(1+\mu_{n,q+\delta+2}(\Sigmama)\frac{\kappa-1}{p}\right) \left(\frac{\sqrt{\ln p}}{\sqrt{p}} + \frac{\tau}{p}\right) p^{1-\frac{q}{\delta}}\\
	 \leq & C_{n,\delta,\gamma}  \left(1+\mu_{n,q+\delta+2}(\Sigmama)\frac{\kappa-1}{p}\right) \left(\frac{\sqrt{\ln p}}{\sqrt{p}} + \frac{\tau}{p}\right) p^{1-\frac{q}{\delta}},
	\end{align*}
	where the second inequality follows from $\theta_1\left(s_0 \sqrt{\ln p}\right)\leq 9+ 4(n-1)=4n+5$ by \eqref{eqn:CEatleadelprenew} and $\tau\leq p/2$; and the last step follows from $q\in \{\delta,\delta+1\}$.
\end{proof}

\begin{proof}[Proof of Proposition \ref{prop:allcloseL2} \ref{item:allcloseL2b} and \ref{item:allcloseL2c}]
(b)
It follows directly from Lemma \ref{lem:NEqNEdelta} with $q=\delta$.\\

\noindent
(c)
By taking square of each terms in Lemma \ref{lem:6quantitiesine}, 
$$
\left(N_{E_\delta}^{(\Pma)}\right)^2 - 2(\delta+1)N_{E_\delta}^{(\Pma)} N_{E_{\delta+1}}^{(\Pma)}   \leq \left(N_{\breve{V}_\delta}^{(\Pma)}\right)^2 \leq \left(N_{V_\delta}^{(\Pma)}\right)^2 \leq \left(N_{E_\delta}^{(\Pma)}\right)^2,
$$
which then implies for $\bar{N}_\delta \in \left\{ N_{\breve{V}_\delta}^{(\Pma)}, N_{V_\delta}^{(\Pma)} \right\} $
$$
\left|\bar{N}_\delta- \left(N_{E_\delta}^{(\Pma)}\right)^2\right|\leq 2(\delta+1)\left(N_{E_{\delta+1}}^{(\Pma)} N_{E_\delta}^{(\Pma)} - N_{E_{\delta+1}}^{(\Rma)} N_{E_\delta}^{(\Rma)} \right) + 2(\delta+1) N_{E_{\delta+1}}^{(\Rma)} N_{E_\delta}^{(\Rma)} .
$$

By Lemma \ref{lem:NEqNEdelta} with $q=\delta+1$,
$$
\E \left|N_{E_{\delta+1}}^{(\Pma)} N_{E_\delta}^{(\Pma)}-N_{E_{\delta+1}}^{(\Rma)} N_{E_\delta}^{(\Rma)}\right|\leq C_{n,\delta,\gamma}  \left(1+\mu_{n,2\delta+3}(\Sigmama)\frac{\kappa-1}{p}\right) \left(\frac{\sqrt{\ln p}}{\sqrt{p}} + \frac{\tau}{p}\right) p^{-\frac{1}{\delta}}.
$$
The proof is then completed by combining the preceding two displays, \eqref{eqn:NEdelta1NEdelta} and the fact that 
$$\left(\frac{\sqrt{\ln p}}{\sqrt{p}} + \frac{\tau}{p}\right)\leq 1.$$
\end{proof}

\subsection{Proofs in Section \ref{sec:limitingcompoundPoisson}}
\label{sec:proofoflem:randgeogra}

\subsubsection{Proof of Lemma \ref{lem:randgeogra} and Lemma \ref{prop:jumpsizedecay}}

\begin{figure}%
	\labellist
	\small \hair 2pt
	\pinlabel $\bm{0}$ at 589 552
	\pinlabel $\tilde{\uve}_{1}$ at 952 547
	\pinlabel $\sqrt{1-\left(\frac{\|\tilde{\uve}_{1}\|_2}{2}\right)^2}$ [l] at 1013 1315
	\pinlabel $\tilde{\uve}_{1}$ at 952 547
	\pinlabel 1  at 336 1005
	\endlabellist
	\centering
	\includegraphics[width=0.45\linewidth]{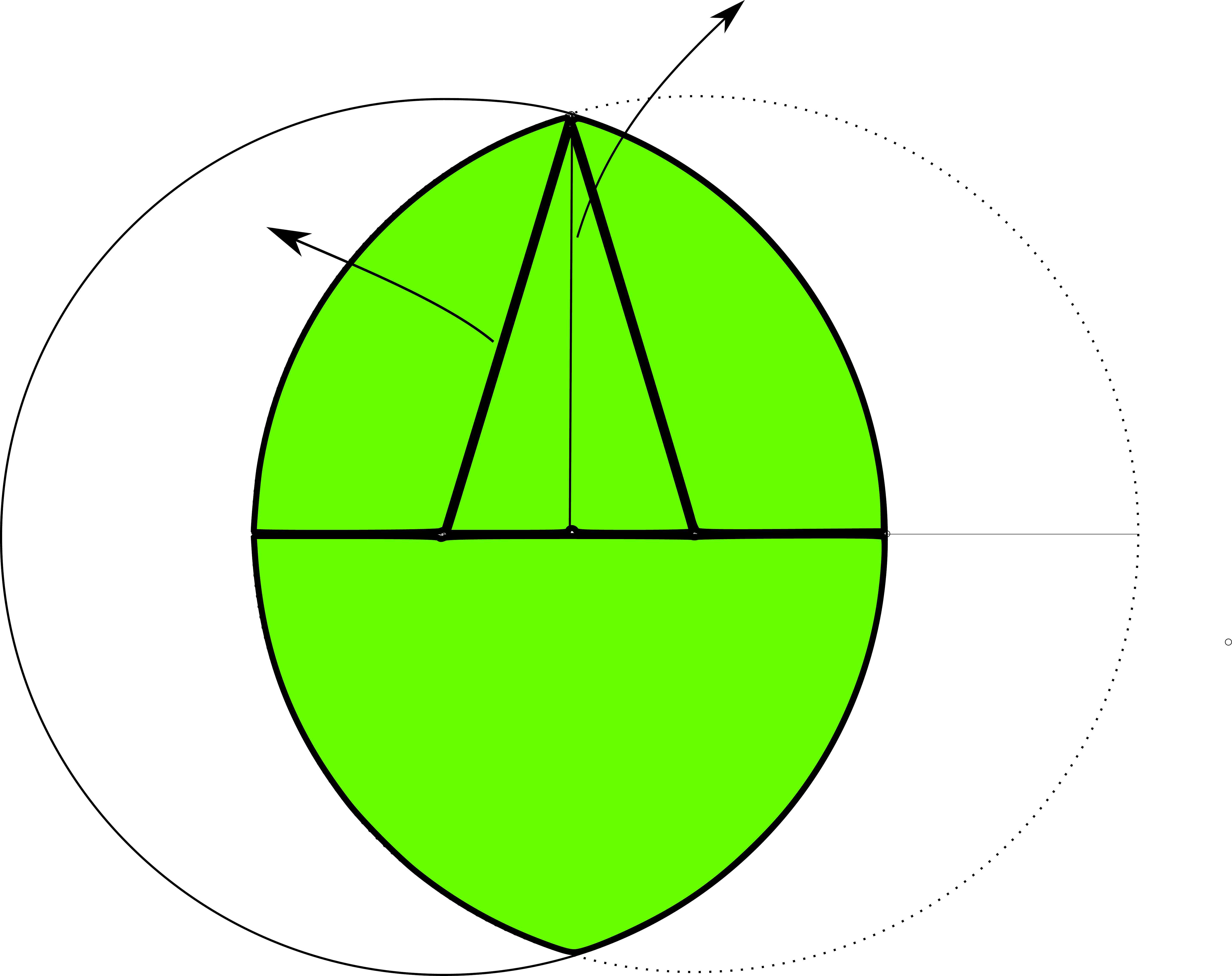}%
	\caption{The solid circle represents the unit Euclidean ball $B_2^n$ in $\R^n$ while the dash circle represents the unit ball centered at $\tilde{\uve}_1$. Their intersection is the green region, which is contained in the ball with center at $\tilde{\uve}_1/2$ and with radius $\sqrt{1-\left(\frac{\|\tilde{\uve}_{1}\|_2}{2}\right)^2}$.
	}.
	\label{fig:ballintersection}
\end{figure}
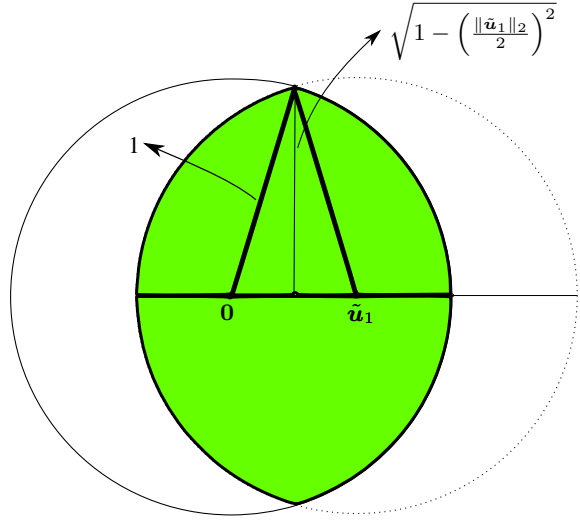

\begin{proof}[Proof of Lemma \ref{lem:randgeogra}]\hfil 
	 (a) 
	 Denote by $\deg(\cdot)$ the degree of a vertex in $\textbf{Ge}\left(\{\tilde{\uve}_i  \}_{i=1}^m, 1  ;  \mathscr{N} \right)$. Then by union bound,
	  \begin{align*}
	  \P\left(\nmd\left(\{\tilde{\uve}_i  \}_{i=1}^{\mgeominus}, 1  ; \Ngeo \right) \geq  1 \right) \leq & \mgeominus \P\left(\deg(\tilde{\uve}_1)=  \mgeominus-1 \right) \\
	  = & \mgeominus \E \left(\P(\deg(\tilde{\uve}_1) = \mgeominus-1 | \tilde{\uve}_1 )\right)\\ 
	  = & \mgeominus \E \left(\P\left( \left.\|\tilde{\uve}_1-\tilde{\uve}_2\|_2\leq 1,\cdots, \|\tilde{\uve}_1-\tilde{\uve}_{\mgeominus}\|_2\leq 1   \right| \tilde{\uve}_1 \right) \right) \\
	  =& \mgeominus \E \left( \P\left( \left.\|\tilde{\uve}_1-\tilde{\uve}_2\|_2\leq 1 \right| \tilde{\uve}_1 \right)^{\mgeominus-1} \right), \numberthis \label{eqn:tep1}
	  \end{align*}
	  where the last equality follows by conditional independence.
	
		As illustrated in Figure \ref{fig:ballintersection}, $\P\left( \|\tilde{\uve}_1-\tilde{\uve}_2\|_2\leq 1 | \tilde{\uve}_1 \right)$ is the ratio between Lebesgue measure of green region and $|B_2^{\Ngeo}|$. Moreover, the Lebesgue measure of the green region is less than $$
		\left(1-\left(\frac{\|\tilde{\uve}_{1}\|_2}{2}\right)^2\right)^{\frac{\Ngeo}{2}}|B_2^{\Ngeo}|.
		$$ 
		Then
	\begin{equation}
	\P\left( \|\tilde{\uve}_1-\tilde{\uve}_2\|_2\leq 1 | \tilde{\uve}_1 \right) \leq \left (1-\left(\frac{\|\tilde{\uve}_{1}\|_2}{2}\right)^2\right)^{\frac{\Ngeo}{2}} \quad a.s.  \label{eqn:u1u21cond}
	\end{equation}
	
	 By combining \eqref{eqn:tep1} and \eqref{eqn:u1u21cond},
	\begin{align*}
	\P\left(\nmd\left(\{\tilde{\uve}_i  \}_{i=1}^{\mgeominus}, 1  ; \Ngeo \right) \geq  1 \right) 
	\leq & \mgeominus \E \left(1-\left(\frac{\|\tilde{\uve}_{1}\|_2}{2}\right)^2\right)^{\frac{(\Ngeo)(\mgeominus-1)}{2}}\\
	= & \mgeominus (\Ngeo)\int_{0}^1 \left(1-\frac{r^2}{4}\right)^{\frac{(\Ngeo)(\mgeominus-1)}{2}} r^{n-3}dr
	\numberthis \label{eqn:inmcompletebetafun} \\
	=& \mgeominus (\Ngeo) 2^{n-3}B\left(\frac{1}{4};\frac{\Ngeo}{2},\frac{(\Ngeo)(\mgeominus-1)}{2}+1\right),
	\end{align*}
	where the first equality follows from expressing the integral in polar coordinates, and the last step follows from changing the variables $r=2\sqrt{y}$.\\

\noindent (b)
	Denote $f\left(r;\alpha,\beta \right)=\left(1-\frac{r^2}{4}\right)^\alpha r^\beta$. Then it is easy to verify that for any $\alpha,\beta>0$,
	\begin{equation}
	\max_{r\in [0,1]} f\left(r;\alpha,\beta\right) = \begin{cases} f\left(1;\alpha,\beta\right) =\left(\frac{3}{4}\right)^\alpha &  \text{ if } 3\beta \geq 2\alpha, \\  f\left(\sqrt{\frac{4\beta}{2\alpha+\beta}};\alpha,\beta\right) = \left(\frac{2\alpha}{2\alpha+\beta}\right)^\alpha \left(\frac{4\beta}{2\alpha+\beta}\right)^{\frac{\beta}{2}} & \text{ if } 3\beta \leq  2\alpha. \end{cases} \label{eqn:eleinetwocase}
	\end{equation}
	Moreover, $f(r;\alpha,\beta)$ is increasing on $[0,1]$ if $3\beta \geq 2\alpha$.
	
	Let $\alpha = \frac{(\Ngeo)(\mgeominus-1)}{2} $ and $\beta=n-3$. If $\mgeominus=2$, then for any $n\geq 4$, $3\beta \geq 2\alpha$ is satisfied. Then since $f(r;\alpha,\beta)$ is increasing on $[0,1]$,
	\begin{equation}
	\int_0^1 f\left(r;\alpha,\beta\right) dr \leq \sqrt{\frac{4}{5}}f\left(\sqrt{\frac{4}{5}};\alpha,\beta\right)+ \left(1-\sqrt{\frac{4}{5}}\right)f\left(1;\alpha,\beta\right). \label{eqn:m2N}
	\end{equation}
	
	If $\mgeominus=3$, then for any $n\geq 5$, $3\beta \geq 2\alpha$ is satisfied and hence \eqref{eqn:m2N} holds. For $n = 4$, $3\beta \leq 2\alpha$ is satisfied and by \eqref{eqn:eleinetwocase},
	\begin{equation}
	\int_0^1 f\left(r;\alpha,\beta\right) dr \leq f\left(\sqrt{\frac{4\beta}{2\alpha+\beta}};\alpha,\beta\right) = f\left(\sqrt{\frac{4}{5}};\alpha,\beta\right). \label{eqn:m3N2}
	\end{equation}
	
	If $\mgeominus \geq 4$, it is easy to see for any $n\geq 4$, $3\beta \leq 2\alpha$ holds. By \eqref{eqn:eleinetwocase}
	\begin{align*}
	&\int_0^1 f\left(r;\alpha,\beta\right) dr \\
	\leq & 
	f\left(\sqrt{\frac{4\beta}{2\alpha+\beta}};\alpha,\beta\right) \\
	= &  \left(\frac{\mgeominus-1}{\mgeominus}\right)^{\frac{(\Ngeo)(\mgeominus-1)}{2}} \left(\frac{4}{\mgeominus}\right)^{\frac{n-3}{2}} \left(\frac{\Ngeo}{\Ngeo-\frac{1}{\mgeominus}}\right)^{\frac{\mgeominus-1}{2}}\left( \left(\frac{\Ngeo}{\Ngeo-\frac{1}{\mgeominus}}\right)^{\mgeominus-1} \left(\frac{n-3}{\Ngeo-\frac{1}{m}}\right) \right)^{\frac{n-3}{2}}\\
	\leq & \exp\left(\frac{1}{4}\right) \left(\frac{\mgeominus-1}{\mgeominus}\right)^{\frac{(\Ngeo)(\mgeominus-1)}{2}} \left(\frac{4}{\mgeominus}\right)^{\frac{n-3}{2}}, \numberthis \label{eqn:inedelta4}	
	\end{align*}
	where the last step follows from $\left(\frac{\Ngeo}{\Ngeo-\frac{1}{\mgeominus}}\right)^{\mgeominus-1} \left(\frac{n-3}{\Ngeo-\frac{1}{\mgeominus}}\right)\leq 1$ and $\left(\frac{\Ngeo}{\Ngeo-\frac{1}{\mgeominus}}\right)^{\frac{\mgeominus-1}{2}}\leq \exp\left(\frac{1}{4}\right)$.

	Then \eqref{eqn:m2N}, \eqref{eqn:m3N2}, \eqref{eqn:inedelta4} and the fact that $f\left(\sqrt{\frac{4}{5}};\alpha,\beta\right) = \left(\frac{4}{5}\right)^{\frac{(\Ngeo)\mgeominus-1}{2}} $ yield the conclusion.
\end{proof} 


\begin{proof}[Proof of Lemma \ref{prop:jumpsizedecay}] 

(a)
Notice that
\begin{align*}
d_{\TV}\left(\bm{\zeta}_{n,\delta}, \Dirac(1) \right) = & \frac{1}{2} \sum_{\ell=1}^{\delta+1}|\bm{\zeta}_{n,\delta}(\ell) - \Dirac(1) \left(\ell\right) | \\ 
= & \sum_{\ell=2}^{\delta+1}\bm{\zeta}_{n,\delta}(\ell) \numberthis \label{eqn:tvatleast2}  \\
=& \frac{\sum_{\ell=2}^{\delta+1}\left(\alpha_\ell / \ell \right) }{\alpha_1+\sum_{\ell=2}^{\delta+1}\left(\alpha_\ell /\ell \right)} \\
\leq &  \frac{\sum_{\ell=2}^{\delta+1}\alpha_\ell   }{\alpha_1+\sum_{\ell=2}^{\delta+1}\alpha_\ell  } \\
= & \sum_{\ell=2}^{\delta+1}\alpha_\ell. \numberthis \label{eqn:tvuppbou}
\end{align*}

\noindent (b)
It follows from that
\begin{align*}
\left| \lambda_{\nn,\delta}(e_{n,\delta}) - \frac{1}{\delta !} \left( e_{n,\delta}\right)^\delta \right| = & \frac{1}{\delta !} \left(e_{n,\delta}\right)^\delta\left| \sum_{\ell=1}^{\delta+1}\left(\alpha_\ell/\ell\right) -1\right|\leq \frac{1}{\delta !} \left( e_{n,\delta}\right)^\delta \frac{3}{2} \sum_{\ell=2}^{\delta+1}\alpha_\ell. 
\end{align*}
(c) It follows directly by part \ref{item:jumpsizedecaya}, part \ref{item:jumpsizedecayb} and  Lemma \ref{eqn:cptv}.
\end{proof}

\subsubsection{Proof of Lemma \ref{lem:randgeograconduppbou}}

\begin{proof}[Proof of Lemma \ref{lem:randgeograconduppbou}]
(a)
Denote
\begin{align*}
I:=\P\left(\nmd\left(\{\uve'_i  \}_{i=1}^\mgeo, r  ; \Ngeo \right) \geq  2 |\deg(\uve'_{\mgeo})=\mgeominus\right). 
\end{align*}
Then by the union bound
\begin{align*}
    I = & \P\left(\bigcup_{i=1}^{\mgeominus}\{\deg(\uve'_i)=\mgeominus\} |\deg(\uve'_{\mgeo})=\mgeominus\right)\\
    \leq & \mgeominus \P\left(\deg(\uve'_1)=\mgeominus |\deg(\uve'_{\mgeo})=\mgeominus\right)\\
    =& \mgeominus \P\left(\deg(\uve'_1)=\mgeominus ,\deg(\uve'_{\mgeo})=\mgeominus\right)/ \P\left(\deg(\uve'_{\mgeo})=\mgeominus \right). \numberthis \label{eqn:Iunionbou}
\end{align*}
Notice that 
\begin{equation}
\P\left(\deg(\uve'_{\mgeo})=\mgeominus \right) = \E \prod_{i=1}^{\mgeominus}\P(\uve'_i\in \SC(r,\uve'_{\mgeo})|\uve'_{\mgeo})=(P_{\Ngeoplusplus}(r))^{\mgeominus}, \label{eqn:denominatordegm-1}
\end{equation}
where $\SC(r,\uve'_{\mgeo})$ and $P_{\Ngeoplusplus}(r)$ are defined in \eqref{eqn:SCdef} and the paragraph after \eqref{eqn:SCdef}. Moreover
\begin{align*}
    &\P\left(\deg(\uve'_1)=\mgeominus ,\deg(\uve'_{\mgeo})=\mgeominus\right) \\
    = &    \E \P\left(\deg(\uve'_1)=\mgeominus ,\deg(\uve'_{\mgeo})=\mgeominus|\uve'_{1} ,\uve'_{\mgeo}\right)\\
    =& \E  1(\|\uve'_1-\uve'_{\mgeo}\|_2\leq r) \prod_{i=2}^{\mgeominus}\P\left(\uve'_i\in \SC(r,\uve'_1)\cap \SC(r,\uve'_{\mgeo})|\uve'_{1} ,\uve'_{\mgeo} \right)\\
    \leq &\E  1(\|\uve'_1-\uve'_{\mgeo}\|_2\leq r) \left(P_{\Ngeoplusplus}(h(r,\|\uve'_1-\uve'_{\mgeo}\|_2))\right)^{\mgeominusminus} \numberthis \label{eqn:jointmaximaldeg}
\end{align*}
where the last inequality follows from Lemma \ref{lem:jointmaximaldgreeuppbou} with 
$$
h(r,d) =\sqrt{2-\frac{2-r^2}{\sqrt{1-(\frac{d}{2})^2}}}.
$$
Observing the random quantity in the expectation of \eqref{eqn:jointmaximaldeg} only depends the distance between $\|\uve'_1-\uve'_{\mgeo}\|_2$, replace $\uve'_{\mgeo}$ with $\vve_0=(1,0,\ldots,0)$ will not change its value. Then
\begin{equation}
\P\left(\deg(\uve'_1)=\mgeominus ,\deg(\uve'_{\mgeo})=\mgeominus\right) \leq \E  1(\|\uve'_{1}-\vve_0\|_2\leq r) \left(P_{\Ngeoplusplus}(h(r,\|\uve'_1-\vve_0\|_2))\right)^{\mgeominusminus}. \label{eqn:jointmaximaluppbou2}
\end{equation}
Use the following coordinate system for each $\uve'_1 = \left(u_{j1}:1\leq j \leq \Ngeoplus \right)^\top$ in the region $\SC(r, \vve_0)$:
\begin{equation*}
\begin{cases}
u_{11}=1-\frac{r^2 r_1^2}{2},\\
u_{21}=r_1r \sqrt{1-\frac{r^2 r_1^2}{4}} \cos(\theta_{2}),\\
\vdots \\
u_{j1}=r_1r \sqrt{1-\frac{r^2 r_1^2}{4}} \cos(\theta_{j})\prod\limits_{m=2}^{j-1}\sin(\theta_{m}),\\
\vdots \\
u_{(n-2)1}=r_1r \sqrt{1-\frac{r^2 r_1^2}{4}}\sin(\theta_{2})\cdots\sin(\theta_{n-3})\cos(\theta_{n-2}),\\
u_{(n-1)1}=r_1r \sqrt{1-\frac{r^2 r_1^2}{4}}\sin(\theta_{2})\cdots\sin(\theta_{n-3})\sin(\theta_{n-2}), 
\end{cases}
\end{equation*}
where
\begin{equation}
 r_1\in [0,1], \theta_{j}\in[0,\pi] \text{ for } 2\leq j\leq n-3 \text{ and } \theta_{n-2}\in[0,2\pi).  \label{eqn:intregionnonprod}
\end{equation}
Then the right hand side of \eqref{eqn:jointmaximaluppbou2} become 
\begin{align*}
    & \E  1(\|\uve'_{1}-\vve_0\|_2\leq r) \left(P_{\Ngeoplusplus}(h(r,\|\uve'_1-\vve_0\|_2))\right)^{\mgeominusminus} \\
    = & \frac{1}{\text{Area}(S^{\Ngeo})} \int_0^1 \left(P_{\Ngeoplusplus}(h(r,r_1r))\right)^{\mgeominusminus}r^{\Ngeo}r_1^{n-3}\left(1-\frac{r^2r_1^2}{4}\right)^{\frac{n-4}{2}}dr_1 \prod_{j=2}^{n-3}\int_0^{\pi}\sin^{\Ngeo-j}(\theta_j)d\theta_j\\
    = & \frac{1}{\int_0^{\pi}\sin^{n-3}(\theta)d\theta} \int_0^1 \left(P_{\Ngeoplusplus}(h(r,r_1r))\right)^{\mgeominusminus} r^{\Ngeo}r_1^{n-3}\left(1-\frac{r^2r_1^2}{4}\right)^{\frac{n-4}{2}}dr_1 \\
    = & \frac{r^{\Ngeo}}{B(\frac{\Ngeo}{2},\frac{1}{2})} \int_0^1 \left(P_{\Ngeoplusplus}(h(r,r_1r))\right)^{\mgeominusminus} r_1^{n-3}\left(1-\frac{r^2r_1^2}{4}\right)^{\frac{n-4}{2}}dr_1\\
    \leq & \frac{r^{\Ngeo}}{B(\frac{\Ngeo}{2},\frac{1}{2})} \int_0^1 \left(P_{\Ngeoplusplus}(h(r,r_1r))\right)^{\mgeominusminus}r_1^{n-3}dr_1. \numberthis \label{eqn:jointmaximaluppbou3}
\end{align*}
Plug \eqref{eqn:denominatordegm-1}, \eqref{eqn:jointmaximaluppbou2} and \eqref{eqn:jointmaximaluppbou3} into \eqref{eqn:Iunionbou} and we obtain
\begin{align*}
I  \leq & \delta \frac{r^{\Ngeo}}{B(\frac{\Ngeo}{2},\frac{1}{2})P_{\Ngeoplusplus}(r)} \int_0^1 \left(\frac{P_{\Ngeoplusplus}(h(r,r_1r))}{P_{\Ngeoplusplus}(r)}\right)^{\mgeominusminus}r_1^{n-3}dr_1 \\
= &\delta (\Ngeo) \frac{a_nr^{\Ngeo}}{P_{\Ngeoplusplus}(r)} \int_0^1 \left(\frac{P_{\Ngeoplusplus}(h(r,r_1r))}{P_{\Ngeoplusplus}(r)}\right)^{\mgeominusminus}r_1^{n-3}dr_1 \numberthis \label{eqn:Iuppbou3int}
\end{align*}
where the equality follows from $a_n =\frac{1}{ (\Ngeo)B(\frac{\Ngeo}{2},\frac{1}{2})}$. By Lemma \ref{pderivative} \ref{item:pderivativea},
\begin{equation}
    \frac{a_nr^{\Ngeo}}{P_{\Ngeoplusplus}(r)}\leq \frac{1}{\left(1-\frac{r^2}{4}\right)^{\frac{n-4}{2}}}. \label{eqn:cnrPnratiouppbou}
\end{equation}

Since when $0<r_1<1$, $0<h(r,d)/r<1$, by Lemma \ref{pderivative} \ref{item:pderivativee},
\begin{align*}
\frac{P_\nn(h(r,r_1 r))}{P_\nn(r)} \leq & \left(\frac{h(r,r_1r)}{r}\right)^{\Ngeo}\left( \frac{ 1-\frac{h^2(r,r_1r)}{4}} { 1-\frac{r^2}{4}}\right )^{\frac{n-4}{2}}\\
\leq & \left(1 -\left(\frac{r_1}{2}\right)^2 \right)^{\frac{n-2}{2}} \frac{1}{\sqrt{1-\left(\frac{r_1r}{2}\right)^2}}\left( \frac{ 1-\frac{h^2(r,r_1r)}{4}} { (1-\frac{r^2}{4})\sqrt{1-\left(\frac{r_1r}{2}\right)^2}}\right )^{\frac{n-4}{2}}\numberthis \label{eqn:Pnratiouppboutemp}
\end{align*}
where the second inequality follows from
$$
\left(\frac{h(r,r_1r)}{r}\right)^2 =  \frac{1}{\sqrt{1-\left(\frac{r_1r}{2}\right)^2}}\left( \frac{-2(\frac{r_1}{2})^2}{1+\sqrt{1-(\frac{r_1r}{2})^2}}+1 \right) \leq \frac{1}{\sqrt{1-\left(\frac{r_1r}{2}\right)^2}}\left( -\left(\frac{r_1}{2}\right)^2+1 \right).
$$
Since $h^2(r,r_1r)$ is decreasing function of $r_1\in [0,1]$,  \eqref{eqn:Pnratiouppboutemp} become
\begin{align*}
\frac{P_\nn(h(r,r_1 r))}{P_\nn(r)}   \leq & \left(1 -\left(\frac{r_1}{2}\right)^2 \right)^{\frac{n-2}{2}} \frac{1}{\sqrt{1-\left(\frac{r}{2}\right)^2}}\left( \frac{ 1-\frac{h^2(r,r)}{4}} { (1-\frac{r^2}{4})\sqrt{1-\left(\frac{r}{2}\right)^2}}\right )^{\frac{n-4}{2}}\\
\leq & \left(1 -\left(\frac{r_1}{2}\right)^2 \right)^{\frac{n-2}{2}} \frac{1}{\sqrt{1-\left(\frac{r}{2}\right)^2}}\left( \frac{ 1} { (1-\frac{r^2}{4})}\right )^{\frac{n-4}{2}}, \numberthis \label{eqn:Pnratiouppbou}
\end{align*}
where the second inequality follows from
$$
1-\frac{h^2(r,r)}{4} \leq \sqrt{1-\frac{r^2}{4}}.
$$
Plugging \eqref{eqn:cnrPnratiouppbou} and \eqref{eqn:Pnratiouppbou} into \eqref{eqn:Iuppbou3int},
\begin{align*}
I\leq &  \delta (\Ngeo) \frac{1}{\left(1-\frac{r^2}{4}\right)^{\frac{n+\delta-5}{2}}}\left( \frac{ 1} { (1-\frac{r^2}{4})}\right )^{\frac{(n-4)(\mgeominusminus)}{2}} \int_0^1 \left(1 -\left(\frac{r_1}{2}\right)^2 \right)^{\frac{(n-2){(\mgeominusminus)}}{2}} r_1^{n-3}dr_1\\
=& \bar{h}\left(\frac{1}{\sqrt{1-r^2/4}},n,\delta\right) \delta (\Ngeo)  \int_0^1 \left(1 -\left(\frac{r_1}{2}\right)^2 \right)^{\frac{(n-2){(\mgeominusminus)}}{2}} r_1^{n-3}dr_1.
\end{align*}

\noindent
(b)
Since  $\frac{1}{\sqrt{1-r^2/4}}$ is decreasing and $\bar{h}(x,n,\delta)$ as a function of $x$ is increasing,
\begin{align*}
\bar{h}\left(\frac{1}{\sqrt{1-r^2/4}},n,\delta\right)
\leq &
\begin{cases}
\bar{h}\left(\left(\frac{5}{4}\right)^{\frac{1}{4}},n,\delta\right), & \delta = 2,3\\
\bar{h}\left(\left(\frac{\delta}{\delta-1}\right)^{\frac{1}{4}},n,\delta\right), & \delta \geq 4
\end{cases} \\
=& 
\begin{cases}
\left(\sqrt{\frac{5}{4}}\right)^{\frac{n+\delta-5}{2}}\left(\sqrt{\frac{5}{4}}\right)^{\frac{(n-2)(\delta-1)}{2}}, & \delta = 2,3\\
\left(\sqrt{\frac{\delta}{\delta-1}}\right)^{\frac{n+\delta-5}{2}}\left(\sqrt{\frac{\delta}{\delta-1}}\right)^{\frac{(n-2)(\delta-1)}{2}}, & \delta \geq 4
\end{cases}. \numberthis \label{eqn:hbaruppbou}
\end{align*}
Then the proof is complete by combining part \ref{item:randgeograconduppboua}, Lemma \ref{lem:randgeogra} \ref{item:randgeograb} and \eqref{eqn:hbaruppbou}.\\

\noindent (c)
Similar to \eqref{eqn:tvuppbou}, we have 
$$
d_{\TV}(\bm{\zeta}_{n,\delta,\rho},\Dirac(1))\leq \sum_{\ell=2}^{\delta+1} \alpha(\ell,r_{\rho})
= \P\left(\nmd\left(\{\uve'_i  \}_{i=1}^\mgeo, r_\rho  ; \Ngeo \right) \geq  2 |\deg(\uve'_{\mgeo})=\mgeominus\right).
$$
where the equality follows from \eqref{eqn:alphalrhorgg}. Then the conclusion follows from part \ref{item:randgeograconduppboua} and part \ref{item:randgeograconduppboub} since $r_\rho$ satisfies the condition there.

\end{proof}

\begin{lem} \label{lem:jointmaximaldgreeuppbou}
Let $n\geq 3$ and $0<r<\sqrt{2}$. If $\wwve_1$ and $\wwve_2$ are two points in $S^{\Ngeo}$ with $\|\wwve_1-\wwve_2\|_2 = d$ satisfying $2-2\sqrt{1-(d/2)^2}<r^2$, then 
    $$
    \P\left(\uve'_1\in \SC(r,\wwve_1)\cap \SC(r,\wwve_{2}) \right) \leq  P_{\Ngeoplusplus}(h(r,d))
    $$
    where $\uve'_1$ has distribution $\unif(S^{\Ngeo})$ and 
    $$
    h(r,d) =\sqrt{2-\frac{2-r^2}{\sqrt{1-(\frac{d}{2})^2}}}.
    $$
\end{lem}
\begin{figure}%
	\labellist
	\small \hair 2pt
	\pinlabel $\bm{0}$ at 20 600 
	\pinlabel $\wwve_1$ at 727 12 
	\pinlabel $\wwve_2$ at 727 1252 
	\pinlabel $\wwve_3$ at 710 600 
	\pinlabel $\wwve_4$ at 780 680 
	\pinlabel $\wwve_5$ at 900 480 
	\pinlabel $r$ at 810 250 
	\pinlabel $r$ at 810 865 
	\pinlabel $\theta$ at 280 640 
	\pinlabel $1$ at 373 266 
	\pinlabel $1$ at 373 1070 
	\endlabellist
	\centering
	\includegraphics[width=7cm, height=4cm]{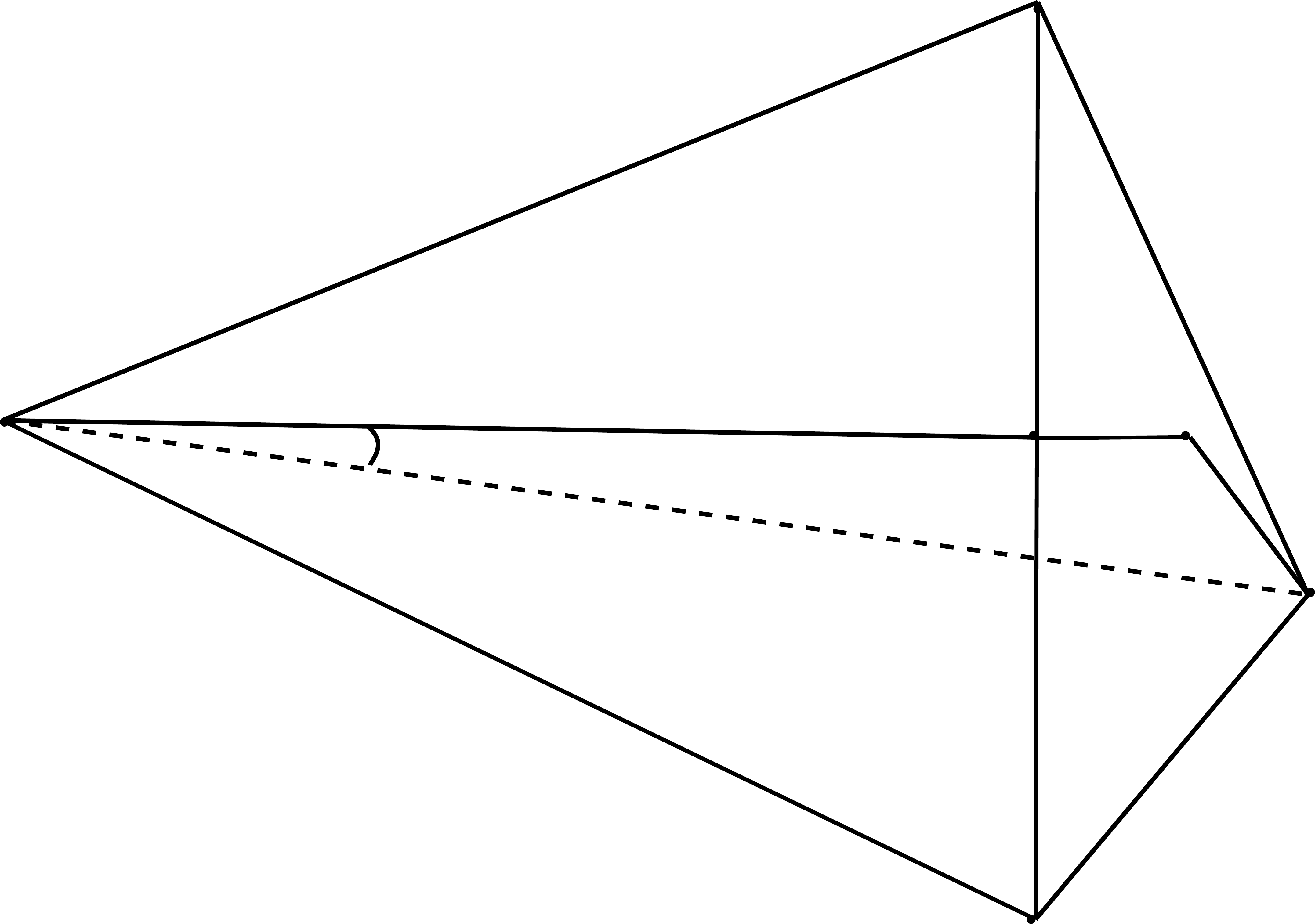}%
	\caption{$\bm{0}$ is the origin in $\R^{\Ngeo}$ and $\wwve_1,\wwve_2,\wwve_4,\wwve_5$ are on $S^{\Ngeo}$. $\wwve_3$ is the midpoint of $z_1$ and $z_2$, while $z_4$ is the midpoint of the shortest arc on $S^{\Ngeo}$ connecting $\wwve_1$ and $\wwve_2$. $\wwve_5$ is one of the two intersection points of the boundary $\SC(r,\wwve_1)$ and the boundary of $\SC(r,\wwve_2)$. The angle between line segment $\bm{0}\wwve_4$ and $\bm{0}\wwve_5$ is $\theta$.
	}
	\label{fig:capintersection}
\end{figure}
\begin{proof} The proof is based on Figure \ref{fig:capintersection} and we use $|\cdot|$ to represent the length of a line segment in this proof. In the right triangle $\bm{0}\wwve_3\wwve_1$, the line segment $\bm{0}\wwve_3$ has length $|\bm{0}\wwve_3|=\sqrt{1-(d/2)^2}$. In the right triangle $\wwve_1\wwve_3\wwve_5$,  $|\wwve_3\wwve_5| =\sqrt{r^2-(d/2)^2}$. In the triangle $\bm{0}\wwve_3\wwve_5$, by the law of Cosines, 
$$
\cos(\theta) = \frac{2-r^2}{2\sqrt{1-(\frac{d}{2})^2}}.
$$
Then in the isosceles triangle, the line segment $\wwve_4\wwve_5$ has length 
$$
|\wwve_4\wwve_5|=2\sin(\theta/2) = \sqrt{2(1-\cos(\theta))} = \sqrt{2-\frac{2-r^2}{\sqrt{1-(\frac{d}{2})^2}}} = h(r,d).
$$
It is easy to deduce that  $|\wwve_1\wwve_4|=\sqrt{2-2\sqrt{1-(d/2)^2}}$. The condition $2-2\sqrt{1-(d/2)^2}<r^2$ entails that $\SC(r,\wwve_1)\cap \SC(r,\wwve_{2})\not= \emptyset$ and that $|\wwve_1\wwve_4|<|\wwve_4\wwve_5|=h(r,d)$. In this case $\SC(r,\wwve_1)\cap \SC(r,\wwve_{2}) \subset \SC(h(r,d),\wwve_4) $. Thus
$$
\P\left(\uve'_1\in \SC(r,\wwve_1)\cap \SC(r,\wwve_{2}) \right) \leq \P\left(\uve'_1\in \SC(h(r,d),\wwve_4)\right) = P_{\Ngeoplusplus}(h(r,d)).
$$
\end{proof}

\subsubsection{Proof of Lemma \ref{lem:alphaldelta2}}

\begin{proof}[Proof of Lemma \ref{lem:alphaldelta2}]
When $\delta=2$, $\alpha_2=0$ since either both vertices have degree $1$ or none. Moreover, 
\begin{align*}
\alpha_3 
= & 
\P(\|\tilde{\uve}_1-\tilde{\uve}_2\|_2\leq 1)\\
= & 
\E \P(\|\tilde{\uve}_1-\tilde{\uve}_2\|_2\leq 1|\tilde{\uve}_1)\\
\overset{(*)}{=} &  
\E \frac{1}{\text{Vol}(B^{n-2})}\times2\times \frac{\pi^{(n-3)/2}}{\Gamma(\frac{n-3}{2}+1)}\int_{0}^{\arccos(\frac{\|\tilde{\uve}_1\|_2}{2})} \sin^{n-2}(\theta)d\theta\\
\overset{(**)}{=}  &  
\frac{1}{\text{Vol}(B^{n-2})}\times2\times \frac{\pi^{(n-3)/2}}{\Gamma(\frac{n-3}{2}+1)}\frac{\text{Area}(S^{n-3})}{\text{Vol}(B^{n-2})} \int_0^1 r^{n-3} \int_{0}^{\arccos(\frac{r}{2})} \sin^{n-2}(\theta)d\theta dr\\
=& \frac{2(n-2)}{B(\frac{n-1}{2},\frac{1}{2})}  \int_0^1 r^{n-3} \int_{0}^{\arccos(\frac{r}{2})} \sin^{n-2}(\theta)d\theta dr \numberthis \label{eqn:alpha3delta2formula1}
\end{align*}
where step $(*)$ follows from the Subsection ``Volume of a hyperspherical cap'' from \cite{li2011concise} and $\text{Vol}(B^{n-2})$ is the volume of $B^{n-2}$, step $(**)$ follows by observing the random quantity only depends on $\tilde{\uve}_1$ through its Euclidean norm, and in the last step $B(\cdot,\cdot)$ is the Beta function. By Fubini's Theorem 
\begin{align*}
&\int_0^1 r^{n-3} \int_{0}^{\arccos(\frac{r}{2})} \sin^{n-2}(\theta)d\theta dr \\
= & \int_{0}^{\frac{\pi}{3}} \int_0^1 r^{n-3}  \sin^{n-2}(\theta)dr d\theta  + \int_{\frac{\pi}{3}}^{\frac{\pi}{2}} \int_0^{2\cos(\theta)} r^{n-3}  \sin^{n-2}(\theta)dr d\theta \\
=& \frac{3}{2(n-2)} \int_{0}^{\frac{\pi}{3}} \sin^{n-2}(\theta) d\theta
\end{align*}
Plugging the preceding formula into \eqref{eqn:alpha3delta2formula1}, $\alpha_3 = \frac{3}{2} I_{\frac{3}{4}}(\frac{n-1}{2},\frac{1}{2})$, where $I_x(a,b)$ is the regularized incomplete Beta function. $\alpha_1=1-\alpha_3$ follows from $\alpha_2=0$.
\end{proof}

\subsection{Auxiliary lemmas}
\label{sec:auxlem}

\begin{lem}
\label{pderivative}
Let $P_\nn(r)$  be defined as in Section \ref{sec:closenessnumedge}. Suppose $\nn\geq 4$.
\begin{enumerate}[label=(\alph*)]
\item \label{item:pderivativea}
Recall $a_n =\frac{b_\nn}{2(\nn-2)}=\frac{\Gamma((\nn-1)/2)}{(n-2)\sqrt{\pi}\Gamma((\nn-2)/2)} \leq 1$. Then
$$
a_nr^{\nn-2}\left(1-\frac{\min\{r^2,4\}}{4}\right)^{\frac{n-4}{2}}\leq P_\nn(r) \leq a_nr^{\nn-2}
$$
 and
	\begin{equation}
	-\frac{n-4}{8} a_nr^{\nn-2}  \min\{r^2,4\}\leq P_\nn(r) - a_nr^{\nn-2}\leq 0. \label{eqn:Pcnrdif}
	\end{equation}
 \item \label{item:pderivativeb}
 $\lim_{r\to 0} P_\nn(r) / \left(a_nr^{\nn-2} \right) =1$.

\item \label{item:pderivativec}
Let $0\leq \beta < 1< \alpha$ and $0<r\leq 2$. Then 
$$
P_\nn(\alpha r)-P_\nn(\beta r)\leq (\nn-2) P_\nn(r) \alpha^{\nn-3} (\alpha - \beta).
$$

\item \label{item:pderivatived}
Consider $\alpha>1$ and $r>0$. Then
$$
P_\nn(\alpha r) \leq \alpha^{n-2}P_\nn(r).
$$

\item \label{item:pderivativee}
Consider $0<\beta<1$ and $0<r<2$. Then
$$
P_\nn(\beta r) \leq  \beta^{n-2}\left( \frac{ 1-\frac{\beta^2 r^2}{4}} { 1-\frac{r^2}{4}}\right )^{\frac{n-4}{2}}P_\nn(r).
$$
\end{enumerate}
\end{lem}

\begin{proof}
(a)
It is easy to verify
\begin{equation}
P_\nn'(x)= \begin{cases} \frac{b_\nn}{2}x^{\nn-3}(1-\frac{x^2}{4})^{\frac{\nn-4}{2}} &  x<{2}, \\
0 & x\geq 2. \end{cases}  \label{eqn:pnrderivative}
\end{equation}
Consider $r>0$. Then
\begin{align}
\frac{P_\nn(r)}{a_nr^{\nn-2}} =  \frac{P'_\nn(\xi)}{(n-2)a_n\xi^{\nn-3}}
 = \left(1-\frac{\xi^2}{4}\right)^{\frac{n-4}{2}}    \label{eqn:pderivativea1}
\end{align}
where in the first equality $\xi \in (0,\min\{r,2\}) $ due to the Cauchy Mean Value Theorem and $P_\nn(r)\not = 0$, and the second equality follows from \eqref{eqn:pnrderivative}. Equation \eqref{eqn:pderivativea1} directly implies
$$
\left(1-\frac{(\min\{r,2\})^2}{4}\right)^{\frac{n-4}{2}} \leq \frac{P_\nn(r)}{a_nr^{\nn-2}} \leq 1.
$$

\noindent (b) 
It follows directly by taking limit $r\to 0^+$ in \eqref{eqn:pderivativea1}.\\
	
\noindent (c)
Since $0\leq \beta < 1< \alpha$ and $0<r\leq 2$, $P_\nn(\alpha r)- P_\nn(\beta r)>0$ and $P_\nn(r)>0$. Then
\begin{align*}
\frac{P_\nn(\alpha r)- P_\nn(\beta r) }{ P_\nn(r) } &= \frac{\left(P_\nn(\alpha r)- P_\nn(\beta r)\right)-\left(P_\nn(\alpha\cdot 0)- P_\nn(\beta \cdot 0)\right) }{ P_\nn(r) - P_\nn(0)}\\
& = \frac{ \left. \frac{d}{dr}\left(P_\nn(\alpha r)- P_\nn(\beta r)\right) \right|_{r=\xi} }{\left. \frac{d}{dr}P_\nn( r) \right|_{r=\xi}} \\
&= \frac{\alpha^{n-2}\left(  1-\frac{\alpha^2\xi^2}{4}\right )^{\frac{n-4}{2}}-\beta^{n-2}\left(  1-\frac{\beta^2\xi^2}{4}\right )^{\frac{n-4}{2}}}{\left(  1-\frac{\xi^2}{4}\right )^{\frac{n-4}{2}}}\\
&  \leq \alpha^{n-2}- \beta^{n-2} \numberthis \label{eqn:Pnderivativetemp}  \\
& \leq (n-2)\alpha^{n-3}(\alpha-\beta),
\end{align*}
where the second equality follows from the Cauchy Mean Value Theorem with $\xi\in (0,r)$, the third equality follows from \eqref{eqn:pnrderivative} together with the fact that the numerator has to be positive, which imply $\alpha\xi<2$, the first inequality follows from $0\leq\beta<1<\alpha$, and the last inequality follows from mean value theorem.\\

\noindent (d)
When $r\geq 2$, 
$
P_\nn(\alpha r) = P_\nn(r) =1 
$
and the conclusion holds trivially. The case $0<r<2$  follows from \eqref{eqn:Pnderivativetemp} with $\beta=0$.\\

\noindent
(e)
Consider $0<\beta<1$ and $0<r<2$. Then
$$
\frac{P_\nn(\beta r)}{P_\nn(r)} = \frac{\beta^{n-2}\left(  1-\frac{\beta^2\xi^2}{4}\right )^{\frac{n-4}{2}}}{\left(  1-\frac{\xi^2}{4}\right )^{\frac{n-4}{2}}}\leq \frac{\beta^{n-2}\left(  1-\frac{\beta^2 r^2}{4}\right )^{\frac{n-4}{2}}}{\left(  1-\frac{r^2}{4}\right )^{\frac{n-4}{2}}},
$$
where the equality follows from Cauchy Mean Value Theorem with $\xi\in (0,r)$. 
\end{proof}

\begin{lem} \label{prop:was1dis}
	Consider $Z_1$ and $Z_2$ be two discrete random variable support on $[\delta]$. Then
	$$
	d_{\text{W}}\left(\mathscr{L}(Z_1), \mathscr{L}(Z_2)\right) \leq \frac{\delta-1}{2} \sum_{\ell=1}^{\delta}\left|\P(Z_1=\ell)-\P(Z_2=\ell)\right|.
	$$
\end{lem}
\begin{proof}
By Remark 2.19 (iii) of Section 2.2 in \cite{villani2003topics}, 
\begin{align*}
d_{\text{W}}\left(\mathscr{L}(Z_1), \mathscr{L}(Z_2)\right) 
=  & 
\sum_{i=1}^{\delta-1}|\P(Z_1\leq i)-\P(Z_2\leq i)| 
\leq 
\sum_{i=1}^{\delta-1} \sum_{j=1}^{i} |\P(Z_1 = j)-\P(Z_2= j)|.
\end{align*}
On the other hand, from the above equality,
\begin{align*}
d_{\text{W}}\left(\mathscr{L}(Z_1), \mathscr{L}(Z_2)\right) 
= &
\sum_{i=1}^{\delta-1}|\P(Z_1\geq i+1)-\P(Z_2\geq i+1)| \\ 
\leq &
\sum_{i=1}^{\delta-1} \sum_{j=i+1}^{\delta} |\P(Z_1 = j)-\P(Z_2= j)|.
\end{align*}
Averaging the above two inequalities yields the desired conclusion.
\end{proof}

\begin{lem}
	\label{prop:tvtomean}
	For any integer-valued random variable $Z_1$ and $Z_2$, 
	$$
	d_{\TV}\left(\mathscr{L}(Z_1), \mathscr{L}(Z_2)\right) \leq \P( Z_1 \not = Z_2) \leq \E \left |Z_1 - Z_2 \right |.
	$$
\end{lem}	
\begin{proof}
	\begin{align*}
	d_{\TV}\left(\mathscr{L}(Z_1), \mathscr{L}(Z_2)\right) & = \max_{A \text{ Borel measurable }} \left| \P(Z_1 \in A) - \P(Z_2 \in A) \right| \\
	& = \max_{A \text{ Borel measurable }} \left| \P(Z_1 \in A, Z_1 \not = Z_2) - \P(Z_2 \in A, Z_1\not = Z_2) \right| \\
	& \leq  \P( Z_1 \not = Z_2) \\
	& = \P( |Z_1 - Z_2| \geq 1)\\
	& \leq \E \left |Z_1 - Z_2 \right |.
	\end{align*}
\end{proof}

\begin{lem} \label{eqn:cptv}
    Consider two compound Poisson distributions $\CP(\lambda_1,\bm{\zeta}_1)$ and $\CP(\lambda_2,\bm{\zeta}_2)$. Then
    \begin{align*}
&    d_{\TV}\left(\CP(\lambda_1,\bm{\zeta}_1),\CP(\lambda_2,\bm{\zeta}_2)\right) \\
    \leq & 
    \min\{\lambda_1,\lambda_2\} d_{\TV}(\bm{\zeta}_1,\bm{\zeta}_2) + d_{\TV}(\Pois(\lambda_1),\Pois(\lambda_2))\\
    \leq  & \min\{\lambda_1,\lambda_2\} d_{\TV}(\bm{\zeta}_1,\bm{\zeta}_2) + \min\left\{|\lambda_1-\lambda_2|,\sqrt{\frac{2}{e}}|\sqrt{\lambda_1}-\sqrt{\lambda_2}|\right\}.
    \end{align*}
\end{lem}
\begin{proof}
By triangular inequality,
    \begin{equation}
    d_{\TV}\left(\CP(\lambda_1,\bm{\zeta}_1),\CP(\lambda_2,\bm{\zeta}_2)\right)\leq d_{\TV}\left(\CP(\lambda_1,\bm{\zeta}_1),\CP(\lambda_1,\bm{\zeta}_2)\right)+d_{\TV}\left(\CP(\lambda_1,\bm{\zeta}_2),\CP(\lambda_2,\bm{\zeta}_2)\right). \label{eqn:totalvariation} 
    \end{equation}
We will bound the two terms in the upper bound separately.

\noindent
\textbf{Step 1}: $d_{\TV}\left(\CP(\lambda_1,\bm{\zeta}_1),\CP(\lambda_1,\bm{\zeta}_2)\right)$\\
Consider $Z_1=\sum_{i=1}^N Y_i$ and $Z_2=\sum_{i=1}^N Y'_i$, where 
$N\sim \Pois(\lambda_1)$, 
$\{Y_i\}\overset{\text{i.i.d.}}{\sim} \bm{\zeta}_1$, 
$\{Y'_i\}\overset{\text{i.i.d.}}{\sim} \bm{\zeta}_2$, and $N$ is independent of $\{Y_i\},\{Y'_i\}$. Then
\begin{align*}
    d_{\TV}\left(\CP(\lambda_1,\bm{\zeta}_1),\CP(\lambda_1,\bm{\zeta}_2)\right) &= d_{\TV}\left(\mathscr{L}(Z_1), \mathscr{L}(Z_2) \right) \\
    & \leq 1-\P(Z_1 = Z_2) \\
    & = 1- \E \P(Z_1 = Z_2|N) \\
    & \leq 1- \E  \left(\P(Y_1 = Y'_1)\right)^N,
\end{align*}
where the first inequality follows from Lemma \ref{prop:tvtomean}. Since the above inequality holds for any coupling $(Y_1,Y'_1)$ of $\zeta_1$ and $\zeta_2$, by Proposition 4.7 in \cite{levin2017markov} taking the infimum of all the coupling then yields 
$$
d_{\TV}\left(\CP(\lambda_1,\bm{\zeta}_1),  \CP(\lambda_1,\bm{\zeta}_2)\right)\leq 1- \E (1-d_{\TV}(\bm{\zeta}_1,\bm{\zeta}_2))^N \leq   \lambda_1 d_{\TV}(\bm{\zeta}_1,\bm{\zeta}_2).
$$

\noindent
\textbf{Step 2}: $d_{\TV}\left(\CP(\lambda_1,\bm{\zeta}_2),\CP(\lambda_2,\bm{\zeta}_2)\right)$\\
Consider $Z_1=\sum_{i=1}^{N_1} Y'_i$ and $Z_2=\sum_{i=1}^N Y'_i$, where 
$N_1\sim \Pois(\lambda_1)$, $N_2\sim \Pois(\lambda_2)$, 
$\{Y'_i\}\overset{\text{i.i.d.}}{\sim} \bm{\zeta}_2$, 
and $\{Y'_i\}$ is independent of $N_1,N_2$. Then
\begin{align*}
    d_{\TV}\left(\CP(\lambda_1,\bm{\zeta}_2),\CP(\lambda_2,\bm{\zeta}_2)\right)
    &= d_{\TV}\left(\mathscr{L}(Z_1), \mathscr{L}(Z_2) \right) \\
    & \leq \P(Z_1 \neq Z_2) \\
    & = \P(N_1 \neq N_2),
\end{align*}
where the first inequality follows from Lemma \ref{prop:tvtomean}. Since the above inequality holds for any coupling $(N_1,N_2)$ of $\Pois(\lambda_1)$ and $\Pois(\lambda_2)$, by Proposition 4.7 in \cite{levin2017markov} taking the infimum of all the coupling then yields $d_{\TV}\left(\CP(\lambda_1,\bm{\zeta}_2),\CP(\lambda_2,\bm{\zeta}_2)\right) \leq d_{\TV}\left(\Pois(\lambda_1), \Pois(\lambda_2) \right) $. 

Plugging step 1 and step 2 into \eqref{eqn:totalvariation},
$$
d_{\TV}\left(\CP(\lambda_1,\bm{\zeta}_1),\CP(\lambda_2,\bm{\zeta}_2)\right) \leq \lambda_1 d_{\TV}(\bm{\zeta}_1,\bm{\zeta}_2) + d_{\TV}\left(\Pois(\lambda_1), \Pois(\lambda_2) \right).
$$
By symmetry property, the first term in the above upper bound can be replace by $\min\{\lambda_1,\lambda_2\} d_{\TV}(\bm{\zeta}_1,\bm{\zeta}_2)$. The proof is then completed by applying equation (2.2) of \cite{adell2006exact} to bound $d_{\TV}\left(\Pois(\lambda_1), \Pois(\lambda_2) \right)$.
\end{proof}

\begin{lem}  \label{basine}
\begin{enumerate}[label=(\alph*)]
  \item \label{item:basineh}
   Let $p,p',m$ be positive integers such that $p\geq p'$. Then $$\prod\limits_{i=0}^m (p-i) - \prod\limits_{i=0}^m (p'-i) \leq (m+1)\left(\prod\limits_{i=0}^{m-1} (p-i)\right)(p-p').$$
  \item \label{item:basinei}
  Let $p,\delta,\kappa$ be positive integers such that $ \delta\leq p-1$. Then
  $$
  \prod\limits_{\ell=1}^\delta (p-\ell) -\prod\limits_{\ell=1}^\delta (p-\ell\kappa) \leq \frac{\delta(\delta+1)}{2} \left(\kappa-1\right) \prod_{\ell=1}^{\delta-1}(p-\ell).
  $$
  \item \label{item:basinej}
  $\left( \frac{1+x}{1-x}\right)^2$ is increasing function on $[0,\frac{1}{2}]$ and $\left( \frac{1+x}{1-x}\right)^2\leq 1+16x$ for $0\leq x\leq \frac{1}{2}$. Then

  \item \label{item:basinek}
  $1+x-\frac{1}{1+x}\leq 2x$ for any $x\geq 0$.
\end{enumerate}
\end{lem}
\begin{proof}
	 \ref{item:basinej} and \ref{item:basinek} are simple quadratic inequalities and hence their proofs are omitted. \\
	\noindent
	(a)
	Let $f(x)=\prod\limits_{i=0}^m (x-i)$. When $p'\geq m$, $f'(x)\leq (m+1) \prod_{i=0}^{m-1}(p-i) $ and the conclusion then follows by the mean value theorem. When $p'\leq m-1$, 
	$$
	f(p)-f(p')\leq f(p) \leq (p-p')\prod_{i=0}^{m-1}(p-i).
	$$
	
	\noindent
	(b)
	Let $f(x)=\prod\limits_{\ell=1}^\delta (p-\ell x)$. When $p<\delta\kappa$, 
	$$
	f(1)-f(\kappa)\leq f(1)\leq (\delta\kappa - \delta) \prod_{\ell=1}^{\delta-1}(p-\ell)\leq \frac{\delta(\delta+1)}{2} \left(\kappa-1\right) \prod_{\ell=1}^{\delta-1}(p-\ell).
	$$
	When $p\geq \delta\kappa$, $f'(x)\geq -\frac{\delta(\delta+1)}{2}\prod_{\ell=1}^{\delta-1}(p-\ell)$ for $x\in [1,\kappa]$. Then the conclusion follows by the mean value theorem.
\end{proof}

\begin{lem}[Perturbation Theory] \label{perthe}
	Consider $\Dma\in \Sbb^n$ and $\Ema\in \Sbb^n$, where $\Sbb^n$ is the set of all real symmetric matrices of dimension $n\times n$. Let $\{\lambda_i(\cdot)\}_{i=1}^n$ be the eigenvalues of corresponding matrix such that $\lambda_1(\cdot)\geq \lambda_2(\cdot)\geq \ldots\geq \lambda_n(\cdot)$. 
	\begin{enumerate}[label=(\alph*)]
		\item \label{item:perthea}
		$$
		|\lambda_i(\Dma+\Ema)-\lambda_i(\Dma)|\leq \|\Ema\|_2 \quad (i=1,2,\ldots,n)
		$$
		\item Assume $\Ema=\omega \xve\xve^\top$, where $\xve\in S^{n-1}$. If $\omega\geq0$, then
		$$
		\lambda_i(\Dma+\Ema)\in [\lambda_i(\Dma),\lambda_{i-1}(\Dma)], \quad ( i=2,3,\ldots,n),
		$$
		while if $\omega\leq0$, then
		$$
		\lambda_i(\Dma+\Ema)\in [\lambda_{i+1}(\Dma),\lambda_{i}(\Dma)], \quad ( i=1,2,\ldots,n-1).
		$$
		In either case, there exist nonnegative $m_1,m_2,\ldots,m_n$ such that
		$$
		\lambda_i(\Dma+\Ema)=\lambda_i(\Dma)+m_i\omega,  \quad (i=1,2,\ldots,n)
		$$
		with $m_1+m_2+\cdots+m_n=1$.
		\item \label{item:perthec}
		Assume $\Ema= \Sum_{i=1}^m\omega_i \xve_i\xve_i^\top$, where $\{\xve_i\}_{i=1}^m \subset S^{n-1}$ and $\omega_i\geq 0$ for all $i$. Then
		$$
		\lambda_n(\Dma+\Ema)\geq \lambda_n(\Dma).
		$$
	\end{enumerate}
\end{lem}
\begin{proof}
	(a) and (b) is Corollary 8.1.6 and Theorem 8.1.8 in \cite{golub2012matrix}. (c) follows by induction on the smallest eigenvalue using part (b) for $\omega \geq 0$.
\end{proof}

\begin{lem}
\label{isoine}
Let $\xve_1,\xve_2$ be two vectors on $S^{n-1}$, and $\Dma\in \R^{n\times n}$ be an invertible matrix. Let $\smax(\Dma)$ and $\smin(\Dma)$ be respectively the largest and smallest singular value of $\Dma$. Define $\bar{\wwve}_i = \Dma \xve_i$ and $\wwve_i=\bar{\wwve}_i/\|\bar{\wwve}_i\|_2$, $(i=1,2)$. Then,
$$  \frac{\smin(\Dma)}{\smax(\Dma)}  \|\xve_1-\xve_2\|_2 \leq \|\wwve_1-\wwve_2\|_2 \leq \frac{\smax(\Dma)}{\smin(\Dma)}\|\xve_1-\xve_2\|_2$$
\end{lem}
\begin{proof}
{\bf Part I (Upper Bound)}

\begin{figure}
\center
\includegraphics[scale=0.3]{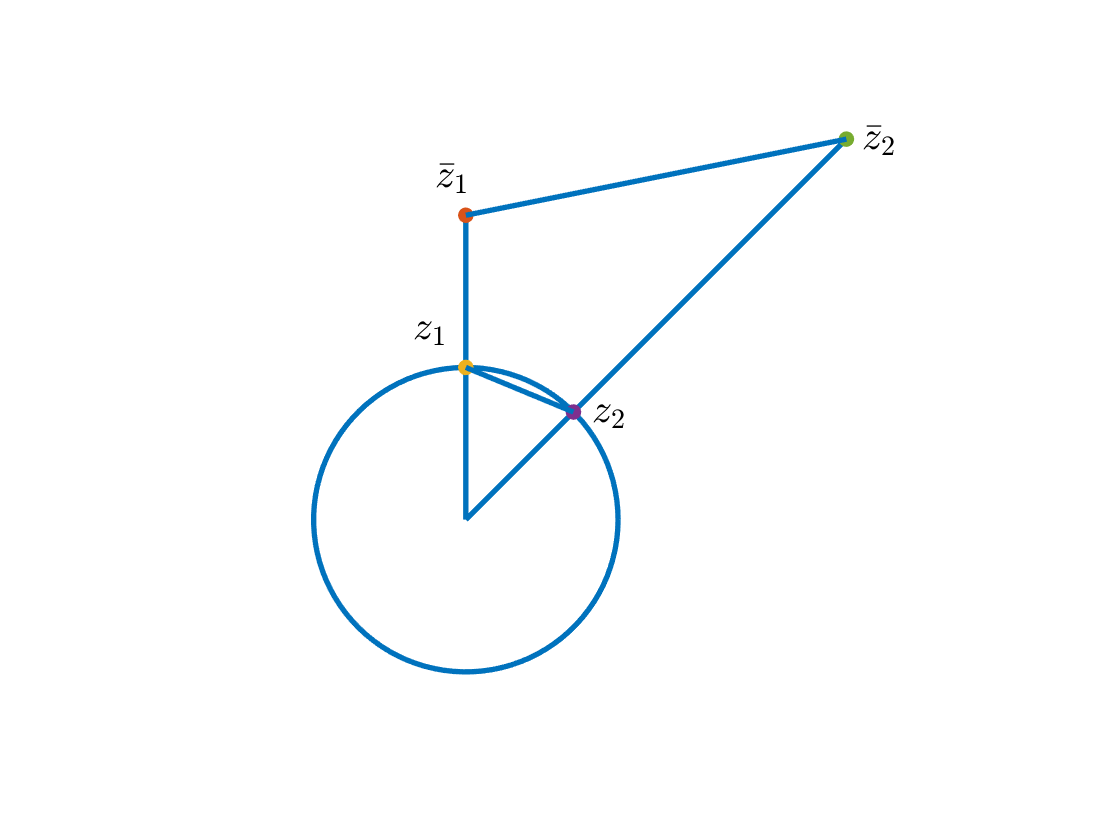}
\caption{$\wwve_1$ and $\wwve_2$ are the normalized vector of $\bar{\wwve}_1$ and $\bar{\wwve}_2$ respectively. }
\label{fig:geogra}
\end{figure}

Denote $\angle(\cdot,\cdot)$ the angle between two vectors. By the Law of Cosines,
$$
\cos(\angle(\wwve_1,\wwve_2))= \frac{\|\wwve_1\|_2^2+\|\wwve_2\|_2^2-\|\wwve_1-\wwve_2\|_2^2}{2\times \|\wwve_1\|_2 \times \|\wwve_2\|_2}= \frac{2-\|\wwve_1-\wwve_2\|_2^2}{2},
$$
and
$$
\cos(\angle(\bar{\wwve}_1,\bar{\wwve}_2))= \frac{\|\bar{\wwve}_1\|_2^2+\|\bar{\wwve}_2\|_2^2-\|\bar{\wwve}_1-\bar{\wwve}_2\|_2^2}{2\times \|\bar{z}_1\|_2 \times \|\bar{\wwve}_2\|_2}.
$$
Observing $\angle(\wwve_1,\wwve_2) = \angle(\bar{\wwve}_1,\bar{\wwve}_2)$, the right hand sides of the above two equations are equal. Solving for $\|\wwve_1-\wwve_2\|_2$, we get
$$
\|\wwve_1-\wwve_2\|_2^2 = \frac{\|\bar{\wwve}_1-\bar{\wwve}_2\|_2^2}{\|\bar{\wwve}_1\|_2\|\bar{\wwve}_2\|_2}+\left( 2 - \frac{\|\bar{\wwve}_2\|_2}{\|\bar{\wwve}_1\|_2} - \frac{\|\bar{\wwve}_1\|_2}{\|\bar{\wwve}_2\|_2}\right) \leq \frac{\|\bar{\wwve}_1-\bar{\wwve}_2\|_2^2}{\|\bar{\wwve}\|_2\|\bar{\wwve}_2\|_2}.
$$
Therefore,
\begin{align*}
\|\wwve_1-\wwve_2\|_2 &\leq \frac{\|\bar{\wwve}_1-\bar{\wwve}_2\|_2}{\sqrt{\|\bar{\wwve}_1\|_2\|\bar{\wwve}_2\|_2}} \\
& \leq \frac{\smax(D)\|\xve_1-\xve_2\|_2}{\sqrt{\smin(\Dma)\|\xve_1\|_2\smin(\Dma)\|\xve_2\|_2}} \\
& =  \frac{\smax(\Dma)}{\smin(\Dma)}\|\xve_1-\xve_2\|_2.
\end{align*}
{\bf Part II(Lower Bound)}\\
Define $\bar{\xve}_i = \Dma^{-1}\wwve_i$, $(i=1,2)$. Note for $\forall i\in \{1,2\}$, $\xve_i$ and $\bar{\xve}_i$ are parallel to each other, since $\xve_i=\Dma^{-1}\bar{\wwve}_i$ and $\bar{\wwve}_i$ is parallel to $\wwve_i$. Thus, we conclude $\xve_i=\bar{\xve}_i/\|\bar{\xve}_i\|_2$, $(i=1,2)$. Reversing the role of $\xve_i$ and $\wwve_i$ in Part I, one has
$$
\|\xve_1-\xve_2\|_2 \leq \frac{\smax(\Dma^{-1})}{\smin(\Dma^{-1})}\|\wwve_1-\wwve_2\|_2.
$$
The lower bound follows from the relation $\frac{\smax(\Dma^{-1})}{\smin(\Dma^{-1})} = \frac{\smax(\Dma)}{\smin(\Dma)}$.
\end{proof}

\begin{lem}
\label{setrelation}
Let $\{\Dc_i\}_{i=1}^m$, $\{\Fc_i\}_{i=1}^m$, $\{\Gc_i\}_{i=1}^m$ and $\{\Hc_i\}_{i=1}^m$ be sets satisfying
$$
\Gc_i \subset \Dc_i\subset \Hc_i,\quad  \Gc_i \subset \Fc_i\subset \Hc_i, \quad (i=1,2,\ldots,m).
$$
Then
\begin{enumerate}[label=(\alph*)]
	\item \label{prop:setinclusiona}
$$
\left(\bigcap_{i=1}^m \Dc_i\right) \bigtriangleup \left(\bigcap_{i=1}^m \Fc_i\right) \subset \bigcup_{i=1}^m \left(\Hc_i \backslash \Gc_i \right)\bigcap \left(\bigcap_{j=1}^m\Hc_j\right)=\bigcup_{i=1}^m \left( \left(\Hc_i \backslash \Gc_i\right) \bigcap \left(\bigcap_{\substack{j=1\\j\not =i}}^m\Hc_j\right) \right) .
$$
\item  \label{item:setinclusionb}
$$
\left(\bigcup_{i=1}^m \Dc_i\right) \bigtriangleup \left(\bigcup_{i=1}^m \Fc_i\right) \subset \bigcup_{i=1}^m \left(\Hc_i \backslash \Gc_i \right).
$$
\end{enumerate}
\end{lem}
\begin{proof}
(a)
Obviously,
\begin{equation}
\bigcap_{i=1}^m \Gc_i \subset \bigcap_{i=1}^m \Dc_i \subset \bigcap_{i=1}^m \Hc_i, \quad  \bigcap_{i=1}^m \Gc_i \subset \bigcap_{i=1}^m \Fc_i \subset \bigcap_{i=1}^m \Hc_i. \label{eqn:setinc1}
\end{equation}
Thus,
\begin{align*}
\left(\bigcap_{i=1}^m \Dc_i\right) \bigtriangleup \left(\bigcap_{i=1}^m \Fc_i\right) &\subset \left(\bigcap_{i=1}^m \Hc_i\right) \backslash \left(\bigcap_{i=1}^m \Gc_i \right).\\
\end{align*}
Take $\forall \omega \in \left(\bigcap_{i=1}^m \Hc_i\right) \backslash \left(\bigcap_{i=1}^m \Gc_i \right)$, we know $\omega\in \bigcap_{i=1}^m \Hc_i$ and $\omega\not\in \bigcap_{i=1}^m \Gc_i$. The later fact shows $\exists j$ (which depends on $\omega$) such that $\omega \not\in \Gc_j$. Then,
\begin{equation}
\omega\in \left(\bigcap_{i=1}^m \Hc_i\right) \backslash \Gc_j\subset \Hc_j\backslash \Gc_j\subset  \bigcup_{i=1}^m \left(\Hc_i \backslash \Gc_i \right). \label{eqn:setinc2}
\end{equation}
The proof is completed by combining \eqref{eqn:setinc1} and \eqref{eqn:setinc2}.\\

\noindent (b)
\begin{align*}
\left(\bigcup_{i=1}^m \Dc_i\right) \bigtriangleup \left(\bigcup_{i=1}^m \Fc_i\right) &= \left(\bigcup_{i=1}^m \Dc_i\right)^c \bigtriangleup \left(\bigcup_{i=1}^m \Fc_i\right)^c \\
& = \left(\bigcap_{i=1}^m \Dc_i^c\right) \bigtriangleup \left(\bigcap_{i=1}^m \Fc_i^c\right) \\
& \subset \bigcup_{i=1}^m \left(\Gc_i^c \backslash \Hc_i^c \right)\\
& = \bigcup_{i=1}^m \left( \Hc_i \backslash \Gc_i  \right),
\end{align*}
where the inclusion step follows from \ref{prop:setinclusiona}.
\end{proof}

\begin{lem}
	\label{proeig}
	Let $\Qma\in \mathbb{R}^{n \times m}(n\leq m)$, with each column $\qve_i$ being i.i.d.  $\unif(\sqrt{n}S^{n-1})$. Let $\lambdamin$ and $\lambdamax$ be respectively the largest and smallest eigenvalue of $\frac{1}{m}\Qma\Qma^\top$. Then with probability at least $ 1 - 2\exp(-ct^2)$,
	\begin{equation}
	\left[1-C\left(\sqrt{\frac{n}{m}}+\frac{t}{\sqrt{m}}\right)\right]^2 \leq \lambdamin \leq \lambdamax \leq \left[1+C\left(\sqrt{\frac{n}{m}}+\frac{t}{\sqrt{m}}\right)\right]^2,    \label{eqn:eigupplowbou}
	\end{equation}
	where $c,C$ are absolute constants.
\end{lem}
\begin{proof}
	Let $\smax$, $\smin$ be respectively the largest and smallest singular value of $\Qma$. Since columns $\{\qve_i\}_{i=1}^m$ are isotropic random vectors with subgaussian norm (or $\psi_2$ norm) being a constant, by applying Theorem 5.39 in \cite{vershynin2012introduction} to $Q^\top$,
	\begin{equation}
	\sqrt{m}-C(\sqrt{n}+t) \leq \smin \leq \smax \leq \sqrt{m}+C(\sqrt{n}+t),
	\end{equation}
	holds with probability at least $ 1 - 2\exp(-ct^2)$, where $c,C$ are absolute constants. The proof is completed by 
	$$
	\lambdamax = \frac{1}{m}\smax^2, \quad \lambdamin=\frac{1}{m} \smin^2.
	$$
\end{proof}

\subsection{Numerical simulations and experiments}
\label{sec:simulation}

\begin{figure}[ht]
\begin{subfigure}{.49\textwidth}
  \centering
  \includegraphics[width=.8\linewidth]{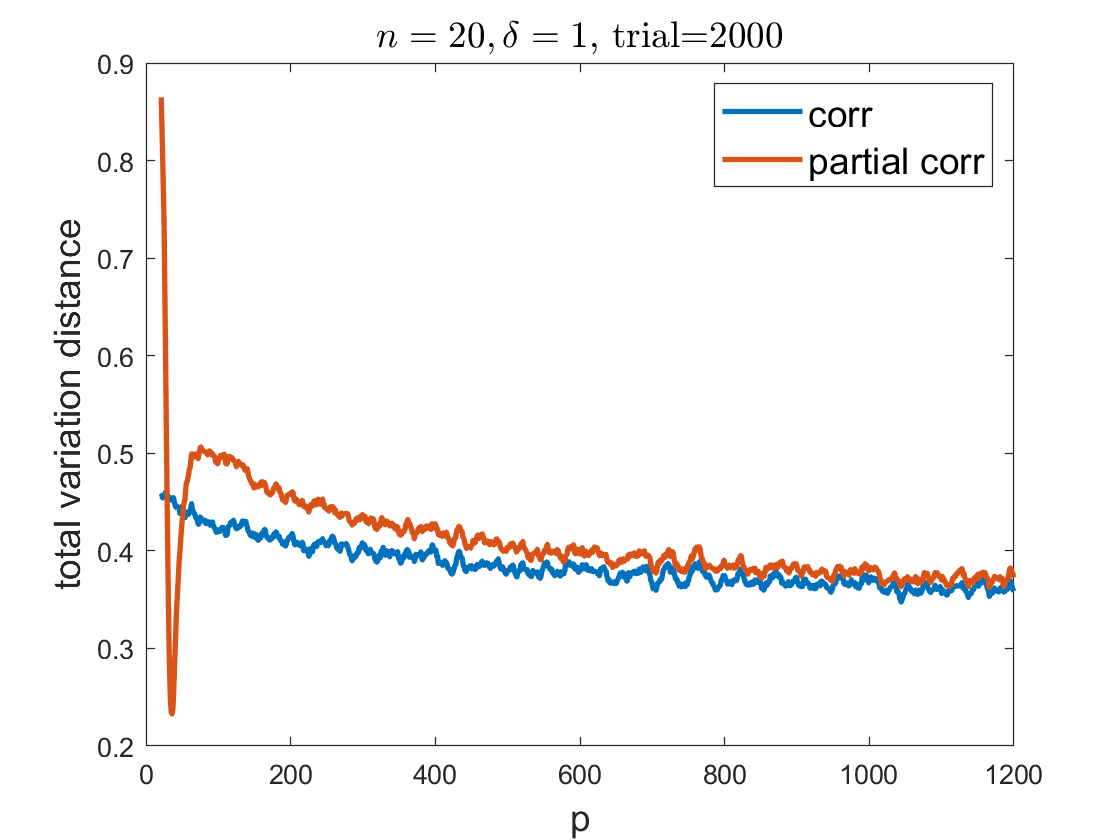}  
  \caption{limiting compound Poisson}
  \label{fig:finiteplimitcomparisona}
\end{subfigure}
\begin{subfigure}{.49\textwidth}
  \centering
  \includegraphics[width=.8\linewidth]{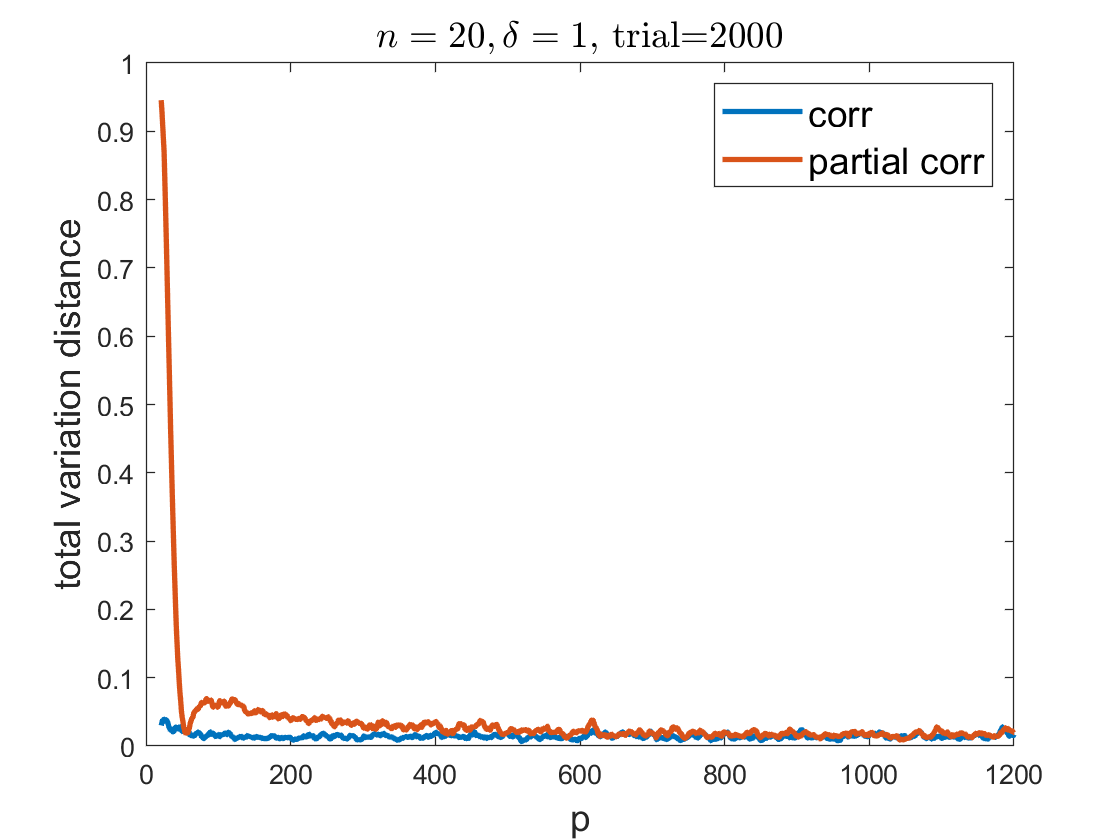}  
  \caption{finite $p$ compound Poisson approximation}
  \label{fig:finiteplimitcomparisonb}
\end{subfigure}
\caption{The vertical axis of (a) is $d_{\TV}(\mathscr{L}\left({N_{V_1}^{(\Psima)}}\right),\CP(\lambda_{20,1}(1),\bm{\zeta}_{20,1}))$ and that of (b) is $d_{\TV}(\mathscr{L}\left(N_{V_1}^{(\Psima)}\right),\CP(\lambda_{p,20,1,\rho},\bm{\zeta}_{20,1,\rho}))$. 
\ONE{The theoretical convergence results have respectively been established in Theorem \ref{cor:Poissonlimit} and Theorem \ref{thm:Poissonultrahigh}.} For both plots the samples are independently generated according to $\mathcal{N}(\bm{0},\Sigmama)$ with $\Sigmama$ being a $(\tau=p^{0.6},\kappa=p^{0.8})$ sparse matrix for each $p$. The parameters are $n=20$, $\delta=1$ and the threshold $\rho$ is chosen according to \eqref{eqn:rhopformula} with $e_{n,\delta}=1$.
The blue curve is for the empirical correlation graph ($\Psima =\Rma$) and the red curve is for the empirical partial correlation graph ($\Psima =\Pma$). Note since $\delta=1$, $\bm{\zeta}_{20,1}=\Dirac(2)=\bm{\zeta}_{20,1,\rho}$, by Example \ref{exa:limitcompoidelta1}. As demonstrated by the plots, for both the empirical correlation and partial correlation graphs, the total variations in (a) decrease very slowly while the total variations in (b) converge to $0$ very fast, which has also been analytically discussed in Remark \ref{rem:comparisonsfinitelimit}. 
\ONE{Our observation that the non-asymptotic compound Poisson distribution provides a better fit to the numerical simulations for small $p$ is a caveat to practitioners who may be tempted to use the Poisson approximation.}
}
    \label{fig:finiteplimitcomparison}
\end{figure}

\begin{figure}[ht]
\begin{subfigure}{.49\textwidth}
  \centering
  \includegraphics[width=.8\linewidth]{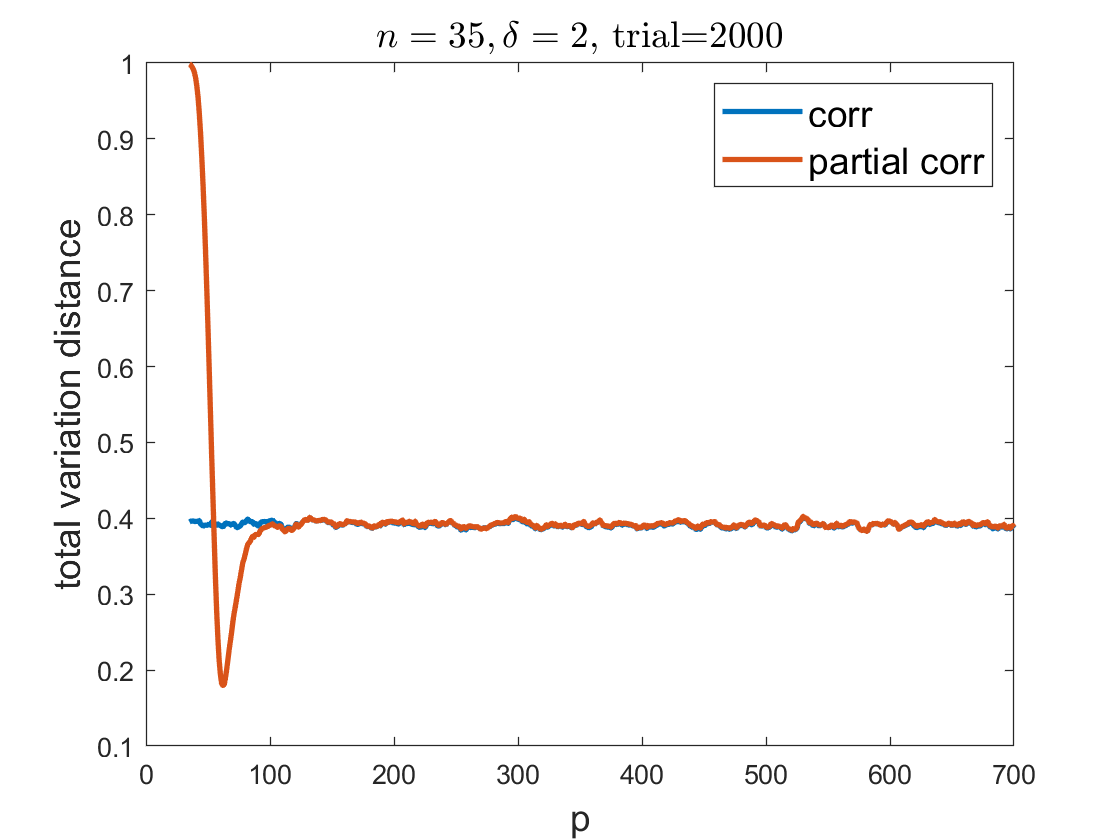}  
  \caption{Poisson limit when $p\to \infty$}
\end{subfigure}
\begin{subfigure}{.49\textwidth}
  \centering
  \includegraphics[width=.8\linewidth]{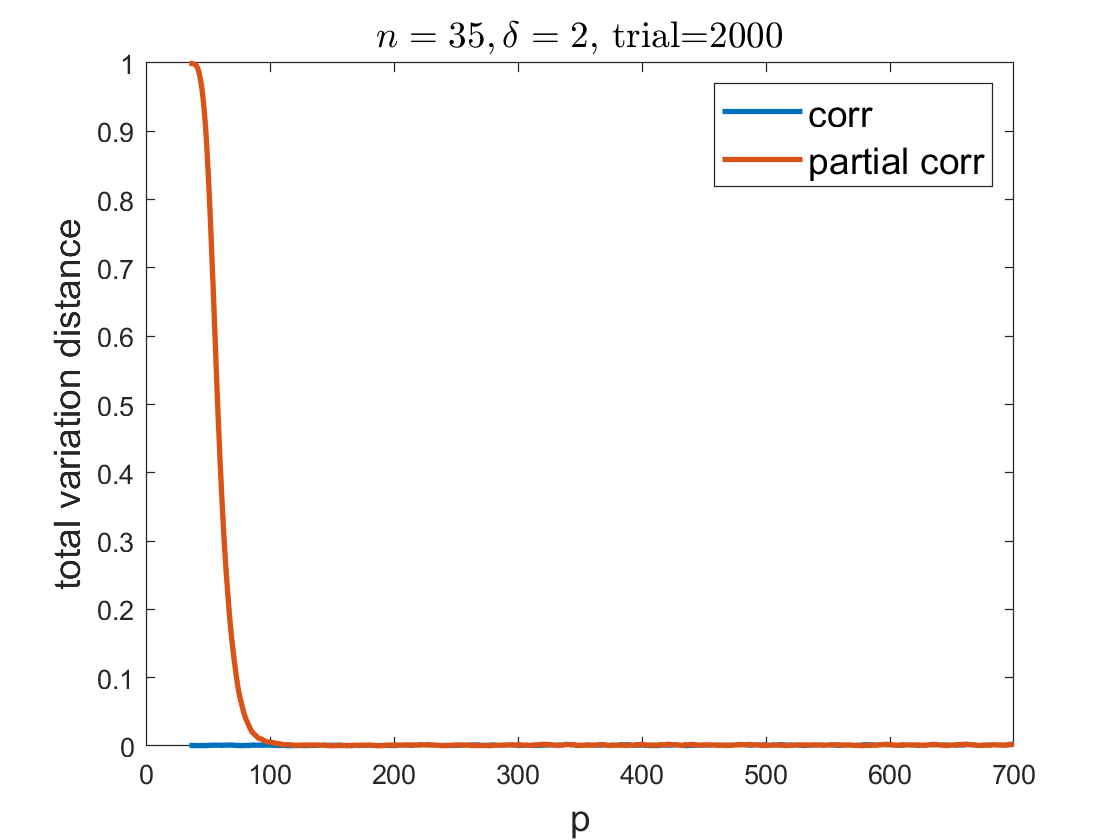}  
  \caption{Poisson approximation for finite $p$}
\end{subfigure}
\caption{ The vertical axis of (a) is $d_{\TV}\left(\mathscr{L}\left(N_{V_\delta}^{(\Psima)}\right),\Pois\left(\frac{(e_{n,\delta})^\delta}{\delta!}\right) \right)$, where we replaced $\CP(\lambda_{n,\delta}(e_{n,\delta}),  \bm{\zeta}_{n,\delta})$ in Theorem \ref{cor:Poissonlimit} by its approximation $\Pois(\frac{(e_{n,\delta})^\delta}{\delta!})$ as in Proposition \ref{prop:poiapplimit}. The vertical axis of (b) is $d_{\TV}(N_{V_\delta}^{(\Psima)},\Pois(\binom{p}{1}\binom{p-1}{\delta}(2P_n(r_\rho))^{\delta}))$, where we replaced $\CP(\lambda_{p,n,\delta,\rho}, \bm{\zeta}_{n,\delta,\rho})$ in Theorem \ref{thm:Poissonultrahigh} by its approximation $\Pois(\binom{p}{1}\binom{p-1}{\delta}(2P_n(r_\rho))^{\delta})$ as in Proposition \ref{prop:poiappfinitep}. For both plots the samples are independently generated according to $\mathcal{N}(\bm{0},\Sigmama)$ with $\Sigmama$ being a $(\tau=p^{0.6},\kappa=p^{0.8})$ sparse matrix for each $p$. The parameters are $n=35$, $\delta=2$ and the threshold $\rho$ is chosen according to \eqref{eqn:rhopformula} with $e_{n,\delta}=1$.
 As demonstrated by the plots, for both the empirical correlation and partial correlation graphs, the total variations in (a) decrease very slowly while the total variations in (b) converge to $0$ very fast. The fast convergence in Figure \ref{fig:Poissonfinitevslimit} (b) verifies the validity of using Poisson distribution $\Pois(\binom{p}{1}\binom{p-1}{\delta}(2P_n(r_\rho))^{\delta})$ to approximate the distribution of random quantities in $\{N_{E_\delta}^{(\Psima)},N_{\breve{V}_\delta}^{(\Psima)},N_{V_\delta}^{(\Psima)}: \Psima \in \{\Rma,\Pma\}  \}$ for large $n$. We now discuss the slow convergence behavior of Figure \ref{fig:Poissonfinitevslimit} (a).
 Note that $n=35$ is large enough to guarantee that $\sum_{\ell=2}^{\delta+1}\alpha_\ell$ is small as indicated by Figure \ref{fig:incrementtodeltauppbou} (b), which implies that $\CP(\lambda_{n,\delta}(e_{n,\delta}),  \bm{\zeta}_{n,\delta})$ is well approximated by $\Pois(\frac{(e_{n,\delta})^\delta}{\delta!})$ by Lemma \ref{prop:jumpsizedecay} \ref{item:jumpsizedecayc}. As a result, the extremely slow decrease in Figure \ref{fig:Poissonfinitevslimit} (a) is not because of using the Poisson approximation, but is
 due to the slow convergence of Theorem \ref{cor:Poissonlimit}, which has been extensively discussed in Remark \ref{rem:comparisonsfinitelimit}. 
This specific example additionally indicates that the slow convergence of Theorem \ref{cor:Poissonlimit} is due to slow convergence of $\lambda_{p,n,\delta,\rho}\to \lambda_{n,\delta}$ since the increment distribution in this large $n$ case are both close to $\Dirac(1)$. } 
\label{fig:Poissonfinitevslimit}
\end{figure}

\end{appendix}

\end{document}